# Some series and integrals involving the Riemann zeta function, binomial coefficients and the harmonic numbers

## Volume IV

Donal F. Connon

18 February 2008

## Abstract

In this series of seven papers, predominantly by means of elementary analysis, we establish a number of identities related to the Riemann zeta function, including the following:

$$\sum_{n=1}^{\infty} \frac{H_n^{(1)}}{n^s} x^n = \int_0^x \frac{Li_s(x) - Li_s(y)}{x - y} \, dy$$

$$\sum_{n=1}^{\infty} \frac{H_n^{(3)}}{n 2^n} = \varsigma(2)\log^2 2 - \frac{7}{8}\varsigma(3)\log 2 + Li_4(1/2) - \frac{1}{6}\log^4 2$$

$$\int_0^{\infty} u^{x-1} e^{-ku} \log^n u \, du = \frac{1}{k^x} \sum_{j=0}^n (-1)^{n-j} \binom{n}{j} \Gamma^{(j)}(x) \log^{n-j} k$$

$$\sum_{n=1}^{\infty} \frac{1}{2^{n+1}} \sum_{k=1}^n \binom{n}{k} \frac{H_k^{(2)}}{k^4} = \varsigma(2)\varsigma(4) + \varsigma^2(3) - \frac{25}{12}\varsigma(6)$$

$$\sum_{n=1}^{\infty} \frac{1}{2^{n+1}} \sum_{k=1}^n \binom{n}{k} \frac{H_k^{(3)}}{k^3} = \frac{1}{2}\varsigma^2(3) + \frac{1}{2}\varsigma(6)$$

$$\frac{1}{3}\sum_{n=0}^{\infty} \frac{1}{n+1} \sum_{k=0}^n \binom{n}{k} \frac{(-1)^k H_{k+1}^{(2)}}{(k+1)^3} = \frac{4}{3}\varsigma^2(3) - \frac{29}{12}\varsigma(6) + \varsigma(2)\varsigma(4)$$

$$\frac{1}{2}\int_0^1 \frac{\log^2(1-x)\log[1 - u(1-x)]}{1-x} \, dx = -Li_4(u)$$

$$\sum_{n=1}^{\infty} \frac{\log n \cos(n\pi/2)}{(\pi n)^2} = \frac{1}{48}\big[\log(2\pi) + \gamma - 1\big] + \frac{1}{4}\varsigma'(-1)$$

Whilst the paper is mainly expository, some of the formulae reported in it are believed to be new, and the paper may also be of interest specifically due to the fact that most of the various identities have been derived by elementary methods.

**CONTENTS OF VOLUMES I TO VI:**　　　　　　　　**Volume/page**

**SECTION:**

**APPENDICES (Volume VI):**

**A**. Some properties of the Bernoulli numbers and the Bernoulli polynomials

**B**. A well-known integral

**C**. Euler's reflection formula for the gamma function and related matters

**D**. A very elementary proof of $\dfrac{\pi^2}{8} = \sum\limits_{n=0}^{\infty} \dfrac{1}{(2n+1)^2}$

**E**. Some aspects of Euler's constant $\gamma$ and the gamma function

**F**. Elementary aspects of Riemann's functional equation for the zeta function

**ACKNOWLEDGEMENTS**

**REFERENCES**

**(xii)** From (4.1.6) we have

$$(4.4.118) \qquad \sum_{k=1}^{n} \binom{n}{k}(-1)^{k+1}\frac{t^k}{k} = \sum_{k=1}^{n}\frac{1-(1-t)^k}{k}$$

Differentiating (4.4.118) gives

$$(4.4.119) \qquad \sum_{k=1}^{n}\binom{n}{k}(-1)^{k+1}t^{k-1} = \sum_{k=1}^{n}(1-t)^{k-1}$$

$$(4.4.120) \qquad = \frac{1-(1-t)^n}{t}$$

Now multiply (4.4.120) by $t$ and integrate over $[0, x]$ to obtain

$$(4.4.121) \qquad \sum_{k=1}^{n}\binom{n}{k}(-1)^{k+1}\frac{x^{k+1}}{k+1} = x + \frac{(1-x)^{n+1}-1}{n+1}$$

Moving the term $x$ to the left-hand side, we have

$$(4.4.122) \qquad \sum_{k=0}^{n}\binom{n}{k}(-1)^{k}\frac{x^{k+1}}{k+1} = \frac{1-(1-x)^{n+1}}{n+1}$$

Letting $x = 1$ we have (as in (4.2.3))

$$(4.4.123) \qquad \sum_{k=0}^{n}\binom{n}{k}\frac{(-1)^k}{k+1} = \frac{1}{n+1}$$

With $x = -1$ we get

$$(4.4.123a) \qquad \sum_{k=0}^{n}\binom{n}{k}\frac{1}{k+1} = \frac{2^{n+1}-1}{n+1}$$

Alternatively, multiply (4.4.120) by $t^u$ and integrate over $[0, y]$ to obtain

$$(4.4.123b) \qquad \sum_{k=1}^{n}\binom{n}{k}(-1)^{k+1}\frac{y^{k+u}}{k+u} = y^u - \int_{0}^{y}t^{u-1}(1-t)^n\,dt$$

Now divide (4.4.122) by $x$ and integrate to produce

$$(4.4.124) \qquad \sum_{k=0}^{n} \binom{n}{k} (-1)^k \frac{t^{k+1}}{(k+1)^2} = \frac{1}{n+1} \int_0^t \frac{1-(1-x)^{n+1}}{x} dx$$

$$(4.4.125) \qquad = \frac{1}{n+1} \sum_{k=1}^{n+1} \binom{n+1}{k} (-1)^{k+1} \frac{t^k}{k}$$

where we have used (4.1.3) to evaluate the integral. Using (4.1.6) this becomes

$$(4.4.126) \qquad \sum_{k=0}^{n} \binom{n}{k} (-1)^k \frac{t^{k+1}}{(k+1)^2} = \frac{1}{n+1} \sum_{k=1}^{n+1} \frac{1-(1-t)^k}{k}$$

Letting $t = 1$ we have (in agreement with (4.2.16))

$$(4.4.127) \qquad \sum_{k=0}^{n} \binom{n}{k} \frac{(-1)^k}{(k+1)^2} = \frac{H_{n+1}}{n+1}$$

As before, now divide (4.4.126) by $t$ and integrate to produce

$$(4.4.128) \qquad \sum_{k=0}^{n} \binom{n}{k} (-1)^k \frac{x^{k+1}}{(k+1)^3} = \frac{1}{n+1} \sum_{k=1}^{n+1} \frac{1}{k} \int_0^x \frac{1-(1-t)^k}{t} dt$$

$$(4.4.129) \qquad = \frac{1}{n+1} \sum_{k=1}^{n+1} \frac{1}{k} \left\{ \sum_{l=1}^{k} \frac{1-(1-x)^l}{l} \right\}$$

Letting $x = 1$ we have

$$(4.4.130) \qquad \sum_{k=0}^{n} \binom{n}{k} \frac{(-1)^k}{(k+1)^3} = \frac{1}{n+1} \sum_{k=1}^{n+1} \frac{H_k}{k} = \frac{1}{2(n+1)} \left\{ \left( H_{n+1}^{(1)} \right)^2 + H_{n+1}^{(2)} \right\}$$

where we have used Adamchik's formula (3.18)

$$(4.4.131) \qquad = \frac{1}{n+1} \left\{ \sum_{k=1}^{n} \frac{H_k}{k} + \frac{H_{n+1}}{n+1} \right\}$$

$$(4.4.132) \qquad = \frac{1}{n+1} \sum_{k=1}^{n} \frac{H_k}{k} + \frac{H_{n+1}}{(n+1)^2}$$

Using (4.4.130) we have

$$(4.4.133) \qquad \sum_{n=0}^{\infty} \frac{1}{(n+1)} \sum_{k=0}^{n} \binom{n}{k} \frac{(-1)^k}{(k+1)^3} = \sum_{n=0}^{\infty} \frac{1}{2(n+1)^2} \left\{ \left(H_{n+1}^{(1)}\right)^2 + H_{n+1}^{(2)} \right\}$$

$$(4.4.134) \qquad = \frac{1}{2} \sum_{n=1}^{\infty} \frac{1}{n^2} \left\{ \left(H_n^{(1)}\right)^2 + H_n^{(2)} \right\}$$

and we have previously seen this in (4.2.20). Some of the above identities are used in Appendix E of Volume VI to evaluate derivatives of the gamma function. Additional proofs of (4.4.123) and (4.4.127) are given in [94aa].

Reference should also be made to the paper by Larcombe et al. [95] where they show that

$$(4.4.135) \qquad m \binom{m+n}{n} \sum_{k=0}^{n} \binom{n}{k} \frac{(-1)^k}{m+k} = 1$$

$$(4.4.135a) \qquad m \binom{m+n}{n} \sum_{k=0}^{n} \binom{n}{k} \frac{(-1)^k}{(m+k)^2} = \sum_{k=m}^{m+n} \frac{1}{k}$$

$$(4.4.135b) \qquad 2m \binom{m+n}{n} \sum_{k=0}^{n} \binom{n}{k} \frac{(-1)^k}{(m+k)^3} = \left(\sum_{k=m}^{m+n} \frac{1}{k}\right)^2 + \sum_{k=m}^{m+n} \frac{1}{k^2}$$

$$(4.4.136) \qquad 6m \binom{m+n}{n} \sum_{k=0}^{n} \binom{n}{k} \frac{(-1)^k}{(m+k)^4} = \left(\sum_{k=m}^{m+n} \frac{1}{k}\right)^3 + 3\left(\sum_{k=m}^{m+n} \frac{1}{k}\right)\left(\sum_{k=m}^{m+n} \frac{1}{k^2}\right) + 2\sum_{k=m}^{m+n} \frac{1}{k^3}$$

for integers $m \geq 1, n \geq 0$. In their paper they employ integrals of the type

$$(4.4.137) \qquad \int_0^{\infty} x^p e^{-mx} (1 - e^{-x})^n \, dx.$$

and, if we use the substitution $t = e^{-x}$, we can immediately see the relationship with (4.4.16) and also note the similarity with Anglesio's identity (4.4.93).

Note that by letting $m = 1$ in (4.4.135a) and (4.4.135b), we obtain (4.4.127) and (4.4.130) respectively.

In addition, Larcombe et al. also show how the above formulae can be derived using the identity from Gould's book, "Combinatorial Identities" [73]

$$(4.4.138) \qquad f(x+y) = y \binom{y+n}{n} \sum_{k=0}^{n} (-1)^k \binom{n}{k} \frac{f(x-k)}{y+k}$$

where $f(t)$ is a polynomial of degree $\leq n$. In [72], Gould attributes the formula (4.4.138) to Z.A. Melzak. Reference should also be made to the paper by Kirschenhofer [86].

The Larcombe identities are also employed in Appendix E in Volume VI of this series of papers.

An alternative proof of the Larcombe identities is set out below.

We have seen in (4.4.13) that

$$g(x) = \frac{n!}{x(1+x)...(n+x)} = \int_0^1 t^{x-1}(1-t)^n\, dt = \sum_{k=0}^{n} \binom{n}{k} \frac{(-1)^k}{k+x}$$

$$= \frac{\Gamma(n+1)\Gamma(x)}{\Gamma(x+n+1)} = B(n+1, x) = B(x, n+1)$$

and the first Larcombe identity becomes evident upon putting $x = m$.

Upon differentiating the above equation we obtain

$$\frac{g'(x)}{n!} = \frac{\Gamma(x+n+1)\Gamma'(x) - \Gamma(x)\Gamma'(x+n+1)}{\Gamma^2(x+n+1)}$$

$$= \frac{\Gamma(x)}{\Gamma(x+n+1)}\left\{\frac{\Gamma'(x)}{\Gamma(x)} - \frac{\Gamma'(x+n+1)}{\Gamma(x+n+1)}\right\}$$

which gives us

(4.4.138a)    $g'(x) = g(x)\left\{\psi(x) - \psi(x+n+1)\right\}$

Hence we have

$$g(x)\left\{\psi(x+n+1) - \psi(x)\right\} = \sum_{k=0}^{n}\binom{n}{k}\frac{(-1)^k}{(k+x)^2}$$

and for $x = m$ we obtain

$$\frac{\Gamma(m)\Gamma(n+1)}{\Gamma(m+n+1)}\left\{\psi(m+n+1) - \psi(m)\right\} = \sum_{k=0}^{n}\binom{n}{k}\frac{(-1)^k}{(k+m)^2}$$

where it is not necessary to have $m$ an integer at this stage.

We see from (E.14) that for $p = 1, 2, 3, ...$

$$\psi(p) = -\gamma + H_{p-1}$$

and hence

$$\psi(m+n+1) - \psi(m) = \sum_{k=m}^{m+n} \frac{1}{k} = H_{m+n} - H_{m-1}$$

Thus we arrive at the second Larcombe identity

$$m \binom{m+n}{n} \sum_{k=0}^{n} \binom{n}{k} \frac{(-1)^k}{(m+k)^2} = \sum_{k=m}^{m+n} \frac{1}{k}$$

Differentiating (4.4.138a) we obtain

$$g''(x) = g(x)\left\{\psi'(x) - \psi'(x+n+1)\right\} + g'(x)\left\{\psi(x) - \psi(x+n+1)\right\}$$

$$(4.4.138b) \qquad = g(x)\left\{\psi'(x) - \psi'(x+n+1)\right\} + g(x)\left\{\psi(x) - \psi(x+n+1)\right\}^2$$

From (E.16) we have $\psi'(p) = \varsigma(2) - H_{p-1}^{(2)}$ and the Larcombe result (4.4.135b) follows immediately. Differentiating (4.4.138b) easily results in a proof of (4.4.136b).

An application of the Larcombe identity is shown below. We refer to (4.4.24a)

$$\sum_{n=0}^{\infty} \frac{1}{2^{n+1}} \sum_{k=0}^{n} \binom{n}{k} \frac{(-1)^k}{(k+x)^s} = \sum_{k=0}^{\infty} \frac{(-1)^k}{(k+x)^s}$$

and for $x = m$ and $s = 2$ we obtain from (4.4.135a)

$$(4.4.138c) \qquad \sum_{n=0}^{\infty} \frac{1}{2^{n+1}} \frac{H_{m+n} - H_{m-1}}{m} \binom{m+n}{n}^{-1} = \sum_{k=0}^{\infty} \frac{(-1)^k}{(k+m)^2}$$

With $m = 2$ we get

$$(4.4.138d) \qquad \sum_{n=0}^{\infty} \frac{1}{2^{n+1}} \frac{H_{n+2} - 1}{(n+1)(n+2)} = \sum_{k=0}^{\infty} \frac{(-1)^k}{(k+2)^2}$$

and this may be written as

$$\sum_{n=0}^{\infty} \frac{1}{2^{n+1}} \frac{H_n}{(n+1)(n+2)} - \sum_{n=0}^{\infty} \frac{1}{2^{n+1}} \frac{n^2+n-1}{(n+1)^2(n+2)^2} = \sum_{k=0}^{\infty} \frac{(-1)^k}{(k+2)^2}$$

We see that

$$\varsigma_a(s) = \frac{1}{1^s} - \frac{1}{2^s} \ldots + \frac{(-1)^N}{(N-1)^s} + \left\{ \frac{(-1)^{N+1}}{N^s} + \frac{(-1)^{N+2}}{(N+1)^s} + \ldots \right\}$$

$$= \frac{1}{1^s} - \frac{1}{2^s} \ldots + \frac{(-1)^N}{(N-1)^s} + (-1)^{N+1} \sum_{n=0}^{\infty} \frac{(-1)^n}{(n+N)^s}$$

With $N = 2$ we have

$$\sum_{n=0}^{\infty} \frac{1}{2^{n+1}} \sum_{k=0}^{n} \binom{n}{k} \frac{(-1)^k}{(k+2)^s} = \sum_{n=0}^{\infty} \frac{(-1)^n}{(n+2)^s} = 1 - \varsigma_a(s)$$

Therefore we obtain

$$\sum_{n=0}^{\infty} \frac{1}{2^{n+1}} \frac{H_n}{(n+1)(n+2)} - \sum_{n=0}^{\infty} \frac{1}{2^{n+1}} \frac{n^2+n-1}{(n+1)^2(n+2)^2} = 1 - \varsigma_a(2) = 1 - \frac{\pi^2}{12}$$

Since we have from (3.30) $\frac{1}{2} \log^2(1-x) = \sum_{n=1}^{\infty} \frac{H_n}{n+1} x^{n+1}$, integration results in

$$\frac{1}{2} \int_0^t \log^2(1-x)\,dx = \int_0^t \sum_{n=1}^{\infty} \frac{H_n}{n+1} x^{n+1} dx = \sum_{n=1}^{\infty} \frac{H_n}{(n+1)(n+2)} t^{n+2}$$

$$\frac{1}{2} \int_0^t \log^2(1-x)\,dx = \frac{1}{2}(1-x)\left[ -1 + \log(1-x) - \frac{1}{2}\log^2(1-x) \right]\Bigg|_0^t$$

$$= \frac{1}{2}(1-t)\left[ -1 + \log(1-t) - \frac{1}{2}\log^2(1-t) \right] + \frac{1}{2}$$

With $t = 1/2$ we have

$$\sum_{n=0}^{\infty} \frac{1}{2^{n+1}} \frac{H_n}{(n+1)(n+2)} = \frac{1}{2}\left[ 1 - \log 2 - \frac{1}{2}\log^2 2 \right]$$

We therefore have

$$\sum_{n=0}^{\infty} \frac{1}{2^{n+1}} \frac{n^2+n-1}{(n+1)^2(n+2)^2} = \frac{\pi^2}{12} - \frac{1}{2}\log 2 - \frac{1}{4}\log^2 2 - \frac{1}{2}$$

Using Euler's series transformation, Knopp [90, p.270] reports that

$$\sum_{n=0}^{\infty} \frac{(-1)^n}{(n+1)(n+2)...(n+p+1)} = \frac{1}{p!} \sum_{n=0}^{\infty} \frac{1}{(p+n+1)2^{n+1}}$$

and with $p=1$ this becomes

$$\sum_{n=0}^{\infty} \frac{(-1)^n}{(n+1)(n+2)} = \frac{1}{2} \sum_{n=0}^{\infty} \frac{1}{(n+2)2^n}$$

We start with the series $S(a,b)$

$$S(a,b) = \sum_{n=0}^{\infty} \frac{(-1)^n}{(n+a)(n+b)} = \frac{1}{a-b} \sum_{n=0}^{\infty} (-1)^n \left[ \frac{1}{n+b} - \frac{1}{n+b} \right]$$

and follow the method of summation employed by Efthimiou [58a].

We have $\dfrac{1}{A} = \displaystyle\int_0^{\infty} e^{-Ax} dx$ and hence

$$S(a,b) = \frac{1}{a-b} \lim_{N \to \infty} \sum_{n=0}^{N} \int_0^{\infty} (-1)^n e^{-nx} (e^{-bx} - e^{-ax}) dx$$

$$= \frac{1}{a-b} \lim_{N \to \infty} \int_0^{\infty} (e^{-bx} - e^{-ax}) \frac{1 + (-1)^N e^{-(N+1)x}}{1 + e^{-x}} dx$$

$$= \frac{1}{a-b} \int_0^{\infty} \frac{(e^{-bx} - e^{-ax})}{1 + e^{-x}} dx$$

Making the substitution $t = e^{-x}$ this becomes

$$S(a,b) = \frac{1}{a-b} \int_0^1 \frac{t^b - t^a}{t(1+t)} dt$$

we therefore have

$$S(2,1) = -\int_0^1 \frac{t^2 - t}{t(1+t)} dt = -\int_0^1 \left(1 - \frac{2}{1+t}\right) dt = 2\log 2 - 1$$

A more direct proof is shown below: we easily see that

$$\sum_{n=0}^{\infty}\frac{1}{(n+2)2^n}=4\sum_{n=0}^{\infty}\frac{1}{(n+2)2^{n+2}}=4\left\{\frac{1}{2.2^2}+\frac{1}{3.2^3}+...\right\}$$

and using

$$\log 2=\frac{1}{1.2}+\left\{\frac{1}{2.2^2}+\frac{1}{3.2^3}+...\right\}$$

we get as before

$$\sum_{n=0}^{\infty}\frac{1}{(n+2)2^n}=4\log 2-2$$

With $m=2$ we may also write (4.4.138b) as

$$\sum_{n=0}^{\infty}\frac{1}{2^{n+1}}\frac{H_{n+2}-1}{(n+1)(n+2)}=\sum_{n=0}^{\infty}\frac{1}{2^{n+1}}\frac{H_{n+2}}{(n+1)(n+2)}-\sum_{n=0}^{\infty}\frac{1}{2^{n+1}}\frac{1}{(n+1)(n+2)}=\sum_{k=0}^{\infty}\frac{(-1)^k}{(k+2)^2}$$

$$=\sum_{n=0}^{\infty}\frac{1}{2^{n+1}}\frac{H_{n+2}}{n+1}-\sum_{n=0}^{\infty}\frac{1}{2^{n+1}}\frac{H_{n+2}}{n+2}-\sum_{n=0}^{\infty}\frac{1}{2^{n+1}}\frac{1}{(n+1)(n+2)}$$

$$\sum_{n=0}^{\infty}\frac{1}{2^{n+1}}\frac{H_{n+2}}{n+1}=\sum_{n=0}^{\infty}\frac{1}{2^{n+1}}\frac{H_{n+1}}{n+1}+\sum_{n=0}^{\infty}\frac{1}{2^{n+1}}\frac{1}{(n+1)^2}$$

$$=\sum_{k=1}^{\infty}\frac{1}{2^k}\frac{H_k}{k}+\sum_{k=1}^{\infty}\frac{1}{k^2 2^k}$$

$$=\sum_{k=1}^{\infty}\frac{1}{2^k}\frac{H_k}{k}+Li_2(1/2)$$

We also have

$$\sum_{n=0}^{\infty}\frac{1}{2^{n+1}}\frac{H_{n+2}}{n+2}=\frac{1}{2}\frac{H_2}{2}+\frac{1}{2^2}\frac{H_3}{3}+\frac{1}{2^3}\frac{H_4}{4}+...$$

$$=2\left[\frac{1}{2^2}\frac{H_2}{2}+\frac{1}{2^3}\frac{H_3}{3}+\frac{1}{2^4}\frac{H_4}{4}+...\right]$$

$$=2\sum_{k=1}^{\infty}\frac{1}{2^k}\frac{H_k}{k}-1$$

We then get

$$\sum_{n=0}^{\infty} \frac{1}{2^{n+1}} \frac{H_{n+2}-1}{(n+1)(n+2)} = 1 - \sum_{k=1}^{\infty} \frac{1}{2^k} \frac{H_k}{k} + Li_2(1/2) - \sum_{n=0}^{\infty} \frac{1}{2^{n+1}} \frac{1}{(n+1)(n+2)} = 1 - \frac{\pi^2}{12}$$

and therefore we have

$$\sum_{n=0}^{\infty} \frac{1}{2^{n+1}} \frac{1}{(n+1)(n+2)} = \frac{\pi^2}{12} - \sum_{k=1}^{\infty} \frac{1}{2^k} \frac{H_k}{k} + Li_2(1/2)$$

Using (3.105a) we see that

$$\sum_{k=1}^{\infty} \frac{1}{2^k} \frac{H_k}{k} = \frac{1}{2} \log^2 2 + Li_2(1/2)$$

and hence we obtain

$$\sum_{n=0}^{\infty} \frac{1}{2^{n+1}} \frac{1}{(n+1)(n+2)} = \frac{\pi^2}{12} - \frac{1}{2} \log^2 2$$

We saw in (4.4.79) in Volume III that

$$\varsigma_a(s,u) = \sum_{n=0}^{\infty} \frac{1}{2^{n+1}} \sum_{k=0}^{n} \binom{n}{k} \frac{(-1)^k}{(k+u)^s}$$

where $\varsigma_a(s,u)$ may be regarded as an alternating Hurwitz zeta function and, using (4.4.24a) in Volume II(b), this may be written as

$$\varsigma_a(s,u) = \sum_{n=0}^{\infty} \frac{(-1)^n}{(n+u)^s}$$

We then have from (4.4.135)

$$\varsigma_a(1,m) = \sum_{n=0}^{\infty} \frac{1}{2^{n+1}} \sum_{k=0}^{n} \binom{n}{k} \frac{(-1)^k}{k+m} = \sum_{n=0}^{\infty} \frac{1}{2^{n+1} m} \binom{m+n}{m}^{-1}$$

and in particular

$$\varsigma_a(1,1) = \sum_{n=0}^{\infty} \frac{1}{2^{n+1}} \sum_{k=0}^{n} \binom{n}{k} \frac{(-1)^k}{k+1} = \sum_{n=0}^{\infty} \frac{1}{n 2^{n+1}} = \log 2$$

which also implies the known result

$$\varsigma_a(1,1) = \varsigma_a(1) = \log 2$$

We also have from (4.4.135a)

$$\varsigma_a(2,m) = \sum_{n=0}^{\infty} \frac{1}{2^{n+1}} \sum_{k=0}^{n} \binom{n}{k} \frac{(-1)^k}{(k+m)^2} = \sum_{n=0}^{\infty} \frac{1}{2^{n+1} m} \binom{m+n}{m}^{-1} \sum_{k=m}^{m+n} \frac{1}{k}$$

With $m = 2$ we get

$$\varsigma_a(2,2) = \sum_{n=0}^{\infty} \frac{1}{2^{n+1}} \sum_{k=0}^{n} \binom{n}{k} \frac{(-1)^k}{(k+2)^2} = \sum_{n=0}^{\infty} \frac{1}{2^{n+1}(n+1)(n+2)} \sum_{k=2}^{n+2} \frac{1}{k}$$

We see that

$$\sum_{k=2}^{n+2} \frac{1}{k} = H_n + \frac{1}{n+1} + \frac{1}{n+2} - 1$$

and therefore we have

$$\varsigma_a(2,2) = \sum_{n=0}^{\infty} \frac{H_n}{2^{n+1}(n+1)(n+2)} + \sum_{n=0}^{\infty} \frac{1}{2^{n+1}(n+1)^2(n+2)}$$

$$+ \sum_{n=0}^{\infty} \frac{1}{2^{n+1}(n+1)(n+2)^2} - \sum_{n=0}^{\infty} \frac{1}{2^{n+1}(n+1)(n+2)}$$

We see that

$$\varsigma_a(2,2) = \sum_{n=0}^{\infty} \frac{(-1)^n}{(n+2)^2} = 1 - \varsigma_a(2)$$

and hence we obtain

$$1 - \varsigma_a(2) = \sum_{n=0}^{\infty} \frac{H_n}{2^{n+1}(n+1)(n+2)} + \sum_{n=0}^{\infty} \frac{1}{2^{n+1}(n+1)^2(n+2)}$$

$$+ \sum_{n=0}^{\infty} \frac{1}{2^{n+1}(n+1)(n+2)^2} - \sum_{n=0}^{\infty} \frac{1}{2^{n+1}(n+1)(n+2)}$$

(**xiii**) In (4.1.1) we considered the integral $J = \int_0^t \frac{1-(1-x)^n}{x} dx$ and evaluated it in two ways to show that

$$\sum_{k=1}^{n}\binom{n}{k}(-1)^{k+1}\frac{1}{k}=\sum_{k=1}^{n}\frac{1}{k}=H_n=H_n^{(1)}.$$

We now employ integration by parts on the slightly modified integral

(4.4.139) $$J_n=\int_a^b\frac{1-(1-t)^n}{t}dt=\left\{1-(1-t)^n\right\}\log t\Big|_a^b-n\int_a^b(1-t)^{n-1}\log t\,dt$$

We have

(4.4.140) $$I_{n-1}=\int_a^b(1-t)^{n-1}\log t\,dt=-\frac{(1-t)^n\log t}{n}\Big|_a^b+\frac{1}{n}\int_a^b\frac{(1-t)^n}{t}dt$$

(4.4.141) $$\int_a^b\frac{(1-t)^n}{t}dt=\int_a^b\sum_{k=0}^{n}\binom{n}{k}(-1)^k t^{k-1}dt$$

$$=\sum_{k=1}^{n}\binom{n}{k}\frac{(-1)^k(b^k-a^k)}{k}+\log b-\log a$$

We therefore have

(4.4.142) $$J_n=\int_a^b\frac{1-(1-t)^n}{t}dt=\sum_{k=1}^{n}\binom{n}{k}\frac{(-1)^{k+1}(b^k-a^k)}{k}$$

Equation (4.4.142) gives the particular identities

(4.4.143) $$\int_0^x\frac{1-(1-t)^n}{t}dt=\sum_{k=1}^{n}\binom{n}{k}\frac{(-1)^{k+1}x^k}{k}$$

(4.4.144) $$\int_0^1\frac{1-(1-t)^n}{t}dt=\sum_{k=1}^{n}\binom{n}{k}\frac{(-1)^{k+1}}{k}$$

Using the obvious substitution $u=1-t$, $I_{n-1}$ can be written as follows

(4.4.145) $$I_{n-1}=\int_a^b(1-t)^{n-1}\log t\,dt=-\int_{1-a}^{1-b}u^{n-1}\log(1-u)\,du$$

One of the first series which any mathematician learns is the geometric series and indeed we used it early on in (2.3) in Volume I. So let's use it again:

(4.4.146)
$$\frac{1-u^n}{1-u} - \sum_{k=1}^{n} u^{k-1} = 0$$

Hence the integral of (4.4.146) is identically equal to zero.

(4.4.147)
$$\int \left\{ \frac{1-u^n}{1-u} - \sum_{k=1}^{n} u^{k-1} \right\} du = 0$$

Therefore we have

(4.4.148)
$$-\log(1-u) - \sum_{k=1}^{n} \frac{u^k}{k} = \int \frac{u^n}{1-u} \, du$$

Integration by parts gives

(4.4.149)
$$\int \frac{u^n}{1-u} \, du = -u^n \log(1-u) + n \int u^{n-1} \log(1-u) \, du$$

and combining (4.4.147) and (4.4.149) we obtain

(4.4.150)
$$n \int u^{n-1} \log(1-u) \, du = -(1-u^n) \log(1-u) - \sum_{k=1}^{n} \frac{u^k}{k}$$

Therefore, we may replace (4.4.139) by

(4.4.151)
$$J_n = \int_a^b \frac{1-(1-t)^n}{t} \, dt = \left\{ 1-(1-t)^n \right\} \log t \Big|_a^b + n \int_{1-a}^{1-b} t^{n-1} \log(1-t) \, dt$$

The integrated terms involving $\log t$ cancel out and we have

(4.4.152)
$$J_n = \int_a^b \frac{1-(1-t)^n}{t} \, dt = \sum_{k=1}^{n} \frac{(1-a)^k - (1-b)^k}{k}$$

together with the particular identities

(4.4.153)
$$\int_0^x \frac{1-(1-t)^n}{t} \, dt = \sum_{k=1}^{n} \frac{1-(1-x)^k}{k}$$

(4.4.154)
$$\int_0^1 \frac{1-(1-t)^n}{t} \, dt = \sum_{k=1}^{n} \frac{1}{k}$$

Equations (4.4.153) and (4.4.154) are of course already familiar to us from (4.1.6) and (4.1.7). Equation (4.4.150) is used later in this paper.

Alternatively, in (4.4.152) using the substitution $u = 1 - t$, we get

$$(4.4.154a) \qquad J_n = \int_{1-a}^{1-b} \frac{1 - u^n}{1 - u} \, du = \int_{1-a}^{1-b} \sum_{k=1}^{n} u^{k-1} \, du$$

$$(4.4.154b) \qquad \qquad = \sum_{k=1}^{n} \frac{(1-a)^k - (1-b)^k}{k}$$

The various identities are summarized below:

$$(4.4.155a) \qquad J_n = \int_a^b \frac{1 - (1-t)^n}{t} \, dt = \sum_{k=1}^{n} \binom{n}{k} (-1)^{k+1} \frac{b^k - a^k}{k}$$

$$= \sum_{k=1}^{n} \frac{(1-a)^k - (1-b)^k}{k}$$

$$(4.4.155b) \qquad \int_0^x \frac{1 - (1-t)^n}{t} \, dt = \sum_{k=1}^{n} \binom{n}{k} (-1)^{k+1} \frac{x^k}{k}$$

$$= \sum_{k=1}^{n} \frac{1 - (1-x)^k}{k}$$

$$(4.4.155c) \qquad \int_0^1 \frac{1 - (1-t)^n}{t} \, dt = \sum_{k=1}^{n} \binom{n}{k} (-1)^{k+1} \frac{1}{k}$$

$$= \sum_{k=1}^{n} \frac{1}{k}$$

From (4.4.7) we have

$$B(x, y) = \int_0^1 t^{x-1} (1-t)^{y-1} \, dt \qquad , (\text{Re}(x) > 0, \, \text{Re}(y) > 0)$$

which, in turn, is related to the gamma function via the identity [115, p.193]

(4.4.155d)     $$B(x, y) = \frac{\Gamma(x)\Gamma(y)}{\Gamma(x+y)}$$

Upon differentiating under the integral sign we obtain

(4.4.155e)     $$\frac{\partial}{\partial x} B(x, y) = \int_0^1 t^{x-1}(1-t)^{y-1} \log t \, dt$$

and specifically we have

$$\left. \frac{\partial}{\partial x} B(x, y) \right|_{x=1} = \int_0^1 (1-t)^{y-1} \log t \, dt$$

Differentiating (4.4.155d) we get

(4.4.155f)     $$\frac{\partial}{\partial x} B(x, y) = \frac{\Gamma(x)\Gamma(y)}{\Gamma(x+y)}[\psi(x) - \psi(x+y)]$$

and with $x = 1$ and $y = n$ we have

$$= \frac{1}{n}\big[\psi(1) - \psi(n+1)\big]$$

Reference to (4.4.139) then shows that

(4.4.155fi)     $$H_n = \int_0^1 \frac{1-(1-t)^n}{t} \, dt = -n \int_0^1 (1-t)^{n-1} \log t \, dt$$

and we therefore get

(4.4.155g)     $$H_n = \psi(n+1) - \psi(1)$$

Equation (4.4.155fi) is also reported by Devoto and Duke [53a, p.30].

By using the binomial expansion in (4.4.155fi) we see that

$$H_n = -n \int_0^1 \sum_{k=0}^n \binom{n}{k} (-1)^k \frac{t^k \log t}{1-t} \, dt$$

and from (4.4.238a) we have

$$\int_0^1 \frac{t^k \log t}{1-t}\, dt = H_k^{(2)} - \varsigma(2)$$

Therefore we obtain

$$\frac{H_n}{n} = \sum_{k=0}^n \binom{n}{k}(-1)^k [H_k^{(2)} - \varsigma(2)] = \sum_{k=0}^n \binom{n}{k}(-1)^k H_k^{(2)} = \sum_{k=1}^n \binom{n}{k}(-1)^k H_k^{(2)} - \varsigma(2)\delta_{n,0}$$

Applying (4.1.26) from Volume II(a) we get

$$\sum_{k=1}^n \frac{H_k}{k^2} = \sum_{k=1}^n \binom{n}{k}(-1)^k \frac{H_k^{(2)}}{k}$$

and accordingly we have

$$\lim_{n \to \infty} \sum_{k=1}^n \binom{n}{k}(-1)^k \frac{H_k^{(2)}}{k} = 2\varsigma(3)$$

We also see from (4.4.155fi) that for $|x| < 1$

$$\sum_{n=1}^\infty \frac{H_n}{n} x^n = -\sum_{n=1}^\infty \int_0^1 \frac{(1-t)^n x^n \log t}{1-t}\, dt$$

$$= -\int_0^1 \frac{x \log t}{1-x(1-t)}\, dt$$

We have by courtesy of the Wolfram Integrator

$$\int \frac{x \log t}{1-x(1-t)}\, dt = \log t \log\left(1 + \frac{xt}{1-x}\right) + Li_2\left(-\frac{xt}{1-x}\right)$$

and therefore we obtain

$$-\int_0^1 \frac{x \log t}{1-x(1-t)}\, dt = -Li_2\left(-\frac{x}{1-x}\right)$$

Hence we have another proof of equation (3.11b) from Volume I

$$\sum_{n=1}^\infty \frac{H_n}{n} x^n = -Li_2\left(-\frac{x}{1-x}\right)$$

This immediately gives us upon integration

$$\sum_{n=1}^{\infty} \frac{H_n}{n^2} u^n = -\int_0^u Li_2\left(-\frac{x}{1-x}\right)\frac{dx}{x}$$

The Wolfram Integrator gives us

$$\int_0^u Li_2\left(-\frac{x}{1-x}\right)\frac{dx}{x} = Li_2(u)\log u + Li_2\left(\frac{-u}{1-u}\right)\log u - \log(1-u)Li_2(1-u)$$

$$+ Li_3(1-u) - Li_3(u) - \varsigma(3)$$

and we therefore see that

$$-\sum_{n=1}^{\infty} \frac{H_n}{n^2} u^n = Li_2(u)\log u + Li_2\left(\frac{-u}{1-u}\right)\log u - \log(1-u)Li_2(1-u)$$

$$+ Li_3(1-u) - Li_3(u) - \varsigma(3)$$

which concurs with (3.111d). The Wolfram Integrator is not able to evaluate the integral $\int_0^u Li_s\left(-\frac{x}{1-x}\right)\frac{dx}{x}$ for $s > 2$.

We also have

$$\sum_{n=1}^{\infty} \frac{H_n}{n^2} u^n = -\int_0^u \frac{dx}{x} \int_0^1 \frac{x\log t}{1-x(1-t)} dt = \int_0^1 \frac{\log t \log[1-u(1-t)]}{1-t} dt$$

The Wolfram Integrator easily evaluates this as

$$\int \frac{\log t \log[1-u(1-t)]}{1-t} dt = -\frac{1}{2}\log u \log^2[1-u(1-t)] + \log u \log t \log[1-u(1-t)]$$

$$-\left(Li_2\frac{[1-u(1-t)]}{ut} - Li_2\frac{[1-u(1-t)]}{t}\right)\log\frac{[1-u(1-t)]}{t}$$

$$-\log t \log(1-t)\log[1-u(1-t)] + \frac{1}{2}\log u \log[1-u(1-t)]\left(\log[1-u(1-t)] - 2\log t\right)$$

$$-\log t\, Li_2(t) - \log[1-u(1-t)]Li_2[1-u(1-t)] + Li_3(t)$$

$$+ Li_3[1-u(1-t)] + Li_3\frac{[1-u(1-t)]}{ut} - Li_3\frac{[1-u(1-t)]}{t}$$

but there appears to be a problem at $t = 0$ that requires further investigation.

We have from (3.111d)

$$-\sum_{n=1}^{\infty}\frac{H_n}{n^2}u^n = Li_2(u)\log u + Li_2\left(\frac{-u}{1-u}\right)\log u - \log(1-u)Li_2(1-u)$$

$$+ Li_3(1-u) - Li_3(u) - \varsigma(3)$$

We may note that the substitution $x = ut$ results in

$$\int_0^1 \frac{\log t \log[1-u(1-t)]}{1-t}\,dt = \int_0^u \frac{\log\{x/u\}\log[1-u(1-\{x/u\})]}{u-x}\,dx$$

and we get the integral identity

$$\int_0^u \frac{\log\{x/u\}\log[1-u(1-\{x/u\})]}{u-x}\,dx = -\int_0^u Li_2\left(-\frac{x}{1-x}\right)\frac{dx}{x}$$

This implies that

$$\int_0^u \left[\frac{\log\{x/u\}\log[1-u+x)]}{u-x} + Li_2\left(-\frac{x}{1-x}\right)\frac{1}{x}\right]dx = 0$$

and the Leibniz differentiation rule

$$\frac{\partial}{\partial u}\int_{a(u)}^{b(u)} f(x,u)\,dx = \int_{a(u)}^{b(u)}\frac{\partial}{\partial u}f(x,u)\,dx + f(b(u),u)\frac{db}{du} - f(a(u),u)\frac{da}{du}$$

then produces an integral representation for $Li_2\left(-\frac{x}{1-x}\right)\frac{1}{x}$.

Using integration by parts we see that

(4.4.155gi) $$\int_0^1 \frac{1-(1-t)^y}{t} \log t \, dt = \left[1-(1-t)^y\right] \frac{\log^2 t}{2} \Bigg|_0^1 - \frac{1}{2} y \int_0^1 (1-t)^{y-1} \log^2 t \, dt$$

$$= -\frac{1}{2} y \int_0^1 (1-t)^{y-1} \log^2 t \, dt$$

From (4.4.155e) we get

$$\frac{\partial^2}{\partial x^2} B(x,y) = \int_0^1 t^{x-1}(1-t)^{y-1} \log^2 t \, dt$$

and we also have by differentiating (4.4.155f)

(4.4.155gii)

$$\frac{\Gamma(x+y)}{\Gamma(y)} \frac{\partial^2}{\partial x^2} B(x,y) = \Gamma''(x) - 2\Gamma'(x)\psi(x+y) - \Gamma(x)\psi'(x+y) + \Gamma(x)\psi^2(x+y)$$

With $x=1$ and $y=n$ we have

$$\frac{\Gamma(n+1)}{\Gamma(n)} \frac{\partial^2}{\partial x^2} B(x,y) \Bigg|_{(1,n)} = \Gamma''(1) - 2\Gamma'(1)\psi(n+1) - \Gamma(1)\psi'(n+1) + \Gamma(1)\psi^2(n+1)$$

$$= \gamma^2 + \varsigma(2) + 2\gamma\left[H_n^{(1)} - \gamma\right] - \left[\varsigma(2) - H_n^{(2)}\right] + \left[H_n^{(1)} - \gamma\right]^2$$

$$= H_n^{(2)} + \left(H_n^{(1)}\right)^2$$

where we have employed (8.57a) and (E.16d). Hence we have (see also (4.4.246))

(4.4.155h) $$n\int_0^1 (1-t)^{n-1} \log^2 t \, dt = H_n^{(2)} + \left(H_n^{(1)}\right)^2$$

The following formula is also reported by Devoto and Duke [53a, p.30]

$$\int_0^1 (1-t)^{n-1} \log^2 t \, dt = \frac{2}{n}\left[H_n^{(2)} + \sum_{k=1}^{n-1} \frac{H_k^{(1)}}{k+1}\right]$$

and the equivalence is readily seen by reference to Adamchik's formula (3.17) contained in Volume I.

The summation of (4.4.155h) gives us

$$\sum_{n=1}^{\infty} \frac{x^n}{n^p} \int_0^1 (1-t)^{n-1} \log^2 t \, dt = \sum_{n=1}^{\infty} \frac{H_n^{(2)}}{n^{p+1}} x^n + \sum_{n=1}^{\infty} \frac{\left(H_n^{(1)}\right)^2}{n^{p+1}} x^n$$

The latter may be written as

(4.4.155i) $$\int_0^1 \frac{Li_p[x(1-t)]\log^2 t}{1-t} \, dt = \sum_{n=1}^{\infty} \frac{H_n^{(2)}}{n^{p+1}} x^n + \sum_{n=1}^{\infty} \frac{\left(H_n^{(1)}\right)^2}{n^{p+1}} x^n$$

We also have

$$\sum_{n=1}^{\infty} x^n \int_0^1 (1-t)^{n-1} \log^2 t \, dt = \sum_{n=1}^{\infty} \frac{H_n^{(2)}}{n} x^n + \sum_{n=1}^{\infty} \frac{\left(H_n^{(1)}\right)^2}{n} x^n$$

and hence using the geometric series this becomes

(4.4.155j) $$x\int_0^1 \frac{\log^2 t}{1-x(1-t)} \, dt = \sum_{n=1}^{\infty} \frac{H_n^{(2)}}{n} x^n + \sum_{n=1}^{\infty} \frac{\left(H_n^{(1)}\right)^2}{n} x^n$$

The Wolfram Integrator provides us with

$$x\int \frac{\log^2 t}{1-x(1-t)} \, dt = \log\left(1+\frac{xt}{1-x}\right)\log^2 t + 2Li_2\left(\frac{-xt}{1-x}\right)\log t - 2Li_3\left(\frac{-xt}{1-x}\right)$$

and hence the definite integral becomes

(4.4.155k) $$x\int_0^1 \frac{\log^2 t}{1-x(1-t)} \, dt = -2Li_3\left(\frac{-x}{1-x}\right)$$

We therefore have for $|x| < 1$

(4.4.155l) $$Li_3\left(\frac{-x}{1-x}\right) = -\frac{1}{2}\left(\sum_{n=1}^{\infty} \frac{H_n^{(2)}}{n} x^n + \sum_{n=1}^{\infty} \frac{\left(H_n^{(1)}\right)^2}{n} x^n\right)$$

This is similar to the formula (3.11b) in Volume I

$$Li_2\left(\frac{-x}{1-x}\right) = -\sum_{n=1}^{\infty}\frac{H_n^{(1)}}{n}x^n$$

$$Li_2\left(y\right) = -\sum_{n=1}^{\infty}\frac{H_n^{(1)}}{n}\left(\frac{y}{y-1}\right)^n$$

$$Li_2\left(1-y\right) = -\sum_{n=1}^{\infty}\frac{H_n^{(1)}}{n}\left(\frac{y-1}{y}\right)^n = -\sum_{n=1}^{\infty}\frac{H_n^{(1)}}{n}\sum_{k=0}^{n}\binom{n}{k}\frac{(-1)^k}{y^k}$$

which we also derived above in an equivalent fashion.

Let us now divide (4.4.155k) by $x$ and integrate to obtain

$$\int_0^u dx \int_0^1 \frac{\log^2 t}{1-x(1-t)}dt = -2\int_0^u Li_3\left(\frac{-x}{1-x}\right)\frac{dx}{x}$$

After integrating with respect to $x$ the left-hand side becomes

$$\int_0^u dx \int_0^1 \frac{\log^2 t}{1-x(1-t)}dt = -\int_0^1 \frac{\log[1-u(1-t)]\log^2 t}{1-t}dt$$

and we therefore see that

$$\int_0^1 \frac{\log[1-u(1-t)]\log^2 t}{1-t}dt = -\left(\sum_{n=1}^{\infty}\frac{H_n^{(2)}}{n^2}u^n + \sum_{n=1}^{\infty}\frac{\left(H_n^{(1)}\right)^2}{n^2}u^n\right)$$

Differentiation results in

$$\int_0^1 \frac{\log^2 t}{1-u(1-t)}dt = -\left(\sum_{n=1}^{\infty}\frac{H_n^{(2)}}{n}u^{n-1} + \sum_{n=1}^{\infty}\frac{\left(H_n^{(1)}\right)^2}{n}u^{n-1}\right)$$

Using (4.4.155h) and (4.4.235) we see that

$$\frac{1}{n}\left[H_n^{(2)} + \left(H_n^{(1)}\right)^2\right] = \int_0^1 \sum_{k=0}^{n}\binom{n}{k}(-1)^k\frac{t^k\log^2 t}{1-t}dt$$

$$= 2\sum_{k=0}^{n}\binom{n}{k}(-1)^{k+1}[H_k^{(3)} - \varsigma(3)]$$

$$= 2\sum_{k=0}^{n}\binom{n}{k}(-1)^{k+1}H_k^{(3)} - 2\varsigma(3)\sum_{k=0}^{n}\binom{n}{k}(-1)^{k+1}$$

$$= 2\sum_{k=1}^{n}\binom{n}{k}(-1)^{k+1}H_k^{(3)}$$

Applying (4.1.26) again we get

$$\sum_{k=1}^{n}\frac{1}{k^2}\left[H_k^{(2)} + \left(H_k^{(1)}\right)^2\right] = 2\sum_{k=1}^{n}\binom{n}{k}(-1)^{k+1}\frac{H_k^{(3)}}{k}$$

Some interesting trigonometric identities may be found by substituting $x = e^{i\theta}$ in (4.4.155l) and employing the usual half-angle formulae.

With $x = 1/2$ in (4.4.155l) we get

(4.4.155li)    $$\varsigma_a(3) = \sum_{n=1}^{\infty}\frac{H_n^{(2)}}{n2^{n+1}} + \sum_{n=1}^{\infty}\frac{\left(H_n^{(1)}\right)^2}{n2^{n+1}}$$

From (4.4.155l) we note that

(4.4.155lii)    $$\int_0^x Li_3\left(\frac{-t}{1-t}\right)\frac{dt}{t} = -\frac{1}{2}\left(\sum_{n=1}^{\infty}\frac{H_n^{(2)}}{n^2}x^n + \sum_{n=1}^{\infty}\frac{\left(H_n^{(1)}\right)^2}{n^2}x^n\right)$$

We also have

$$\int_0^{-x} Li_3\left(\frac{-t}{1-t}\right)\frac{dt}{t} = -\frac{1}{2}\left(\sum_{n=1}^{\infty}(-1)^n\frac{H_n^{(2)}}{n^2}x^n + \sum_{n=1}^{\infty}(-1)^n\frac{\left(H_n^{(1)}\right)^2}{n^2}x^n\right)$$

and hence

$$\int_0^x Li_3\left(\frac{t}{1+t}\right)\frac{dt}{t} = -\frac{1}{2}\left(\sum_{n=1}^{\infty}(-1)^n\frac{H_n^{(2)}}{n^2}x^n + \sum_{n=1}^{\infty}(-1)^n\frac{\left(H_n^{(1)}\right)^2}{n^2}x^n\right)$$

Unsurprisingly, the Wolfram Integrator was not able to evaluate the above integral.

Later in (4.4.156a) we will see that

$$\sum_{n=1}^{\infty} \frac{H_n^{(2)}}{n} x^n = 2Li_3(1-x) - 2\varsigma(3) + Li_3(x) - \log(1-x)\left[Li_2(1-x) + \varsigma(2)\right]$$

and hence we obtain

(4.4.155liii)

$$-2Li_3\left(\frac{-x}{1-x}\right) = 2Li_3(1-x) - 2\varsigma(3) + Li_3(x) - \log(1-x)\left[Li_2(1-x) + \varsigma(2)\right] + \sum_{n=1}^{\infty} \frac{\left(H_n^{(1)}\right)^2}{n} x^n$$

An expression for $\displaystyle\sum_{n=1}^{\infty} \frac{\left(H_n^{(1)}\right)^2}{n2^n}$ may be easily obtained by letting $x = 1/2$.

We will also see in (4.4.247b) that

$$Li_3(x) = \sum_{n=1}^{\infty} \frac{H_n^{(2)}}{n} x^n + \log(1-x)\varsigma(2) + \sum_{n=1}^{\infty} \frac{1}{n^2} \sum_{k=1}^{n} \frac{x^k}{k}$$

and this could also be substituted in (4.4.155l). We may also recall Landen's identity (3.115) from Volume I

$$Li_3\left(\frac{-x}{1-x}\right) = \varsigma(2)\log(1-x) - \frac{1}{2}\log x \log^2(1-x) - Li_3(1-x) + \varsigma(3) + \frac{1}{6}\log^3(1-x) - Li_3(x)$$

which shows that

$$-\frac{1}{2}\left(\sum_{n=1}^{\infty} \frac{H_n^{(2)}}{n} x^n + \sum_{n=1}^{\infty} \frac{\left(H_n^{(1)}\right)^2}{n} x^n\right)$$

$$= \varsigma(2)\log(1-x) - \frac{1}{2}\log x \log^2(1-x) - Li_3(1-x) + \varsigma(3) + \frac{1}{6}\log^3(1-x) - Li_3(x)$$

From (3.46a) we see that

$$\sum_{n=1}^{\infty} \frac{\left(H_n^{(1)}\right)^2}{n} x^n = -\frac{1}{3}\log^3(1-x) + Li_3(x) - Li_2(x)\log(1-x)$$

Dividing the above by $x$ and integrating we get

$$\sum_{n=1}^{\infty} \frac{\left(H_n^{(1)}\right)^2}{n^2} x^n = -\frac{1}{3}\int_0^x \frac{\log^3(1-x)}{x} dx + Li_4(x) - \int_0^x \frac{Li_2(x)\log(1-x)}{x} dx$$

We have from (3.122) in Volume I

$$\int_0^x \frac{\log^3(1-x)}{x}\,dx =$$

$$\log x \log^3(1-x) + 3Li_2(1-x)\log^2(1-x) - 6Li_3(1-x)\log(1-x) + 6Li_4(1-x) - 6\varsigma(4)$$

$$\int_0^x \frac{Li_2(x)\log(1-x)}{x}\,dx = -\frac{1}{2}\big[Li_2(x)\big]^2$$

Accordingly, as mentioned in (3.110b), we obtain

$$\sum_{n=1}^{\infty} \frac{\left(H_n^{(1)}\right)^2}{n^2}x^n = -\frac{1}{3}\log x \log^3(1-x) - Li_2(1-x)\log^2(1-x) + 2Li_3(1-x)\log(1-x)$$

$$-2Li_4(1-x) + 2\varsigma(4) + Li_4(x) + \frac{1}{2}\big[Li_2(x)\big]^2$$

From (3.34) we see that

$$\frac{Li_2(x)}{1-x} = \sum_{n=1}^{\infty} H_n^{(2)} x^n \qquad , x \in [0,1)$$

and therefore

$$\int_0^x \frac{Li_2(x)}{x(1-x)}\,dx = \sum_{n=1}^{\infty} \frac{H_n^{(2)}}{n}x^n$$

$$= Li_3(x) - Li_2(x)\log(1-x) - \log^2(1-x)\log x - 2\log(1-x)Li_2(1-x) + 2Li_3(1-x) - 2\varsigma(3)$$

Thus a further integration results in

$$\sum_{n=1}^{\infty} \frac{H_n^{(2)}}{n^2}x^n = Li_4(x) + \frac{1}{2}\big[Li_2(x)\big]^2 - \int_0^x \frac{\log^2(1-x)\log x}{x}\,dx$$

$$-2\int_0^x \frac{\log(1-x)Li_2(1-x)}{x}\,dx - 2\int_0^x \frac{\varsigma(3) - Li_3(1-x)}{x}\,dx$$

We note from (3.320) that

$$\int_0^x \frac{\log^2(1-x)\log x}{x}\,dx =$$

$$-Li_4(x) + Li_2(x)Li_2(1-x) - \varsigma(2)Li_2(x) + \frac{1}{2}\big[Li_2(x)\big]^2 + \sum_{n=1}^\infty \frac{H_n^{(2)}}{n^2} x^n$$

$$-\log^2 x \log^2(1-x) - Li_2(x)\log x \log(1-x) + 2\varsigma(2)\log x \log(1-x) - 2\log x\, Li_3(1-x) + 2\varsigma(3)\log x$$

Also, we note from (4.4.168j) that

$$\int_0^x \frac{\varsigma(3) - Li_3(1-x)}{x}\,dx = \big[\varsigma(3) - Li_3(1-x)\big]\log x - \frac{1}{2}\big[Li_2(1-x)\big]^2 + \frac{1}{2}\varsigma^2(2)$$

We then obtain

$$2\int_0^x \frac{\log(1-x)Li_2(1-x)}{x}\,dx = Li_4(x) + \frac{1}{2}\big[Li_2(x)\big]^2 - Li_4(x) - Li_2(x)Li_2(1-x) + \varsigma(2)Li_2(x)$$

$$-\frac{1}{2}\big[Li_2(x)\big]^2 - \log^2 x \log^2(1-x) + Li_2(x)\log x \log(1-x)$$

$$-2\varsigma(2)\log x \log(1-x) + 4\log x\big[Li_3(1-x) - \varsigma(3)\big]$$

$$-2\sum_{n=1}^\infty \frac{H_n^{(2)}}{n^2} x^n + \big[Li_2(1-x)\big]^2 - \varsigma^2(2)$$

Therefore we may obtain a relationship involving $\int_0^x \frac{\log(1-x)Li_2(1-x)}{x}\,dx$ and $\int_0^x Li_3\left(\frac{-x}{1-x}\right)\frac{dx}{x}$ (undoubtedly connected with (3.119)).

Integrating (4.4.155j) we obtain

$$\int_0^x dx \int_0^1 \frac{\log^2 t}{1-x(1-t)}\,dt = \sum_{n=1}^\infty \frac{H_n^{(2)}}{n^2} x^n + \sum_{n=1}^\infty \frac{\left(H_n^{(1)}\right)^2}{n^2} x^n$$

$$\int\limits_{0}^{x} dx \int\limits_{0}^{1} \frac{\log^2 t}{1 - x(1-t)}\, dt = -\int\limits_{0}^{1} \frac{\log\left[1 - x(1-t)\right]\log^2 t}{1-t}\, dt$$

and hence we get (as we have already seen in Volume I)

(4.4.155m) $\qquad -\int\limits_{0}^{1} \frac{\log\left[1 - x(1-t)\right]\log^2 t}{1-t}\, dt = \sum\limits_{n=1}^{\infty} \frac{H_n^{(2)}}{n^2} x^n + \sum\limits_{n=1}^{\infty} \frac{\left(H_n^{(1)}\right)^2}{n^2} x^n$

The Wolfram Integrator was not able to evaluate the above integral. Letting $x = 1$ we get

$$-\int\limits_{0}^{1} \frac{\log^3 t}{1 - t}\, dt = 6\varsigma(4) = \sum\limits_{n=1}^{\infty} \frac{H_n^{(2)}}{n^2} + \sum\limits_{n=1}^{\infty} \frac{\left(H_n^{(1)}\right)^2}{n^2}$$

which concurs with (4.4.167s) and (4.4.168).

Dividing (4.4.155m) by $x$ and integrating results in

(4.4.155n) $\qquad \int\limits_{0}^{1} \frac{Li_2\left[x(1-t)\right]\log^2 t}{1-t}\, dt = \sum\limits_{n=1}^{\infty} \frac{H_n^{(2)}}{n^3} x^n + \sum\limits_{n=1}^{\infty} \frac{\left(H_n^{(1)}\right)^2}{n^3} x^n$

Further integrations give us back (4.4.155i)

(4.4.155o) $\qquad \int\limits_{0}^{1} \frac{Li_p\left[x(1-t)\right]\log^2 t}{1-t}\, dt = \sum\limits_{n=1}^{\infty} \frac{H_n^{(2)}}{n^{p+1}} x^n + \sum\limits_{n=1}^{\infty} \frac{\left(H_n^{(1)}\right)^2}{n^{p+1}} x^n$

Differentiating (4.4.155gii) we obtain

$$\frac{\Gamma(x+y)}{\Gamma(y)} \frac{\partial^3}{\partial x^3} B(x, y) + \frac{\Gamma'(x+y)}{\Gamma(y)} \frac{\partial^2}{\partial x^2} B(x, y) =$$

$$\Gamma^{(3)}(x) - 3\Gamma'(x)\psi'(x+y) - 2\Gamma''(x)\psi(x+y) - \Gamma(x)\psi''(x+y)$$

$$+ 2\Gamma(x)\psi(x+y)\psi'(x+y) + \Gamma'(x)\psi^2(x+y)$$

With $x = 1$ and $y = n$ we have

$$\frac{\Gamma(x+y)}{\Gamma(y)}\frac{\partial^3}{\partial x^3}B(x,y)\Bigg|_{(1,n)} + \frac{\Gamma'(x+y)}{\Gamma(y)}\frac{\partial^2}{\partial x^2}B(x,y)\Bigg|_{(1,n)} =$$

$$\Gamma^{(3)}(1) - 3\Gamma'(1)\psi'(n+1) - 2\Gamma''(1)\psi(n+1) - \Gamma(1)\psi''(n+1)$$

$$+2\Gamma(1)\psi(n+1)\psi'(n+1) + \Gamma'(1)\psi^2(n+1)$$

and, using (E.16e) this is equal to

$$= -\gamma^3 - 3\gamma\varsigma(2) - 2\varsigma(3) + 3\gamma\left[\varsigma(2) - H_n^{(2)}\right] - 2\left[\gamma^2 + \varsigma(2)\right]\left[H_n^{(1)} - \gamma\right] - 2\left[H_n^{(3)} - \varsigma(3)\right]$$

$$+2\left[H_n^{(1)} - \gamma\right]\left[\varsigma(2) - H_n^{(2)}\right] - \gamma\left[H_n^{(1)} - \gamma\right]^2$$

This nicely simplifies to

$$= -2H_n^{(1)}H_n^{(2)} - 2H_n^{(3)} - \gamma H_n^{(2)} - \gamma\left[H_n^{(1)}\right]^2$$

Since

$$\frac{\Gamma(n+1)}{\Gamma(n)}\frac{\partial^2}{\partial x^2}B(x,y)\Bigg|_{(1,n)} = H_n^{(2)} + \left[H_n^{(1)}\right]^2$$

we have

$$\frac{\Gamma'(x+y)}{\Gamma(y)}\frac{\partial^2}{\partial x^2}B(x,y)\Bigg|_{(1,n)} = \frac{\Gamma'(n+1)}{n\Gamma(n)}\left(H_n^{(2)} + \left[H_n^{(1)}\right]^2\right)$$

In (E.20) we have shown that

$$\frac{\Gamma'(m)}{\Gamma(m)} = -\gamma - \sum_{k=0}^{\infty}\left(\frac{1}{m+k} - \frac{1}{k+1}\right) = H_{m-1}^{(1)} - \gamma$$

and hence we have

$$\frac{\Gamma'(x+y)}{\Gamma(y)}\frac{\partial^2}{\partial x^2}B(x,y)\Bigg|_{(1,n)} = \left[H_n^{(1)} - \gamma\right]\left[H_n^{(2)} + \left[H_n^{(1)}\right]^2\right]$$

Accordingly we see that

$$\frac{\Gamma(x+y)}{\Gamma(y)}\frac{\partial^3}{\partial x^3}B(x,y)\bigg|_{(1,n)} = -2H_n^{(1)}H_n^{(2)} - 2H_n^{(3)} - \gamma H_n^{(2)} - \gamma\left[H_n^{(1)}\right]^2$$

$$-\left[H_n^{(1)} - \gamma\right]\left[H_n^{(2)} + \left[H_n^{(1)}\right]^2\right]$$

(4.4.155p)
$$= -3H_n^{(1)}H_n^{(2)} - 2H_n^{(3)} - \left[H_n^{(1)}\right]^3$$

$$= -6\left(\frac{1}{6}\left[H_n^{(1)}\right]^3 + \frac{1}{2}H_n^{(1)}H_n^{(2)} + \frac{1}{3}H_n^{(3)}\right)$$

where we have written the result in a form reminiscent of equation (3.16c) in Volume I.

$$\sum_{k=1}^n \binom{n}{k}\frac{(-1)^k}{k^3} = -6\left(\frac{1}{6}\left[H_n^{(1)}\right]^3 + \frac{1}{2}H_n^{(1)}H_n^{(2)} + \frac{1}{3}H_n^{(3)}\right)$$

Hence, we have

(4.4.155q)
$$n\int_0^1 (1-t)^{n-1}\log^3 t\,dt = -6\left(\frac{1}{6}\left[H_n^{(1)}\right]^3 + \frac{1}{2}H_n^{(1)}H_n^{(2)} + \frac{1}{3}H_n^{(3)}\right)$$

The following formula is also reported by Devoto and Duke [53a, p.30]

$$\int_0^1 (1-t)^{n-1}\log^3 t\,dt = -\frac{6}{n}\left[H_n^{(3)} + \sum_{k=1}^{n-1}\frac{H_k^{(1)}}{(k+1)^2} + \sum_{k=1}^{n-1}\frac{H_k^{(2)}}{k+1} + \sum_{k=1}^{n-1}\frac{1}{k+1}\sum_{j=1}^{k-2}\frac{H_j^{(1)}}{j+1}\right]$$

and the equivalence is readily seen by reference to Adamchik's formula (3.18) contained in Volume I.

The summation of (4.4.155q) gives us

$$\sum_{n=1}^\infty \frac{x^n}{n^p}\int_0^1 (1-t)^{n-1}\log^3 t\,dt = -6\left(\frac{1}{6}\sum_{n=1}^\infty\frac{\left[H_n^{(1)}\right]^3}{n^{p+1}}x^n + \frac{1}{2}\sum_{n=1}^\infty\frac{H_n^{(1)}H_n^{(2)}}{n^{p+1}}x^n + \frac{1}{3}\sum_{n=1}^\infty\frac{H_n^{(3)}}{n^{p+1}}x^n\right)$$

The latter may be written as
(4.4.155r)

$$\int_0^1 \frac{Li_p[x(1-t)]\log^3 t}{1-t}\,dt = -6\left(\frac{1}{6}\sum_{n=1}^\infty\frac{\left[H_n^{(1)}\right]^3}{n^{p+1}}x^n + \frac{1}{2}\sum_{n=1}^\infty\frac{H_n^{(1)}H_n^{(2)}}{n^{p+1}}x^n + \frac{1}{3}\sum_{n=1}^\infty\frac{H_n^{(3)}}{n^{p+1}}x^n\right)$$

With $x = 1$ we get

$$\int_0^1 \frac{Li_p(1-t)\log^3 t}{1-t}\,dt = \sum_{n=1}^\infty \frac{1}{n^p}\sum_{k=1}^n \binom{n}{k}\frac{(-1)^k}{k^3}$$

We also have

$$\sum_{n=1}^\infty x^n \int_0^1 (1-t)^{n-1}\log^3 t\,dt = -6\left(\frac{1}{6}\sum_{n=1}^\infty \frac{\left[H_n^{(1)}\right]^3}{n}x^n + \frac{1}{2}\sum_{n=1}^\infty \frac{H_n^{(1)}H_n^{(2)}}{n}x^n + \frac{1}{3}\sum_{n=1}^\infty \frac{H_n^{(3)}}{n}x^n\right)$$

and hence using the geometric series this becomes

(4.4.155s)

$$x\int_0^1 \frac{\log^3 t}{1-x(1-t)}\,dt = -6\left(\frac{1}{6}\sum_{n=1}^\infty \frac{\left[H_n^{(1)}\right]^3}{n}x^n + \frac{1}{2}\sum_{n=1}^\infty \frac{H_n^{(1)}H_n^{(2)}}{n}x^n + \frac{1}{3}\sum_{n=1}^\infty \frac{H_n^{(3)}}{n}x^n\right)$$

The Wolfram Integrator provides us with

$$x\int \frac{\log^3 t}{1-x(1-t)}\,dt = \log\left(1+\frac{xt}{1-x}\right)\log^3 t + 3Li_2\left(\frac{-xt}{1-x}\right)\log^2 t - 6Li_3\left(\frac{-xt}{1-x}\right)\log t + 6Li_4\left(\frac{-xt}{1-x}\right)$$

and hence the definite integral becomes

(4.4.155t) $\qquad x\displaystyle\int_0^1 \frac{\log^3 t}{1-x(1-t)}\,dt = 6Li_4\left(\frac{-x}{1-x}\right)$

We therefore have for $|x| < 1$

(4.4.155u) $\qquad -Li_4\left(\frac{-x}{1-x}\right) = \frac{1}{6}\sum_{n=1}^\infty \frac{\left[H_n^{(1)}\right]^3}{n}x^n + \frac{1}{2}\sum_{n=1}^\infty \frac{H_n^{(1)}H_n^{(2)}}{n}x^n + \frac{1}{3}\sum_{n=1}^\infty \frac{H_n^{(3)}}{n}x^n$

Letting $x = 1/2$ in (4.4.155t) gives us

$$\int_0^1 \frac{\log^3 t}{1+t}\,dt = 6Li_4(-1)$$

Integrating (4.4.155s) we obtain

$$\int_0^x dx \int_0^1 \frac{\log^3 t}{1-x(1-t)}\, dt = -6\left(\frac{1}{6}\sum_{n=1}^{\infty}\frac{\left[H_n^{(1)}\right]^3}{n^2}x^n + \frac{1}{2}\sum_{n=1}^{\infty}\frac{H_n^{(1)}H_n^{(2)}}{n^2}x^n + \frac{1}{3}\sum_{n=1}^{\infty}\frac{H_n^{(3)}}{n^2}x^n\right)$$

$$\int_0^x dx \int_0^1 \frac{\log^3 t}{1-x(1-t)}\, dt = -\int_0^1 \frac{\log\left[1-x(1-t)\right]\log^3 t}{1-t}\, dt$$

and hence we get

(4.4.155v)   $$\int_0^1 \frac{\log\left[1-x(1-t)\right]\log^3 t}{1-t}\, dt =$$

$$6\left(\frac{1}{6}\sum_{n=1}^{\infty}\frac{\left[H_n^{(1)}\right]^3}{n^2}x^n + \frac{1}{2}\sum_{n=1}^{\infty}\frac{H_n^{(1)}H_n^{(2)}}{n^2}x^n + \frac{1}{3}\sum_{n=1}^{\infty}\frac{H_n^{(3)}}{n^2}x^n\right)$$

Not surprisingly, the Wolfram Integrator was unable to evaluate the above integral. With $x = 1$ we have

(4.4.155vi)   $$\int_0^1 \frac{\log^4 t}{1-t}\, dt = 6\left(\frac{1}{6}\sum_{n=1}^{\infty}\frac{\left[H_n^{(1)}\right]^3}{n^2} + \frac{1}{2}\sum_{n=1}^{\infty}\frac{H_n^{(1)}H_n^{(2)}}{n^2} + \frac{1}{3}\sum_{n=1}^{\infty}\frac{H_n^{(3)}}{n^2}\right)$$

We have

$$\int \frac{\log^4 t}{1-t}\, dt = -24\left[\begin{array}{l}\frac{1}{24}\log(1-x)\log^4 x + \frac{1}{6}Li_2(x)\log^3 x \\[2ex] + \frac{1}{2}Li_3(x)\log^2 x + Li_4(x)\log x - Li_5(x)\end{array}\right]$$

and therefore we see that

$$\int_0^1 \frac{\log^4 t}{1-t}\, dt = 24\varsigma(5) = 6\left(\frac{1}{6}\sum_{n=1}^{\infty}\frac{\left[H_n^{(1)}\right]^3}{n^2} + \frac{1}{2}\sum_{n=1}^{\infty}\frac{H_n^{(1)}H_n^{(2)}}{n^2} + \frac{1}{3}\sum_{n=1}^{\infty}\frac{H_n^{(3)}}{n^2}\right)$$

Integrating (4.4.155v) results in

(4.4.155w)   $$\int_0^1 \frac{Li_2\left[x(1-t)\right]\log^3 t}{1-t}\, dt =$$

$$-6\left(\frac{1}{6}\sum_{n=1}^{\infty}\frac{\left[H_n^{(1)}\right]^3}{n^3}x^n+\frac{1}{2}\sum_{n=1}^{\infty}\frac{H_n^{(1)}H_n^{(2)}}{n^3}x^n+\frac{1}{3}\sum_{n=1}^{\infty}\frac{H_n^{(3)}}{n^3}x^n\right)$$

Further integrations give us (4.4.155r)

(4.4.155x)    $$\int_0^1\frac{Li_p\left[x(1-t)\right]\log^3 t}{1-t}dt=$$

$$-6\left(\frac{1}{6}\sum_{n=1}^{\infty}\frac{\left[H_n^{(1)}\right]^3}{n^{p+1}}x^n+\frac{1}{2}\sum_{n=1}^{\infty}\frac{H_n^{(1)}H_n^{(2)}}{n^{p+1}}x^n+\frac{1}{3}\sum_{n=1}^{\infty}\frac{H_n^{(3)}}{n^{p+1}}x^n\right)$$

With $x=1$ we have

(4.4.155xi)   $$\int_0^1\frac{Li_p(1-t)\log^3 t}{1-t}dt=-6\left(\frac{1}{6}\sum_{n=1}^{\infty}\frac{\left[H_n^{(1)}\right]^3}{n^{p+1}}+\frac{1}{2}\sum_{n=1}^{\infty}\frac{H_n^{(1)}H_n^{(2)}}{n^{p+1}}+\frac{1}{3}\sum_{n=1}^{\infty}\frac{H_n^{(3)}}{n^{p+1}}\right)$$

Any reader with the requisite energy may continue the above process to infinity and beyond.

Integrating (4.4.155u) gives us

$$-\int_0^x Li_4\left(\frac{-x}{1-x}\right)\frac{dx}{x}=\frac{1}{6}\sum_{n=1}^{\infty}\frac{\left[H_n^{(1)}\right]^3}{n^2}x^n+\frac{1}{2}\sum_{n=1}^{\infty}\frac{H_n^{(1)}H_n^{(2)}}{n^2}x^n+\frac{1}{3}\sum_{n=1}^{\infty}\frac{H_n^{(3)}}{n^2}x^n$$

Differentiating (4.4.155gi) we get

$$\int_0^1\frac{(1-t)^y\log t\log(1-t)}{t}dt=\frac{1}{2}y\int_0^1(1-t)^{y-1}\log^2 t\log(1-t)\,dt+\frac{1}{2}\int_0^1(1-t)^{y-1}\log^2 t\,dt$$

With $y=n$ and reference to (4.4.155h) we get

$$\int_0^1\frac{(1-t)^n\log t\log(1-t)}{t}dt=\frac{1}{2}n\int_0^1(1-t)^{n-1}\log^2 t\log(1-t)\,dt+\frac{1}{2n}\left[H_n^{(2)}+\left(H_n^{(1)}\right)^2\right]$$

Differentiating (4.4.155gii) with respect to $y$ we have

$$\frac{\Gamma(x+y)}{\Gamma(y)} \frac{\partial}{\partial y} \frac{\partial^2}{\partial x^2} B(x,y) + \frac{\Gamma(y)\Gamma'(x+y) - \Gamma(x+y)\Gamma'(y)}{\Gamma^2(y)} \frac{\partial^2}{\partial x^2} B(x,y) =$$

$$-2\Gamma'(x)\psi'(x+y) - \Gamma(x)\psi''(x+y) + 2\Gamma(x)\psi(x+y)\psi'(x+y)$$

With $y = n$ and $x = 1$ we have

$$n \frac{\partial}{\partial y} \frac{\partial^2}{\partial x^2} B(x,y) \bigg|_{(1,n)} + \frac{\Gamma(n)\Gamma'(n+1) - \Gamma(n+1)\Gamma'(n)}{\Gamma^2(n)} \frac{1}{n} \left[ H_n^{(2)} + \left( H_n^{(1)} \right)^2 \right] =$$

$$2\gamma\psi'(n+1) - \psi''(n+1) + 2\psi(n+1)\psi'(n+1)$$

which may be written as

$$n \frac{\partial}{\partial y} \frac{\partial^2}{\partial x^2} B(x,y) \bigg|_{(1,n)} + \left[ \psi(n+1) - \psi(n) \right] \left[ H_n^{(2)} + \left( H_n^{(1)} \right)^2 \right] =$$

$$2\gamma\psi'(n+1) - \psi''(n+1) + 2\psi(n+1)\psi'(n+1)$$

We therefore obtain using (8.57a)

$$n \frac{\partial}{\partial y} \frac{\partial^2}{\partial x^2} B(x,y) \bigg|_{(1,n)} + \frac{1}{n} \left[ H_n^{(2)} + \left( H_n^{(1)} \right)^2 \right] =$$

$$2\gamma \left[ \varsigma(2) - H_n^{(2)} \right] + 2 \left[ \varsigma(3) - H_n^{(3)} \right] + 2 \left[ H_n^{(1)} - \gamma \right] \left[ \varsigma(2) - H_n^{(2)} \right]$$

$$= 2 \left[ \varsigma(3) - H_n^{(3)} \right] + 2 H_n^{(1)} \left[ \varsigma(2) - H_n^{(2)} \right]$$

We see that

$$\frac{\partial}{\partial y} \frac{\partial^2}{\partial x^2} B(x,y) = \int_0^1 t^{x-1} (1-t)^{y-1} \log^2 t \log(1-y) \, dt$$

and therefore we get

$$\frac{1}{2} n \int_0^1 (1-t)^{n-1} \log^2 t \log(1-t) \, dt = -\frac{1}{2n} \left[ H_n^{(2)} + \left( H_n^{(1)} \right)^2 \right] + \left[ \varsigma(3) - H_n^{(3)} \right] + H_n^{(1)} \left[ \varsigma(2) - H_n^{(2)} \right]$$

Accordingly we determine that

(4.4.155y) $\qquad \int\limits_0^1 \dfrac{(1-t)^n \log t \log(1-t)}{t} dt = \left[ \varsigma(3) - H_n^{(3)} \right] + H_n^{(1)} \left[ \varsigma(2) - H_n^{(2)} \right]$

Completing the summation we obtain

(4.4.155z)

$$\int\limits_0^1 \dfrac{Li_p(1-t) \log t \log(1-t)}{t} dt = \varsigma(3)\varsigma(p) - \sum_{n=1}^{\infty} \dfrac{H_n^{(3)}}{n^p} + \varsigma(2) \sum_{n=1}^{\infty} \dfrac{H_n^{(1)}}{n^p} - \sum_{n=1}^{\infty} \dfrac{H_n^{(1)} H_n^{(2)}}{n^p}$$

We have using integration by parts

$$\int\limits_0^1 \dfrac{\left[ 1-(1-t)^n \right] \log^{p-1} t}{t} dt = \left[ 1-(1-t)^n \right] \dfrac{\log^p t}{p} \Bigg|_0^1 - \dfrac{n}{p} \int\limits_0^1 \dfrac{(1-t)^{n-1} \log^p t}{t} dt$$

$$= -\dfrac{n}{p} \int\limits_0^1 \dfrac{(1-t)^{n-1} \log^p t}{t} dt$$

The binomial theorem gives us

$$\int\limits_0^1 \dfrac{\left[ 1-(1-t)^n \right] \log^{p-1} t}{t} dt = \int\limits_0^1 \sum_{k=1}^n \binom{n}{k} (-1)^{k-1} t^{k-1} \log^{p-1} t \, dt$$

$$= (-1)^{p+1}(p-1)! \sum_{k=1}^n \binom{n}{k} \dfrac{(-1)^{k-1}}{k^p}$$

where, in the last step we have used (3.86). We therefore have shown that

(4.4.155zi) $\qquad (-1)^{p+1} n \int\limits_0^1 (1-t)^{n-1} \log^p t \, dt = p! \sum_{k=1}^n \binom{n}{k} \dfrac{(-1)^k}{k^p}$

We note from (4.4.155fi) that

(4.4.155zi) $\quad \left[ H_n^{(1)} \right]^2 = n^2 \int\limits_0^1 (1-u)^{n-1} \log u \, du \int\limits_0^1 (1-v)^{n-1} \log v \, dv$

$$= n^2 \int\limits_0^1 \int\limits_0^1 (1-u)^{n-1} (1-v)^{n-1} \log u \log v \, du \, dv$$

and making the summation we get

$$(4.4.155\text{zii}) \quad \sum_{n=1}^{\infty} \frac{\left[H_n^{(1)}\right]^2}{n^{s+2}} x^n = \int_0^1 \int_0^1 \frac{Li_s[(1-u)(1-v)x]}{(1-u)(1-v)} \log u \log v \, du \, dv$$

The Wolfram Integrator can only evaluate the integral $\int_0^1 \frac{Li_s[(1-u)(1-v)x]}{(1-u)(1-v)} \log du$ in the simplest case where $s = 1$.

It may be noted from (4.4.155zi) that we have a double integral representation for $\left[H_n^{(1)}\right]^2$ and, with the aid of (4.4.155h), we then have integral representation for $H_n^{(2)}$. Multiplying gives us an integral representation for $H_n^{(1)} H_n^{(2)}$ and we may also determine a triple integral representation for $\left[H_n^{(1)}\right]^3$. Using (4.4.155q) then enables us to determine a triple integral representation for $H_n^{(3)}$.

We may also consider the gamma function

$$\Gamma(x) = \int_0^{\infty} t^{x-1} e^{-t} dt$$

and with $t = au$ we get

$$\Gamma(x) = a^x \int_0^{\infty} u^{x-1} e^{-au} du$$

With $x = 2$ and $a = n$ we have

$$\frac{1}{n^2} = \int_0^{\infty} u \, e^{-nu} du = -\int_0^1 (1-x)^{n-1} \log(1-x) dx$$

and we obtain the summation

$$H_n^{(2)} = -\int_0^1 \frac{1-(1-x)^n}{x} \log(1-x) dx$$

Another summation gives us

$$\sum_{n=1}^{\infty} \frac{H_n^{(2)}}{n^s} u^n = \int_0^1 \frac{Li_s[(1-x)u] - Li_s(u)}{x} \log(1-x) dx$$

We shall see an equivalent formula later in (4.4.230). The Wolfram Integrator is not able to evaluate the integral $\int_0^1 \frac{Li_2[(1-x)]}{x} \log(1-x) dx$.

The formula (4.4.155a) can be employed to derive a well-known combinatorial identity involving the reciprocal of the binomial numbers.

We note from (4.4.1) that

$$\binom{n}{k}^{-1} = \frac{\Gamma(k+1)\Gamma(n-k+1)}{\Gamma(n+1)}$$

$$= (n+1)B(k+1, n-k+1)$$

$$= (n+1)\int_0^1 x^k (1-x)^{n-k} dx$$

Therefore we get

$$\sum_{k=0}^n \binom{n}{k}^{-1} = (n+1)\int_0^1 \sum_{k=0}^n x^k (1-x)^{n-k} dx$$

$$= (n+1)\int_0^1 \frac{(1-x)^{n+1} - x^{n+1}}{1-2x} dx$$

We then make the substitution $x = (1-y)/2$ to obtain

$$\sum_{k=0}^n \binom{n}{k}^{-1} = \frac{(n+1)}{2^{n+1}} \int_{-1}^1 \frac{(1+y)^{n+1} - (1-y)^{n+1}}{2y} dy$$

$$= \frac{(n+1)}{2^{n+1}} \int_{-1}^1 \frac{(1+y)^{n+1} - 1 + 1 - (1-y)^{n+1}}{2y} dy$$

$$= \frac{(n+1)}{2^{n+1}} \int_{-1}^1 \frac{1 - (1-y)^{n+1}}{2y} dy - \frac{(n+1)}{2^{n+1}} \int_{-1}^1 \frac{1 - (1+y)^{n+1}}{2y} dy$$

$$= \frac{(n+1)}{2^{n+1}} \int\limits_{-1}^{1} \frac{1-(1-y)^{n+1}}{y} dy$$

Using (4.4.155a) we have

$$\int\limits_{-1}^{1} \frac{1-(1-t)^{n+1}}{t} dt = \sum_{k=1}^{n+1} \binom{n+1}{k} \frac{2^k}{k} = \sum_{k=1}^{n+1} \frac{2^k}{k}$$

and therefore we obtain

(4.4.155ziii)
$$\sum_{k=0}^{n} \binom{n}{k}^{-1} = \frac{(n+1)}{2^{n+1}} \sum_{k=1}^{n+1} \binom{n+1}{k} \frac{2^k}{k}$$

(4.4.155ziv)
$$= \frac{(n+1)}{2^{n+1}} \sum_{k=1}^{n+1} \frac{2^k}{k}$$

Several alternative proofs of this result are given in [32a] which contains several references to other sources for this identity. Mansour [101aa] also derived (4.4.155ziii) in 2001.

Incidentally, it may be noted that Borwein et al. [30] showed that

$$\sum_{n=1}^{\infty} \frac{1}{n^2} \sum_{k=1}^{n-1} \frac{(-1)^{k+1}}{k} = \sum_{n=1}^{\infty} \frac{1}{n2^n} \sum_{k=1}^{n-1} \frac{2^k}{k}$$

**(xiv) Theorem 4.9:**

(4.4.156a) $\quad \sum_{n=1}^{\infty} \frac{H_n^{(2)}}{n} x^n = 2Li_3(1-x) - 2\varsigma(3) + Li_3(x) - \log(1-x)\left[Li_2(1-x) + \varsigma(2)\right]$

(4.4.156b) $\quad \sum_{n=1}^{\infty} \frac{H_n^{(2)}}{n2^n} = \frac{5}{8}\varsigma(3)$

(4.4.156c) $\quad \sum_{n=1}^{\infty} \frac{H_n^{(1)}}{n^s} x^n = \int\limits_{0}^{x} \frac{Li_s(x) - Li_s(y)}{x-y} dy$

(4.4.156d) $\quad \sum_{n=1}^{\infty} \frac{H_n^{(1)}}{n^3} = \frac{1}{2}\varsigma^2(2) = \frac{\pi^4}{72} = \frac{5}{4}\varsigma(4)$

(4.4.156e)  $\displaystyle\sum_{n=1}^{\infty}\frac{\left(H_n^{(1)}\right)^2}{n^2}=\frac{17}{360}\pi^4=\frac{17}{4}\varsigma(4)$

(4.4.156f)  $\displaystyle\sum_{n=1}^{\infty}\frac{H_n^{(3)}}{n2^n}=\varsigma(2)\log^2 2-\frac{7}{8}\varsigma(3)\log 2+Li_4(1/2)-\frac{1}{6}\log^4 2$

**Proof:**

Let us now revisit one of the very first formulas quoted in Volume I, namely Euler's identity for $\varsigma(2)$ in (1.2)

(4.4.157)  $\varsigma(2)=\log x\log(1-x)+\displaystyle\sum_{n=1}^{\infty}\frac{x^n}{n^2}+\sum_{n=1}^{\infty}\frac{(1-x)^n}{n^2}$

$=\log x\log(1-x)+Li_2(x)+Li_2(1-x)$

Now (i) transfer the term involving $Li_2(1-x)$ to the left-hand side, (ii) divide by $x$ and (iii) integrate over the interval $[a,b]$. We then obtain

(4.4.158)  $\displaystyle\int_a^b\frac{\varsigma(2)-Li_2(1-x)}{x}\,dx=\int_a^b\frac{\log x\log(1-x)}{x}\,dx+\int_a^b\frac{Li_2(x)}{x}\,dx$

Using L'Hôpital's rule, we can show that all of the integrands in (4.4.158) are finite at $x=0$. The left-hand side of (4.4.158) may be written as

(4.4.159)  $\displaystyle\int_a^b\sum_{n=1}^{\infty}\frac{1-(1-x)^n}{n^2x}\,dx=\sum_{n=1}^{\infty}\frac{1}{n^2}\int_a^b\frac{1-(1-x)^n}{x}\,dx=\sum_{n=1}^{\infty}\frac{J_n}{n^2}$

$=\displaystyle\sum_{n=1}^{\infty}\frac{1}{n^2}\sum_{k=1}^{n}\frac{(1-a)^k-(1-b)^k}{k}$

where we have used (4.4.155a). Integration by parts gives

(4.4.160)  $\displaystyle\int_a^b\frac{\log x\log(1-x)}{x}\,dx=-\log x\,Li_2(x)\Big|_a^b+\int_a^b\frac{Li_2(x)}{x}\,dx$

$=Li_3(b)-Li_3(a)-\left[\log b\,Li_2(b)-\log a\,Li_2(a)\right]$

Therefore we have

(4.4.161)
$$\sum_{n=1}^{\infty}\frac{1}{n^2}\sum_{k=1}^{n}\frac{(1-a)^k-(1-b)^k}{k}=2Li_3(b)-2Li_3(a)-\left[\log b\,Li_2(b)-\log a\,Li_2(a)\right]$$

Letting $a=0$ and $b=1$ in (4.4.161) we obtain

(4.4.162)
$$\sum_{n=1}^{\infty}\frac{1}{n^2}\sum_{k=1}^{n}\frac{1}{k}=2Li_3(1)$$

because $\lim_{a\to 0}\left[\log a\,Li_2(a)\right]=0$ by L'Hôpital's rule. By definition $Li_3(1)=\varsigma(3)$ and hence we have

(4.4.163)
$$\sum_{n=1}^{\infty}\frac{1}{n^2}\sum_{k=1}^{n}\frac{1}{k}=\sum_{n=1}^{\infty}\frac{H_n}{n^2}=2\varsigma(3)$$

and this was previously demonstrated in (4.2.33).

Alternatively, letting $a=0$ and $b=1/2$ in (4.4.161) the left-hand side becomes

(4.4.164)
$$\sum_{n=1}^{\infty}\frac{1}{n^2}\sum_{k=1}^{n}\left(\frac{1}{k}-\frac{1}{k2^k}\right)=\sum_{n=1}^{\infty}\frac{1}{n^2}\sum_{k=1}^{n}\frac{1}{k}-\sum_{n=1}^{\infty}\frac{1}{n^2}\sum_{k=1}^{n}\frac{1}{k2^k}$$

Reversing the order of summation this becomes

$$=2\varsigma(3)-\sum_{n=1}^{\infty}\frac{1}{n2^n}\sum_{k=n}^{\infty}\frac{1}{k^2}$$

$$=2\varsigma(3)-\sum_{n=1}^{\infty}\frac{1}{n2^n}\left(\varsigma(2)-H_{n-1}^{(2)}\right)$$

$$=2\varsigma(3)-\log 2\varsigma(2)+\sum_{n=1}^{\infty}\frac{1}{n2^n}H_{n-1}^{(2)}$$

$$=2\varsigma(3)-\log 2\varsigma(2)+\sum_{n=1}^{\infty}\frac{1}{n2^n}\left(H_n^{(2)}-\frac{1}{n^2}\right)$$

(4.4.165)
$$=2\varsigma(3)-\log 2\varsigma(2)+\sum_{n=1}^{\infty}\frac{H_n^{(2)}}{n2^n}-Li_3(1/2)$$

The right hand side of (4.4.161) for $a=0$ and $b=1/2$ is

(4.4.166)                    $= 2Li_3(1/2) + \log 2 Li_2(1/2)$

Therefore, using the Euler/Landen identities (3.43a) and (3.43b), we obtain

(4.4.167)              $\displaystyle\sum_{n=1}^{\infty} \frac{H_n^{(2)}}{n 2^n} = \frac{5}{8}\varsigma(3)$

The identity (4.4.167) can also be obtained by using the identity (3.34), where dividing by $x$ and integrating we have

(4.4.167a)            $\displaystyle\int_0^{\frac{1}{2}} \frac{Li_2(x)}{x(1-x)}\,dx = \int_0^{\frac{1}{2}} \sum_{n=1}^{\infty} H_n^{(2)} x^{n-1}\,dx$

The integral on the left hand side has previously been evaluated in (3.45), and we end up with the result (4.4.167).

We may generalise the above result by letting $a = 0$ and $b = 1 - x$ in (4.4.161): the left-hand side becomes

$$\sum_{n=1}^{\infty} \frac{1}{n^2} \sum_{k=1}^{n} \left( \frac{1}{k} - \frac{x^k}{k} \right) = \sum_{n=1}^{\infty} \frac{1}{n^2} \sum_{k=1}^{n} \frac{1}{k} - \sum_{n=1}^{\infty} \frac{1}{n^2} \sum_{k=1}^{n} \frac{x^k}{k}$$

Reversing the order of summation this becomes

$$= 2\varsigma(3) - \sum_{n=1}^{\infty} \frac{x^n}{n} \sum_{k=n}^{\infty} \frac{1}{k^2}$$

$$= 2\varsigma(3) - \sum_{n=1}^{\infty} \frac{x^n}{n} \left( \varsigma(2) - H_{n-1}^{(2)} \right)$$

Provided $-1 \le x < 1$ we may write this as

$$= 2\varsigma(3) + \log(1-x)\varsigma(2) + \sum_{n=1}^{\infty} \frac{x^n}{n} H_{n-1}^{(2)}$$

$$= 2\varsigma(3) + \log(1-x)\varsigma(2) + \sum_{n=1}^{\infty} \frac{x^n}{n} \left( H_n^{(2)} - \frac{1}{n^2} \right)$$

$$= 2\varsigma(3) + \log(1-x)\varsigma(2) + \sum_{n=1}^{\infty} \frac{H_n^{(2)}}{n} x^n - Li_3(x)$$

The right-hand side of (4.4.161) with $a = 0$ and $b = 1 - x$ is

$$= 2Li_3(1-x) - \log(1-x)Li_2(1-x)$$

Therefore we obtain

(4.4.167b) $\quad \sum_{n=1}^{\infty} \frac{H_n^{(2)}}{n} x^n = 2Li_3(1-x) - 2\varsigma(3) + Li_3(x) - \log(1-x)\big[Li_2(1-x) + \varsigma(2)\big]$

The identity (4.4.167b) can also be obtained by using the identity (3.34) where, dividing by $x$ and integrating, we obtain

(4.4.167c) $\quad \int_0^x \frac{Li_2(t)}{t(1-t)} dt = \int_0^x \sum_{n=1}^{\infty} H_n^{(2)} t^{n-1} dt = \sum_{n=1}^{\infty} \frac{H_n^{(2)}}{n} x^n$

$$\int_0^x \frac{Li_2(t)}{t(1-t)} dt = \int_0^x \frac{Li_2(t)}{t} dt + \int_0^x \frac{Li_2(t)}{(1-t)} dt$$

The second integral on the right-hand side was previously evaluated in (3.45), and we end up with

$$\int_0^x \frac{Li_2(t)}{t(1-t)} dt = 2Li_3(1-x) - 2\varsigma(3) + Li_3(x) - \log(1-x)\big[2Li_2(1-x) + Li_2(x) + \log(1-x)\log x\big]$$

and, by equating these results, we simply end up with a restatement of Euler's identity for the dilogarithm.

With reference to (4.4.158) we also have

$$\int_0^t \frac{\varsigma(2) - Li_2(1-x)}{x} dx = \big[\varsigma(2) - Li_2(1-t)\big]\log t + \int_0^t \frac{\log^2 x}{1-x} dx$$

A further integration by parts gives us

$$\int_0^t \frac{\log^2 x}{1-x} dx = -\log^2 t \log(1-t) + 2\int_0^t \frac{\log x \log(1-x)}{x} dx$$

and therefore we get

$$\int_0^t \frac{\varsigma(2) - Li_2(1-x)}{x} dx = \big[\varsigma(2) - Li_2(1-t)\big]\log t - \log^2 t \log(1-t) + 2\int_0^t \frac{\log x \log(1-x)}{x} dx$$

Comparing this with (4.4.158) we have

$$\int_0^t \frac{\log x \log(1-x)}{x} dx = -\left[\varsigma(2) - Li_2(1-t)\right]\log t + \log^2 t \log(1-t) + Li_3(t)$$

Reference to (3.40) shows that

$$\int \frac{\log(1-x)\log x}{1-x} dx = \log(1-x) Li_2(1-x) - Li_3(1-x)$$

and we have

$$\int_0^t \frac{\log x \log(1-x)}{x} dx = \int_{1-t}^1 \frac{\log(1-u)\log u}{1-u} du = \int_0^1 \frac{\log(1-u)\log u}{1-u} du - \int_0^{1-t} \frac{\log(1-u)\log u}{1-u} du$$

$$= -\log t\, Li_2(t) + Li_3(t)$$

Therefore

$$-\left[\varsigma(2) - Li_2(1-t)\right]\log t + \log^2 t \log(1-t) + Li_3(t) = -\log t\, Li_2(t) + Li_3(t)$$

and hence again we simply regain Euler's identity (1.2).

Now let us divide equation (4.4.167b) by $x$ and integrate to obtain

$$\sum_{n=1}^\infty \frac{H_n^{(2)}}{n^2} t^n = \int_0^t \left(2Li_3(1-x) - 2\varsigma(3) + Li_3(x) - \log(1-x)\left[Li_2(1-x) + \varsigma(2)\right]\right)/x\, dx$$

$$= -2\int_0^t \frac{\varsigma(3) - Li_3(1-x)}{x} dx + \int_0^t \frac{Li_3(x)}{x} dx - \int_0^t \frac{\log(1-x) Li_2(1-x)}{x} dx - \varsigma(2)\int_0^t \frac{\log(1-x)}{x} dx$$

(4.4.167d)

$$= -2\int_0^t \frac{\varsigma(3) - Li_3(1-x)}{x} dx + Li_4(t) - \int_0^t \frac{\log(1-x) Li_2(1-x)}{x} dx + \varsigma(2) Li_2(t)$$

because $\int_0^t \frac{\log(1-x)}{x} dx = -Li_2(x)\big|_0^t = -Li_2(t)$

The first integral in (4.4.167d) may be written as

$$(4.4.167e) \qquad \int_0^t \frac{\varsigma(3) - Li_3(1-x)}{x} dx = \int_0^t \sum_{n=1}^\infty \frac{1-(1-x)^n}{n^3 x} dx$$

$$= \sum_{n=1}^\infty \frac{1}{n^3} \int_0^t \frac{1-(1-x)^n}{x} dx$$

$$= \sum_{n=1}^\infty \frac{1}{n^3} \sum_{k=1}^n \frac{1-(1-t)^k}{k}$$

where we have used (4.4.154b). The above is equivalent to

$$\sum_{n=1}^\infty \frac{1}{n^3} \sum_{k=1}^n \left( \frac{1}{k} - \frac{(1-t)^k}{k} \right) = \sum_{n=1}^\infty \frac{1}{n^3} \sum_{k=1}^n \frac{1}{k} - \sum_{n=1}^\infty \frac{1}{n^3} \sum_{k=1}^n \frac{(1-t)^k}{k}$$

Using (3.23) to reverse the order of summation we have

$$\sum_{n=1}^\infty \frac{1}{n^3} \sum_{k=1}^n \frac{(1-t)^k}{k} = \sum_{n=1}^\infty \frac{(1-t)^n}{n} \sum_{k=n}^\infty \frac{1}{k^3}$$

$$= \sum_{n=1}^\infty \frac{(1-t)^n}{n} \left( \varsigma(3) - H_{n-1}^{(3)} \right)$$

Provided $0 < t \le 2$ we may write this as

$$= -\varsigma(3) \log t - \sum_{n=1}^\infty \frac{(1-t)^n}{n} H_{n-1}^{(3)}$$

$$= -\varsigma(3) \log t - \sum_{n=1}^\infty \frac{(1-t)^n}{n} \left( H_n^{(3)} - \frac{1}{n^3} \right)$$

Hence we obtain for $0 < t \le 2$

$$(4.4.167f) \qquad \sum_{n=1}^\infty \frac{1}{n^3} \sum_{k=1}^n \frac{(1-t)^k}{k} = -\varsigma(3) \log t - \sum_{n=1}^\infty \frac{H_n^{(3)}}{n} (1-t)^n + Li_4(1-t)$$

With $1 - t = x$ we have

$$(4.4.167fi) \qquad \sum_{n=1}^\infty \frac{1}{n^3} \sum_{k=1}^n \frac{x^k}{k} = -\varsigma(3) \log(1-x) - \sum_{n=1}^\infty \frac{H_n^{(3)}}{n} x^n + Li_4(x)$$

Therefore we have for $0 < t \leq 1$

(4.4.167g)  $\displaystyle\int\limits_0^t \frac{\varsigma(3) - Li_3(1-x)}{x}\, dx = \sum_{n=1}^{\infty} \frac{H_n^{(1)}}{n^3} + \varsigma(3)\log t + \sum_{n=1}^{\infty} \frac{H_n^{(3)}}{n}(1-t)^n - Li_4(1-t)$

(this integral is considered further in (4.4.168j)). With $t = 1$ we have

(4.4.167ga)  $\displaystyle\int\limits_0^1 \frac{\varsigma(3) - Li_3(1-x)}{x}\, dx = \sum_{n=1}^{\infty} \frac{H_n^{(1)}}{n^3}$

and this is a particular case of the general identity (4.4.167j) considered below.

In his book, "The Art of Computer Programming", D.E. Knuth [90b, p.78] states that if $f(x) = \sum_{n=1}^{\infty} a_n x^n$ is convergent, then

(4.4.167h)  $\displaystyle\sum_{n=1}^{\infty} a_n H_n^{(1)} x^n = \int\limits_0^1 \frac{f(x) - f(tx)}{1-t}\, dt$

With $f(x) = Li_s(x) = \sum_{n=1}^{\infty} \dfrac{x^n}{n^s}$ we therefore obtain

$\displaystyle\sum_{n=1}^{\infty} \frac{H_n^{(1)}}{n^s} x^n = \int\limits_0^1 \frac{Li_s(x) - Li_s(tx)}{1-t}\, dt$

Letting $y = tx$ in the above integral we obtain the symmetrical form

(4.4.167i)  $\displaystyle\sum_{n=1}^{\infty} \frac{H_n^{(1)}}{n^s} x^n = \int\limits_0^x \frac{Li_s(x) - Li_s(y)}{x-y}\, dy$

and with $x = \pm 1$ we have

(4.4.167j)  $\displaystyle\sum_{n=1}^{\infty} \frac{H_n^{(1)}}{n^s} = \int\limits_0^1 \frac{\varsigma(s) - Li_s(y)}{1-y}\, dy$

(4.4.167k)  $\displaystyle\sum_{n=1}^{\infty} (-1)^n \frac{H_n^{(1)}}{n^s} = \int\limits_0^1 \frac{(2^{1-s}-1)\varsigma(s) - Li_s(-y)}{1-y}\, dy$

We also see that

$$\sum_{n=1}^{\infty} \frac{H_n^{(1)}}{n^{s+1}} u^n = \int_0^u \int_0^x \frac{Li_s(x) - Li_s(y)}{x(x-y)} dx dy$$

Integrating (4.4.167fi) we obtain

$$\sum_{n=1}^{\infty} \frac{1}{n^3} \sum_{k=1}^{n} \frac{x^k}{k^2} = \varsigma(3) Li_2(x) - \sum_{n=1}^{\infty} \frac{H_n^{(3)}}{n^2} x^n + Li_5(x)$$

and letting $x = 1$ we have

$$\sum_{n=1}^{\infty} \frac{H_n^{(2)}}{n^3} + \sum_{n=1}^{\infty} \frac{H_n^{(3)}}{n^2} = \varsigma(2)\varsigma(3) + \varsigma(5)$$

in agreement with (4.4.232a).

Further integrations result in

$$\sum_{n=1}^{\infty} \frac{1}{n^3} \sum_{k=1}^{n} \frac{x^k}{k^{p+2}} = \varsigma(3) Li_{p+2}(x) - \sum_{n=1}^{\infty} \frac{H_n^{(3)}}{n^{p+2}} x^n + Li_{p+5}(x)$$

In 1991 de Doelder [55] showed that (as corrected by Coffey [45d])

$$\sum_{n=1}^{\infty} (-1)^n \frac{H_n^{(1)}}{n^3} = \varsigma(4) - \frac{1}{2} \int_0^1 \frac{\log x \log^2(1+x)}{x} dx$$

and

$$\sum_{n=1}^{\infty} (-1)^n \frac{H_n^{(1)}}{n^3} = -\frac{11}{4}\varsigma(4) + \frac{7}{4}\varsigma(3)\log 2 - \frac{1}{12}\pi^2 \log^2 2 + \frac{1}{12}\log^4 2 + 2Li_4(1/2)$$

With $s = 3$ in (4.4.167k) we obtain

$$\sum_{n=1}^{\infty} (-1)^n \frac{H_n^{(1)}}{n^3} = \int_0^1 \frac{\varsigma_a(3) - Li_3(-y)}{1-y} dy$$

We get using integration by parts

$$\int_0^1 \frac{\varsigma_a(3) - Li_3(-y)}{1-y} dy = -\left[\varsigma_a(3) - Li_3(-y)\right] \log(1-y) \Big|_0^1 - \int_0^1 \frac{\log(1-y) Li_2(-y)}{y} dy$$

$$= -\int_0^1 \frac{\log(1-y)Li_2(-y)}{y}\,dy$$

and

$$\int_0^1 \frac{\log(1-y)Li_2(-y)}{y}\,dy = -Li_2(-y)Li_2(y)\Big|_0^1 + \int_0^1 \frac{\log(1+y)Li_2(y)}{y}\,dy$$

$$= -\varsigma_a(3)\varsigma(3) + \int_0^1 \frac{\log(1+y)Li_2(y)}{y}\,dy$$

Hence we obtain

$$\int_0^1 \frac{\log(1+y)Li_2(y)}{y}\,dy = \varsigma_a(3)\varsigma(3) - \sum_{n=1}^{\infty}(-1)^n \frac{H_n^{(1)}}{n^3}$$

As shown by Sitaramachandrarao [120a] in 1987, and also noted in [69], we have

$$\sum_{n=0}^{\infty}(-1)^n \frac{H_n^{(1)}}{(2n+1)^3} = 3\sum_{n=0}^{\infty}\frac{(-1)^n}{(2n+1)^4} - \frac{7}{16}\pi\varsigma(3) - \frac{1}{16}\pi^3\log 2$$

Equation (4.4.167g) was valid for $0 < t \leq 1$ and as $t \to 0$ this suggests that

$$\lim_{t \to 0}\left[\varsigma(3)\log t + \sum_{n=1}^{\infty}\frac{H_n^{(3)}}{n}(1-t)^n\right] = Li_4(1) - \sum_{n=1}^{\infty}\frac{H_n^{(1)}}{n^3}$$

and using (4.4.167u) this becomes

$$\lim_{t \to 0}\left[\varsigma(3)\log t + \sum_{n=1}^{\infty}\frac{H_n^{(3)}}{n}(1-t)^n\right] = -\frac{1}{4}\varsigma(4)$$

Referring to (3.33) we see that

$$\sum_{n=1}^{\infty}H_n^{(3)}(1-t)^n = \frac{Li_3(1-t)}{t} \qquad , t \in (0,1]$$

We see from (3.26) that

$$f(x) = \frac{1}{2}\log^2(1-x) + Li_2(x) = \sum_{n=1}^{\infty} \frac{H_n^{(1)}}{n} x^n$$

and using Knuth's formula (4.4.167h) we obtain

$$2\sum_{n=1}^{\infty} \frac{\left[H_n^{(1)}\right]^2}{n} x^n = \int_0^1 \frac{\log^2(1-x) + 2Li_2(x) - \log^2(1-tx) - 2Li_2(tx)}{1-t} dt$$

The above integral may be evaluated with the Wolfram Integrator in terms of polylogarithms. One could also employ, for example, (3.33).

After that minor diversion, let us now return to the main issue.

Using Euler's identity (1.6c) for the dilogarithm, we have for the second integral in (4.4.167d)

(4.4.167l)

$$\int_0^t \frac{\log(1-x)Li_2(1-x)}{x} dx = \int_0^t \frac{\log(1-x)\left[\varsigma(2) - Li_2(x) - \log x \log(1-x)\right]}{x} dx$$

$$= -\varsigma(2)Li_2(t) - \int_0^t \frac{\log(1-x)Li_2(x)}{x} dx - \int_0^t \frac{\log x \log^2(1-x)}{x} dx$$

Integration by parts gives us

$$\int_0^t \frac{\log(1-x)Li_2(x)}{x} dx = -\left[Li_2(x)\right]^2 \Big|_0^t - \int_0^t \frac{Li_2(x)\log(1-x)}{x} dx$$

and therefore we have (which we have shown previously in (3.110e))

(4.4.167m)    $$\int_0^t \frac{\log(1-x)Li_2(x)}{x} dx = -\frac{1}{2}\left[Li_2(t)\right]^2$$

This then gives us

(4.4.167ma)    $$\int_0^t \frac{\log(1-x)Li_2(1-x)}{x} dx + \int_0^t \frac{\log x \log^2(1-x)}{x} dx = -\varsigma(2)Li_2(t) + \frac{1}{2}\left[Li_2(t)\right]^2$$

From (3.42) we have

(4.4.167n)    $\int \dfrac{\log^2(1-x)dx}{x} = \log x \log^2(1-x) + 2\log(1-x)Li_2(1-x) - 2Li_3(1-x)$

Using integration by parts we obtain with the use of (4.4.167n) (where, for the moment, we assume that $a \neq 0$ due to convergence issues)

$$\int\limits_a^t \dfrac{\log x \log^2(1-x)}{x}\, dx = \log x \Big[ \log x \log^2(1-x) + 2\log(1-x)Li_2(1-x) - 2Li_3(1-x) \Big]\Big|_a^t$$

$$-\int\limits_a^t \dfrac{\log x \log^2(1-x) + 2\log(1-x)Li_2(1-x) - 2Li_3(1-x)}{x}\, dx$$

Therefore we obtain

(4.4.167o)

$$2\int\limits_a^t \dfrac{\log x \log^2(1-x)}{x}\, dx = \log x \Big[ \log^2(1-x)\log x + 2\log(1-x)Li_2(1-x) - 2Li_3(1-x) \Big]\Big|_a^t$$

$$-2\int\limits_a^t \dfrac{\log(1-x)Li_2(1-x)}{x}\, dx + 2\int\limits_a^t \dfrac{Li_3(1-x)}{x}\, dx$$

The convergence issue may be remedied by writing (4.4.167o) as

$$2\int\limits_0^t \dfrac{\log x \log^2(1-x)}{x}\, dx = \log x \Big[ \log^2(1-x)\log x + 2\log(1-x)Li_2(1-x) - 2Li_3(1-x) + 2\varsigma(3) \Big]\Big|_a^t$$

$$-2\int\limits_0^t \dfrac{\log(1-x)Li_2(1-x)}{x}\, dx + 2\int\limits_0^t \dfrac{Li_3(1-x) - \varsigma(3)}{x}\, dx$$

With integration by parts we obtain

$$\int\limits_0^t \dfrac{\log(1-x)Li_2(1-x)}{x}\, dx = -Li_2(t)Li_2(1-t) + \int\limits_0^t \dfrac{Li_2(x)\log x}{1-x}\, dx$$

since $\dfrac{d}{dx} Li_2(1-x) = \dfrac{\log x}{1-x}$

I tried long and hard to complete the evaluation of the integral

$$F(t) = \int_0^t \frac{\log x \log^2(1-x)}{x}\, dx$$

but failed miserably and, with some reluctance, I resorted to the Wolfram Integrator: when I saw the result I fully comprehended my inability to derive a human proof! The answer is:

(4.4.167p)

$$\int_0^x \frac{\log t \log^2(1-t)}{t}\, dt = \frac{1}{2}\log^2(1-x)\log^2 x + \frac{1}{12}\log^4 x - \log^2(1-x)\log^2 x + \frac{2}{3}\log(1-x)\log^3 x$$

$$-\left[\log(1-x) + \frac{1}{3}\log x\right]\log^2 x \log\left(\frac{x}{1-x}\right) + \frac{1}{2}\log^2 x\left[\log\left(\frac{x}{1-x}\right)\right]^2$$

$$-\frac{1}{4}\left[\log\left(\frac{x}{1-x}\right)\right]^4 + \log^2(1-x)Li_2(1-x) - \log^2 x\, Li_2(x)$$

$$-\left[\log\left(\frac{x}{1-x}\right)\right]^2 Li_2\left(\frac{-x}{1-x}\right) - 2\log(1-x)Li_3(1-x) + 2\log x\, Li_3(x)$$

$$+2\log\left(\frac{x}{1-x}\right)Li_3\left(\frac{-x}{1-x}\right) + 2\left[Li_4(1-x) - Li_4(x) - Li_4\left(\frac{-x}{1-x}\right)\right] - 2\zeta(4)$$

The machine generated proof may of course be easily verified by differentiating $F(x)$ where $F(x)$ is the expression given by (4.4.167p).

The presence of so many terms involving $\log x$ initially made me think that the integral was not convergent at $x = 0$ but closer inspection reveals the net sum of the terms only involving $\log^4 x$ is zero (the factor of $\log(1-x)$ elsewhere ensures convergence for the other powers involving $\log x$). However, it is not immediately apparent to me how the integrated part behaves as $x \to 1$, having regard to the polylogarithmic terms involving $(-x/1-x)$: perhaps the Wolfram Integrator output is only valid for $x < 1$ because the Maclaurin series expansion for $\log(1-x)$ is not convergent at $x = 1$.

In fact, the following definite integral is well known; see, for example, [69a] and also the proof given in (3.226).

(4.4.167q) $\qquad \displaystyle\int_0^1 \frac{\log x \log^2(1-x)}{x}\, dx = \int_0^1 \frac{\log(1-x)\log^2 x}{1-x}\, dx = -\frac{1}{2}\zeta(4)$

(a further proof is given in (4.4.168f)).

In addition, in 1934 Rutledge and Douglass [116aa] also showed that

(4.4.167qa) $$\int_0^1 \frac{\log x \log^2(1-x)}{x} \, dx = -2\sum_{n=1}^{\infty} \frac{H_n^{(1)}}{(n+1)^3} = -\frac{1}{2}\varsigma(4)$$

(a similar integral is given in (4.4.245a)). Accordingly, reference to (4.4.167l), (4.4.167m) and (4.4.167q) gives us

(4.4.167r) $$\int_0^1 \frac{\log(1-x)Li_2(1-x)}{x} \, dx = \int_0^1 \frac{\log x \, Li_2(x)}{1-x} \, dx$$

$$= -\varsigma^2(2) + \frac{1}{2}\varsigma^2(2) + \frac{1}{2}\varsigma(4)$$

$$= -\frac{3}{4}\varsigma(4) = -\frac{1}{120}\pi^4$$

Similarly, reference to (4.4.167d) gives us

$$\sum_{n=1}^{\infty} \frac{H_n^{(2)}}{n^2} = -2\sum_{n=1}^{\infty} \frac{H_n^{(1)}}{n^3} + \frac{7}{4}\varsigma(4) + \varsigma^2(2)$$

Using (4.2.43) we have

(4.4.167s) $$\sum_{n=1}^{\infty} \frac{H_n^{(2)}}{n^2} = \frac{7}{4}\varsigma(4) = \frac{7\pi^4}{360}$$

which is a particular case of the general formula which is derived in (4.4.232a)

(4.4.167t) $$\sum_{n=1}^{\infty} \frac{H_n^{(r)}}{n^r} = \frac{1}{2}\left[\varsigma^2(r) + \varsigma(2r)\right]$$

Therefore we obtain

(4.4.167u) $$\sum_{n=1}^{\infty} \frac{H_n^{(1)}}{n^3} = \frac{1}{2}\varsigma^2(2) = \frac{\pi^4}{72} = \frac{5}{4}\varsigma(4)$$

which was derived by de Doelder [55] in 1991 and subsequently by Borwein and Borwein [27] in 1995. Employing this in (4.4.167g) we obtain

(4.4.167v) $\qquad \int_0^1 \dfrac{\varsigma(3) - Li_3(1-x)}{x} dx = \int_0^1 \dfrac{\varsigma(3) - Li_3(x)}{1-x} dx = \sum_{n=1}^{\infty} \dfrac{H_n^{(1)}}{n^3} = \dfrac{\pi^4}{72} = \dfrac{5}{4}\varsigma(4)$

From (3.108b) we have

(4.4.167w) $\qquad 6\sum_{n=1}^{\infty} \dfrac{H_n^{(1)}}{n^3} t^n + 3\sum_{n=1}^{\infty} \dfrac{H_n^{(2)}}{n^2} t^n - 3\sum_{n=1}^{\infty} \dfrac{\left(H_n^{(1)}\right)^2}{n^2} t^n =$

$$\log^3(1-t)\log t + 3\log^2(1-t)Li_2(1-t) - 6\log(1-t)Li_3(1-t) + 6Li_4(1-t) + 6Li_4(t) - 6\varsigma(4)$$

and therefore with $t = 1$ we have

$$2\sum_{n=1}^{\infty} \dfrac{H_n^{(1)}}{n^3} + \sum_{n=1}^{\infty} \dfrac{H_n^{(2)}}{n^2} - \sum_{n=1}^{\infty} \dfrac{\left(H_n^{(1)}\right)^2}{n^2} = 0$$

Hence we obtain

(4.4.168) $\qquad \sum_{n=1}^{\infty} \dfrac{\left(H_n^{(1)}\right)^2}{n^2} = \dfrac{17}{360}\pi^4 = \dfrac{17}{4}\varsigma(4)$

This result was first obtained by de Doelder [55] in 1991 and subsequently by Borwein and Borwein [27] in 1995, by Flajolet and Salvy [69] in 1996, W. Chu in 1997, Basu and Apostol in 2000 and most recently by Freitas [69a] in 2004 (see [69a] for further references). The method employed by Freitas is very elegant in its simplicity.

We recall (3.46a) from Volume I

$$\sum_{n=1}^{\infty} \dfrac{\left(H_n^{(1)}\right)^2}{n} x^n = -\dfrac{1}{3}\log^3(1-x) + Li_3(x) - Li_2(x)\log(1-x)$$

and we now divide this by $x$ and integrate to obtain (3.110ea)

$$\sum_{n=1}^{\infty} \dfrac{\left(H_n^{(1)}\right)^2}{n^2} x^n = Li_4(x) + \dfrac{1}{2}\left[Li_2(x)\right]^2$$

$$-\dfrac{1}{3}\left[\log^3(1-x)\log x + 3\log^2(1-x)Li_2(1-x) - 6\log(1-x)Li_3(1-x) + 6Li_4(1-x) - 6\varsigma(4)\right]$$

and with $x = 1$ we obtain (4.4.168).

From (4.4.167d) we have

(4.4.168a)

$$\sum_{n=1}^{\infty} \frac{H_n^{(2)}}{n^2} t^n = -2\int_0^t \frac{\varsigma(3) - Li_3(1-x)}{x} dx + Li_4(t) - \int_0^t \frac{\log(1-x)Li_2(1-x)}{x} dx + \varsigma(2)Li_2(t)$$

Equation (4.4.167ma) gives us

$$\int_0^t \frac{\log(1-x)Li_2(1-x)}{x} dx = -\int_0^t \frac{\log x \log^2(1-x)}{x} dx - \varsigma(2)Li_2(t) + \frac{1}{2}\left[Li_2(t)\right]^2$$

Therefore, using (4.4.167g)

$$\int_0^t \frac{\varsigma(3) - Li_3(1-x)}{x} dx = \sum_{n=1}^{\infty} \frac{H_n^{(1)}}{n^3} + \varsigma(3)\log t + \sum_{n=1}^{\infty} \frac{H_n^{(3)}}{n}(1-t)^n - Li_4(1-t)$$

and (4.4.167o)

$$\sum_{n=1}^{\infty} \frac{H_n^{(1)}}{n^3} = \frac{1}{2}\varsigma^2(2)$$

we end up with

(4.4.168b)

$$\sum_{n=1}^{\infty} \frac{H_n^{(2)}}{n^2} t^n = -\varsigma^2(2) - 2\varsigma(3)\log t - 2\sum_{n=1}^{\infty} \frac{H_n^{(3)}}{n}(1-t)^n + 2Li_4(1-t) + Li_4(t) - \frac{1}{2}\left[Li_2(t)\right]^2 + 2\varsigma(2)Li_2(t) + F(t)$$

where $F(t)$ is the expression given by (4.4.167p). In view of my propensity to make arithmetical slips, I checked the validity of the above identity at $t = 1$ and was pleased to see that it all worked out correctly! Reference should also be made to (4.4.168m) for a different representation.

In principle, one could divide (4.4.168b) by $t$ and integrate as before to obtain an expression for $\sum_{n=1}^{\infty} \frac{H_n^{(2)}}{n^3} t^n$: however, some of the resulting integrals involving $F(t)$ appear to be quite horrendous (can they be evaluated using Mathematica?). The difficulties here may well be linked with the comments made in this paper following equation (4.2.50). Freitas [69a] records the following integral which would be relevant to the exercise

$(4.4.168c)$ 
$$\int_0^1 \frac{\left[Li_2(t)\right]^2}{t}\,dt = 2\varsigma(2)\varsigma(3) - 3\varsigma(5)$$

From (4.4.167f) we have

$$\sum_{n=1}^{\infty} \frac{H_n^{(3)}}{n}(1-t)^n = Li_4(1-t) - \varsigma(3)\log t - \sum_{n=1}^{\infty} \frac{1}{n^3}\sum_{k=1}^{n}\frac{(1-t)^k}{k}$$

It is easily seen that

$$-\sum_{n=1}^{\infty} \frac{1}{n^3}\sum_{k=1}^{n}\frac{(1-t)^k}{k} = \sum_{n=1}^{\infty} \frac{1}{n^3}\sum_{k=1}^{n}\frac{1-(1-t)^k}{k} - \sum_{n=1}^{\infty} \frac{H_n^{(1)}}{n^3}$$

Therefore upon dividing by $t$ and integrating we have

$$-\int_a^x \sum_{n=1}^{\infty} \frac{1}{n^3}\sum_{k=1}^{n}\frac{(1-t)^k}{kt}\,dt = \int_a^x \sum_{n=1}^{\infty} \frac{1}{n^3}\sum_{k=1}^{n}\frac{1-(1-t)^k}{kt}\,dt - \frac{1}{2}\varsigma^2(2)(\log x - \log a)$$

where we have used (4.4.167u).Using (4.4.155a) we obtain

$(4.4.168d)$

$$\int_a^x \sum_{n=1}^{\infty} \frac{1}{n^3}\sum_{k=1}^{n}\frac{1-(1-t)^k}{kt}\,dt = \sum_{n=1}^{\infty} \frac{1}{n^3}\sum_{k=1}^{n}\frac{1}{k}\int_a^x \frac{1-(1-t)^k}{t}\,dt = \sum_{n=1}^{\infty} \frac{1}{n^3}\sum_{k=1}^{n}\frac{1}{k}\sum_{l=1}^{k}\frac{(1-a)^l-(1-x)^l}{l}$$

We are now in the realm of multiple Euler-Zagier sums and I shall sensibly leave any further analysis in this area to the experts.

A long time after I wrote the above, I noticed an integral provided by Freitas [69a], which I will use later in (4.4.239a), and realised that it provided the missing link in my work on (4.4.167o):

$(4.4.168e)$ 
$$\int_0^1 \frac{\log x\, Li_2(x)}{1-x}\,dx = \int_0^1 \frac{\log(1-x)\, Li_2(1-x)}{x}\,dx = -\frac{3}{4}\varsigma(4)$$

So let's revisit that problem after our little diversion!

From (4.4.167o) we have

$$2\int_a^t \frac{\log x\log^2(1-x)}{x}\,dx = \log x\left[\log^2(1-x)\log x + 2\log(1-x)Li_2(1-x) - 2Li_3(1-x)\right]\Big|_a^t$$

$$-2\int_a^t \frac{Li_2(1-x)\log(1-x)}{x}\,dx + 2\int_a^t \frac{Li_3(1-x)}{x}\,dx$$

Therefore, by judiciously adding and subtracting a factor of $2\varsigma(3)(\log t - \log a)$, we have

$$2\int_a^t \frac{\log x \log^2(1-x)}{x}\,dx = \log x\Big[\log^2(1-x)\log x + 2\log(1-x)Li_2(1-x) - 2Li_3(1-x)\Big]\Big|_a^t$$

$$-2\int_a^t \frac{Li_2(1-x)\log(1-x)}{x}\,dx - 2\int_a^t \frac{\varsigma(3) - Li_3(1-x)}{x}\,dx + 2\varsigma(3)(\log t - \log a)$$

We now take the limit as $a \to 0$ and noting that

$$\lim_{a\to 0}(2Li_3(1-a)\log a - 2\varsigma(3)\log a) = 0$$

we obtain

(4.4.168f)

$$\int_0^t \frac{\log x \log^2(1-x)}{x}\,dx = \frac{1}{2}\log t\Big[\log^2(1-t)\log t + 2\log(1-t)Li_2(1-t) - 2Li_3(1-t) + 2\varsigma(3)\Big]$$

$$-\int_0^t \frac{Li_2(1-x)\log(1-x)}{x}\,dx - \int_0^t \frac{\varsigma(3) - Li_3(1-x)}{x}\,dx$$

An alternative proof of (4.4.168f) was given in Volume I and is repeated below for ease of reference. We have from (4.4.100gii) in Volume III

$$\int_0^t \frac{\log^2(1-x)}{x}\,dx = \log t \log^2(1-t) + 2\log(1-t)Li_2(1-t) - 2Li_3(1-t) + 2\varsigma(3)$$

and integration by parts gives us

$$\int_0^t \frac{\log x \log^2(1-x)}{x}\,dx = \log t\Big[\log t \log^2(1-t) + 2\log(1-t)Li_2(1-t) - 2Li_3(1-t) + 2\varsigma(3)\Big]$$

$$-\int_0^t \Big[\log t \log^2(1-t) + 2\log(1-t)Li_2(1-t) - 2Li_3(1-t) + 2\varsigma(3)\Big]\frac{dx}{x}$$

Therefore we have

$$\int_0^t \frac{\log x \log^2(1-x)}{x}\,dx = \frac{1}{2}\log t\left[\log t\log^2(1-t) + 2\log(1-t)Li_2(1-t) - 2Li_3(1-t) + 2\varsigma(3)\right]$$

$$-\int_0^t \frac{\log(1-t)Li_2(1-t)}{x}\,dx - \int_0^t \frac{\varsigma(3) - Li_3(1-t)}{x}\,dx$$

□

Letting $t = 1$, using (4.4.167v) and (4.4.168e) we have

(4.4.168g)    $$\int_0^1 \frac{\varsigma(3) - Li_3(1-x)}{x}\,dx = \int_0^1 \frac{\varsigma(3) - Li_3(x)}{1-x}\,dx = \sum_{n=1}^\infty \frac{H_n^{(1)}}{n^3} = \frac{\pi^4}{72} = \frac{5}{4}\varsigma(4)$$

and we therefore rediscover the Rutledge and Douglass [116aa] formula (4.4.167qa)

$$\int_0^1 \frac{\log x \log^2(1-x)}{x}\,dx = \frac{3}{4}\varsigma(4) - \frac{5}{4}\varsigma(4) = -\frac{1}{2}\varsigma(4)$$

Somewhat more generally, using (4.4.167g) for $0 < t \leq 1$

$$\int_0^t \frac{\varsigma(3) - Li_3(1-x)}{x}\,dx = \sum_{n=1}^\infty \frac{H_n^{(1)}}{n^3} + \varsigma(3)\log t + \sum_{n=1}^\infty \frac{H_n^{(3)}}{n}(1-t)^n - Li_4(1-t)$$

we have

$$\int_0^t \frac{\log x \log^2(1-x)}{x}\,dx + \int_0^t \frac{Li_2(1-x)\log(1-x)}{x}\,dx =$$

$$\frac{1}{2}\log^2(1-t)\log^2 t + Li_2(1-t)\log(1-t)\log t - Li_3(1-t)\log t + Li_4(1-t)$$

$$-\sum_{n=1}^\infty \frac{H_n^{(1)}}{n^3} - \sum_{n=1}^\infty \frac{H_n^{(3)}}{n}(1-t)^n$$

Using (4.4.167ma) we have

$$\int_0^t \frac{\log(1-x)Li_2(1-x)}{x}\,dx + \int_0^t \frac{\log x \log^2(1-x)}{x}\,dx = -\varsigma(2)Li_2(t) + \frac{1}{2}\left[Li_2(t)\right]^2$$

and we therefore get

$(4.4.168h)$ $\quad \dfrac{1}{2}\log^2(1-t)\log^2 t + Li_2(1-t)\log(1-t)\log t - Li_3(1-t)\log t + Li_4(1-t)$

$$-\sum_{n=1}^{\infty}\frac{H_n^{(1)}}{n^3} - \sum_{n=1}^{\infty}\frac{H_n^{(3)}}{n}(1-t)^n = -\varsigma(2)Li_2(t) + \frac{1}{2}\left[Li_2(t)\right]^2$$

Letting $t = 1$ we obtain as before $\displaystyle\sum_{n=1}^{\infty}\frac{H_n^{(1)}}{n^3} = \frac{1}{2}\varsigma^2(2)$. Letting $t = 1/2$ we get

$$\sum_{n=1}^{\infty}\frac{H_n^{(3)}}{n2^n} = \frac{1}{2}\log^4 2 + Li_2(1/2)\log^2 2 - Li_3(1/2)\log 2 + Li_4(1/2)$$

and using the Euler/Landen identities this may be written as

$(4.4.168i)$ $\quad \displaystyle\sum_{n=1}^{\infty}\frac{H_n^{(3)}}{n2^n} = \varsigma(2)\log^2 2 - \frac{7}{8}\varsigma(3)\log 2 + Li_4(1/2) - \frac{1}{6}\log^4 2$

In $(4.4.247cii)$ we will show that

$$Li_4(x) = \sum_{n=1}^{\infty}\frac{H_n^{(3)}}{n}x^n + \varsigma(3)\log(1-x) + \sum_{n=1}^{\infty}\frac{1}{n^3}\sum_{k=1}^{n}\frac{x^k}{k}$$

and hence we have

$$\sum_{n=1}^{\infty}\frac{H_n^{(3)}}{n2^n} = Li_4(1/2) + \varsigma(3)\log 2 - \sum_{n=1}^{\infty}\frac{1}{n^3}\sum_{k=1}^{n}\frac{1}{k2^k}$$

This then results in

$(4.4.168ii)$ $\quad \displaystyle\sum_{n=1}^{\infty}\frac{1}{n^3}\sum_{k=1}^{n}\frac{1}{k2^k} = \frac{15}{8}\varsigma(3)\log 2 - \varsigma(2)\log^2 2 + \frac{1}{6}\log^4 2$

From $(4.4.167g)$ we have for $0 < t \le 1$

$$\int_0^t \frac{\varsigma(3) - Li_3(1-x)}{x}dx = \sum_{n=1}^{\infty}\frac{H_n^{(1)}}{n^3} + \varsigma(3)\log t + \sum_{n=1}^{\infty}\frac{H_n^{(3)}}{n}(1-t)^n - Li_4(1-t)$$

I didn't realise until I had nearly completed this paper that one could directly evaluate the above integral until I spotted it in one of Blümlein's papers "Analytic Continuation of Mellin Transforms up to two-loop Order" [24aa]. Using integration by parts we have

$$\int_0^t \frac{\varsigma(3) - Li_3(1-x)}{x}\, dx = \left[\varsigma(3) - Li_3(1-t)\right]\log t - \int_0^t \frac{\log x\, Li_2(1-x)}{1-x}\, dx$$

and hence we get for $0 \le t \le 1$

(4.4.168j) $\quad \int_0^t \frac{\varsigma(3) - Li_3(1-x)}{x}\, dx = \left[\varsigma(3) - Li_3(1-t)\right]\log t - \frac{1}{2}\left[Li_2(1-t)\right]^2 + \frac{1}{2}\varsigma^2(2)$

Hence we obtain

$$\sum_{n=1}^{\infty} \frac{H_n^{(1)}}{n^3} + \sum_{n=1}^{\infty} \frac{H_n^{(3)}}{n}(1-t)^n - Li_4(1-t) = -Li_3(1-t)\log t - \frac{1}{2}\left[Li_2(1-t)\right]^2 + \frac{1}{2}\varsigma^2(2)$$

Letting $t = 1$ we very easily obtain a further proof of (4.4.167u)

$$\sum_{n=1}^{\infty} \frac{H_n^{(1)}}{n^3} = \frac{1}{2}\varsigma^2(2)$$

and we then obtain for $0 < t \le 1$

(4.4.168k) $\quad \sum_{n=1}^{\infty} \frac{H_n^{(3)}}{n}(1-t)^n = Li_4(1-t) - Li_3(1-t)\log t - \frac{1}{2}\left[Li_2(1-t)\right]^2$

With $t = 1/2$ we get

(4.4.168l) $\quad \sum_{n=1}^{\infty} (-1)^n \frac{H_n^{(3)}}{n2^n} = Li_4(1/2) + Li_3(1/2)\log 2 - \frac{1}{2}\left[Li_2(1/2)\right]^2$

Substituting (4.4.168j) in (4.4.168a) we get

(4.4.168m) $\quad \sum_{n=1}^{\infty} \frac{H_n^{(2)}}{n^2}t^n = -2\left[\varsigma(3) - Li_3(1-t)\right]\log t + \left[Li_2(1-t)\right]^2 - \varsigma^2(2)$

$$+ 2Li_4(1-t) + Li_4(t) - \frac{1}{2}\left[Li_2(t)\right]^2 + 2\varsigma(2)Li_2(t) + F(t)$$

It may be noted that by substituting (4.4.168j) in (4.4.158) we simply reproduce Euler's dilogarithm identity.

Let us now divide (4.4.168k) by $(1-t)$ and integrate to obtain

(4.4.168n)

$$\sum_{n=1}^{\infty} \frac{H_n^{(3)}}{n^2} - \sum_{n=1}^{\infty} \frac{H_n^{(3)}}{n^2}(1-x)^n = \varsigma(5) - Li_5(1-x) - \int_0^x \frac{Li_3(1-t)\log t}{1-t}dt - \frac{1}{2}\int_0^x \frac{\left[Li_2(1-t)\right]^2}{1-t}dt$$

Integration by parts gives us

$$\int_0^x \frac{\left[Li_2(1-t)\right]^2}{1-t}dt = -Li_3(1-x)Li_2(1-x) + \varsigma(3)\varsigma(2) + \int_0^x \frac{Li_3(1-t)\log t}{1-t}dt$$

and we therefore obtain

(4.4.168o)

$$\sum_{n=1}^{\infty} \frac{H_n^{(3)}}{n^2} - \sum_{n=1}^{\infty} \frac{H_n^{(3)}}{n^2}(1-x)^n + Li_5(1-x) - \varsigma(5) + Li_2(1-x)Li_3(1-x) - \varsigma(2)\varsigma(3) + \frac{3}{2}\int_0^x \frac{\left[Li_2(1-t)\right]^2}{1-t}dt = 0$$

When $x = 1$ we obtain

$$\sum_{n=1}^{\infty} \frac{H_n^{(3)}}{n^2} - \varsigma(5) - \varsigma(2)\varsigma(3) + \frac{3}{2}\int_0^1 \frac{\left[Li_2(1-t)\right]^2}{1-t}dt = 0$$

Since $\int_0^1 \frac{\left[Li_2(1-t)\right]^2}{1-t}dt = \int_0^1 \frac{\left[Li_2(t)\right]^2}{t}dt$ and from (4.4.168c) we have

$$\int_0^1 \frac{\left[Li_2(t)\right]^2}{t}dt = 2\varsigma(2)\varsigma(3) - 3\varsigma(5)$$

we conclude that (as previously seen in (3.211f) in Volume I)

(4.4.168p)  $$\sum_{n=1}^{\infty} \frac{H_n^{(3)}}{n^2} = \frac{11}{2}\varsigma(5) - 2\varsigma(2)\varsigma(3)$$

With reference to (4.4.168n) we note that employing integration by parts results in

$$\int_a^x \frac{Li_3(1-t)\log t}{1-t}dt =$$

$$-Li_4(1-x)\log x + Li_4(1-a)\log a + \varsigma(4)\big[\log x - \log a\big] + \int\limits_a^x \frac{Li_4(1-t) - \varsigma(4)}{t}\,dt$$

$$= -Li_4(1-x)\log x + [Li_4(1-a) - \varsigma(4)]\log a + \varsigma(4)\log x + \int\limits_a^x \frac{Li_4(1-t) - \varsigma(4)}{t}\,dt$$

Therefore as $a \to 0$ we obtain

$$\int\limits_0^x \frac{Li_3(1-t)\log t}{1-t}\,dt = -Li_4(1-x)\log x + \varsigma(4)\log x + \int\limits_0^x \frac{Li_4(1-t) - \varsigma(4)}{t}\,dt$$

From (3.42) and (3.45) we have

$$\int \frac{Li_2(x)}{1-x}\,dx = -Li_2(x)\log(1-x) - \int \frac{\log^2(1-x)}{x}\,dx$$

$$\int \frac{\log^2(1-x)}{x}\,dx = \log^2(1-x)\log x + 2\log(1-x)Li_2(1-x) - 2Li_3(1-x)$$

and therefore we get

$$\int \frac{Li_2(x)}{1-x}\,dx = -Li_2(x)\log(1-x) - \log^2(1-x)\log x - 2\log(1-x)Li_2(1-x) + 2Li_3(1-x)$$

Hence, employing integration by parts, we obtain

$$\int\limits_a^t \frac{Li_2(x)\log x}{1-x}\,dx =$$

$$\log x\Big[ -Li_2(x)\log(1-x) - \log^2(1-x)\log x - 2\log(1-x)Li_2(1-x) + 2Li_3(1-x) - 2\varsigma(3)\Big]\Big|_a^t$$

$$-\int\limits_a^t \frac{-Li_2(x)\log(1-x) - \log^2(1-x)\log x - 2\log(1-x)Li_2(1-x) + 2Li_3(1-x) - 2\varsigma(3)}{x}\,dx$$

As before, this then can be written as

$$\int\limits_a^t \frac{Li_2(x)\log x}{1-x}\,dx - \int\limits_a^t \frac{\log x \log^2(1-x)}{x}\,dx =$$

$$\log x\Big[-Li_2(x)\log(1-x)-\log^2(1-x)\log x-2\log(1-x)Li_2(1-x)+2Li_3(1-x)-2\varsigma(3)\Big]\Big|_a^t$$

$$-\int_a^t \frac{-Li_2(x)\log(1-x)-2\log(1-x)Li_2(1-x)+2Li_3(1-x)-2\varsigma(3)}{x}\,dx$$

Therefore we have

$$\int_a^t \frac{Li_2(x)\log x}{1-x}\,dx-\int_a^t \frac{\log x\log^2(1-x)}{x}\,dx=$$

$$=\log x\Big[-Li_2(x)\log(1-x)-\log^2(1-x)\log x-2\log(1-x)Li_2(1-x)+2Li_3(1-x)-2\varsigma(3)\Big]\Big|_a^t$$

$$+\frac{1}{2}\Big[Li_2(1-t)\Big]^2-\frac{1}{2}\varsigma^2(2)+2\int_a^t \frac{\log(1-x)Li_2(1-x)}{x}\,dx-2\int_a^t \frac{Li_3(1-x)-\varsigma(3)}{x}\,dx$$

and, in the limit as $a\to 0$, we get

(4.4.168q)

$$\int_0^t \frac{Li_2(x)\log x}{1-x}\,dx-\int_0^t \frac{\log x\log^2(1-x)}{x}\,dx=$$

$$=\log t\Big[-Li_2(t)\log(1-t)-\log^2(1-t)\log t-2\log(1-t)Li_2(1-t)+2Li_3(1-t)-2\varsigma(3)\Big]$$

$$+\frac{1}{2}\Big[Li_2(1-t)\Big]^2-\frac{1}{2}\varsigma^2(2)+2\int_0^t \frac{\log(1-x)Li_2(1-x)}{x}\,dx-2\int_0^t \frac{Li_3(1-x)-\varsigma(3)}{x}\,dx$$

We now refer back to (4.4.168f)

(4.4.168r)

$$\int_0^t \frac{\log x\log^2(1-x)}{x}\,dx=\frac{1}{2}\log t\Big[\log^2(1-t)\log t+2\log(1-t)Li_2(1-t)-2Li_3(1-t)+2\varsigma(3)\Big]$$

$$-\int_0^t \frac{Li_2(x)\log x}{1-x}\,dx-\int_0^t \frac{\varsigma(3)-Li_3(1-x)}{x}\,dx$$

Combining (4.4.168q) and (4.4.168r) we get

$$\int_0^t \frac{Li_2(x)\log x}{1-x}\,dx =$$

$$\log t\Big[-Li_2(t)\log(1-t)-\log^2(1-t)\log t-2\log(1-t)Li_2(1-t)+2Li_3(1-t)-2\varsigma(3)\Big]$$

$$+\frac{1}{2}\big[Li_2(1-t)\big]^2-\frac{1}{2}\varsigma^2(2)+2\int_0^t \frac{\log(1-x)Li_2(1-x)}{x}\,dx-2\int_0^t \frac{Li_3(1-x)-\varsigma(3)}{x}\,dx$$

$$+\frac{1}{2}\log t\Big[\log^2(1-t)\log t+2\log(1-t)Li_2(1-t)-2Li_3(1-t)+2\varsigma(3)\Big]$$

$$-\int_0^t \frac{Li_2(x)\log x}{1-x}\,dx-\int_0^t \frac{\varsigma(3)-Li_3(1-x)}{x}\,dx$$

This then becomes

$$2\int_0^t \frac{Li_2(x)\log x}{1-x}\,dx =$$

$$\log t\Big[-Li_2(t)\log(1-t)-\log^2(1-t)\log t-2\log(1-t)Li_2(1-t)+2Li_3(1-t)-2\varsigma(3)\Big]$$

$$+\frac{1}{2}\big[Li_2(1-t)\big]^2-\frac{1}{2}\varsigma^2(2)+2\int_0^t \frac{\log(1-x)Li_2(1-x)}{x}\,dx-\int_0^t \frac{Li_3(1-x)-\varsigma(3)}{x}\,dx$$

$$+\frac{1}{2}\log t\Big[\log^2(1-t)\log t+2\log(1-t)Li_2(1-t)-2Li_3(1-t)+2\varsigma(3)\Big]$$

Letting $t=1$ we get

$$2\int_0^1 \frac{Li_2(x)\log x}{1-x}\,dx =$$

$$-\frac{1}{2}\varsigma^2(2)+2\int_0^1 \frac{\log(1-x)Li_2(1-x)}{x}\,dx-\int_0^1 \frac{Li_3(1-x)-\varsigma(3)}{x}\,dx$$

and since $\displaystyle\int_0^1 \frac{Li_2(x)\log x}{1-x}\,dx=\int_0^1 \frac{\log(1-x)Li_2(1-x)}{x}\,dx$ we see that

$$\int_0^1 \frac{Li_3(1-x) - \varsigma(3)}{x} \, dx = -\frac{5}{4}\varsigma(4) = -\frac{1}{2}\varsigma^2(2)$$

Employing Euler's dilogarithm identity (1.6c) we have

$$\int_0^t \frac{Li_2(x)\log x}{1-x} \, dx = \varsigma(2)\int_0^t \frac{\log x}{1-x} \, dx - \int_0^t \frac{\log^2 x \log(1-x)}{1-x} \, dx - \int_0^t \frac{Li_2(1-x)\log x}{1-x} \, dx$$

$$= \varsigma(2)\big[Li_2(1-t) - \varsigma(2)\big] - \int_0^t \frac{\log^2 x \log(1-x)}{1-x} \, dx - \frac{1}{2}\big[Li_2(1-t)\big]^2 + \frac{1}{2}\varsigma^2(2)$$

Using the Wolfram Integrator we obtain

$$-\int \frac{\log(1-x)\log^2 x}{1-x} \, dx =$$

$$\frac{1}{2}\log^2(1-x)\log^2 x + \frac{1}{12}\log^4(1-x) - \log^2(1-x)\log^2 x + \frac{2}{3}\log^3(1-x)\log x$$

$$-\left[\log x + \frac{1}{3}\log(1-x)\right]\log^2(1-x)\log\left(\frac{1-x}{x}\right) + \frac{1}{2}\log^2(1-x)\left[\log\left(\frac{1-x}{x}\right)\right]^2$$

$$-\frac{1}{4}\left[\log\left(\frac{1-x}{x}\right)\right]^4 + \log^2 x \, Li_2(x) - \log^2(1-x)Li_2(1-x)$$

$$-\left[\log\left(\frac{1-x}{x}\right)\right]^2 Li_2\left(\frac{x-1}{x}\right) - 2\log x \, Li_3(x) + 2\log(1-x)\, Li_3(1-x)$$

$$+2\log\left(\frac{1-x}{x}\right)Li_3\left(\frac{x-1}{x}\right) + 2\left[Li_4(x) - Li_4(1-x) - Li_4\left(\frac{x-1}{x}\right)\right]$$

However, it is not immediately apparent to me how the integrated part behaves as $x \to 0$, having regard to the polylogarithmic terms involving $(x-1/x)$. It appears to me that a combination of the above identities would produce a mixed functional equation involving, inter alia,

$$Li_4(x), \; Li_4(1-x), \; Li_4\left(-\frac{x}{1-x}\right) \text{ and } \; Li_4\left(-\frac{1-x}{x}\right)$$

combined with similar expressions for the two polylogarithms of lower order.

**(xv)** For completeness, proofs by induction of Adamchik's formulae (3.17), (3.18) and (3.19) in Volume I are set out below.

**Theorem 4.10(i):**

$$(4.4.169) \qquad \sum_{k=1}^{n} \frac{H_k^{(1)}}{k} = \frac{1}{2}\left(H_n^{(1)}\right)^2 + \frac{1}{2}H_n^{(2)}$$

**Proof:**

The theorem is obviously true for $n = 1$. Assume that it is valid for $n = N$. Then we have

$$\sum_{k=1}^{N+1} \frac{H_k^{(1)}}{k} = \sum_{k=1}^{N} \frac{H_k^{(1)}}{k} + \frac{H_{N+1}^{(1)}}{N+1}$$

$$(4.4.170) \qquad = \frac{1}{2}\left(H_N^{(1)}\right)^2 + \frac{1}{2}H_N^{(2)} + \frac{H_{N+1}^{(1)}}{N+1}$$

Now write the right hand side of (4.4.169) as

$$\frac{1}{2}\left(H_{N+1}^{(1)}\right)^2 + \frac{1}{2}H_{N+1}^{(2)} = \frac{1}{2}\left(H_N^{(1)} + \frac{1}{N+1}\right)^2 + \frac{1}{2}H_N^{(2)} + \frac{1}{2}\frac{1}{(N+1)^2}$$

$$= \frac{1}{2}\left(H_N^{(1)}\right)^2 + \frac{1}{2}H_N^{(2)} + \frac{H_N^{(1)}}{N+1} + \frac{1}{(N+1)^2}$$

$$= \frac{1}{2}\left(H_N^{(1)}\right)^2 + \frac{1}{2}H_N^{(2)} + \frac{H_{N+1}^{(1)}}{N+1}$$

$$= \sum_{k=1}^{N+1} \frac{H_k^{(1)}}{k} \quad, \text{ using (4.4.170)}$$

Therefore the theorem has been proved by mathematical induction.

**Theorem 4.10(ii):**

$$(4.4.171) \qquad \sum_{k=1}^{n} \frac{H_k^{(2)}}{k} + \sum_{k=1}^{n} \frac{H_k^{(1)}}{k^2} = H_n^{(3)} + H_n^{(1)}H_n^{(2)}$$

**Proof:**

The theorem is obviously true for $n = 1$. Assume that it is valid for $n = N$. Then we have

$$H_{N+1}^{(3)} + H_{N+1}^{(1)} H_{N+1}^{(2)} = H_N^{(3)} + \frac{1}{(N+1)^3} + \left( H_N^{(1)} + \frac{1}{N+1} \right)\left( H_N^{(2)} + \frac{1}{(N+1)^2} \right)$$

$$= H_N^{(3)} + H_N^{(1)} H_N^{(2)} + \frac{H_N^{(1)}}{(N+1)^2} + \frac{H_N^{(2)}}{N+1} + \frac{2}{(N+1)^2}$$

Now write the left hand side of (4.4.171) as

$$\sum_{k=1}^{N+1} \frac{H_k^{(2)}}{k} + \sum_{k=1}^{N+1} \frac{H_k^{(1)}}{k^2} = \sum_{k=1}^{N} \frac{H_k^{(2)}}{k} + \frac{H_{N+1}^{(2)}}{N+1} + \sum_{k=1}^{N} \frac{H_k^{(1)}}{k^2} + \frac{H_{N+1}^{(1)}}{(N+1)^2}$$

$$= \sum_{k=1}^{N} \frac{H_k^{(2)}}{k} + \sum_{k=1}^{N} \frac{H_k^{(1)}}{k^2} + \frac{H_N^{(2)}}{N+1} + \frac{1}{(N+1)^3} + \frac{H_N^{(1)}}{(N+1)^2} + \frac{1}{(N+1)^3}$$

$$= H_N^{(3)} + H_N^{(1)} H_N^{(2)} + \frac{H_N^{(1)}}{(N+1)^2} + \frac{H_N^{(2)}}{N+1} + \frac{2}{(N+1)^2}$$

$$= H_{N+1}^{(3)} + H_{N+1}^{(1)} H_{N+1}^{(2)}$$

Hence the identity is also valid for $n = N+1$ and the theorem is proved.

**Theorem 4.10(iii):**

$$(4.4.172) \qquad \sum_{k=1}^{n} \frac{\left( H_k^{(1)} \right)^2}{k} + \sum_{k=1}^{n} \frac{H_k^{(2)}}{k} = \frac{1}{3} \left\{ \left( H_n^{(1)} \right)^3 + 3 H_n^{(1)} H_n^{(2)} + 2 H_n^{(3)} \right\}$$

**Proof:**

The theorem is obviously true for $n = 1$. Assume that it is valid for $n = N$. Then we have

$$\frac{1}{3} \left\{ \left( H_{N+1}^{(1)} \right)^3 + 3 H_{N+1}^{(1)} H_{N+1}^{(2)} + 2 H_{N+1}^{(3)} \right\} =$$

$$\frac{1}{3} \left\{ \left( H_n^{(1)} + \frac{1}{n+1} \right)^3 + 3 \left( H_n^{(1)} + \frac{1}{n+1} \right)\left( H_n^{(2)} + \frac{1}{(n+1)^2} \right) + 2 \left( H_n^{(3)} + \frac{1}{(n+1)^3} \right) \right\}$$

$$= \frac{1}{3}\left\{\left(H_N^{(1)}\right)^3 + 3H_N^{(1)}H_N^{(2)} + 2H_N^{(3)}\right\} + \frac{\left(H_N^{(1)}\right)^2}{N+1} + 2\frac{H_N^{(1)}}{(N+1)^2} + \frac{H_N^{(2)}}{N+1} + 2\frac{1}{(N+1)^3}$$

We also have

$$\sum_{k=1}^{N+1}\frac{\left(H_k^{(1)}\right)^2}{k} + \sum_{k=1}^{N+1}\frac{H_k^{(2)}}{k} = \sum_{k=1}^{N}\frac{\left(H_k^{(1)}\right)^2}{k} + \sum_{k=1}^{N}\frac{H_k^{(2)}}{k} + \frac{\left(H_{N+1}^{(1)}\right)^2}{N+1} + \frac{H_{N+1}^{(2)}}{N+1}$$

$$= \sum_{k=1}^{N}\frac{\left(H_k^{(1)}\right)^2}{k} + \sum_{k=1}^{N}\frac{H_k^{(2)}}{k} + \frac{\left(H_N^{(1)} + \frac{1}{N+1}\right)^2}{N+1} + \frac{\left(H_N^{(2)} + \frac{1}{(N+1)^2}\right)}{N+1}$$

$$= \sum_{k=1}^{N}\frac{\left(H_k^{(1)}\right)^2}{k} + \sum_{k=1}^{N}\frac{H_k^{(2)}}{k} + \frac{\left(H_N^{(1)}\right)^2}{N+1} + 2\frac{H_N^{(1)}}{(N+1)^2} + \frac{H_N^{(2)}}{N+1} + \frac{2}{(N+1)^3}$$

It can therefore be seen that the proofs of these rather complicated looking identities are actually fairly trivial.

We also see that

(4.4.172a) $\qquad \sum_{k=1}^{n}\frac{\left(H_k^{(1)}\right)^2}{k} + \sum_{k=1}^{n}\frac{H_k^{(2)}}{k} = \sum_{k=1}^{n}\frac{\left(H_k^{(1)}\right)^2 + H_k^{(2)}}{k} = 2\sum_{k=1}^{n}\frac{1}{k}\sum_{j=1}^{k}\frac{H_j^{(1)}}{j} = 2\sum_{k=1}^{n}\frac{1}{k}\sum_{j=1}^{k}\frac{1}{j}\sum_{l=1}^{j}\frac{1}{l}$

and therefore we get

(4.4.172b) $\qquad \sum_{k=1}^{n}\frac{1}{k}\sum_{j=1}^{k}\frac{1}{j}\sum_{l=1}^{j}\frac{1}{l} = \frac{1}{6}\left[\left(H_n^{(1)}\right)^3 + 3H_n^{(1)}H_n^{(2)} + 2H_n^{(3)}\right]$

**(xvi)** In connection with Lemma 3.4 in Volume I we stated that $P_s(1)$ was not convergent for $s = 1$. The proof of this assertion is set out below.

$P_s(1)$ was defined as:

(4.4.173) $$P_s(1) = \sum_{n=1}^{\infty}\frac{1}{2^n}\sum_{k=1}^{n}\binom{n}{k}\frac{1}{k^s} = 2\varsigma(s)$$

Letting $t = -1$ in (4.1.6) we obtain

(4.4.174) $$\sum_{k=1}^{n}\binom{n}{k}\frac{1}{k} = \sum_{k=1}^{n}\frac{2^k - 1}{k}$$

$$> \sum_{k=1}^{n} \frac{2^k - 1}{n}$$

$$= \frac{1}{n}\left(2^{n+1} - 2 - n\right)$$

and hence we have

(4.4.175)
$$\sum_{n=1}^{N} \frac{1}{2^n} \sum_{k=1}^{n} \binom{n}{k} \frac{1}{k} > \sum_{n=1}^{N} \frac{1}{2^n} \frac{1}{n}\left(2^{n+1} - 2 - n\right)$$

(4.4.176)
$$= 2\sum_{n=1}^{N} \frac{1}{n} - \sum_{n=1}^{N} \frac{1}{n2^n} - \sum_{n=1}^{N} \frac{1}{2^n}$$

As $N \to \infty$ the harmonic series in (4.4.176) clearly diverges, whilst the other two series there converge to $\log 2$ and 1 respectively. This therefore proves that the series $P_1(1)$ is divergent. It should be noted that I carried out this exercise before I had realised that $P_s(x) = 2Li_s(x)$ and hence that $P_1(x) = -2\log(1-x)$.

**(xvii)** Often in mathematics we look for divine inspiration but we do not usually expect to obtain it from a canonised saint. This is indeed the source of the next remark. The calculation of $g^{(m)}(x)$, as defined in (4.2.1), effectively involves the derivative of a composite function $g(f(t))$ and the general formula for this was discovered by Francesco Faà di Bruno (1825-1888) who was declared a Saint by Pope John Paul II in St. Peter's Square in Rome on 25 September 1988 [127a].

Di Bruno [82] showed that

(4.4.177)
$$\frac{d^{(m)}}{dt^{(m)}} g(f(t)) = \sum \frac{m!}{b_1! b_2! \ldots b_m!} g^{(k)}(f(t)) \left(\frac{f^{(1)}(t)}{1!}\right)^{b_1} \left(\frac{f^{(2)}(t)}{2!}\right)^{b_2} \cdots \left(\frac{f^{(m)}(t)}{m!}\right)^{b_m}$$

where the sum is over all different solutions in non-negative integers $b_1, \ldots, b_m$ of $b_1 + 2b_2 + \ldots + mb_m$, and $k = b_1 + \ldots + b_m$. In our case, the composite function was of the form $g(f(t))$ where $g(t) = 1/t$ and $f(t) = t(t+1)\ldots(t+n)$ and we therefore have an explanation why the Flajolet and Sedgewick formulation of (3.14) and the di Bruno formula both involve the use of Bell polynomials.

In [123c] Gould reminded the mathematical community of the "not well-known" formula for the $n$th derivative of a composite function $f(z)$ where $z$ is a function of $x$, namely

$$(4.4.177a) \qquad D_x^{(n)} f(z) = \sum_{k=1}^{n} D_z^{(k)} f(z) \frac{(-1)^k}{k!} \sum_{j=1}^{k} (-1)^j \binom{k}{j} z^{k-j} \qquad , \text{ for } n \geq 1$$

This expression is frequently easier to handle than the di Bruno algorithm.

**(xviii)** From (4.4.25) we have

$$(4.4.178) \qquad Li_s(x) = \frac{x}{\Gamma(s)} \int_0^\infty \frac{u^{s-1}}{e^u - x} du$$

In (4.4.178) let $x = 1$, $u = 2\pi t$ and $s = 2n$ to obtain

$$(4.4.179) \qquad \frac{\varsigma(2n)\Gamma(2n)}{(2\pi)^{2n}} = \int_0^\infty \frac{t^{2n-1}}{e^{2\pi t} - 1} dt$$

Using Euler's identity (1.7) we therefore have an integral expression for the Bernoulli numbers

$$(4.4.180) \qquad B_{2n} = 4n(-1)^{n+1} \int_0^\infty \frac{t^{2n-1}}{e^{2\pi t} - 1} dt$$

and this corrects a typographical error in [25, p.223]. With the substitution $u = 2\pi t$ we have

$$B_{2n} = \frac{4n(-1)^{n+1}}{(2\pi)^n} \int_0^\infty \frac{u^{2n-1}}{e^u - 1} du$$

We have $\int \frac{1}{e^u - 1} du = \int \frac{e^{-u}}{1 - e^{-u}} du = \log(1 - e^{-u})$ and using integration by parts we obtain

$$\int_0^\infty \frac{u^k}{e^u - 1} du = u^k \log(1 - e^{-u}) \Big|_0^\infty - k \int_0^\infty u^{k-1} \log(1 - e^{-u}) du$$

Hence we have

$$(4.4.181) \qquad \int_0^\infty \frac{u^k}{e^u - 1} du = -k \int_0^\infty u^{k-1} \log(1 - e^{-u}) du$$

and we therefore get

$$(4.4.182) \qquad B_{2n} = \frac{2(-1)^n 2n(2n-1)}{(2\pi)^n} \int\limits_0^\infty u^{2n-2} \log(1-e^{-u}) du$$

With the substitution $y = e^{-u}$ we have

$$\int \log(1-e^{-u}) du = -\int \frac{\log(1-y)}{y} dy = Li_2(e^{-u})$$

and we therefore obtain

$$\int\limits_0^\infty u^k \log(1-e^{-u}) du = u^k Li_2(e^{-u})\Big|_0^\infty - k\int\limits_0^\infty u^{k-1} Li_2(e^{-u}) du$$

$$= -k\int\limits_0^\infty u^{k-1} Li_2(e^{-u}) du$$

Therefore we have

$$(4.4.183) \qquad \int\limits_0^\infty u^{2n-2} \log(1-e^{-u}) du = -(2n-2)\int\limits_0^\infty u^{2n-3} Li_2(e^{-u}) du$$

Hence we have $\qquad B_{2n} = \frac{2(-1)^{n+1} 2n(2n-1)(2n-2)}{(2\pi)^n} \int\limits_0^\infty u^{2n-3} Li_2(e^{-u}) du$

Using the series expansion of the dilogarithm we obtain

$$\int\limits_0^\infty u^{2n-3} Li_2(e^{-u}) du = \sum_{k=1}^\infty \frac{1}{k^2} \int\limits_0^\infty u^{2n-3} e^{-ku} du$$

Since by definition $\Gamma(\nu) = \int\limits_0^\infty u^{\nu-1} e^{-u} du$ we have $\int\limits_0^\infty u^{\nu-1} e^{-ku} du = \frac{\Gamma(\nu)}{k^\nu}$ and hence

$$\int\limits_0^\infty u^{2n-3} e^{-ku} du = \frac{\Gamma(2n-2)}{k^{2n-2}}$$

Completing the summation we have

$$\sum_{k=1}^\infty \frac{1}{k^2} \int\limits_0^\infty u^{2n-3} e^{-ku} du = \Gamma(2n-2)\varsigma(2n)$$

and therefore

$$B_{2n} = \frac{2(-1)^{n+1}2n(2n-1)(2n-2)}{(2\pi)^n}\,\Gamma(2n-2)\varsigma(2n)$$

$$= \frac{2(-1)^{n+1}(2n)!\varsigma(2n)}{(2\pi)^n}$$

The above analysis is somewhat circular but it highlighted the following identity

$$f(\nu,s) = \int_0^\infty u^{\nu-1} Li_s(e^{-u})du = \sum_{k=1}^\infty \frac{1}{k^s} \int_0^\infty u^{\nu-1}e^{-ku}du$$

$$= \Gamma(\nu)\varsigma(s+\nu)$$

Reversing the roles of $\nu$ and $s$ we have

$$f(s,\nu) = \int_0^\infty u^{s-1} Li_\nu(e^{-u})du = \Gamma(s)\varsigma(s+\nu)$$

and it is easily seen that

$$f(s,\nu) = \frac{\Gamma(s)}{\Gamma(\nu)} f(\nu,s)$$

Upon differentiation with respect to $s$ we have

$$\frac{\partial}{\partial s} f(s,\nu) = \int_0^\infty u^{s-1} \log u\, Li_\nu(e^{-u})du = \Gamma'(s)\varsigma(s+\nu)+\Gamma(s)\varsigma'(s+\nu)$$

where, for ease of notation, $\varsigma'(s+\nu) = \dfrac{\partial}{\partial s}\varsigma(s+\nu)$. Alternatively, we have

$$\frac{\partial}{\partial s} f(s,\nu) = \frac{\partial}{\partial s}\left[\frac{\Gamma(s)}{\Gamma(\nu)} f(\nu,s)\right]$$

$$= \frac{\Gamma'(s)}{\Gamma(\nu)} f(\nu,s) + \frac{\Gamma(s)}{\Gamma(\nu)}\frac{\partial}{\partial s} f(\nu,s)$$

$$= \frac{\Gamma'(s)}{\Gamma(\nu)} f(\nu,s) + \frac{\Gamma(s)}{\Gamma(\nu)}\frac{\partial}{\partial s}\int_0^\infty u^{\nu-1} Li_s(e^{-u})du$$

$$= \frac{\Gamma'(s)}{\Gamma(\nu)} f(\nu,s) + \frac{\Gamma(s)}{\Gamma(\nu)} \int_0^\infty u^{\nu-1} \frac{\partial}{\partial s} Li_s(e^{-u}) du$$

We have

$$\frac{\partial}{\partial s} Li_s(e^{-u}) = \frac{\partial}{\partial s} \sum_{k=1}^\infty \frac{e^{-ku}}{k^s} = -\sum_{k=1}^\infty \frac{e^{-ku} \log k}{k^s}$$

and therefore

(4.4.184) $$\int_0^\infty u^{\nu-1} \frac{\partial}{\partial s} Li_s(e^{-u}) du = -\int_0^\infty \sum_{k=1}^\infty \frac{u^{\nu-1} e^{-ku} \log k}{k^s} du$$

$$= -\sum_{k=1}^\infty \frac{\Gamma(\nu) \log k}{k^\nu k^s} = -\Gamma(\nu) \sum_{k=1}^\infty \frac{\log k}{k^{\nu+s}}$$

$$= \Gamma(\nu) \varsigma'(s+\nu)$$

Accordingly we have

$$\frac{\partial}{\partial s} f(s,\nu) = \frac{\Gamma'(s)}{\Gamma(\nu)} f(\nu,s) + \frac{\Gamma(s)}{\Gamma(\nu)} \Gamma(\nu) \varsigma'(s+\nu)$$

$$= \frac{\Gamma'(s)}{\Gamma(\nu)} f(\nu,s) + \Gamma(s) \varsigma'(s+\nu)$$

Therefore, as before, we have

$$\frac{\partial}{\partial s} f(s,\nu) = \Gamma'(s) \varsigma(s+\nu) + \Gamma(s) \varsigma'(s+\nu)$$

From the above we see that we have an integral expression involving the derivative of the Riemann zeta function

(4.4.185) $$\int_0^\infty u^{s-1} \log u \, Li_\nu(e^{-u}) du = \Gamma'(s) \varsigma(s+\nu) + \Gamma(s) \varsigma'(s+\nu)$$

$$= \frac{\partial}{\partial s} \big[ \Gamma(s) \varsigma(s+\nu) \big]$$

For example, with $\nu = s = 3/2$, we have

$$\int_0^\infty \sqrt{u} \, \log u \, Li_{3/2}(e^{-u}) \, du = \Gamma'(3/2)\varsigma(3) + \Gamma(3/2)\varsigma'(3)$$

With the substitution $u \to ku$ in the definition of the gamma function $\Gamma(x) = \int_0^\infty u^{x-1} e^{-u} du$

we have

$$\Gamma(x) = k \int_0^\infty (ku)^{x-1} e^{-ku} du = k^x \int_0^\infty u^{x-1} e^{-ku} du$$

Hence

(4.4.186) $$\Gamma'(x) = k \int_0^\infty (ku)^{x-1} e^{-ku} \log(ku) \, du$$

$$\Gamma'(1) = k \int_0^\infty \log(ku) \, e^{-ku} du = k \int_0^\infty \log k \, e^{-ku} du + k \int_0^\infty \log u \, e^{-ku} du$$

Therefore we obtain

$$\int_0^\infty \log u \, e^{-ku} du = -\frac{1}{k}\gamma - \log k \int_0^\infty e^{-ku} du$$

$$= -\frac{1}{k}\gamma - \frac{\log k}{k}$$

If we let $k = \log t$ with $t > 1$ in the above equation we obtain

$$\int_0^\infty \log u \, e^{-ku} du = \int_0^\infty \frac{\log u}{t^u} du = -\frac{1}{\log t}(\gamma + \log(\log t))$$

and we learn from the Wolfram Integrator that

$$\int \frac{\log u}{t^u} du = \frac{Ei(-u \log t) - t^{-u} \log t}{\log t}$$

where $Ei(x) = -\int_{-x}^\infty \frac{e^{-t}}{t} dt$ is the exponential integral.

We can therefore make the following summation

$$\sum_{k=1}^{\infty} \frac{1}{k^s} \int_0^{\infty} \log u \, e^{-ku} du = -\gamma \varsigma(s+1) - \sum_{k=1}^{\infty} \frac{\log k}{k^{s+1}}$$

which implies that

(4.4.187)
$$\int_0^{\infty} \log u \, Li_s(e^{-u}) \, du = \varsigma'(s+1) - \gamma \varsigma(s+1)$$

or alternatively

(4.4.187a)
$$\int_0^1 \frac{\log \log(1/x) Li_s(x)}{x} \, dx = \varsigma'(s+1) - \gamma \varsigma(s+1)$$

A further proof is given in (C.70).

Differentiating the above with respect to $s$, and using (4.4.184) we obtain

$$-\int_0^{\infty} \log u \sum_{k=1}^{\infty} \frac{e^{-ku} \log k}{k^s} du = \varsigma''(s+1) - \gamma \varsigma'(s+1)$$

$$-\sum_{k=1}^{\infty} \int_0^{\infty} \log u \frac{e^{-ku} \log k}{k^s} du = \varsigma''(s+1) - \gamma \varsigma'(s+1)$$

Therefore we have

$$\sum_{k=1}^{\infty} (\gamma + \log k) \frac{\log k}{k^{s+1}} = \varsigma''(s+1) - \gamma \varsigma'(s+1)$$

as one would expect from the definition of the zeta function (even dead ends are sometimes worth reporting!).

We also have

$$\Gamma^{(m)}(x) = k \int_0^{\infty} (ku)^{x-1} \log^m(ku) \, e^{-ku} du$$

$$= k^x \int_0^{\infty} u^{x-1} \left( \log k + \log u \right)^m e^{-ku} du$$

$$= k^x \int_0^\infty u^{x-1} e^{-ku} \sum_{j=0}^m \binom{m}{j} \log^{m-j} k \, \log^j u \, du$$

In particular we have

$$\Gamma''(x) = k^x \int_0^\infty u^{x-1} (\log^2 k + 2\log k \log u + \log^2 u) e^{-ku} du$$

$$= k^x \log^2 k \int_0^\infty u^{x-1} e^{-ku} du + 2k^x \log k \int_0^\infty u^{x-1} e^{-ku} \log u \, du + k^x \int_0^\infty u^{x-1} e^{-ku} \log^2 u \, du$$

$$= \Gamma(x) \log^2 k + 2k^x \log k \int_0^\infty u^{x-1} e^{-ku} \log u \, du + \Gamma(x) \Big[ (\psi(x) - \log k)^2 + \varsigma(2, x) \Big]$$

where we have used the following integral which is contained in [59] (which was originally derived by Kölbig in 1983)

(4.4.188)  $$\int_0^\infty u^{x-1} e^{-ku} \log^2 u \, du = \frac{\Gamma(x)}{k^x} \Big[ (\psi(x) - \log k)^2 + \varsigma(2, x) \Big]$$

Therefore for $k \neq 1$ we have

(4.4.188a)
$$\int_0^\infty u^{x-1} e^{-ku} \log u \, du = \frac{1}{2k^x} \frac{\Gamma''(x)}{\log k} - \frac{1}{2k^x} \Gamma(x) \log k - \frac{1}{2k^x} \frac{\Gamma(x)}{\log k} \Big[ (\psi(x) - \log k)^2 + \varsigma(2, x) \Big]$$

From the series definition of the polylogarithm function it is easily seen that

$$\sum_{k=1}^\infty \frac{1}{k^s} \int_0^\infty u^{x-1} e^{-ku} \log u \, du = \int_0^\infty u^{x-1} Li_s(e^{-u}) \log u \, du$$

Starting the summation at $k = 2$ gives us

$$\sum_{k=2}^\infty \frac{1}{k^s} \int_0^\infty u^{x-1} e^{-ku} \log u \, du = \int_0^\infty u^{x-1} Li_s(e^{-u}) \log u \, du - \int_0^\infty u^{x-1} e^{-u} \log u \, du$$

Alternatively we get from (4.4.188a)

$$\sum_{k=2}^\infty \frac{1}{k^s} \int_0^\infty u^{x-1} e^{-ku} \log u \, du =$$

$$= \frac{\Gamma''(x)}{2} \sum_{k=2}^{\infty} \frac{1}{k^{x+s} \log k} - \frac{\Gamma(x)}{2} \sum_{k=2}^{\infty} \frac{\log k}{k^{x+s}} - \frac{1}{2} \sum_{k=2}^{\infty} \frac{\Gamma(x)}{k^{x+s} \log k} \Big[ (\psi(x) - \log k)^2 + \varsigma(2,x) \Big]$$

(4.4.188b)

$$= \frac{\Gamma''(x)}{2} \sum_{k=2}^{\infty} \frac{1}{k^{x+s} \log k} - \frac{\Gamma(x)}{2} \sum_{k=2}^{\infty} \frac{\log k}{k^{x+s}} - \frac{1}{2} \Gamma(x) \Big[ \psi^2(x) + \varsigma(2,x) \Big] \sum_{k=2}^{\infty} \frac{1}{k^{x+s} \log k}$$

$$+ \Gamma(x)\psi(x)\varsigma(x+s) - \frac{1}{2}\Gamma(x) \sum_{k=2}^{\infty} \frac{\log k}{k^{x+s}}$$

$$= \frac{1}{2} \Big\{ \Gamma''(x) - \Gamma(x) \Big[ \psi^2(x) + \varsigma(2,x) \Big] \Big\} \sum_{k=2}^{\infty} \frac{1}{k^{x+s} \log k} + \Gamma(x)\psi(x)\varsigma(x+s) - \Gamma(x) \sum_{k=2}^{\infty} \frac{\log k}{k^{x+s}}$$

We have the well-known formula (see for example [119] and also Appendix E)

(4.4.189) $$\psi'(x) = \frac{d}{dx} \frac{\Gamma'(x)}{\Gamma(x)} = \frac{\Gamma(x)\Gamma''(x) - (\Gamma'(x))^2}{\Gamma^2(x)} = \varsigma(2,x)$$

or alternatively

(4.4.190) $$\Gamma''(x) = \Gamma(x) \Big[ \psi^2(x) + \varsigma(2,x) \Big]$$

Consequently, it is seen that the series in (4.4.188b), with a troublesome factor of $\log k$ in the denominator, therefore vanishes and we obtain

$$\sum_{k=2}^{\infty} \frac{1}{k^s} \int_0^{\infty} u^{x-1} e^{-ku} \log u \, du = \Gamma(x)\psi(x)\varsigma(x+s) + \Gamma(x)\varsigma'(x+s)$$

$$= \Gamma'(x)\varsigma(x+s) + \Gamma(x)\varsigma'(x+s)$$

Hence, since $\int_0^{\infty} u^{x-1} e^{-u} \log u \, du = \Gamma'(x)$ we obtain for $s > 1$

$$\int_0^{\infty} u^{x-1} Li_s(e^{-u}) \log u \, du = \Gamma'(x)\varsigma(x+s) + \Gamma(x)\varsigma'(x+s)$$

and we have previously met this in (4.4.185).

I must have had a temporary mental aberration when I did this work because I simply failed to notice at the time that it was in fact very easy to carry out this analysis much more directly. From (4.4.186) we have

$$\Gamma'(x) = k \int_0^\infty (ku)^{x-1} e^{-ku} \log(ku)\, du$$

$$= k^x \log k \int_0^\infty u^{x-1} e^{-ku} du + k^x \int_0^\infty u^{x-1} e^{-ku} \log u\, du$$

Therefore we have

$$k^x \int_0^\infty u^{x-1} e^{-ku} \log u\, du = \Gamma'(x) - \Gamma(x) \log k$$

As noted previously we have

$$\Gamma''(x) = k^x \log^2 k \int_0^\infty u^{x-1} e^{-ku} du + 2k^x \log k \int_0^\infty u^{x-1} e^{-ku} \log u\, du + k^x \int_0^\infty u^{x-1} e^{-ku} \log^2 u\, du$$

$$= \Gamma(x) \log^2 k + 2 \log k \left[ \Gamma'(x) - \Gamma(x) \log k \right] + k^x \int_0^\infty u^{x-1} e^{-ku} \log^2 u\, du$$

Therefore, in a considerably more straightforward manner, we obtain

$$(4.4.191) \qquad k^x \int_0^\infty u^{x-1} e^{-ku} \log^2 u\, du = \Gamma''(x) - \Gamma(x) \log^2 k - 2 \log k \left[ \Gamma'(x) - \Gamma(x) \log k \right]$$

and, as demonstrated above, this is equivalent to Kölbig's integral

$$(4.4.192) \qquad \int_0^\infty u^{x-1} e^{-ku} \log^2 u\, du = \frac{\Gamma(x)}{k^x} \left[ (\psi(x) - \log k)^2 + \varsigma(2, x) \right]$$

In turn, we may differentiate (4.4.192) with respect to $x$ to obtain (and noting that $\psi'(x) = \varsigma(2, x)$ )

$$\int_0^\infty u^{x-1} e^{-ku} \log^3 u\, du = \frac{\Gamma(x)}{k^x} \left[ 2(\psi(x) - \log k)\psi'(x) + \psi''(x) \right]$$

$$+ \frac{\Gamma'(x) - \Gamma(x) \log k}{k^x} \left[ (\psi(x) - \log k)^2 + \psi'(x) \right]$$

Hence we have

$$(4.4.193) \qquad k^x \int_0^\infty u^{x-1} e^{-ku} \log^3 u \, du = \Gamma(x) \left[ 2(\psi(x) - \log k) \psi'(x) + \psi''(x) \right]$$

$$+ \left[ \Gamma'(x) - \Gamma(x) \log k \right] \left[ (\psi(x) - \log k)^2 + \psi'(x) \right]$$

It is clear that similar results may be obtained for $\int_0^\infty u^{x-1} e^{-ku} \log^m u \, du$ by successive differentiations.

With $k = 1$ we have

$$\int_0^\infty u^{x-1} e^{-u} \log^3 u \, du = \Gamma(x) \left[ 2\psi(x)\psi'(x) + \psi''(x) \right] + \Gamma'(x) \left[ \psi^2(x) + \psi'(x) \right]$$

and with $x = 1$ this becomes

$$(4.4.194) \qquad \Gamma^{(3)}(1) = \int_0^\infty e^{-u} \log^3 u \, du = -\gamma^3 - \gamma \frac{\pi^2}{2} - 2\varsigma(3)$$

(where we have used $\psi'(x) = \sum_{j=0}^\infty \frac{1}{(j+x)^2}$ and $\psi''(x) = -2\sum_{j=0}^\infty \frac{1}{(j+x)^3}$ and accordingly

$\psi'(1) = \varsigma(2)$ and $\psi''(1) = -2\varsigma(3)$). This is in agreement with Levenson's result in [99]. See also Choi's recent paper [45acii].

Since $\Gamma(x) = k^x \int_0^\infty u^{x-1} e^{-ku} \, du$ we have

$$\frac{d^n}{dx^n} \left[ \frac{\Gamma(x)}{k^x} \right] = \int_0^\infty u^{x-1} e^{-ku} \log^n u \, du$$

Using Leibniz's rule we have

$$\frac{d^n}{dx^n} \left[ \Gamma(x) k^{-x} \right] = \sum_{j=0}^n \binom{n}{j} \Gamma^{(j)}(x) \frac{d^{n-j}}{dx^{n-j}} (k^{-x})$$

and using

$$\frac{d^{n-j}}{dx^{n-j}} (k^{-x}) = (-1)^{n-j} \frac{\log^{n-j} k}{k^x}$$

we therefore get

$$\frac{d^n}{dx^n}\left[\Gamma(x)k^{-x}\right] = \sum_{j=0}^{n}\binom{n}{j}\Gamma^{(j)}(x)\ (-1)^{n-j}\frac{\log^{n-j}k}{k^x}$$

Hence we obtain a relatively closed form expression for the integral

(4.4.195)  $$\int_0^\infty u^{x-1}e^{-ku}\log^n u\,du = \frac{1}{k^x}\sum_{j=0}^{n}\binom{n}{j}\Gamma^{(j)}(x)(-\log k)^{n-j}$$

$$= \frac{1}{k^x}\sum_{j=0}^{n}(-1)^{n-j}\binom{n}{j}\Gamma^{(j)}(x)\log^{n-j}k$$

This is a slight improvement over Levenson's result in [99] because (4.4.195) gives an explicit formula whereas Levenson only gives a recursion relation. With $k=1$ we have

$$\frac{d^n}{dx^n}\int_0^\infty u^{x-1}e^{-u}\,du = \int_0^\infty u^{x-1}e^{-u}\log^n u\,du = \sum_{j=0}^{n}(-1)^{n-j}\binom{n}{j}\Gamma^{(j)}(x)\log^{n-j}1 = \Gamma^{(n)}(x)$$

We have the geometric series

$$\sum_{k=1}^{\infty}e^{-ku} = \frac{1}{e^u-1}$$

Hence for $n \geq 0$ we have

$$\sum_{k=1}^{\infty}\int_0^\infty u^{x-1}e^{-ku}\log^n u\,du = \int_0^\infty \frac{u^{x-1}\log^n u}{e^u-1}\,du$$

$$= \sum_{k=1}^{\infty}\frac{1}{k^x}\sum_{j=0}^{n}(-1)^{n-j}\binom{n}{j}\Gamma^{(j)}(x)\log^{n-j}k$$

$$= \sum_{j=0}^{n}(-1)^{n-j}\binom{n}{j}\Gamma^{(j)}(x)\sum_{k=1}^{\infty}\frac{\log^{n-j}k}{k^x}$$

$$= \sum_{j=0}^{n}(-1)^{n-j}\binom{n}{j}\Gamma^{(j)}(x)(-1)^{n-j}\varsigma^{(n-j)}(x)$$

$$= \sum_{j=0}^{n} \binom{n}{j} \Gamma^{(j)}(x) \varsigma^{(n-j)}(x)$$

$$= \frac{d^n}{dx^n} \big[ \Gamma(x) \varsigma(x) \big]$$

where in the final part we have used the Leibniz rule for the nth derivative of a product.

We have therefore shown that

$$\frac{d^n}{dx^n} \left[ \int_0^\infty \frac{u^{x-1} \log^n u}{e^u - 1} \, du \right] = \int_0^\infty \frac{u^{x-1} \log^n u}{e^u - 1} \, du = \frac{d^n}{dx^n} \big[ \Gamma(x) \varsigma(x) \big]$$

and hence we have an alternative proof of formula (4.4.38)

$$\int_0^\infty \frac{u^{x-1}}{e^u - 1} \, du = \Gamma(x) \varsigma(x)$$

Letting $x = 1$ in (4.4.188) we see that

$$\int_0^\infty e^{-ku} \log^2 u \, du = \frac{1}{k} \Big[ (\gamma + \log k)^2 + \varsigma(2) \Big]$$

and making the summation we obtain

$$\sum_{k=1}^\infty \frac{1}{k} \int_0^\infty e^{-ku} \log^2 u \, du = \sum_{k=1}^\infty \frac{1}{k^2} \Big[ (\gamma + \log k)^2 + \varsigma(2) \Big]$$

This then gives us

$$-\int_0^\infty \log[1 - e^{-u}] \log^2 u \, du = \Big[ \gamma^2 + \varsigma(2) \Big] \varsigma(2) - 2\gamma \varsigma'(2) + \varsigma''(2)$$

With the substitution $x = e^{-u}$ this becomes

(4.4.195a) $\qquad -\int_0^1 \log[1 - x] \left[ \log \log \left( \frac{1}{x} \right) \right]^2 \frac{dx}{x} = \Big[ \gamma^2 + \varsigma(2) \Big] \varsigma(2) - 2\gamma \varsigma'(2) + \varsigma''(2)$

Letting $x = 1$ in (4.4.188) we get

$$\int_0^\infty e^{-ku} \log^2 u \, du = \frac{1}{k}\Big[ (\gamma + \log k)^2 + \varsigma(2) \Big]$$

and summation gives us

$$\sum_{k=1}^\infty \frac{1}{k^s} \int_0^\infty e^{-ku} \log^2 u \, du = \sum_{k=1}^\infty \frac{1}{k^{s+1}} \Big[ (\gamma + \log k)^2 + \varsigma(2) \Big]$$

Hence we obtain

$$\int_0^\infty Li_s(e^{-u}) \log^2 u \, du = [\gamma^2 + \varsigma(2)]\varsigma(s+1) - 2\gamma\varsigma'(s+1) + \varsigma''(s+1)$$

### SOME INTEGRALS INVOLVING $\log\Gamma(x)$

**(xix)** In his 1998 paper "Polygamma functions of Negative Order" [4] Adamchik gave a proof of Gosper's formula [71]

(4.4.196) $$\int_0^q \log\Gamma(x)dx = \frac{q(1-q)}{2} + \frac{q}{2}\log 2\pi - \varsigma'(-1) + \varsigma'(-1,q)$$

where $\varsigma(s,a)$ is the Hurwitz zeta function defined by

(4.4.197) $$\varsigma(s,a) = \sum_{n=0}^\infty \frac{1}{(n+a)^s} \qquad (\operatorname{Re}(s) > 1; \operatorname{Re}(a) > 0)$$

and

(4.4.198) $$\varsigma'(s,a) = \frac{\partial}{\partial s}\varsigma(s,a)$$

If $q=1$, then we have $\int_0^1 \log\Gamma(x)dx = \frac{1}{2}\log 2\pi$ (see (C.43b) of Volume VI)

The functional equation for the Riemann zeta function [81] was employed in the proof of (4.4.196).

(4.4.199) $$\varsigma(1-s) = 2(2\pi)^{-s}\Gamma(s)\cos(\pi s/2)\varsigma(s)$$

For $q=1/2$ in (4.4.196) we have

$$(4.4.200) \qquad \int_0^{1/2} \log \Gamma(x) dx = \frac{1}{8} + \frac{1}{4} \log 2\pi - \varsigma'(-1) + \varsigma'\left(-1, \frac{1}{2}\right)$$

A couple of pages later, in the same volume of the Journal of Computation and Applied Mathematics, in a paper entitled "Derivatives of the Hurwitz Zeta Function for Rational Arguments", Miller and Adamchik [103] show that

$$(4.4.201) \qquad \varsigma'\left(-2k+1, \frac{1}{2}\right) = -\frac{B_{2k} \log 2}{4^k k} - \frac{\left(2^{2k-1}-1\right) \varsigma'(-2k+1)}{2^{2k-1}}$$

so that for $k = 1$ we have

$$(4.4.202) \qquad \varsigma'\left(-1, \frac{1}{2}\right) = -\frac{B_2 \log 2}{4} - \frac{\varsigma'(-1)}{2}$$

$$(4.4.203) \qquad\qquad = -\frac{\log 2}{24} - \frac{\varsigma'(-1)}{2}$$

where we have substituted the numerical value of the Bernoulli number (see (A.8) in Appendix A). Combining (4.4.200) and (4.4.203) we obtain

$$(4.4.204) \qquad \int_0^{1/2} \log \Gamma(x) dx = \frac{1}{8} + \frac{1}{4} \log 2\pi - \frac{\log 2}{24} - \frac{3}{2} \varsigma'(-1)$$

Later in 2002, in their paper "On some integrals involving the Hurwitz zeta function. Part 2", Espinosa and Moll [60] note that

$$(4.4.205) \qquad \varsigma'(-1) = \frac{1}{12}(1 - \gamma - \log 2\pi) - \frac{1}{2\pi^2} \sum_{n=1}^{\infty} \frac{\log n}{n^2}$$

$$\qquad\qquad = \frac{1}{12}(1 - \gamma - \log 2\pi) + \frac{1}{2\pi^2} \varsigma'(2)$$

The result (4.4.205) may be obtained by the logarithmic differentiation of Riemann's functional equation (4.4.199) to obtain

$$(4.4.205a) \qquad \log \varsigma(1-s) = \log 2 - s \log(2\pi) + \log \Gamma(s) + \log \cos\left(\frac{\pi s}{2}\right) + \log \varsigma(s)$$

$$(4.4.205b) \qquad -\frac{\varsigma'(1-s)}{\varsigma(1-s)} = -\log(2\pi) + \frac{\Gamma'(s)}{\Gamma(s)} - \frac{\pi}{2} \tan\left(\frac{\pi s}{2}\right) + \frac{\varsigma'(s)}{\varsigma(s)}$$

We now let $s = 2$ and obtain (4.4.205): in the process we use the following

(4.4.205c) $$\varsigma(-1) = -\frac{1}{12}$$

(4.4.205d) $$\psi(x+1) = \psi(x) + \frac{1}{x} \quad \Rightarrow \quad \psi(2) = 1 - \gamma$$

Equation (4.4.205c) was obtained with the use of (3.11b). Riemann's functional equation also enables us to easily see the trivial zeros of the zeta function by letting $s = 2n + 1$

(4.4.205e) $$\varsigma(-2n) = 0 \text{ because } \cos\frac{(2n+1)\pi}{2} = 0$$

Hence we have

(4.4.206) $$\int_0^{\frac{1}{2}} \log\Gamma(x)\,dx = \frac{1}{8}\gamma + \frac{3}{8}\log 2\pi - \frac{\log 2}{24} + \frac{3}{4\pi^2}\sum_{n=1}^{\infty}\frac{\log n}{n^2}$$

and this is in agreement with equation (6.16) in the earlier paper by Espinosa and Moll [59]. See also (4.4.212a).

(4.4.207) $$\int_0^{\frac{1}{2}} \log\Gamma(x+1)\,dx = \frac{1}{8}\gamma + \frac{3}{4}\log\sqrt{2\pi} - \frac{13\log 2}{24} - \frac{3\varsigma'(2)}{4\pi^2} - \frac{1}{2}$$

Using the formulae

(4.4.208) $$\log\Gamma(x+1) = \log x + \log\Gamma(x)$$

(4.4.209) $$\varsigma'(2) = -\sum_{n=1}^{\infty}\frac{\log n}{n^2}$$

it is easily seen that (4.4.206) is equivalent to (4.4.207). Formula (4.4.207) was originally obtained by Gosper [71].

In 1985 Berndt [20] gave an elementary proof of the Fourier series expansion for $\log\Gamma(x)$

(4.4.210) $$\log\Gamma(x) = \frac{1}{2}\log\pi - \frac{1}{2}\log\sin\pi x + \sum_{n=1}^{\infty}\frac{(\gamma + \log 2\pi n)\sin 2\pi nx}{\pi n} \quad (0 < x < 1)$$

(this formula was originally derived by Kummer in 1847 [94]). Reference to (7.8) confirms that (4.4.210) is properly described as a Fourier series expansion for $\log \Gamma(x)$. Using (7.5) we may also write Kummer's formula as (cf. Nielsen [104a, p.79])

(4.4.210a)     $$\log \Gamma(x) = \frac{1}{2} \log \frac{\pi}{\sin \pi x} + \left(\frac{1}{2} - x\right)(\gamma + \log 2\pi) + \frac{1}{\pi} \sum_{n=1}^{\infty} \frac{\log n}{n} \sin 2\pi n x$$

A proof of Kummer's formula is given in Appendix E in Volume VI.

Assuming that it is valid to integrate (4.4.210) over the interval [0, 1/2] (both sides of equation (4.4.210) $\to \infty$ as $x \to 0$), we obtain

$$\int_0^{\frac{1}{2}} \log \Gamma(x)\, dx =$$

$$\frac{1}{4} \log \pi + \frac{1}{4} \log 2 + \frac{\gamma}{2\pi^2}\{\varsigma_a(2) + \varsigma(2)\} + \frac{1}{2\pi^2} \sum_{n=1}^{\infty} \frac{(-1)^{n+1} \log 2\pi n}{n^2} + \frac{1}{2\pi^2} \sum_{n=1}^{\infty} \frac{\log 2\pi n}{n^2}$$

(4.4.211)     $$= \frac{1}{4} \log 2\pi + \frac{1}{4} \log 2 + \frac{1}{8}\gamma + \frac{\log 2\pi}{8} + \frac{1}{2\pi^2} \sum_{n=1}^{\infty} \frac{(-1)^{n+1} \log n}{n^2} + \frac{1}{2\pi^2} \sum_{n=1}^{\infty} \frac{\log n}{n^2}$$

From (1.1) we have

$$\varsigma(s)\left(1 - 2^{1-s}\right) = \sum_{n=1}^{\infty} \frac{(-1)^{n+1}}{n^s}$$

and upon differentiation we obtain

$$\varsigma'(s)\left(1 - 2^{1-s}\right) + 2^{1-s} \varsigma(s) \log 2 = -\sum_{n=1}^{\infty} \frac{(-1)^{n+1} \log n}{n^s}$$

which becomes with $s = 2$

(4.4.212)     $$\frac{1}{2}\varsigma'(2) + \frac{1}{2}\varsigma(2) \log 2 = -\sum_{n=1}^{\infty} \frac{(-1)^{n+1} \log n}{n^2} = \varsigma_a'(2)$$

Combining (4.4.209), (4.4.211) and (4.4.212) we obtain (4.4.206)

(4.4.212a)
$$\int_0^{1/2} \log \Gamma(x)\, dx = \frac{1}{8}\gamma + \frac{3}{8}\log 2\pi - \frac{\log 2}{24} + \frac{3}{4\pi^2}\sum_{n=1}^{\infty}\frac{\log n}{n^2}$$

This method can also be employed to evaluate integrals of the form $\displaystyle\int_0^{p/q} \log \Gamma(x)\, dx$.

In [59], using the Hurwitz zeta function, Espinosa and Moll showed that

(4.4.213a)
$$\int_0^1 B_{2n}(x)\log \Gamma(x)\, dx = (-1)^{n+1}\frac{(2n)!\,\varsigma(2n+1)}{2(2\pi)^{2n}} = -\varsigma'(-2n)$$

(4.4.213b)
$$\int_0^1 B_{2n-1}(x)\log \Gamma(x)\, dx = \frac{B_{2n}}{2n}\left[\frac{\varsigma'(2n)}{\varsigma(2n)} - A\right]$$

where

$$A = 2\log\sqrt{2\pi} + \gamma = -2\frac{d}{dx}\big(\varsigma(x)\Gamma(1-x)\big)\Big|_{x=0}$$

(the identity for $A$ follows directly from (E.16), (E.16) and (F.6)).
In particular, for $n=1$ we have

(4.4.213c)
$$\int_0^1 \left(x^2 - x + \frac{1}{6}\right)\log \Gamma(x)\, dx = \frac{\varsigma(3)}{4\pi^2}$$

(4.4.213d)
$$\int_0^1 \left(x - \frac{1}{2}\right)\log \Gamma(x)\, dx = \frac{1}{12}\left(\frac{6\varsigma'(2)}{\pi^2} - 2\log\sqrt{2\pi} - \gamma\right)$$

A further proof of (4.4.213d) is shown in Appendix E of Volume VI.

Using Kummer's identity (4.4.210) we have

$$\int_0^1\left(x-\frac{1}{2}\right)\log\Gamma(x)\,dx = \frac{1}{2}\log\pi\int_0^1\left(x-\frac{1}{2}\right)dx - \frac{1}{2}\int_0^1\left(x-\frac{1}{2}\right)\log\sin\pi x\,dx + \frac{(\gamma+\log 2\pi)}{\pi}\sum_{n=1}^{\infty}\int_0^1\frac{\left(x-\frac{1}{2}\right)\sin 2\pi n x}{n}\,dx$$

$$+\frac{1}{\pi}\sum_{n=1}^{\infty}\int_0^1\frac{\log n\left(x-\frac{1}{2}\right)\sin 2\pi n x\,dx}{n}$$

We have from (3.2)

$$\int_0^{\pi/2} \log \sin x \, dx = -\frac{\pi}{2} \log 2$$

and it is clear that

$$\int_0^{\pi} \log \sin x \, dx = \int_0^{\pi/2} \log \sin x \, dx + \int_{\pi/2}^{\pi} \log \sin x \, dx$$

With the substitution $x \rightarrow \pi - x$ in the last integral, we deduce that

$$\int_0^{\pi} \log \sin x \, dx = 2 \int_0^{\pi/2} \log \sin x \, dx = -\pi \log 2$$

Therefore we have

$$\int_0^1 \log \sin \pi x \, dx = \frac{1}{\pi} \int_0^{\pi} \log \sin y \, dy = -\log 2$$

As shown below, we also have

$$\int_0^1 x \log \sin \pi x \, dx = \frac{1}{\pi^2} \int_0^{\pi} y \log \sin y \, dy = -\frac{1}{2} \log 2$$

To prove this, we consider a different integral

$$\int_0^1 x^2 \log \sin \pi x \, dx = \frac{1}{\pi^3} \int_0^{\pi} y^2 \log \sin y \, dy$$

and using the substitution $y = \pi - t$ we have

$$\int_0^{\pi} y^2 \log \sin y \, dy = \int_0^{\pi} (\pi - t)^2 \log \sin t \, dt$$

The last integral therefore proves that $\int_0^{\pi} (\pi^2 - 2\pi t) \log \sin t \, dt = 0$ and hence we have

$\int_0^{\pi} t \log \sin t \, dt = -\frac{\pi^2}{2} \log 2$. Therefore we get $\int_0^1 x \log \sin \pi x \, dx = -\frac{1}{2} \log 2$ and accordingly we have

(4.4.213e) $$\int_0^1 \left( x - \frac{1}{2} \right) \log \sin \pi x \, dx = 0$$

In fact (4.4.213d) can be shown much more directly by using the substitution $y = x - 1/2$ and then noting that the integrand of the resulting integral is an odd function.

By integration by parts we easily obtain

$$\int_0^1 \left( x - \frac{1}{2} \right) \sin 2\pi n x \, dx = \frac{\cos 2\pi n x}{4\pi n} - \frac{x \cos 2\pi n x}{2\pi n} + \frac{\sin 2\pi n x}{4\pi^2 n^2} \bigg|_0^1$$

$$= -\frac{1}{2\pi n}$$

and we therefore get

$$\int_0^1 \left( x - \frac{1}{2} \right) \log \Gamma(x) \, dx = -\frac{(\gamma + \log 2\pi)}{2\pi^2} \sum_{n=1}^\infty \frac{1}{n^2} - \frac{1}{2\pi^2} \sum_{n=1}^\infty \frac{\log n}{n^2}$$

$$= -\frac{(\gamma + \log 2\pi)\varsigma(2)}{2\pi^2} - \frac{1}{2\pi^2} \sum_{n=1}^\infty \frac{\log n}{n^2}$$

$$= -\frac{(\gamma + \log 2\pi)}{12} + \frac{\varsigma'(2)}{2\pi^2}$$

(and this is in agreement with the Espinosa and Moll result (4.4.213d)).

More generally, using Kummer's Fourier series expansion (4.4.210) for $\log \Gamma(x)$ we obtain

(4.4.213f)

$$\int_0^1 B_{2n}(x) \log \Gamma(x) \, dx = \frac{1}{2} \log \pi \int_0^1 B_{2n}(x) \, dx - \frac{1}{2} \int_0^1 B_{2n}(x) \log \sin \pi x \, dx + \sum_{n=1}^\infty \int_0^1 \frac{(\gamma + \log 2\pi n) \sin 2\pi n x}{\pi n} B_{2n}(x) \, dx$$

We have the elementary integral (where we use integration by parts and the identity (A.14b) for the Bernoulli polynomials)

$$\int_0^1 B_{2n}(x) \, dx = \frac{B_{2n+1}(x)}{2n+1} \bigg|_0^1 = 0$$

In addition, using a simple substitution we see that

$$\int_0^1 B_{2n}(x)\sin 2\pi nx\, dx = \int_0^1 B_{2n}(1-y)\sin 2\pi n(1-y)\, dy$$

$$= -\int_0^1 B_{2n}(y)\sin 2\pi ny\, dy$$

and hence we get

$$\int_0^1 B_{2n}(x)\sin 2\pi nx\, dx = 0$$

Accordingly, (4.4.213f) simplifies to

(4.4.213g) $\quad \int_0^1 B_{2n}(x)\log\Gamma(x)\, dx = -\dfrac{1}{2}\int_0^1 B_{2n}(x)\log\sin\pi x\, dx$

Using integration by parts we have

$$\int_0^1 B_{2n+1}(x)\cot\pi x\, dx = B_{2n+1}(x)\left.\frac{\log\sin\pi x}{\pi}\right|_0^1 - \frac{2n+1}{\pi}\int_0^1 B_{2n}(x)\log\sin\pi x\, dx$$

$$= -\frac{2n+1}{\pi}\int_0^1 B_{2n}(x)\log\sin\pi x\, dx$$

In (6.34) we will see that

$$\int_0^1 B_{2n}(x)\log\sin\pi x\, dx = (-1)^n\frac{(2n)!\,\varsigma(2n+1)}{(2\pi)^{2n}}$$

and hence we have obtained the Espinosa and Moll formula (4.4.213a) by an alternate route

$$\int_0^1 B_{2n}(x)\log\Gamma(x)\, dx = (-1)^{n+1}\frac{(2n)!\,\varsigma(2n+1)}{2(2\pi)^{2n}}$$

In [59] and [60] Espinosa and Moll show that

(4.4.214)

$$\int x^n \log \Gamma(x)\, dx = -\varsigma'(0) \frac{x^{n+1}}{n+1} + n! \sum_{k=1}^{n+1} (-1)^{k+1} \frac{x^{n+1-k}}{k!(n+1-k)!} \left[ \frac{A_{k+1}(x)}{k+1} - \frac{H_k}{k+1} B_{k+1}(x) \right]$$

where $H_k$ and $B_k(x)$ are the harmonic numbers and the Bernoulli polynomials respectively, and

(4.4.215a)
$$\frac{A_{k+1}(x)}{k+1} = \frac{\partial}{\partial z} \varsigma(z,x)\Big|_{z=-k}$$

(4.4.215b)
$$A_1(x) = \varsigma'(0,x) = \log \Gamma(x) + \varsigma'(0)$$

(4.4.215c)
$$\varsigma'(0) = -\log \sqrt{2\pi}$$

With $n = 1$ we have

(4.4.216)

$$\int_0^{\frac{1}{2}} x \log \Gamma(x)\, dx = -\frac{1}{8} \varsigma'(0) - \frac{1}{6} A_k(0) + \sum_{k=1}^{2} (-1)^{k+1} \frac{(1/2)^{2-k}}{k!(2-k)!} \left[ \frac{A_{k+1}(1/2)}{k+1} - \frac{H_k}{k+1} B_{k+1}(1/2) \right]$$

From [60] we have

(4.4.217)
$$A_k(1/2) = (-1)^{k-1} B_k 2^{1-k} \log 2 - \left(1 - 2^{1-k}\right) k \varsigma'(1-k)$$

and hence using (6.60) we get

(4.4.218a)
$$A_1(1/2) = \log \Gamma(1/2) + \varsigma'(0) = -\frac{1}{2} \log 2$$

(4.4.218b)
$$A_2(1/2) = -\frac{1}{12} \log 2 - \varsigma'(-1)$$

(4.4.218c)
$$A_3(1/2) = -\frac{9}{4} \varsigma'(-2) \quad (\text{since } B_{2n+1} = 0 \ \ \forall \ n \geq 1)$$

Using [59] we have (easily obtained by differentiating Riemann's functional equation (4.4.199) for $\varsigma(s)$)

(4.4.219)
$$\varsigma'(-2n) = \frac{(-1)^n (2n)! \varsigma(2n+1)}{2(2\pi)^{2n}}$$

and therefore we get

(4.4.220) $\qquad \varsigma'(-2) = -\dfrac{1}{4\pi^2}\varsigma(3) \quad \text{and} \quad A_3(1/2) = \dfrac{9\varsigma(3)}{16\pi^2}$

From Appendix A we have

(4.4.221) $\qquad B_2(x) = x^2 - x + \dfrac{1}{6}$

and therefore $B_2(1/2) = -\dfrac{1}{12}$: we also have $B_3(1/2) = 0$ because $B_{2n+1}(1/2) = 0$.

We therefore obtain from (4.4.216)

$$\int\limits_0^{1/2} x\log\Gamma(x)\,dx = \frac{1}{8}\log\sqrt{2\pi} + \frac{1}{2}\left[ -\frac{1}{24}\log 2 - \frac{1}{2}\varsigma'(-1) + \frac{1}{24}\right] - \frac{1}{2}\left[\frac{3\varsigma(3)}{16\pi^2}\right]$$

(4.4.222) $\qquad = \dfrac{1}{16}\log\pi + \dfrac{1}{24}\log 2 - \dfrac{1}{4}\varsigma'(-1) + \dfrac{1}{48} - \dfrac{7\varsigma(3)}{32\pi^2}$

Multiplying Kummer's identity (4.4.210) by $x$ and integrating, we have

$$\int\limits_0^{1/2} x\log\Gamma(x)\,dx = \frac{1}{16}\log\pi - \frac{1}{2}\int\limits_0^{1/2} x\log\sin\pi x\,dx + \frac{(\gamma+\log 2\pi)}{\pi}\sum_{n=1}^{\infty}\int\limits_0^{1/2}\frac{x\sin 2\pi nx}{n}\,dx + \frac{1}{\pi}\sum_{n=1}^{\infty}\int\limits_0^{1/2}\frac{(\log n)\,x\sin 2\pi nx}{n}\,dx$$

Since $\int\limits_0^{1/2}\dfrac{x\sin 2\pi nx\,dx}{n} = \dfrac{(-1)^{n+1}}{4\pi n^2}$ we have

$$\int\limits_0^{1/2} x\log\Gamma(x)\,dx = \frac{1}{16}\log\pi - \frac{1}{2}\int\limits_0^{1/2} x\log\sin\pi x\,dx + \frac{(\gamma+\log 2\pi)}{4\pi^2}\varsigma_a(2) + \frac{1}{4\pi^2}\sum_{n=1}^{\infty}\frac{(-1)^{n+1}\log n}{n^2}$$

Using (4.4.212) we have

$$\frac{1}{4\pi^2}\sum_{n=1}^{\infty}\frac{(-1)^{n+1}\log n}{n^2} = -\frac{1}{8\pi^2}\varsigma'(2) - \frac{1}{48}\log 2$$

and hence we get

$$\int\limits_0^{1/2} x \log \Gamma(x)\, dx = \frac{1}{16}\log \pi - \frac{1}{2}\int\limits_0^{1/2} x \log \sin \pi x\, dx + \frac{(\gamma + \log 2\pi)}{48} - \frac{1}{8\pi^2}\varsigma'(2) - \frac{1}{48}\log 2$$

$$= \frac{1}{12}\log \pi - \frac{1}{2}\int\limits_0^{1/2} x \log \sin \pi x\, dx + \frac{\gamma}{48} - \frac{1}{8\pi^2}\varsigma'(2)$$

We have using Riemann's integral identity (1.11) from Volume I

$$\int\limits_0^{1/2} x \log \sin \pi x\, dx = \frac{1}{\pi^2}\int\limits_0^{\pi/2} x \log \sin x\, dx = \frac{7}{16\pi^2}\varsigma(3) - \frac{1}{8}\log 2$$

and deduce

(4.4.223)     $$\int\limits_0^{1/2} x \log \Gamma(x)\, dx = \frac{1}{12}\log \pi - \frac{7}{32\pi^2}\varsigma(3) + \frac{1}{16}\log 2 + \frac{\gamma}{48} - \frac{1}{8\pi^2}\varsigma'(2)$$

We now equate this with (4.4.222) as shown below

$$\int\limits_0^{1/2} x \log \Gamma(x)\, dx = \frac{1}{16}\log \pi + \frac{1}{24}\log 2 - \frac{1}{4}\varsigma'(-1) + \frac{1}{48} - \frac{7\varsigma(3)}{32\pi^2}$$

and obtain

(4.4.224)     $$\frac{1}{12} - \varsigma'(-1) = \frac{1}{12}\gamma + \frac{1}{12}\log 2\pi - \frac{1}{2}\frac{\varsigma'(2)}{\pi^2}$$

We also have the identity from Adamchik's paper [2a]

(4.4.225)     $$\log A = \frac{1}{12} - \varsigma'(-1)$$

where $A$ is the Glaisher-Kinkelin constant defined by (see also (6.83))

(4.4.226)     $$A = \lim_{n \to \infty}\left[ \sum_{k=1}^n k \log k - \left( \frac{n^2}{2} + \frac{n}{2} + \frac{1}{12} \right)\log n + \frac{n^2}{4} \right]$$

Equation (4.4.225) also follows from Adamchik's formula for the generalised Glaisher-Kinkelin constants derived in his 1998 paper "PolyGamma functions of negative order" [4]

$$(4.4.227) \qquad \log A_n = \frac{B_{n+1} H_n}{n+1} - \varsigma'(-n) \quad (\text{where } A_1 = A)$$

Using Riemann's functional equation for $\varsigma(s)$ we will also show in (F.7) that

$$(4.4.228) \qquad \frac{1}{12} - \varsigma'(-1) = \frac{1}{12}\gamma + \frac{1}{12}\log 2\pi - \frac{1}{2}\frac{\varsigma'(2)}{\pi^2}$$

<div align="center">ANOTHER DETERMINATION OF $\log G(1/2)$</div>

Kummer's Fourier series expansion (4.4.210) may be written as

$$\log \Gamma(x) = \frac{1}{2}\log\frac{\pi}{\sin\pi x} + \sum_{n=1}^{\infty}\frac{(\gamma + \log 2\pi n)\sin 2\pi nx}{\pi n} \quad (0 < x < 1)$$

and hence using Euler's reduction formula $\Gamma(x)\Gamma(1-x) = \dfrac{\pi}{\sin\pi x}$ we obtain

$$(4.4.228a) \qquad \frac{1}{2}\log\frac{\Gamma(x)}{\Gamma(1-x)} = \sum_{n=1}^{\infty}\frac{(\gamma + \log 2\pi n)\sin 2\pi nx}{\pi n}$$

Multiplying by $x$ and integrating we get for $t < 1$

$$(4.4.228b) \qquad \int_0^t x\log\frac{\Gamma(x)}{\Gamma(1-x)}dx = 2\sum_{n=1}^{\infty}\int_0^t\frac{(\gamma + \log 2\pi n)x\sin 2\pi nx}{\pi n}dx$$

It may be noted that

$$\lim_{x\to 0}\big[x\log\Gamma(x)\big] = \lim_{x\to 0}\big[x\log\Gamma(x+1) - x\log x\big] = 0$$

and hence Kummer's Fourier series expansion (when multiplied by $x$) may be validly used at $x = 0$.

We have (and this is also valid at $t = 1$)

$$(4.4.228c) \qquad \int_0^t x\log\frac{\Gamma(x)}{\Gamma(1-x)}dx = \int_0^t x\log\Gamma(x)\,dx + \int_1^{1-t}(1-x)\log\Gamma(x)\,dx$$

$$= \int_0^t x \log \Gamma(x)\,dx + \int_0^1 x \log \Gamma(x)\,dx - \int_0^{1-t} x \log \Gamma(x)\,dx - \int_0^1 \log \Gamma(x)\,dx + \int_0^{1-t} \log \Gamma(x)\,dx$$

With $t = 1/2$ we have

$$\int_0^{1/2} x \log \frac{\Gamma(x)}{\Gamma(1-x)}\,dx$$

$$= \int_0^{1/2} x \log \Gamma(x)\,dx + \int_0^1 x \log \Gamma(x)\,dx - \int_0^{1/2} x \log \Gamma(x)\,dx - \int_0^1 \log \Gamma(x)\,dx + \int_0^{1/2} \log \Gamma(x)\,dx$$

$$= \int_0^1 x \log \Gamma(x)\,dx - \int_0^1 \log \Gamma(x)\,dx + \int_0^{1/2} \log \Gamma(x)\,dx$$

In equation (6.116) in Volume V we will see that

$$\int_0^1 x \log \Gamma(1+x)\,dx = \frac{1}{4}\big[\log(2\pi) - 1\big] - \log A$$

and we deduce that

$$\int_0^1 x \log \Gamma(x)\,dx = \int_0^1 x \log \Gamma(1+x)\,dx - \int_0^1 x \log x\,dx$$

$$= \frac{1}{4}\log(2\pi) - \log A$$

We see from (C.43b) in Volume VI that

$$\int_0^1 \log \Gamma(x)\,dx = \frac{1}{2}\log(2\pi)$$

and we accordingly obtain

$$\int_0^{1/2} x \log \frac{\Gamma(x)}{\Gamma(1-x)}\,dx = \frac{1}{4}\log(2\pi) - \log A - \frac{1}{2}\log(2\pi) + \int_0^{1/2} \log \Gamma(x)\,dx$$

$$= -\frac{1}{4}\log(2\pi) - \log A + \int_0^{1/2} \log \Gamma(x)\,dx$$

Integration by parts shows that

$$\int_0^{1/2} x \sin 2n\pi x\, dx = \frac{(-1)^{n+1}}{4n\pi}$$

and therefore

$$2\sum_{n=1}^{\infty} \int_0^{1/2} \frac{(\gamma + \log 2\pi n)x \sin 2\pi nx}{\pi n}\, dx = \frac{1}{2\pi^2}(\gamma + \log 2\pi)\varsigma_a(2) - \frac{1}{2\pi^2}\varsigma_a'(2)$$

Hence we get from (4.4.228b)

$$-\frac{1}{4}\log(2\pi) - \log A + \int_0^{1/2} \log \Gamma(x)\, dx = \frac{1}{2\pi^2}(\gamma + \log 2\pi)\varsigma_a(2) - \frac{1}{2\pi^2}\varsigma_a'(2)$$

We have $\varsigma_a(2) = \frac{1}{2}\varsigma(2)$ and from (F.8h) in Volume VI we have

$$\varsigma_a'(2) = 6\varsigma(2)\left[\varsigma'(-1) - \frac{1}{12}(1 - \gamma - \log 2\pi)\right] + \frac{1}{2}\varsigma(2)\log 2$$

and using (4.4.225) $\log A = \frac{1}{12} - \varsigma'(-1)$ we get

(4.4.228ci) $\qquad \varsigma_a'(2) = \pi^2\left[\frac{1}{12}(\gamma + \log 2\pi + \log 2) - \log A\right]$

which results in

$$\frac{1}{2\pi^2}(\gamma + \log 2\pi)\varsigma_a(2) - \frac{1}{2\pi^2}\varsigma_a'(2) = \frac{1}{2}\log A - \frac{1}{24}\log 2$$

We therefore obtain [126, p.35]

(4.4.228d) $\qquad \int_0^{1/2} \log \Gamma(t)\, dt = \frac{5}{24}\log 2 + \frac{1}{4}\log \pi + \frac{3}{2}\log A$

where $A$ is the Glaisher-Kinkelin constant.

We have from (4.3.85) in Volume II(a) and also [126, p.32]

$$\int_0^z \log \Gamma(1+t)\,dt = \frac{1}{2}\big[\log(2\pi)-1\big]z - \frac{z^2}{2} + z\log\Gamma(1+z) - \log G(1+z)$$

and therefore

$$\int_0^{1/2} \log \Gamma(1+t)\,dt = \frac{1}{4}\big[\log(2\pi)-1\big] - \frac{1}{8} + \frac{1}{2}\log\Gamma(3/2) - \log G(3/2)$$

Using $G(3/2) = \Gamma(1/2)G(1/2)$ and $\Gamma(3/2) = \frac{1}{2}\Gamma(1/2)$ this becomes

$$\int_0^{1/2} \log \Gamma(1+t)\,dt = -\frac{1}{4}\log 2 - \frac{3}{8} - \log G(1/2)$$

We see that

$$\int_0^{1/2} \log \Gamma(1+t)\,dt = \int_0^{1/2} \log \Gamma(t)\,dt + \int_0^{1/2} \log t\,dt = \int_0^{1/2} \log \Gamma(t)\,dt - \frac{1}{2}\log 2 - \frac{1}{2}$$

and therefore we obtain

$$\int_0^{1/2} \log \Gamma(t)\,dt = \frac{1}{4}\log 2 + \frac{1}{8} - \log G(1/2)$$

Comparing this with (4.4.228d) we deduce that

(4.4.228e)     $\log G(1/2) = -\frac{3}{2}\log A - \frac{1}{4}\log\pi + \frac{1}{8} + \frac{1}{24}\log 2$

As we shall see later in (6.127) in Volume V, the value of $G(1/2)$ was originally determined by Barnes [17aa] in 1899 as (see also [126, p.26])

$$G(1/2) = A^{-\frac{3}{2}}\pi^{-\frac{1}{4}}e^{\frac{1}{8}}2^{\frac{1}{24}}$$

We may also write the above integrals as

$$\int_0^{1/2} \log \Gamma(t)\,dt = \frac{5}{24}\log 2 + \frac{1}{4}\log\pi + \frac{3}{2}\log A$$

$$\int\limits_0^{1/2} \log \Gamma\left(1+t\right) dt = \frac{1}{4}\left[\log(2\pi)-1\right] - \frac{1}{8} + \frac{1}{2}\log\frac{\sqrt{\pi}}{2} + \frac{3}{2}\log A + \frac{1}{4}\log\pi - \frac{1}{8} - \frac{1}{24}\log 2 - \log\sqrt{\pi}$$

$$= -\frac{1}{2} - \frac{7}{24}\log 2 + \frac{1}{4}\log\pi + \frac{3}{2}\log A$$

$$\int\limits_0^{1/2} x\log\frac{\Gamma(x)}{\Gamma(1-x)} dx = -\frac{1}{2} - \frac{7}{24}\log 2 + \frac{1}{2}\log A$$

## AN APPLICATION OF KUMMER'S FOURIER SERIES FOR $\log\Gamma(x)$

The above analysis may be generalised as follows. From (4.3.87ai) we have

$$2\int\limits_0^t x\log\Gamma(1+x)dx = \left(\frac{1}{4}-2\log A\right)t + \frac{1}{2}\log(2\pi)t^2 - \frac{1}{4}t^2 - \frac{1}{2}t^3 + t^2\log\Gamma(1+t)$$

$$-\log G(1+t) - 2\log\Gamma_3(1+t)$$

Accordingly we have

$$2\int\limits_0^{1-t} x\log\Gamma(x)dx = \left(\frac{1}{4}-2\log A\right)(1-t) + \frac{1}{2}\log(2\pi)(1-t)^2 + \frac{1}{4}(1-t)^2 - \frac{1}{2}(1-t)^3$$

$$+(1-t)^2\log\Gamma(1-t) - \log G(2-t) - 2\log\Gamma_3(2-t)$$

and hence we see that

$$\int\limits_0^t x\log\Gamma(x)dx - \int\limits_0^{1-t} x\log\Gamma(x)dx =$$

$$\frac{1}{2}\left(\frac{1}{4}-2\log A\right)t - \frac{1}{2}\left(\frac{1}{4}-2\log A\right)(1-t) + \frac{1}{4}\log(2\pi)\left[t^2-(1-t)^2\right] + \frac{1}{8}\left[t^2-(1-t)^2\right]$$

$$-\frac{1}{4}\left[t^3-(1-t)^3\right] + \frac{1}{2}t^2\log\Gamma(t) - \frac{1}{2}(1-t)^2\log\Gamma(1-t) - \frac{1}{2}\log G(1+t) + \frac{1}{2}\log G(2-t)$$

$$-\log\Gamma_3(1+t) + \log\Gamma_3(2-t)$$

From (4.4.228c) we have
$$\int_0^t x \log \frac{\Gamma(x)}{\Gamma(1-x)}\,dx =$$

$$\int_0^t x \log \Gamma(x)\,dx - \int_0^{1-t} x \log \Gamma(x)\,dx + \int_0^1 x \log \Gamma(x)\,dx - \int_0^1 \log \Gamma(x)\,dx + \int_0^{1-t} \log \Gamma(x)\,dx$$

and in (4.3.85) we saw that

$$\int_0^t \log \Gamma(x)\,dx = \frac{1}{2}t(1-t) + \frac{1}{2}t\log(2\pi) - \log G(t+1) + t\log \Gamma(t)$$

This gives us

$$\int_0^{1-t} \log \Gamma(x)\,dx = \frac{1}{2}t(1-t) + \frac{1}{2}(1-t)\log(2\pi) - \log G(2-t) + (1-t)\log \Gamma(1-t)$$

and hence we obtain

$$\int_0^t x \log \Gamma(x)\,dx - \int_0^{1-t} x \log \Gamma(x)\,dx + \int_0^1 x \log \Gamma(x)\,dx - \int_0^1 \log \Gamma(x)\,dx + \int_0^{1-t} \log \Gamma(x)\,dx =$$

$$\frac{1}{2}\left(\frac{1}{4} - 2\log A\right)t - \frac{1}{2}\left(\frac{1}{4} - 2\log A\right)(1-t) + \frac{1}{4}\log(2\pi)\left[t^2 - (1-t)^2\right] + \frac{1}{8}\left[t^2 - (1-t)^2\right]$$

$$-\frac{1}{4}\left[t^3 - (1-t)^3\right] + \frac{1}{2}t^2\log \Gamma(t) - \frac{1}{2}(1-t)^2\log \Gamma(1-t) - \frac{1}{2}\log G(1+t) + \frac{1}{2}\log G(2-t)$$

$$-\log \Gamma_3(1+t) + \log \Gamma_3(2-t) - \frac{1}{4}\log(2\pi) - \log A$$

$$+\frac{1}{2}t(1-t) + \frac{1}{2}(1-t)\log(2\pi) - \log G(2-t) + (1-t)\log \Gamma(1-t)$$

(4.4.228f)

$$= \frac{1}{2}\left(\frac{1}{4} - 2\log A\right)(2t-1) + \frac{1}{4}\log(2\pi)\left[t^2 - (1-t)^2\right] + \frac{1}{8}\left[t^2 - (1-t)^2\right]$$

$$-\frac{1}{4}\left[t^3 - (1-t)^3\right] + \frac{1}{2}t^2\log \Gamma(t) - \frac{1}{2}(1-t)^2\log \Gamma(1-t) - \frac{1}{2}\log G(1+t) + \frac{1}{2}\log G(2-t)$$

$$-\log \Gamma_3(1+t) + \log \Gamma_3(2-t) - \frac{1}{4}\log(2\pi) - \log A$$

$$+\frac{1}{2}t(1-t) + \frac{1}{2}(1-t)\log(2\pi) - \log G(2-t) + (1-t)\log \Gamma(1-t)$$

Noting that $\log G(2-t) = \log G(1-t) + \log \Gamma(1-t)$ and similarly
$\log \Gamma_3(2-t) = \log \Gamma_3(1-t) + \log G(1-t)$, we therefore obtain with a little bit of algebra

(4.4.228g)

$$\int_0^t x \log \frac{\Gamma(x)}{\Gamma(1-x)}\,dx = \frac{1}{4}t(1-t)(2t+1) - 2t\log A + \frac{1}{2}t^2 \log\frac{\Gamma(t)}{\Gamma(1-t)} - \frac{1}{2}\log\frac{G(1+t)}{G(1-t)} - \log\frac{\Gamma_3(1+t)}{\Gamma_3(1-t)}$$

Let us now turn to the right-hand side of (4.4.228b). Integration by parts shows that

$$\int_0^t x \sin 2n\pi x\,dx = -\frac{t\cos 2n\pi t}{2n\pi} + \frac{\sin 2n\pi t}{(2n\pi)^2}$$

and hence we have

$$S(t) = 2\sum_{n=1}^{\infty}\int_0^t \frac{(\gamma + \log 2\pi n)x\sin 2\pi nx}{\pi n}\,dx$$

$$= -4\sum_{n=1}^{\infty}\frac{(\gamma + \log 2\pi n)t\cos 2n\pi t}{(2\pi n)^2} + 4\sum_{n=1}^{\infty}\frac{(\gamma + \log 2\pi n)\sin 2n\pi t}{(2\pi n)^3}$$

We may write this as

$$S(t) = -4(\gamma + \log 2\pi)t\sum_{n=1}^{\infty}\frac{\cos 2n\pi t}{(2\pi n)^2} - 4t\sum_{n=1}^{\infty}\frac{\log n \cos 2n\pi t}{(2\pi n)^2}$$

$$+4(\gamma + \log 2\pi)\sum_{n=1}^{\infty}\frac{\sin 2n\pi t}{(2\pi n)^3} + 4\sum_{n=1}^{\infty}\frac{\log n \sin 2n\pi t}{(2\pi n)^3}$$

Using (6.25b) and (6.25c) from Volume V

$$B_{2N}(t) = (-1)^{N+1} 2(2N)! \sum_{n=1}^{\infty} \frac{\cos 2n\pi t}{(2\pi n)^{2N}} \qquad ,(N = 1, 2, \ldots)$$

$$B_{2N+1}(t) = (-1)^{N+1} 2(2N+1)! \sum_{n=1}^{\infty} \frac{\sin 2n\pi t}{(2\pi n)^{2N+1}} \qquad ,(N = 0, 1, 2, \ldots)$$

this may in turn be written in terms of Bernoulli polynomials

(4.4.228h)

$$S(t) = (\gamma + \log 2\pi)\left[\frac{1}{3}B_3(t) - tB_2(t)\right] - 4t\sum_{n=1}^{\infty} \frac{\log n \cos 2n\pi t}{(2\pi n)^2} + 4\sum_{n=1}^{\infty} \frac{\log n \sin 2n\pi t}{(2\pi n)^3}$$

In Titchmarsh's treatise [129, p.37] we have Hurwitz's formula for the Fourier expansion of $\varsigma(p,t)$

(4.4.228i) $\quad \varsigma(p,t) = 2\Gamma(1-p)\left[\sin\left(\frac{\pi p}{2}\right)\sum_{n=1}^{\infty} \frac{\cos 2n\pi t}{(2\pi n)^{1-p}} + \cos\left(\frac{\pi p}{2}\right)\sum_{n=1}^{\infty} \frac{\sin 2n\pi t}{(2\pi n)^{1-p}}\right]$

where $\text{Re}(p) < 0$ and $0 < t \leq 1$. In 2000, Boudjelkha [30c] showed that this formula also applies in the region $\text{Re}(p) < 1$. It may be noted that when $t = 1$ this reduces to Riemann's functional equation for $\varsigma(p)$. Letting $s = 1 - p$ we may write this as

(4.4.228j) $\quad \dfrac{\varsigma(1-s,t)}{\Gamma(s)} = 2\left[\cos\left(\frac{\pi s}{2}\right)\sum_{n=1}^{\infty} \frac{\cos 2n\pi t}{(2\pi n)^s} + \sin\left(\frac{\pi s}{2}\right)\sum_{n=1}^{\infty} \frac{\sin 2n\pi t}{(2\pi n)^s}\right]$

This may also be written as [126, p.89]

(4.4.228ji) $\quad \varsigma(1-s,t) = \dfrac{\Gamma(s)}{(2\pi)^s}\left[e^{-\frac{1}{2}\pi i t}L(t,s) + e^{\frac{1}{2}\pi i t}L(-t,s)\right]$

where the periodic (or Lerch) zeta function $L(t,s)$ is defined by

$$L(t,s) = \sum_{n=1}^{\infty} \frac{e^{2\pi i t}}{n^s}$$

Differentiating (4.4.228j) with respect to $s$ gives us

(4.4.228k)

$$-\frac{\varsigma'(1-s,t)+\varsigma(1-s,t)\psi(s)}{\Gamma(s)}=$$

$$-2\left[\cos\left(\frac{\pi s}{2}\right)\sum_{n=1}^{\infty}\frac{\log(2\pi n)\cos 2n\pi t}{(2\pi n)^s}+\sin\left(\frac{\pi s}{2}\right)\sum_{n=1}^{\infty}\frac{\log(2\pi n)\sin 2n\pi t}{(2\pi n)^s}\right]$$

$$+\pi\left[-\sin\left(\frac{\pi s}{2}\right)\sum_{n=1}^{\infty}\frac{\cos 2n\pi t}{(2\pi n)^s}+\cos\left(\frac{\pi s}{2}\right)\sum_{n=1}^{\infty}\frac{\sin 2n\pi t}{(2\pi n)^s}\right]$$

Let us now consider the case where $s$ is either even or odd. First of all, with $s=2k$ we get

(4.4.228l)

$$\frac{\varsigma'(1-2k,t)+\varsigma(1-2k,t)\psi(2k)}{\Gamma(2k)}=2(-1)^k\sum_{n=1}^{\infty}\frac{\log(2\pi n)\cos 2n\pi t}{(2\pi n)^{2k}}-\pi(-1)^k\sum_{n=1}^{\infty}\frac{\sin 2n\pi t}{(2\pi n)^{2k}}$$

where $\psi(2k)=H_{2k-1}-\gamma$.

Using (4.3.110) in Volume II(a) we have $\varsigma(1-m,t)=-\frac{B_m(t)}{m}$ and hence

$$\varsigma(1-2k,t)=-\frac{B_{2k}(t)}{2k}\quad\text{and}\quad\varsigma(-2k,t)=-\frac{B_{2k+1}(t)}{2k+1}$$

and therefore we have

$$\varsigma'(1-2k,t)-\frac{B_{2k}(t)\left(H_{2k-1}-\gamma\right)}{2k}=$$

$$2(-1)^k\Gamma(2k)\sum_{n=1}^{\infty}\frac{\log(2\pi n)\cos 2n\pi t}{(2\pi n)^{2k}}-\pi(-1)^k\Gamma(2k)\sum_{n=1}^{\infty}\frac{\sin 2n\pi t}{(2\pi n)^{2k}}$$

$$=2(-1)^k\Gamma(2k)\log(2\pi)\sum_{n=1}^{\infty}\frac{\cos 2n\pi t}{(2\pi n)^{2k}}+2(-1)^k\Gamma(2k)\sum_{n=1}^{\infty}\frac{\log n\cos 2n\pi t}{(2\pi n)^{2k}}$$

$$-\pi(-1)^k\Gamma(2k)\sum_{n=1}^{\infty}\frac{\sin 2n\pi t}{(2\pi n)^{2k}}$$

The above may be represented in terms of the Clausen function $\text{Cl}_n(x)$ defined by

$$\text{Cl}_{2n+1}(x) = \sum_{n=1}^{\infty} \frac{\cos kx}{k^{2n+1}}$$

$$\text{Cl}_{2n}(x) = \sum_{n=1}^{\infty} \frac{\sin kx}{k^{2n}}$$

and we therefore have

$$\varsigma'(1-2k,t) - \frac{B_{2k}(t)\left(H_{2k-1}-\gamma\right)}{2k} =$$

$$-\log(2\pi)\frac{1}{2k}B_{2k}(t) + 2(-1)^k \Gamma(2k)\sum_{n=1}^{\infty} \frac{\log n \cos 2n\pi t}{(2\pi n)^{2k}} - \frac{\pi(-1)^k}{(2\pi)^{2k}}\Gamma(2k)\text{Cl}_{2k}(2\pi t)$$

Hence we obtain

(4.4.228m) $\quad 2k\varsigma'(1-2k,t) - B_{2k}(t)\left(H_{2k-1}-\gamma-\log(2\pi)\right) =$

$$2(-1)^k(2k)!\sum_{n=1}^{\infty} \frac{\log n \cos 2n\pi t}{(2\pi n)^{2k}} - \frac{\pi(-1)^k(2k)!}{(2\pi)^{2k}}\text{Cl}_{2k}(2\pi t)$$

With $k=1$ we get

(4.4.228n) $\quad 2\varsigma'(-1,t) - B_2(t)\left[1-\gamma-\log(2\pi)\right] = -4\sum_{n=1}^{\infty} \frac{\log n \cos 2n\pi t}{(2\pi n)^2} + \frac{1}{2\pi}\text{Cl}_2(2\pi t)$

which may also be written as

(4.4.228o)

$$2\varsigma'(-1,t) - B_2(t)\left[1-\gamma-\log(2\pi)\right] = -4\sum_{n=1}^{\infty} \frac{\log n \cos 2n\pi t}{(2\pi n)^2} - \frac{1}{2\pi}\int_0^{2\pi t} \log\left[2\sin(x/2)\right]dx$$

With $t=1$ in (4.4.228m) we have

$$2\varsigma'(1-2k) - B_{2k}\left(H_{2k-1}-\gamma-\log(2\pi)\right) = 2(-1)^k \frac{(2k)!}{(2\pi)^{2k}}\sum_{n=1}^{\infty} \frac{\log n}{n^{2k}}$$

We therefore obtain the well-known result given by Miller and Adamchik [103] (see also (F.8) in Volume VI)

$$\varsigma'(2k) = \frac{(-1)^{k+1}(2\pi)^{2k}}{2(2k)!}\left(2k\varsigma'(1-2k) - B_{2k}\left[\psi(2k) - \log(2\pi)\right]\right)$$

and with $k=1$ we have the familiar result (see (F.7) in Volume VI)

$$\varsigma'(-1) = \frac{1}{2\pi^2}\varsigma'(2) - \frac{1}{12}\left[\log(2\pi) + \gamma - 1\right]$$

Letting $t = 1/2$ in (4.4.228n) we obtain

(4.4.228p) $\quad \varsigma'\left(-1, \frac{1}{2}\right) - \frac{1}{2}B_2\left(\frac{1}{2}\right)\left[1 - \gamma - \log(2\pi)\right] = \frac{1}{2\pi^2}\sum_{n=1}^{\infty}\frac{(-1)^{n+1}\log n}{n^2}$

$$= -\frac{1}{2\pi^2}\varsigma'_a\left(2\right)$$

and using (4.4.228ci) we get

$$\varsigma'\left(-1, \frac{1}{2}\right) - \frac{1}{2}B_2\left(\frac{1}{2}\right)\left[1 - \gamma - \log(2\pi)\right] = -\frac{1}{2}\left[\frac{1}{12}(\gamma + \log 2\pi + \log 2) - \log A\right]$$

After some cancellation this gives us

(4.4.228pi) $\quad \varsigma'\left(-1, \frac{1}{2}\right) = -\frac{1}{24}\log 2 - \frac{1}{2}\varsigma'(-1)$

which we have seen before in (4.3.140) in Volume II(a).

Letting $t = 1/4$ in (4.4.228n) we obtain

$$2\varsigma'\left(-1, \frac{1}{4}\right) - B_2\left(\frac{1}{4}\right)\left[1 - \gamma - \log(2\pi)\right] = -\frac{1}{\pi^2}\sum_{n=1}^{\infty}\frac{\log n \cos(n\pi/2)}{n^2} + \frac{1}{2\pi}\mathrm{Cl}_2\left(\frac{\pi}{2}\right)$$

and since $\mathrm{Cl}_2(\pi/2) = G$ we may write this as

$$\varsigma'\left(-1, \frac{1}{4}\right) - \frac{1}{96}\left[1 - \gamma - \log(2\pi)\right] = -\frac{1}{2\pi^2}\sum_{n=1}^{\infty}\frac{\log n \cos(n\pi/2)}{n^2} + \frac{G}{4\pi}$$

By inspection we see that

$$\sum_{n=1}^{\infty}\frac{\log n \cos(n\pi/2)}{n^2} = \sum_{n=1}^{\infty}(-1)^n\frac{\log(2n)}{(2n)^2} = \frac{1}{4}\log 2\sum_{n=1}^{\infty}\frac{(-1)^n}{n^2} + \frac{1}{4}\sum_{n=1}^{\infty}(-1)^n\frac{\log n}{n^2}$$

$$= -\frac{1}{4}\log 2\,\varsigma_a(2) + \frac{1}{4}\varsigma_a'(2)$$

and using (4.4.228ci)

$$\varsigma_a'(2) = \pi^2 \left[ \frac{1}{12}(\gamma + \log 2\pi + \log 2) - \log A \right]$$

we obtain

$$\sum_{n=1}^{\infty} \frac{\log n \cos(n\pi/2)}{n^2} = \frac{1}{4}\pi^2 \left[ \frac{1}{12}(\gamma + \log 2\pi) - \log A \right]$$

We shall see this again in (4.4.229j). We then deduce that

$$(4.4.228\text{pii}) \qquad \varsigma'\left(-1, \frac{1}{4}\right) = \frac{G}{4\pi} - \frac{1}{8}\varsigma'(-1)$$

and this formula is contained in Adamchik's paper [2a].

Similarly, with $s = 2k+1$ in (4.4.228k) we get

$$\frac{\varsigma'(-2k,t) + \varsigma(-2k,t)\psi(2k+1)}{\Gamma(2k+1)} = 2(-1)^k \sum_{n=1}^{\infty} \frac{\log(2\pi n)\sin 2n\pi t}{(2\pi n)^{2k+1}} + \pi(-1)^k \sum_{n=1}^{\infty} \frac{\cos 2n\pi t}{(2\pi n)^{2k+1}}$$

and we note that $\psi(2k+1) = H_{2k} - \gamma$. This may be written as

$$\varsigma'(-2k,t) + \varsigma(-2k,t)\psi(2k+1) = 2(-1)^k(2k)! \sum_{n=1}^{\infty} \frac{\log(2\pi n)\sin 2n\pi t}{(2\pi n)^{2k+1}} + \pi(-1)^k(2k)! \sum_{n=1}^{\infty} \frac{\cos 2n\pi t}{(2\pi n)^{2k+1}}$$

and we then obtain

$$(4.4.228\text{q}) \qquad \varsigma'(-2k,t) + \varsigma(-2k,t)\psi(2k+1) =$$

$$-\frac{1}{2k}\log(2\pi)B_{2k+1}(t) + 2(-1)^k(2k)! \sum_{n=1}^{\infty} \frac{\log n \sin 2n\pi t}{(2\pi n)^{2k+1}} + \frac{(-1)^k(2k)!}{2(2\pi)^{2k}}\text{Cl}_{2k+1}(2\pi t)$$

Since $\varsigma(-2k) = 0$, with $t = 1$ we obtain the well-known result (see (F.8a) in Volume VI)

$$\varsigma'(-2k) = \frac{(-1)^k(2k)!\varsigma(2k+1)}{2(2\pi)^{2k}}$$

With $k = 1$ we deduce that

$$\varsigma'(-2, t) + \varsigma(-2, t)\left(\frac{3}{2} - \gamma\right) = -\frac{1}{2}\log(2\pi)B_3(t) - 4\sum_{n=1}^{\infty}\frac{\log n \sin 2n\pi t}{(2\pi n)^3} - \frac{1}{(2\pi)^2}\text{Cl}_3(2\pi t)$$

and using (4.3.110) from Volume II(a) this may be written as

$$(4.4.228r) \quad \varsigma'(-2, t) - B_3(t)\left[\frac{1}{2} - \frac{1}{3}\gamma - \frac{1}{2}\log(2\pi)\right] = -4\sum_{n=1}^{\infty}\frac{\log n \sin 2n\pi t}{(2\pi n)^3} - \frac{1}{(2\pi)^2}\text{Cl}_3(2\pi t)$$

which may be contrasted with (4.4.228n) derived above

$$2\varsigma'(-1, t) - B_2(t)\left(1 - \gamma - \log(2\pi)\right) = -4\sum_{n=1}^{\infty}\frac{\log n \cos 2n\pi t}{(2\pi n)^2} + \frac{1}{2\pi}\text{Cl}_2(2\pi t)$$

Letting $t = 1$ in (4.4.228r) confirms that

$$(4.4.228s) \quad \varsigma'(-2) = -\frac{1}{(2\pi)^2}\text{Cl}_3(2\pi) = -\frac{\varsigma(3)}{4\pi^2}$$

in agreement with (F.8b) in Volume VI. Further identities may be obtained by integrating (4.4.228r).

Using (4.4.228n) and (4.4.228r) we see that

$$2t\varsigma'(-1, t) - t\,B_2(t)\left[1 - \gamma - \log(2\pi)\right] - \frac{1}{2\pi}t\,\text{Cl}_2(2\pi t) = -4t\sum_{n=1}^{\infty}\frac{\log n \cos 2n\pi t}{(2\pi n)^2}$$

$$-\varsigma'(-2, t) + B_3(t)\left[\frac{1}{2} - \frac{1}{3}\gamma - \frac{1}{2}\log(2\pi)\right] - \frac{1}{(2\pi)^2}\text{Cl}_3(2\pi t) = 4\sum_{n=1}^{\infty}\frac{\log n \sin 2n\pi t}{(2\pi n)^3}$$

and we see from (4.4.228g) that

$$\frac{1}{4}t(1-t)(2t+1) - 2t\log A + \frac{1}{2}t^2\log\frac{\Gamma(t)}{\Gamma(1-t)} - \frac{1}{2}\log\frac{G(1+t)}{G(1-t)} - \log\frac{\Gamma_3(1+t)}{\Gamma_3(1-t)} =$$

$$(\gamma + \log 2\pi)\left[\frac{1}{3}B_3(t) - t B_2(t)\right] - 4t\sum_{n=1}^{\infty}\frac{\log n \cos 2n\pi t}{(2\pi n)^2} + 4\sum_{n=1}^{\infty}\frac{\log n \sin 2n\pi t}{(2\pi n)^3}$$

$$= (\gamma + \log 2\pi)\left[\frac{1}{3}B_3(t) - tB_2(t)\right] + 2t\varsigma'(-1,t) - t\,B_2(t)\big[1 - \gamma - \log(2\pi)\big] - \frac{1}{2\pi}t\,\mathrm{Cl}_2(2\pi t)$$

$$-\varsigma'(-2,t) + B_3(t)\left[\frac{1}{2} - \frac{1}{3}\gamma - \frac{1}{2}\log(2\pi)\right] - \frac{1}{(2\pi)^2}\mathrm{Cl}_3(2\pi t)$$

This may be simplified to give a reflection formula for the triple Barnes gamma function

(4.4.228t)

$$\frac{1}{4}t(1-t)(2t+1) - 2t\log A + \frac{1}{2}t^2\log\frac{\Gamma(t)}{\Gamma(1-t)} - \frac{1}{2}\log\frac{G(1+t)}{G(1-t)} - \log\frac{\Gamma_3(1+t)}{\Gamma_3(1-t)} =$$

$$2t\varsigma'(-1,t) - \varsigma'(-2,t) - t\,B_2(t) + B_3(t)\left[\frac{1}{2} - \frac{1}{6}\log(2\pi)\right] - \frac{1}{(2\pi)^2}\mathrm{Cl}_3(2\pi t) - \frac{1}{2\pi}t\,\mathrm{Cl}_2(2\pi t)$$

With $t = 1/2$ in (4.4.228t) we get

$$\frac{1}{8} - \log A - \frac{1}{2}\log\frac{G(3/2)}{G(1/2)} - \log\frac{\Gamma_3(3/2)}{\Gamma_3(1/2)} = \frac{1}{8} - \log A - \frac{1}{4}\log\pi - \log G(1/2)$$

$$= \varsigma'\left(-1,\frac{1}{2}\right) + \frac{1}{24} - \varsigma'\left(-2,\frac{1}{2}\right) + \frac{3\varsigma(3)}{16\pi^2}$$

and using (4.3.168d) this becomes

$$= \frac{1}{24} - \frac{1}{24}\log 2 - \frac{1}{2}\varsigma'(-1)$$

We then see that

$$\frac{1}{8} - \log A - \frac{1}{4}\log\pi - \log G(1/2) = \frac{1}{24} - \frac{1}{24}\log 2 - \frac{1}{2}\varsigma'(-1)$$

and we deduce the well-known result (see 4.3.126bi) and [126, p.26])

(4.4.228ti)     $$\log G(1/2) = \frac{1}{8} + \frac{1}{24}\log 2 - \frac{3}{2}\log A - \frac{1}{4}\log\pi$$

Using (6.69c) from Volume V

$$\text{Cl}_2(2\pi t) = \sum_{n=1}^{\infty} \frac{\sin 2\pi n t}{n^2} = -2\pi t \log\left[\frac{\sin \pi t}{\pi}\right] - 2\pi \log\frac{G(1+t)}{G(1-t)}$$

we may express (4.4.228t) as follows

(4.4.228tii)

$$\frac{1}{4}t(1-t)(2t+1) - 2t\log A + \frac{1}{2}t^2 \log\frac{\Gamma(t)}{\Gamma(1-t)} - \left(\frac{1}{2}+t\right)\log\frac{G(1+t)}{G(1-t)} - \log\frac{\Gamma_3(1+t)}{\Gamma_3(1-t)} =$$

$$2t\varsigma'(-1,t) - \varsigma'(-2,t) - t\,B_2(t) + B_3(t)\left[\frac{1}{2} - \frac{1}{6}\log(2\pi)\right] - \frac{1}{(2\pi)^2}\text{Cl}_3(2\pi t) + t^2 \log\frac{\sin \pi t}{\pi}$$

At first glance, the left-hand side of (4.4.228tii) does not appear to be finite at $t=1$ but from (4.4.228c) we see that

$$\int_0^1 x\log\frac{\Gamma(x)}{\Gamma(1-x)}\,dx = \int_0^1 x\log\Gamma(x)\,dx + \int_1^0 (1-x)\log\Gamma(x)\,dx$$

$$= 2\int_0^1 x\log\Gamma(x)\,dx - \int_0^1 \log\Gamma(x)\,dx$$

and from (6.126) in Volume V we see that

$$\int_0^1 x\log\Gamma(1+x)\,dx = \frac{1}{4}\log(2\pi) - \frac{1}{4} - \log A$$

and therefore we have

$$\int_0^1 x\log\Gamma(x)\,dx = \frac{1}{4}\log(2\pi) - \log A$$

Hence we see that

$$\int_0^1 x\log\frac{\Gamma(x)}{\Gamma(1-x)}\,dx = -2\log A$$

With $t=1$ in the right-hand side of (4.4.228t) we get

$$= -\frac{1}{6}(\gamma + \log 2\pi) + 2\varsigma'(-1) - \frac{1}{6}\left[1 - \gamma - \log(2\pi)\right] - \varsigma'(-2) - \frac{\varsigma(3)}{(2\pi)^2}$$

$$= 2\varsigma'(-1) - \frac{1}{6} = -2\log A$$

and the equality is a useful check on my algebra.

□

With $t = 1/4$ in (4.4.228t) we get

(4.4.228u)

$$\frac{9}{128} - \frac{1}{2}\log A + \frac{1}{32}\log\frac{\Gamma(1/4)}{\Gamma(3/4)} - \frac{1}{2}\log\frac{G(5/4)}{G(3/4)} - \log\frac{\Gamma_3(5/4)}{\Gamma_3(3/4)} =$$

$$\frac{1}{2}\varsigma'\left(-1,\frac{1}{4}\right) - \varsigma'\left(-2,\frac{1}{4}\right) - \frac{1}{192} + \frac{3}{64}\left[\frac{1}{2} - \frac{1}{6}\log(2\pi)\right] - \frac{1}{(2\pi)^2}\text{Cl}_3\left(\frac{\pi}{2}\right) - \frac{1}{8\pi}\text{Cl}_2\left(\frac{\pi}{2}\right)$$

With Euler's reflection formula we see that

$$\log\Gamma(3/4) = \log\pi + \frac{1}{2}\log 2 - \log\Gamma(1/4)$$

Using a duplication formula for the G-function, Choi and Srivastava [45aa] have shown that

$$\log G\left(\frac{1}{4}\right) = \frac{3}{32} - \frac{G}{4\pi} - \frac{3}{4}\log\Gamma\left(\frac{1}{4}\right) - \frac{9}{8}\log A$$

$$\log G\left(\frac{3}{4}\right) = -\frac{1}{8}\log 2 - \frac{1}{4}\log\pi + \frac{3}{32} + \frac{G}{4\pi} - \frac{9}{8}\log A + \frac{1}{4}\log\Gamma\left(\frac{1}{4}\right)$$

$$\log G\left(\frac{5}{4}\right) - \log G\left(\frac{3}{4}\right) = \frac{1}{8}\log 2 + \frac{1}{4}\log\pi - \frac{G}{2\pi}$$

Choi, Cho and Srivastava [45ac] show that

$$\log\Gamma_3(1+x) = -\frac{1}{24} + \frac{1}{2}\log A + \frac{\varsigma(3)}{8\pi^2} + \left(\frac{1}{12} - \log A\right)x - \frac{1}{2}x(1-x)\log\Gamma(1+x)$$

$$+ \left(\frac{1}{2} - x\right)\varsigma'(-1,1+x) + \frac{1}{2}\varsigma'(-2,1+x)$$

and from this it follows with $x = -1/2$ that

$$\Gamma_3\left(\frac{1}{2}\right) = 2^{-\frac{1}{24}} \pi^{\frac{3}{16}} \exp\left[-\frac{1}{8} + \frac{7\varsigma(3)}{32\pi^2}\right] A^{\frac{1}{2}}$$

Adamchik [6c] also reports that

$$\varsigma'(-2, x) - \varsigma'(-2) = 2\log\Gamma_3(x) + (3 - 2x)\log G(x) + (1 - x)^2 \log\Gamma(x)$$

We have from (4.3.161a) and (4.3.161b)

$$\varsigma'\left(-1, \frac{1}{4}\right) = \frac{G}{4\pi} - \frac{1}{8}\varsigma'(-1)$$

$$\varsigma'\left(-1, \frac{3}{4}\right) = -\frac{G}{4\pi} - \frac{1}{8}\varsigma'(-1)$$

Hence, using the above data, we are able to derive an expression for $\log\dfrac{\Gamma_3(5/4)}{\Gamma_3(3/4)}$.

We recall Adamchik's result [2a] from (4.3.167)

$$\varsigma'(-2k, t) + \varsigma'(-2k, 1 - t) = (-1)^k \frac{(2k)!}{(2\pi)^{2k}} \mathrm{Cl}_{2k+1}(2\pi t)$$

$$\varsigma'(-2k - 1, t) - \varsigma'(-2k - 1, 1 - t) = \frac{(2k+1)!}{(2\pi)^{2k+1}} \mathrm{Cl}_{2k+2}(2\pi t)$$

where with $k = 1$ we get

$$\varsigma'(-2, t) + \varsigma'(-2, 1 - t) = -\frac{1}{2\pi^2} \mathrm{Cl}_3(2\pi t)$$

$$\varsigma'(-1, t) - \varsigma'(-1, 1 - t) = \frac{1}{2\pi} \mathrm{Cl}_2(2\pi t)$$

$\square$

Upon integrating (4.4.228n) we obtain

$$2\int_0^x \varsigma'(-1, t)\,dt - \frac{1}{3} B_3(x)\left(1 - \gamma - \log(2\pi)\right) = -4\sum_{n=1}^{\infty} \frac{\log n \sin 2n\pi x}{(2\pi n)^3} + \frac{1}{2\pi}\int_0^x \mathrm{Cl}_2(2\pi t)\,dt$$

and we recall Adamchik's formula (4.3.131) from Volume II(a)

$$n\int_0^x \varsigma'(1-n,t)\,dt = \frac{B_{n+1} - B_{n+1}(x)}{n(n+1)} + \varsigma'(-n,x) - \varsigma'(-n)$$

This gives us

$$2\int_0^x \varsigma'(-1,t)\,dt = -\frac{1}{6}B_3(x) + \varsigma'(-2,x) - \varsigma'(-2)$$

From the definition of the Clausen function it is readily seen that [126, p.115]

$$\int_0^x \mathrm{Cl}_{2n}(at)\,dt = \frac{1}{a}\big[\varsigma(2n+1) - \mathrm{Cl}_{2n+1}(ax)\big]$$

$$\int_0^x \mathrm{Cl}_{2n+1}(at)\,dt = \frac{1}{a}\mathrm{Cl}_{2n+2}(ax)$$

Hence we have

$$-\frac{1}{6}B_3(x) + \varsigma'(-2,x) - \varsigma'(-2) - \frac{1}{3}B_3(x)\big[1 - \gamma - \log(2\pi)\big] =$$

$$-4\sum_{n=1}^{\infty} \frac{\log n \sin 2n\pi x}{(2\pi n)^3} + \frac{1}{4\pi^2}\big[\varsigma(3) - \mathrm{Cl}_3(2\pi x)\big]$$

which, using (4.4.228s), may be written as

$$\varsigma'(-2,x) - B_3(x)\left[\frac{1}{2} - \frac{1}{3}\gamma - \frac{1}{3}\log(2\pi)\right] = -4\sum_{n=1}^{\infty} \frac{\log n \sin 2n\pi x}{(2\pi n)^3} - \frac{1}{4\pi^2}\mathrm{Cl}_3(2\pi x)$$

We note that this corresponds with (4.4.228r) which was derived above and hence this could also serve as an alternative derivation of Adamchik's formula (4.3.131).

Hence with $x = 1/2$ we obtain

$$\varsigma'\left(-2,\frac{1}{2}\right) = -\frac{1}{4\pi^2}\mathrm{Cl}_3(\pi)$$

With $x = 1/4$ we get

$$\varsigma'\left(-2,\frac{1}{4}\right) - B_3\left(\frac{1}{4}\right)\left[\frac{1}{2} - \frac{1}{3}\gamma - \frac{1}{3}\log(2\pi)\right] = -4\sum_{n=1}^{\infty}\frac{\log n \sin(n\pi/2)}{(2\pi n)^3} - \frac{1}{4\pi^2}\mathrm{Cl}_3\left(\frac{\pi}{2}\right)$$

and this may be written as

$$\varsigma'\left(-2,\frac{1}{4}\right) - B_3\left(\frac{1}{4}\right)\left[\frac{1}{2} - \frac{1}{3}\gamma - \frac{1}{3}\log(2\pi)\right] = -\frac{4}{\pi^3}\sum_{n=0}^{\infty}(-1)^n\frac{\log(2n+1)}{(2n+1)^3} - \frac{\varsigma(3)}{128\pi^2}$$

It is easily seen from the definition of the Clausen function that

$$\mathrm{Cl}_{2n}(\pi) = \mathrm{Cl}_{2n}(2\pi) = 0$$

$$\mathrm{Cl}_{2n+1}(\pi) = (2^{-2n}-1)\varsigma(2n+1) = -\varsigma_a(2n+1)$$

$$\mathrm{Cl}_{2n+1}(2\pi) = \varsigma(2n+1)$$

and we therefore deduce that

$$\varsigma'\left(-2,\frac{1}{2}\right) = \frac{3\varsigma(3)}{16\pi^2}$$

which we have seen previously in (4.3.168d).

We also have

$$\mathrm{Cl}_2(\pi/2) = G = -\mathrm{Cl}_2(3\pi/2)$$

$$\frac{1}{2}\mathrm{Cl}_2(2x) = \mathrm{Cl}_2(x) - \mathrm{Cl}_2(\pi - x)$$

which implies that

$$\mathrm{Cl}_2\left(\frac{2\pi}{3}\right) = \frac{2}{3}\mathrm{Cl}_2\left(\frac{\pi}{3}\right)$$

Adamchik [6a] has shown that

$$\mathrm{Cl}_2\left(\frac{2\pi}{3}\right) = \frac{\psi^{(1)}\left(\frac{1}{3}\right)}{3\sqrt{3}} - \frac{2\pi^2}{9\sqrt{3}}$$

and this may be derived by letting $t = 1/3$ in (4.4.228n).

The Clausen function may be expressed in closed form in at least three other cases and from Lewin's book [100, p.198] we have

$$\text{Cl}_{2n+1}\left(\frac{\pi}{2}\right) = -2^{-2n-1}(1-2^{-2n})\varsigma(2n+1)$$

$$\text{Cl}_{2n+1}\left(\frac{\pi}{3}\right) = \frac{1}{2}(1-2^{-2n})(1-3^{-2n})\varsigma(2n+1)$$

$$\text{Cl}_{2n+1}\left(\frac{2\pi}{3}\right) = -\frac{1}{2}(1-3^{-2n})\varsigma(2n+1)$$

For example, we see from the definition that

$$\text{Cl}_{2n+1}\left(\frac{\pi}{2}\right) = -\frac{1}{2^{2n+1}} + \frac{1}{4^{2n+1}} - \frac{1}{6^{2n+1}} + \dots$$

$$= -\frac{1}{2^{2n+1}}\left[\frac{1}{1^{2n+1}} - \frac{1}{2^{2n+1}} + \dots\right] = -\frac{1}{2^{2n+1}}\varsigma_a(2n+1)$$

We also have a number of formulae involving $\text{Cl}_2\left(\dfrac{p\pi}{q}\right)$ in Browkin's paper in Lewin's survey [101, p.244], including for example

$$\text{Cl}_2\left(\frac{\pi}{6}\right) + \text{Cl}_2\left(\frac{5\pi}{6}\right) = \frac{4}{3}G$$

Using PSLQ, Bailey et al. [16a] discovered experimentally that

$$6\text{Cl}_2\left(\alpha\right) - 6\text{Cl}_2\left(2\alpha\right) + 2\text{Cl}_2\left(3\alpha\right) = 7\text{Cl}_2\left(\frac{2\pi}{7}\right) + 7\text{Cl}_2\left(\frac{4\pi}{7}\right) - 7\text{Cl}_2\left(\frac{6\pi}{7}\right)$$

where $\alpha = 2\tan^{-1}\sqrt{7}$. Using this in conjunction with (4.4.228m) does look rather fearsome.

With reference to (4.3.162)

$$-\log\frac{G(1+t)}{G(1-t)} = \varsigma'(-1,t) - \varsigma'(-1,1-t) + t\log\frac{\sin\pi t}{\pi}$$

we may write (4.4.228t) as

$$\frac{1}{4}t(1-t)(2t+1) - 2t\log A + \frac{1}{2}t^2\log\frac{\Gamma(t)}{\Gamma(1-t)} - \log\frac{\Gamma_3(1+t)}{\Gamma_3(1-t)} =$$

$$\left(2t - \frac{1}{2}\right)\varsigma'(-1,t) + \frac{1}{2}\varsigma'(-1,1-t) - \varsigma'(-2,t) - \frac{1}{2}t\log\frac{\sin\pi t}{\pi}$$

$$-t\,B_2(t) + B_3(t)\left[\frac{1}{2} - \frac{1}{6}\log(2\pi)\right] - \frac{1}{(2\pi)^2}\operatorname{Cl}_3(2\pi t) - \frac{1}{2\pi}t\operatorname{Cl}_2(2\pi t)$$

As mentioned by Adamchik [5b] we also have the reflection formula for the G-function which is valid for $0 < t < 1$

$$\log\frac{G(1+t)}{G(1-t)} = -\frac{1}{2\pi}\operatorname{Cl}_2(2\pi t) - t\log\frac{\sin\pi t}{\pi}$$

□

Multiplying (4.4.228a) by $x^2$ and integrating we get

$$\int_0^t x^2\log\frac{\Gamma(x)}{\Gamma(1-x)}\,dx = 2\sum_{n=1}^\infty\int_0^t\frac{(\gamma + \log 2\pi n)x^2\sin 2\pi nx}{\pi n}\,dx$$

With $t = 1/2$ we obtain

$$\int_0^{1/2} x^2\log\frac{\Gamma(x)}{\Gamma(1-x)}\,dx = -\frac{1}{4\pi^2}\sum_{n=1}^\infty(\gamma + \log 2\pi n)\frac{1}{n^2} + \frac{1}{2\pi^4}\sum_{n=1}^\infty(\gamma + \log 2\pi n)\left[\frac{(-1)^n}{n^4} - \frac{1}{n^4}\right]$$

$$= -\frac{\gamma + \log 2\pi}{24} + \frac{1}{4\pi^2}\varsigma'(2) - \frac{\gamma + \log 2\pi}{2\pi^4}\left[\varsigma_a(4) + \varsigma'(4)\right] + \frac{1}{2\pi^4}\left[\varsigma_a'(4) + \varsigma'(4)\right]$$

Choi and Srivastava [45ab] have evaluated the integral $\int_0^t x^2\log\Gamma(x+a)\,dx$ in terms involving the multiple gamma functions $\Gamma_n(x+a)$, with $n = 1, 2$ and $3$, and the integral $\int_0^t\log\Gamma_3(x+a)\,dx$. It may be worthwhile extending the analysis here to, for example, $\int_0^t x^p\log\frac{\Gamma(x)}{\Gamma(1-x)}\,dx$ or $\int_0^t B_p(x)\log\frac{\Gamma(x)}{\Gamma(1-x)}\,dx$.

# SOME FOURIER SERIES CONNECTIONS

We now recall Hurwitz's formula (4.4.228i) which is valid for $\operatorname{Re}(p) < 1$

$$\varsigma(p,t) = 2\Gamma(1-p)\left[\sin\left(\frac{\pi p}{2}\right)\sum_{n=1}^{\infty}\frac{\cos 2n\pi t}{(2\pi n)^{1-p}} + \cos\left(\frac{\pi p}{2}\right)\sum_{n=1}^{\infty}\frac{\sin 2n\pi t}{(2\pi n)^{1-p}}\right]$$

and letting $t \to 1-t$ we get

$$\varsigma(p,1-t) = 2\Gamma(1-p)\left[\sin\left(\frac{\pi p}{2}\right)\sum_{n=1}^{\infty}\frac{\cos 2n\pi t}{(2\pi n)^{1-p}} - \cos\left(\frac{\pi p}{2}\right)\sum_{n=1}^{\infty}\frac{\sin 2n\pi t}{(2\pi n)^{1-p}}\right]$$

We therefore see that

$$\varsigma(p,t) + \varsigma(p,1-t) = 4\Gamma(1-p)\sin\left(\frac{\pi p}{2}\right)\sum_{n=1}^{\infty}\frac{\cos 2n\pi t}{(2\pi n)^{1-p}}$$

$$\varsigma(p,t) - \varsigma(p,1-t) = 4\Gamma(1-p)\cos\left(\frac{\pi p}{2}\right)\sum_{n=1}^{\infty}\frac{\sin 2n\pi t}{(2\pi n)^{1-p}}$$

From the above are immediately found the corresponding Fourier coefficients

(4.4.229a) $\qquad \int_0^1 \varsigma(p,t)\sin 2n\pi t\, dt = \frac{\Gamma(1-p)}{(2\pi n)^{1-p}}\cos\left(\frac{\pi p}{2}\right)$

(4.4.229b) $\qquad \int_0^1 \varsigma(p,t)\cos 2n\pi t\, dt = \frac{\Gamma(1-p)}{(2\pi n)^{1-p}}\sin\left(\frac{\pi p}{2}\right)$

Using Euler's reflection formula these may be written as (see [59])

(4.4.229c) $\qquad \int_0^1 \varsigma(p,t)\sin 2n\pi t\, dt = \frac{(2\pi)^p n^{p-1}}{4\Gamma(p)}\csc\left(\frac{\pi p}{2}\right)$

(4.4.229d) $\qquad \int_0^1 \varsigma(p,t)\cos 2n\pi t\, dt = \frac{(2\pi)^p n^{p-1}}{4\Gamma(p)}\sec\left(\frac{\pi p}{2}\right)$

Having regard to (4.4.229a) we consider

$$f(p) = \frac{\Gamma(1-p)}{(2\pi n)^{1-p}} \cos\left(\frac{\pi p}{2}\right)$$

and logarithmic differentiation gives us

$$\frac{f'(p)}{f(p)} = \log(2\pi n) - \psi(1-p) - \frac{\pi}{2}\tan\left(\frac{\pi p}{2}\right)$$

Therefore differentiating (4.4.229a) results in

(4.4.229e)

$$\int_0^1 \varsigma'(p,t) \sin 2n\pi t \, dt = \left[\log(2\pi n) - \psi(1-p) - \frac{\pi}{2}\tan\left(\frac{\pi p}{2}\right)\right]\frac{\Gamma(1-p)}{(2\pi n)^{1-p}}\cos\left(\frac{\pi p}{2}\right)$$

Applying Lerch's formula (4.3.116) we then have as $p \to 0$

$$\int_0^1 \left[\log\Gamma(t) - \frac{1}{2}\log(2\pi)\right] \sin 2n\pi t \, dt = \frac{\log(2\pi n) + \gamma}{2\pi n}$$

which becomes the well-known result

(4.4.229f)    $$\int_0^1 \log\Gamma(t) \sin 2n\pi t \, dt = \frac{\log(2\pi n) + \gamma}{2\pi n}$$

Similarly differentiating (4.4.229b) gives us

(4.4.229g)

$$\int_0^1 \varsigma'(p,t) \cos 2n\pi t \, dt = \left[\log(2\pi n) - \psi(1-p) + \frac{\pi}{2}\cot\left(\frac{\pi p}{2}\right)\right]\frac{\Gamma(1-p)}{(2\pi n)^{1-p}}\sin\left(\frac{\pi p}{2}\right)$$

and as $p \to 0$ we obtain

$$\int_0^1 \left[\log\Gamma(t) - \frac{1}{2}\log(2\pi)\right] \cos 2n\pi t \, dt = \frac{1}{4n}$$

This then gives us

(4.4.229h)    $$\int_0^1 \log\Gamma(t) \cos 2n\pi t \, dt = \frac{1}{4n}$$

We have thereby rediscovered Kummer's Fourier series (4.4.210) for $\log \Gamma(t)$.

With $p = -1$ in (4.4.229e) and (4.4.229g) we obtain

(4.4.229hi)     $\displaystyle \int_0^1 \varsigma'(-1,t)\sin 2n\pi t\, dt = \frac{1}{8\pi n^2}$

(4.4.229hii)     $\displaystyle \int_0^1 \varsigma'(-1,t)\cos 2n\pi t\, dt = -\frac{1}{(2\pi n)^2}\big[\log(2\pi n) + \gamma - 1\big]$

and from (4.3.145) we see that

$$\int_0^1 \varsigma'(-1,t)\, dt = 0.$$

We then deduce from these Fourier coefficients that for $0 \le t \le 1$

(4.4.229i)     $\displaystyle \varsigma'(-1,t) = -2\sum_{n=1}^{\infty}\frac{\log(2n\pi) + \gamma - 1}{(2n\pi)^2}\cos 2n\pi t + \pi\sum_{n=1}^{\infty}\frac{\sin 2n\pi t}{(2n\pi)^2}$

and this may be written as

$$\varsigma'(-1,t) = -\frac{1}{2}[\log(2\pi) + \gamma - 1]B_2(t) - 2\sum_{n=1}^{\infty}\frac{\log n}{(2n\pi)^2}\cos 2n\pi t + \frac{1}{4\pi}\mathrm{Cl}_2(2\pi t)$$

More generally, we have the Fourier series

(4.4.229ii)

$$\varsigma'(p,t) = 2\sum_{n=1}^{\infty}\left[\log(2\pi n) - \psi(1-p) + \frac{\pi}{2}\cot\left(\frac{\pi p}{2}\right)\right]\frac{\Gamma(1-p)}{(2\pi n)^{1-p}}\sin\left(\frac{\pi p}{2}\right)\cos 2n\pi t$$

$$+ 2\sum_{n=1}^{\infty}\left[\log(2\pi n) - \psi(1-p) - \frac{\pi}{2}\tan\left(\frac{\pi p}{2}\right)\right]\frac{\Gamma(1-p)}{(2\pi n)^{1-p}}\cos\left(\frac{\pi p}{2}\right)\sin 2n\pi t$$

We may note from (4.4.228i) that for $\mathrm{Re}(p) < 1$

$$\int_0^1 \varsigma(p,t)\, dt = 0$$

and from (4.4.229ii) that

$$\int_0^1 \varsigma'(p,t)\,dt = 0$$

which we also saw in (4.3.145) in Volume II(a).

With $t=1$ or $t=0$ in (4.4.229i) we obtain

(4.4.229iii)　$\varsigma'(-1) = -2\sum_{n=1}^{\infty} \frac{\log(2\pi)+\gamma-1+\log n}{(2\pi n)^2}$

which may be written in its more familiar form

$$\varsigma'(-1) = \frac{1-\log(2\pi)-\gamma}{12} + \frac{\varsigma'(2)}{2\pi^2}$$

Letting $t \to 1-t$ we obtain

$$\varsigma'(-1,1-t) = -2\sum_{n=1}^{\infty} \frac{\log(2\pi n)+\gamma-1}{(2\pi n)^2}\cos 2n\pi t - \sum_{n=1}^{\infty} \frac{\sin 2n\pi t}{4\pi n^2}$$

and we therefore have

(4.4.229iv)　$\varsigma'(-1,t) - \varsigma'(-1,1-t) = \frac{1}{2\pi}\sum_{n=1}^{\infty}\frac{\sin 2n\pi t}{n^2} = \frac{1}{2\pi}\mathrm{Cl}_2(2\pi t)$

which is in accordance with Adamchik's result [2a].

With $t=1/4$ in (4.4.229i) we obtain

$$\varsigma'\left(-1,\frac{1}{4}\right) = -2\sum_{n=1}^{\infty}\frac{\log(2\pi)+\log n+\gamma-1}{(2\pi n)^2}\cos(n\pi/2) + \sum_{n=1}^{\infty}\frac{\sin(n\pi/2)}{4\pi n^2}$$

$$= \frac{1}{96}\big[\log(2\pi)+\gamma-1\big] - 2\sum_{n=1}^{\infty}\frac{\log n\cos(n\pi/2)}{(2\pi n)^2} + \frac{G}{4\pi}$$

where $G$ is Catalan's constant

$$G = \sum_{n=0}^{\infty}\frac{(-1)^n}{(2n+1)^2} = \sum_{n=1}^{\infty}\frac{\sin(n\pi/2)}{n^2}$$

and

$$\sum_{n=1}^{\infty} \frac{\cos(n\pi/2)}{n^2} = -\frac{1}{2^2} + \frac{1}{4^2} - ... = -\frac{1}{2^2}\left[\frac{1}{1^2} - \frac{1}{2^2} + ...\right] = -\frac{1}{2^2}\varsigma_a(2) = -\frac{\pi^2}{48}$$

$G$ is a particular case of the Dirichlet beta function $\beta(s)$ defined by

$$\beta(s) = \sum_{n=0}^{\infty} \frac{(-1)^n}{(2n+1)^s} = \frac{1}{2^{2s}}\left[\varsigma\left(s,\frac{1}{4}\right) - \varsigma\left(s,\frac{3}{4}\right)\right]$$

We recall Adamchik's formula (4.3.161a)

$$\varsigma'\left(-1,\frac{1}{4}\right) = \frac{G}{4\pi} - \frac{1}{8}\varsigma'(-1)$$

and we therefore obtain

(4.4.229j) $$\sum_{n=1}^{\infty} \frac{\log n \cos(n\pi/2)}{(\pi n)^2} = \frac{1}{48}\left[\log(2\pi) + \gamma - 1\right] + \frac{1}{4}\varsigma'(-1)$$

which we have seen earlier in this paper.

From (4.3.129c) we have

$$\varsigma'(-1, N+1) - \varsigma'(-1,1) = \sum_{k=1}^{N} k \log k$$

but we cannot let $t = N$ in (4.4.229i) because of its limited range of validity.

Differentiating (4.4.229i) and making use of (4.3.142)

$$\frac{\partial}{\partial t}\frac{\partial}{\partial s}\varsigma(s,t) = -\varsigma(s+1,t) - s\frac{\partial}{\partial s}\varsigma(s+1,t)$$

we obtain

$$\varsigma'(0,t) - \varsigma(0,t) = 2\sum_{n=1}^{\infty} \frac{\log(2\pi n) + \gamma - 1}{2\pi n}\sin 2n\pi t + \sum_{n=1}^{\infty} \frac{\cos 2n\pi t}{2n}$$

and, upon using Lerch's identity, we see that

$$\log\Gamma(t) - \frac{1}{2}\log(2\pi) - \varsigma(0,t) = \sum_{n=1}^{\infty} \frac{\log(2\pi n) + \gamma}{\pi n}\sin 2n\pi t - \frac{1}{\pi}\sum_{n=1}^{\infty} \frac{\sin 2n\pi t}{n} + \sum_{n=1}^{\infty} \frac{\cos 2n\pi t}{2n}$$

which, by reference to (7.5) and (7.5), is equivalent to Kummer's Fourier series expansion (4.4.210).

Upon integrating (4.4.229i) we obtain

$$(4.4.229k) \quad \int_0^x \varsigma'(-1,t)dt = -\frac{1}{12}B_3(x) + \frac{1}{2}[\varsigma'(-2,x) - \varsigma'(-2)]$$

$$= -2\sum_{n=1}^{\infty} \frac{\log(2\pi n) + \gamma - 1}{(2\pi n)^3}\sin 2n\pi x - \frac{1}{8\pi^2}\sum_{n=1}^{\infty}\frac{\cos 2n\pi x - 1}{n^3}$$

where we have employed Adamchik's integral (4.3.131) from Volume II(a). This integration exercise obviously may be continued indefinitely.

Letting $x \to 1 - x$ we obtain

$$-\frac{1}{12}B_3(1-x) + \frac{1}{2}[\varsigma'(-2,1-x) - \varsigma'(-2)] =$$

$$+2\sum_{n=1}^{\infty}\frac{\log(2\pi n) + \gamma - 1}{(2\pi n)^3}\sin 2n\pi x - \frac{1}{8\pi^2}\sum_{n=1}^{\infty}\frac{\cos 2n\pi x - 1}{n^3}$$

and hence we get

$$-\frac{1}{12}B_3(x) + \frac{1}{2}[\varsigma'(-2,x) - \varsigma'(-2)] - \frac{1}{12}B_3(1-x) + \frac{1}{2}[\varsigma'(-2,1-x) - \varsigma'(-2)] =$$

$$-\frac{1}{4\pi^2}\sum_{n=1}^{\infty}\frac{\cos 2n\pi x - 1}{n^3}$$

Using (A.14) of Volume VI

$$B_n(1-x) = (-1)^n B_n(x)$$

this is easily simplified to

(4.4.229l)

$$\varsigma'(-2,x) + \varsigma'(-2,1-x) = 2\varsigma'(-2) - \frac{1}{2\pi^2}\sum_{n=1}^{\infty}\frac{\cos 2n\pi x - 1}{n^3} = -\frac{1}{2\pi^2}\sum_{n=1}^{\infty}\frac{\cos 2n\pi x}{n^3}$$

as previously noted by Adamchik [2a]. See also (4.3.166) in Volume II(a).

Letting $x = 1/2$ in (4.4.229k) and using (F.8b) we immediately obtain the well-known result (4.3.168d)

$$\varsigma'\left(-2, \frac{1}{2}\right) = \frac{1}{6} B_3\left(\frac{1}{2}\right) + \varsigma'(-2) + \frac{\varsigma_a(3) + \varsigma(3)}{8\pi^2} = \frac{3\varsigma(3)}{16\pi^2}$$

With $x = 1/4$ and $x = 1/4$ in (4.4.229k) we have

$$\varsigma'\left(-2, \frac{1}{4}\right) = \frac{1}{6} B_3\left(\frac{1}{4}\right) - \frac{\varsigma(3)}{(2\pi)^2} - 4\sum_{n=1}^{\infty} \frac{\log(2\pi) + \log n + \gamma - 1}{(2\pi n)^3} \sin(n\pi/2) - \frac{1}{8\pi^2} \sum_{n=1}^{\infty} \frac{\cos(n\pi/2) - 1}{n^3}$$

$$\varsigma'\left(-2, \frac{3}{4}\right) = \frac{1}{6} B_3\left(\frac{3}{4}\right) - \frac{\varsigma(3)}{(2\pi)^2} - 4\sum_{n=1}^{\infty} \frac{\log(2\pi) + \log n + \gamma - 1}{(2\pi n)^3} \sin(3n\pi/2) - \frac{1}{8\pi^2} \sum_{n=1}^{\infty} \frac{\cos(3n\pi/2) - 1}{n^3}$$

Since $\sin(n\pi/2) = -\sin(3n\pi/2)$ and $\cos(n\pi/2) = \cos(3n\pi/2)$, upon combining the above two equations we find that

$$\varsigma'\left(-2, \frac{1}{4}\right) + \varsigma'\left(-2, \frac{3}{4}\right) = \frac{3\varsigma(3)}{64\pi^2}$$

as previously discovered by Adamchik [2a].

Reference to (4.3.152a)

$$\log G(t+1) = t[\varsigma'(0,t) + \log A_0] - [\varsigma(-1,t) + \log A_1] - \varsigma'(-1,t)$$

would then get us part of the way to obtaining a Fourier series involving the Barnes multiple gamma $\log G(t+1)$. The pathway is completed by referring to a result by Espinosa and Moll in [59, Eq.(2.8)] where it is demonstrated how a Fourier series may be determined. See also (6.117k) in Volume V.

It may also be noted that by letting $p = 1 - 2N$ and $p = 2 - 2N$ in Hurwitz's formula (4.4.228i) we would obtain the Fourier series (6.25b) and (6.25b) for the Bernoulli polynomials.

We recall Elizalde's formula (6.117ai) from Volume V

$$\varsigma'(-1, t) =$$

$$-\varsigma(-1,t)\log t - \frac{1}{4}t^2 + \frac{1}{12} + \sum_{n=1}^{\infty}\frac{\sin 2n\pi t}{4\pi n^2} - \frac{1}{2\pi^2}\sum_{n=1}^{\infty}\frac{1}{n^2}[\cos(2n\pi t)Ci(2n\pi t) + \sin(2n\pi t)Si(2n\pi t)]$$

and in (4.4.229i) we showed that

$$\varsigma'(-1,t) = -2\sum_{n=1}^{\infty}\frac{\log(2\pi) + \log n + \gamma - 1}{(2n\pi)^2}\cos 2n\pi t + \sum_{n=1}^{\infty}\frac{\sin 2n\pi t}{4\pi n^2}$$

Therefore we deduce that (see also (4.4.228n))

(4.4.229m)

$$-\varsigma(-1,t)\log t - \frac{1}{4}t^2 + \frac{1}{12} - \frac{1}{2\pi^2}\sum_{n=1}^{\infty}\frac{1}{n^2}[\cos(2n\pi t)Ci(2n\pi t) + \sin(2n\pi t)Si(2n\pi t)] =$$

$$-\frac{1}{2\pi^2}\sum_{n=1}^{\infty}\frac{\log(2\pi) + \log n + \gamma - 1}{n^2}\cos 2n\pi t$$

Integration results in

(4.4.229n)

$$\frac{t}{72}\big[t(9-4t) + 6(t-1)(2t-1)\log t - 6\big] - \frac{1}{12}t^3 + \frac{1}{12}t$$

$$= \frac{1}{4\pi^3}\sum_{n=1}^{\infty}\frac{1}{n^3}[\sin(2n\pi t)Ci(2n\pi t) - \cos(2n\pi t)Si(2n\pi t)] - \frac{1}{4\pi^3}\sum_{n=1}^{\infty}\frac{\log(2\pi) + \log n + \gamma - 1}{n^3}\sin 2n\pi t$$

and additional identities may be obtained from successive integrations. Letting $t = 1$ results in

$$\sum_{n=1}^{\infty}\frac{Si(2n\pi)}{n^3} = \frac{1}{18}\pi^3$$

which we shall also determine in (6.117e). We may also compare (4.4.229n) with (6.117d).

□

We shall see in (6.5) and (6.5a) from Volume V that as a direct consequence of the Riemann-Lebesgue lemma we have

$$\frac{1}{2}\int_a^b p(x)\,dx = \sum_{n=0}^\infty \int_a^b p(x)\cos\alpha nx\,dx$$

$$\frac{1}{2}\int_a^b p(x)\cot(\alpha x/2)\,dx = \sum_{n=0}^\infty \int_a^b p(x)\sin\alpha nx\,dx$$

The above equations are valid provided (i) $p(x)$ is twice continuously differentiable on $[a,b]$ and either (ii) $\sin(\alpha x/2)\neq 0 \ \forall\, x\in[a,b]$ or, alternatively, (iii) if $\sin(\alpha a/2)=0$ then $p(a)=0$ also.

We now consider the case where $p(x)=\varsigma'(-1,x)-\varsigma'(-1)$, $[a,b]=[0,1]$ and $\alpha=2\pi$. We then have $p(0)=p(1)=0$ and therefore obtain

$$-\frac{1}{2}\int_0^1 [\varsigma'(-1,x)-\varsigma'(-1)]\,dx = \sum_{n=1}^\infty \int_0^1 [\varsigma'(-1,x)-\varsigma'(-1)]\cos 2\pi nx\,dx$$

Since from (4.3.145) $\int_0^1 \varsigma'(-1,x)\,dx=0$ and using (4.4.229i) we deduce that

$$\varsigma'(-1)=-2\sum_{n=1}^\infty \frac{\log(2\pi n)+\gamma-1}{(2\pi n)^2}$$

which is confirmed by (4.4.229iii). Unfortunately, as explained below, this proof, as presently constituted, is not rigorous.

If our choice of $p(x)$ was twice continuously differentiable, we would also have

(4.4.229ni) $\quad \dfrac{1}{2}\displaystyle\int_0^1 [\varsigma'(-1,x)-\varsigma'(-1)]\cot(\pi x)\,dx = \sum_{n=1}^\infty \int_0^1 [\varsigma'(-1,x)-\varsigma'(-1)]\sin 2\pi nx\,dx$

$$=\sum_{n=1}^\infty \int_0^1 \varsigma'(-1,x)\sin 2\pi nx\,dx$$

$$=\sum_{n=1}^\infty \frac{1}{8\pi n^2}=\frac{\pi}{48}$$

where we have employed (4.4.229hi). However, as will be seen later, this is **not** correct.

Now, using integration by parts we see that

$$I = \int_0^1 [\varsigma'(-1,x) - \varsigma'(-1)]\cot(\pi x)\,dx = \frac{1}{\pi}[\varsigma'(-1,x) - \varsigma'(-1)]\log\sin(\pi x)\Big|_0^1$$

$$-\frac{1}{\pi}\int_0^1 \log\sin(\pi x)\frac{\partial}{\partial x}\varsigma'(-1,x)\,dx$$

The integrated part may be represented by

$$c = [\varsigma'(-1,x) - \varsigma'(-1)]\log\sin(\pi x)\Big|_0^1 = \lim_{a\to 1}[\varsigma'(-1,a) - \varsigma'(-1)]\log\sin(\pi a)$$

$$-\lim_{a\to 0}[\varsigma'(-1,a) - \varsigma'(-1)]\log\sin(\pi a)$$

We see that

$$\lim_{a\to 1}[\varsigma'(-1,a) - \varsigma'(-1)]\log\sin(\pi a) = \lim_{a\to 0}[\varsigma'(-1,1-a) - \varsigma'(-1)]\log\sin\pi(1-a)$$

$$= \lim_{a\to 0}[\varsigma'(-1,1-a) - \varsigma'(-1)]\log\sin(\pi a)$$

and therefore we have

$$c = \lim_{a\to 0}[\varsigma'(-1,1-a) - \varsigma'(-1,a)]\log\sin(\pi a)$$

Using (4.4.229iv) this becomes

$$c = -\frac{1}{2\pi}\lim_{a\to 0}[\mathrm{Cl}(2\pi a)\log\sin(\pi a)]$$

and we now consider the limit

$$\lim_{a\to 0}[\sin(2\pi ka)\log\sin(\pi a)] = \lim_{a\to 0}\left[\frac{\sin(2\pi ka)}{2\pi ka}2\pi ka\log\sin(\pi a)\right] = 0$$

Hence it is easily seen that $c = 0$ and we therefore obtain

$$I = -\frac{1}{\pi}\int_0^1 \log\sin(\pi x)\frac{\partial}{\partial x}\varsigma'(-1,x)\,dx$$

From (4.3.126di) we saw that

$$\frac{\partial}{\partial x}\varsigma'(-1,x) = -\varsigma(-1,x) + \varsigma'(0,x)$$

whereupon, using Lerch's identity (4.3.116), this becomes

(4.4.229o) $$\frac{\partial}{\partial x}\varsigma'(-1,x) = -\varsigma(-1,x) + \log\Gamma(x) - \frac{1}{2}\log(2\pi)$$

and we then obtain

$$I = -\frac{1}{\pi}\int_0^1 \log\sin(\pi x)\left[-\varsigma(-1,x) + \log\Gamma(x) - \frac{1}{2}\log(2\pi)\right]dx$$

Then using (4.3.110a) $\varsigma(-1,x) = -\frac{1}{2}B_2(x)$ we get

$$I = -\frac{1}{2\pi}\int_0^1 B_2(x)\log\sin(\pi x)\,dx - \frac{1}{\pi}\int_0^1 \log\Gamma(x)\log\sin(\pi x)\,dx + \frac{1}{2\pi}\log(2\pi)\int_0^1 \log\sin(\pi x)\,dx$$

From the interesting paper by Espinosa and Moll [59] we learn that

(4.4.229p) $$\int_0^1 B_2(x)\log\sin(\pi x)\,dx = \frac{\varsigma(3)}{4\pi^2}$$

which we have also proved in Volume V. Espinosa and Moll [59] also show that

(4.4.229q) $$\int_0^1 \log\Gamma(x)\log\sin(\pi x)\,dx = -\frac{1}{2}\log^2 2 - \frac{1}{2}\log 2\log\pi - \frac{\pi^2}{24}$$

and of course, more than 300 years ago, Euler told us that (see (3.2) from Volume I)

(4.4.229r) $$\int_0^1 \log\sin(\pi x)\,dx = -\log 2$$

This gives us

$$I = -\frac{\varsigma(3)}{8\pi^3} + \frac{1}{2\pi}\log^2 2 + \frac{1}{2\pi}\log 2\log\pi + \frac{\pi}{24} - \frac{1}{2\pi}\log(2\pi)\log 2$$

and hence we see that

(4.4.229t) $$\int_0^1 [\varsigma'(-1,x) - \varsigma'(-1)]\cot(\pi x)\,dx = -\frac{\varsigma(3)}{8\pi^3} + \frac{\pi}{24}$$

which does **not** agree with the previous result (4.4.229ni).

Using the Gosper/Vardi functional equation (4.3.126)

$$\log G(x+1) - x\log\Gamma(x) = \varsigma'(-1) - \varsigma'(-1,x)$$

we obtain

(4.4.229u) $$\int_0^1 [\log G(x+1) - x\log\Gamma(x)]\cot(\pi x)\,dx = \frac{\varsigma(3)}{8\pi^3} - \frac{\pi}{24}$$

It should however be noted that we are not permitted to apply (6.5a) because (4.4.229o) shows that the derivative of $p(x) = [\varsigma'(-1,x) - \varsigma'(-1)]$ is not finite at $x = 0$. On the other hand, the functions $x^2[\varsigma'(-1,x) - \varsigma'(-1)]$ and $[\varsigma'(-1,x) - \varsigma'(-1)] - x\log\Gamma(x)$ do satisfy that condition. Differentiating (4.4.229i) also indicates the nature of the difficulty at $x = 0$

$$\frac{\partial}{\partial x}\varsigma'(-1,x) = 2\sum_{n=1}^{\infty}\frac{\log(2\pi n) + \gamma - 1}{2\pi n}\sin 2n\pi x + \sum_{n=1}^{\infty}\frac{\cos 2n\pi x}{2n}$$

because we end up with a divergent series.

It may be possible to extend this method with a little more ingenuity!!

We could perhaps refer again to Hurwitz's formula (4.4.228i)

$$\varsigma(p,x) = 2\Gamma(1-p)\left[\sin\left(\frac{\pi p}{2}\right)\sum_{n=1}^{\infty}\frac{\cos 2n\pi x}{(2\pi n)^{1-p}} + \cos\left(\frac{\pi p}{2}\right)\sum_{n=1}^{\infty}\frac{\sin 2n\pi x}{(2\pi n)^{1-p}}\right]$$

and from this find the Fourier coefficients for $x\varsigma(p,x)$.

We may however safely apply equations (6.5) and (6.5a) of Volume V to $p(x) = \varsigma'(-2m,x) - \varsigma'(-2m)$. For example we have

$$-\frac{1}{2}\int_0^1 [\varsigma'(-2,x) - \varsigma'(-2)]\,dx = \sum_{n=1}^{\infty}\int_0^1 [\varsigma'(-2,x) - \varsigma'(-2)]\cos 2\pi nx\,dx$$

and by reference to (4.4.229g) we have

$$\int_0^1 \varsigma'(-2,x)\cos 2\pi nx\,dx = -\frac{\Gamma(3)}{16\pi^2 n^3}$$

We then obtain the familiar result (4.4.228s)

$$\varsigma'(-2) = -\frac{\varsigma(3)}{8\pi^2}$$

More generally from (4.4.229g) we have where $m$ is an integer

$$\int_0^1 \varsigma'(-2m,t)\cos 2n\pi t\,dt = \frac{\pi}{2}\cos(m\pi)\frac{\Gamma(1+2m)}{(2\pi n)^{2m+1}}$$

and hence we obtain

$$\frac{1}{2}\varsigma'(-2m) = \sum_{n=1}^{\infty}\frac{\pi}{2}(-1)^m\frac{\Gamma(1+2m)}{(2\pi n)^{2m+1}}$$

and this clearly results in (F.8a) of Volume V

$$\varsigma'(-2m) = (-1)^m\frac{(2m)!}{2(2\pi)^{2m}}\varsigma(2m+1)$$

We now consider the case where $p(x) = \varsigma(1-2m,x) - \varsigma(1-2m)$, $[a,b] = [0,1]$ and $\alpha = 2\pi$. We then have $p(0) = p(1) = 0$ and we therefore obtain

$$-\frac{1}{2}\int_0^1 [\varsigma(1-2m,x) - \varsigma(1-2m)]\,dx = \sum_{n=1}^{\infty}\int_0^1 [\varsigma(1-2m,x) - \varsigma(1-2m)]\cos 2\pi nx\,dx$$

$$= \sum_{n=1}^{\infty}\int_0^1 \varsigma(1-2m,x)\cos 2\pi nx\,dx$$

Since from (4.3.134a) $\int_0^1 \varsigma(1-2m,x)\,dx = 0$ and using (4.4.229b)

$$\int_0^1 \varsigma(p,x)\cos 2n\pi x\,dx = \frac{\Gamma(1-p)}{(2\pi n)^{1-p}}\sin\left(\frac{\pi p}{2}\right)$$

we have

$$\int_0^1 \varsigma(1-2m,x)\cos 2n\pi x\,dx = \frac{\Gamma(2m)}{(2\pi n)^{2m}}\sin\left(\frac{\pi(1-2m)}{2}\right) = (-1)^m\frac{\Gamma(2m)}{(2\pi n)^{2m}}$$

We have by reference to (4.3.110) in Volume II(a)

$$\int_0^1 [\varsigma(1-2m,x)-\varsigma(1-2m)]\,dx = -\frac{1}{2m}\int_0^1 [B_{2m}(x)-B_{2m}]\,dx$$

$$= \frac{B_{2m}}{2m}$$

and we therefore obtain

$$-\frac{B_{2m}}{4m} = (-1)^m\sum_{n=1}^{\infty}\frac{\Gamma(2m)}{(2\pi n)^{2m}} = \frac{(-1)^m}{(2\pi)^{2m}}\Gamma(2m)\varsigma(2m)$$

This then gives us Euler's formula (1.7)

$$\varsigma(2m) = (-1)^{m+1}\frac{(2\pi)^{2m}B_{2m}}{2.(2m)!}$$

We also have

$$\frac{1}{2}\int_0^1 [\varsigma(1-2m,x)-\varsigma(1-2m)]\cot(\pi x)\,dx = \sum_{n=1}^{\infty}\int_0^1 [\varsigma(1-2m,x)-\varsigma(1-2m)]\sin 2\pi nx\,dx$$

$$= \sum_{n=1}^{\infty}\int_0^1 \varsigma(1-2m,x)\sin 2\pi nx\,dx$$

and using (4.4.229a) we have

$$\int_0^1 \varsigma(1-2m,x)\sin 2\pi nx\,dx = \frac{\Gamma(2m)}{(2\pi n)^{2m}}\cos\left(\frac{\pi(1-2m)}{2}\right) = 0$$

Therefore we see that

$$\frac{1}{2m}\int_0^1 [B_{2m}(x)-B_{2m}]\cot(\pi x)\,dx = 0$$

Integration by parts gives us

(4.4.229v)

$$\int_0^1 [B_{2m}(x) - B_{2m}] \cot(\pi x)\, dx = \frac{1}{\pi} [B_{2m}(x) - B_{2m}] \log \sin(\pi x)\Big|_0^1 - \frac{2m}{\pi} \int_0^1 B_{2m-1}(x) \log \sin(\pi x)\, dx$$

Since $B_{2m}(0) = B_{2m}(1) = B_{2m}$ we have

$$= -\frac{2m}{\pi} \int_0^1 B_{2m-1}(x) \log \sin(\pi x)\, dx$$

and we have therefore proved that the final integral is equal to zero (which will be seen in (6.35) in Volume V). Similarly, considering the integral

$$\int_0^1 [\varsigma(-2m, x) - \varsigma(-2m)] \cot(\pi x)\, dx$$

provides an alternative derivation of (6.34)

$$\int_0^1 B_{2n}(x) \log \sin \pi x\, dx = (-1)^n \frac{(2n)!\, \varsigma(2n+1)}{(2\pi)^{2n}}$$

FURTHER APPEARANCES OF THE RIEMANN FUNCTIONAL EQUATION

Apostol [14b] has shown that for $-1 < s < 0$

$$\lim_{x \to 0+} \varsigma(s, x) = \varsigma(s)$$

and we therefore have $p(0) = p(1) = 0$ where $p(x) = \varsigma(s, x) - \varsigma(s)$. Therefore, assuming that (6.5) applies to this function $p(x)$, we have

$$-\frac{1}{2} \int_0^1 [\varsigma(s, x) - \varsigma(s)]\, dx = \sum_{n=1}^\infty \int_0^1 [\varsigma(s, x) - \varsigma(s)] \cos 2\pi n x\, dx$$

and therefore we see that

$$\frac{1}{2} \varsigma(s) = \sum_{n=1}^\infty \int_0^1 \varsigma(s, x) \cos 2\pi n x\, dx$$

Using (4.4.229b) this becomes

$$\frac{1}{2}\varsigma(s) = \sum_{n=1}^{\infty} \frac{\Gamma(1-s)}{(2\pi n)^{1-s}} \sin\left(\frac{\pi s}{2}\right) = \frac{\Gamma(1-s)\varsigma(1-s)}{(2\pi)^{1-s}} \sin\left(\frac{\pi s}{2}\right)$$

and hence we recover the Riemann functional equation.

Similarly we have using (6.5a)

$$\frac{1}{2}\int_0^1 [\varsigma(s,x) - \varsigma(s)] \cot \pi x \, dx = \sum_{n=1}^{\infty} \int_0^1 [\varsigma(s,x) - \varsigma(s)] \sin 2\pi n x \, dx$$

and, with integration by parts, the left-hand side becomes

$$= \frac{1}{2\pi} [\varsigma(s,x) - \varsigma(s)] \log \sin \pi x \Big|_0^1 - \frac{1}{2\pi} \int_0^1 \frac{\partial}{\partial x} \varsigma(s,x) \log \sin \pi x \, dx$$

$$= \frac{s}{2\pi} \int_0^1 \varsigma(s+1,x) \log \sin \pi x \, dx$$

Espinosa and Moll [59] have shown that for $s < 0$

(4.4.229w) $$\int_0^1 \log \sin \pi x \, \varsigma(s,x) \, dx = -\frac{1}{2} \frac{\varsigma(s)\varsigma(2-s)}{\varsigma(1-s)}$$

and we then have for $s < -1$

$$\frac{1}{2}\int_0^1 [\varsigma(s,x) - \varsigma(s)] \cot \pi x \, dx = -\frac{1}{4\pi} \frac{s\varsigma(s+1)\varsigma(1-s)}{\varsigma(-s)}$$

For the right-hand side we have

$$\sum_{n=1}^{\infty} \int_0^1 [\varsigma(s,x) - \varsigma(s)] \sin 2\pi n x \, dx = \sum_{n=1}^{\infty} \int_0^1 \varsigma(s,x) \sin 2\pi n x \, dx$$

and using (4.4.229a) this becomes

$$= \sum_{n=1}^{\infty} \frac{\Gamma(1-s)}{(2\pi n)^{1-s}} \cos\left(\frac{\pi s}{2}\right) = \frac{\Gamma(1-s)}{(2\pi)^{1-s}} \cos\left(\frac{\pi s}{2}\right) \varsigma(1-s)$$

Therefore we have

$$-\frac{1}{4\pi}\frac{s\varsigma(s+1)\varsigma(1-s)}{\varsigma(-s)}=\frac{\Gamma(1-s)}{(2\pi)^{1-s}}\cos\left(\frac{\pi s}{2}\right)\varsigma(1-s)$$

and letting $s \to -p$ this becomes

$$\frac{1}{4\pi}p\varsigma(1-p)=\frac{\Gamma(1+p)}{(2\pi)^{1+p}}\cos\left(\frac{\pi p}{2}\right)\varsigma(p)$$

which again produces the Riemann functional equation (F.1)

$$\varsigma(1-p)=2(2\pi)^{-p}\Gamma(p)\cos(\pi s/2)\varsigma(p)$$

$\square$

We have seen in (4.3.132) that

$$\lim_{x\to0}[\varsigma'(-n,x)-\varsigma'(-n)]=0$$

and we now conjecture that for $p<0$

$$\lim_{x\to0}[\varsigma'(p,x)-\varsigma'(p)]=0$$

We have using integration by parts

$$\int_0^1\cot\pi x[\varsigma'(p,x)-\varsigma'(p)]dx=$$

$$\frac{1}{\pi}\log\sin\pi x[\varsigma'(p,x)-\varsigma'(p)]\Big|_0^1-\frac{1}{\pi}\int_0^1\log\sin\pi x\frac{\partial}{\partial x}[\varsigma'(p,x)-\varsigma'(p)]dx$$

$$=-\frac{1}{\pi}\int_0^1\log\sin\pi x\frac{\partial}{\partial x}\varsigma'(p,x)\,dx$$

As noted previously we have

$$\frac{\partial}{\partial x}\frac{\partial}{\partial p}\varsigma(p,x)=\frac{\partial}{\partial p}\frac{\partial}{\partial x}\varsigma(p,x)=-\frac{\partial}{\partial p}[p\varsigma(p+1,x)]$$

$$=-\varsigma(p+1,x)-p\frac{\partial}{\partial p}\varsigma(p+1,x)$$

$$= -\varsigma(p+1,x) - p\varsigma'(p+1,x)$$

Hence we have

$$\int_0^1 \cot \pi x \, [\varsigma'(p,x) - \varsigma'(p)]dx = \frac{1}{\pi}\int_0^1 \log \sin \pi x \, [\varsigma(p+1,x) + p\varsigma'(p+1,x)]dx$$

As mentioned earlier, Espinosa and Moll [59] have shown that for $p < 0$

$$\int_0^1 \log \sin \pi x \, \varsigma(p,x)\,dx = -\frac{1}{2}\frac{\varsigma(p)\varsigma(2-p)}{\varsigma(1-p)}$$

and differentiation results in

(4.4.229x)
$$\int_0^1 \log \sin \pi x \, \varsigma'(p,x)\,dx = -\frac{1}{2}\frac{\varsigma(1-p)[-\varsigma(p)\varsigma'(2-p)+\varsigma'(p)\varsigma(2-p)]+\varsigma(p)\varsigma(2-p)\varsigma'(1-p)}{\varsigma^2(1-p)}$$

We then have for $p < -1$

$$\int_0^1 \log \sin \pi x \, \varsigma(p+1,x)\,dx = -\frac{1}{2}\frac{\varsigma(p+1)\varsigma(1-p)}{\varsigma(-p)}$$

$$\int_0^1 \log \sin \pi x \, \varsigma'(p+1,x)\,dx = -\frac{1}{2}\frac{\varsigma(-p)[-\varsigma(p+1)\varsigma'(1-p)+\varsigma'(p+1)\varsigma(1-p)]+\varsigma(p+1)\varsigma(1-p)\varsigma'(-p)}{\varsigma^2(-p)}$$

We then obtain

$$\pi\int_0^1 \cot \pi x \, [\varsigma'(p,x) - \varsigma'(p)]dx = -\frac{1}{2}\frac{\varsigma(p+1)\varsigma(1-p)}{\varsigma(-p)}$$

$$-\frac{p}{2}\frac{\varsigma(-p)[-\varsigma(p+1)\varsigma'(1-p)+\varsigma'(p+1)\varsigma(1-p)]+\varsigma(p+1)\varsigma(1-p)\varsigma'(-p)}{\varsigma^2(-p)}$$

Then, assuming that (6.5a) may be validly employed, we obtain

$$\frac{1}{2}\pi\int_0^1 \cot \pi x \, [\varsigma'(p,x) - \varsigma'(p)]dx = \pi\sum_{n=0}^{\infty}\int_0^1 [\varsigma'(p,x) - \varsigma'(p)]\sin 2\pi nx \, dx$$

$$= \pi \sum_{n=1}^{\infty} \int_0^1 \varsigma'(p,x) \sin 2\pi nx \, dx$$

Then using (4.4.229e)

$$\int_0^1 \varsigma'(p,x) \sin 2n\pi x \, dx = \left[ \log(2\pi n) - \psi(1-p) + \frac{\pi}{2} \tan\left(\frac{\pi p}{2}\right) \right] \frac{\Gamma(1-p)}{(2\pi n)^{1-p}} \cos\left(\frac{\pi p}{2}\right)$$

we get

$$\frac{1}{2} \pi \int_0^1 \cot \pi x \left[ \varsigma'(p,x) - \varsigma'(p) \right] dx = \pi \sum_{n=1}^{\infty} \left[ \log(2\pi n) - \psi(1-p) + \frac{\pi}{2} \tan\left(\frac{\pi p}{2}\right) \right] \frac{\Gamma(1-p)}{(2\pi n)^{1-p}} \cos\left(\frac{\pi p}{2}\right)$$

This results in

$$-\frac{1}{2} \frac{\varsigma(p+1)\varsigma(1-p)}{\varsigma(-p)} - \frac{p}{2} \frac{\varsigma(-p)[-\varsigma(p+1)\varsigma'(1-p) + \varsigma'(p+1)\varsigma(1-p)] + \varsigma(p+1)\varsigma(1-p)\varsigma'(-p)}{\varsigma^2(-p)} =$$

$$2\pi \sum_{n=1}^{\infty} \left[ \log(2\pi n) - \psi(1-p) + \frac{\pi}{2} \tan\left(\frac{\pi p}{2}\right) \right] \frac{\Gamma(1-p)}{(2\pi n)^{1-p}} \cos\left(\frac{\pi p}{2}\right)$$

which may be written as a functional equation

(4.4.229x)

$$\varsigma(p+1)\varsigma(1-p)\varsigma(-p) + p\varsigma(-p)[-\varsigma(p+1)\varsigma'(1-p) + \varsigma'(p+1)\varsigma(1-p)] + p\varsigma(p+1)\varsigma(1-p)\varsigma'(-p) =$$

$$-4\pi \varsigma^2(-p) \left[ \log(2\pi) - \psi(1-p) + \frac{\pi}{2} \tan\left(\frac{\pi p}{2}\right) \right] \frac{\Gamma(1-p)}{(2\pi)^{1-p}} \cos\left(\frac{\pi p}{2}\right) \varsigma(1-p)$$

$$+4\pi \varsigma^2(-p) \frac{\Gamma(1-p)}{(2\pi)^{1-p}} \cos\left(\frac{\pi p}{2}\right) \varsigma'(1-p)$$

For example, with $p = -3$ and using (3.11c) $\varsigma(-2m) = 0$ we get the well-known result (F.8b) from Volume VI

$$\varsigma'(-2) = -\frac{\varsigma(3)}{4\pi^2}$$

and with $p = -2$ we obtain

$$\varsigma(-1)\varsigma(3)\varsigma(2) - 2\varsigma(2)[-\varsigma(-1)\varsigma'(3) + \varsigma'(-1)\varsigma(3)] - 2\varsigma(-1)\varsigma(3)\varsigma'(2) =$$

$$-\frac{\pi^2}{36}[\varsigma'(3) - \varsigma(3)\log(2\pi) + \varsigma(3)\psi(3)]$$

It is easily seen that the terms involving $\varsigma'(3)$ cancel and we have

$$\varsigma(-1)\varsigma(3)\varsigma(2) - 2\varsigma(2)\varsigma'(-1)\varsigma(3) - 2\varsigma(-1)\varsigma(3)\varsigma'(2) = -\frac{\pi^2}{36}[-\varsigma(3)\log(2\pi) + \varsigma(3)\psi(3)]$$

and with some simple algebra we simply obtain another well-known result

$$\varsigma'(-1) = \frac{1}{12}(1 - \gamma - \log 2\pi) + \frac{1}{2\pi^2}\varsigma'(2)$$

When $p = -4$ we have

$$\varsigma(-3)\varsigma(5)\varsigma(4) - 4\varsigma(4)[-\varsigma(-3)\varsigma'(5) + \varsigma'(-3)\varsigma(5)] - 4\varsigma(-3)\varsigma(5)\varsigma'(4) =$$

$$-4\pi\varsigma^2(4)\left[\log(2\pi) - \psi(5)\right]\frac{4!}{(2\pi)^5}\varsigma(5) + 4\pi\varsigma^2(4)\frac{4!}{(2\pi)^5}\varsigma'(5)$$

We know from (F.5) that

$$\varsigma'(1-p) = \left[\log(2\pi) - \psi(p) + \frac{\pi}{2}\tan\left(\frac{\pi p}{2}\right) - \frac{\varsigma'(p)}{\varsigma(p)}\right]\varsigma(1-p)$$

and with $p = 4$ we obtain

$$\varsigma'(-3) = \left[\log(2\pi) - \psi(4) - \frac{\varsigma'(4)}{\varsigma(4)}\right]\varsigma(-3)$$

From (F.4a) we have $\varsigma(1-2n) = -\frac{B_{2n}}{2n}$ and thus $\varsigma(-3) = -\frac{B_4}{4} = \frac{1}{120}$. Therefore we have

$$\varsigma'(-3) = \frac{1}{120}\left[\log(2\pi) - \frac{11}{6} + \gamma - \frac{\varsigma'(4)}{\varsigma(4)}\right]$$

We then have

$$\frac{1}{120}\varsigma(5)\varsigma(4) - 4\varsigma(4)\left\{-\frac{1}{120}\varsigma'(5) + \frac{1}{120}\left[\log(2\pi) - \frac{11}{6} + \gamma - \frac{\varsigma'(4)}{\varsigma(4)}\right]\varsigma(5)\right\} - \frac{4}{120}\varsigma(5)\varsigma'(4) =$$

$$4\pi\,\varsigma^2(4)\frac{4!}{(2\pi)^5}\left[\varsigma'(5) - \log(2\pi)\varsigma(5) + \psi(5)\varsigma(5)\right]$$

and with a little more algebra we simply find that everything cancels out! I initially thought that (4.4.229x) might be new, but in fact it seems to be a heavily disguised version of the derivative of the Riemann functional equation.

□

Using (6.5) gives us for $p < 0$

$$-\frac{1}{2}\int_0^1[\varsigma'(p,x) - \varsigma'(p)]dx = \sum_{n=1}^{\infty}\int_0^1[\varsigma'(p,x) - \varsigma'(p)]\cos 2\pi nx\,dx$$

$$= \sum_{n=1}^{\infty}\int_0^1\varsigma'(p,x)\cos 2\pi nx\,dx$$

and the left-hand side is equal to $\frac{1}{2}\varsigma'(p)$. For the right-hand side we use (4.4.229g)

$$\int_0^1\varsigma'(p,t)\cos 2n\pi t\,dt = \left[\log(2\pi n) - \psi(1-p) + \frac{\pi}{2}\cot\left(\frac{\pi p}{2}\right)\right]\frac{\Gamma(1-p)}{(2\pi n)^{1-p}}\sin\left(\frac{\pi p}{2}\right)$$

and we therefore obtain

$$\frac{1}{2}\varsigma'(p) = \left[\log(2\pi) - \psi(1-p) + \frac{\pi}{2}\cot\left(\frac{\pi p}{2}\right)\right]\frac{\Gamma(1-p)\varsigma(1-p)}{(2\pi)^{1-p}}\sin\left(\frac{\pi p}{2}\right)$$

$$+ \frac{\Gamma(1-p)}{(2\pi)^{1-p}}\sin\left(\frac{\pi p}{2}\right)\sum_{n=1}^{\infty}\frac{\log n}{n^{1-p}}$$

or alternatively

$$\varsigma'(p) = \left[\log(2\pi) - \psi(1-p) + \frac{\pi}{2}\cot\left(\frac{\pi p}{2}\right)\right]\frac{2\Gamma(1-p)\varsigma(1-p)}{(2\pi)^{1-p}}\sin\left(\frac{\pi p}{2}\right)$$

$$- \frac{2\Gamma(1-p)\varsigma'(1-p)}{(2\pi)^{1-p}}\sin\left(\frac{\pi p}{2}\right)$$

Employing (F.1a) this becomes

$$\varsigma'(p) = \left[ \log(2\pi) - \psi(1-p) + \frac{\pi}{2} \cot\left(\frac{\pi p}{2}\right) \right] \varsigma(p) - \frac{2\Gamma(1-p)\varsigma(p)\varsigma'(1-p)}{(2\pi)^{1-p}\varsigma(1-p)}$$

and if $\varsigma(p) \neq 0$ this may be written as

$$\frac{\varsigma'(p)}{\varsigma(p)} = \log(2\pi) - \psi(1-p) + \frac{\pi}{2} \cot\left(\frac{\pi p}{2}\right) - \frac{\varsigma'(1-p)}{\varsigma(1-p)}$$

Using $\psi(1-p) - \psi(p) = \pi \cot \pi p$ we simply find that

$$\frac{\varsigma'(1-p)}{\varsigma(1-p)} = \log(2\pi) - \psi(p) + \frac{\pi}{2} \tan\left(\frac{\pi p}{2}\right) - \frac{\varsigma'(p)}{\varsigma(p)}$$

as noted in (F.5).

□

We now multiply (4.4.229i) by $t$ and integrate to obtain

$$\int_0^x t\varsigma'(-1,t)\,dt = -2x\sum_{n=1}^{\infty} \frac{\log(2\pi n) + \gamma - 1}{(2\pi n)^3} \sin 2n\pi x - 2\sum_{n=1}^{\infty} \frac{\log(2\pi n) + \gamma - 1}{(2\pi n)^4} \cos 2n\pi x$$

$$+ 2\sum_{n=1}^{\infty} \frac{\log(2\pi n) + \gamma - 1}{(2\pi n)^4} - \pi x\sum_{n=1}^{\infty} \frac{\cos 2n\pi x}{(2n\pi)^3} + \pi \sum_{n=1}^{\infty} \frac{\sin 2n\pi x}{(2n\pi)^4}$$

and in particular with $x = 1$ we have

$$(4.4.229y) \qquad \int_0^1 t\varsigma'(-1,t)\,dt = -\pi \sum_{n=1}^{\infty} \frac{1}{(2n\pi)^3} = -\frac{\varsigma(3)}{8\pi^2}$$

Alternatively, using integration by parts and Adamchik's formula (4.3.131)

$$n\int_0^x \varsigma'(1-n,u)\,du = \frac{B_{n+1} - B_{n+1}(x)}{n(n+1)} + \varsigma'(-n,x) - \varsigma'(-n)$$

$$\int_0^x \varsigma'(-1,u)\,du = -\frac{B_3(x)}{12} + \frac{1}{2}\varsigma'(-2,x) - \frac{1}{2}\varsigma'(-2)$$

we have

$$\int\limits_0^1 t\varsigma'(-1,t)dt = t\left[-\frac{B_3(t)}{12}+\frac{1}{2}\varsigma'(-2,t)-\frac{1}{2}\varsigma'(-2)\right]\Bigg|_0^1 - \int\limits_0^1\left[-\frac{B_3(t)}{12}+\frac{1}{2}\varsigma'(-2,t)-\frac{1}{2}\varsigma'(-2)\right]dt$$

$$=-\int\limits_0^1\left[-\frac{B_3(t)}{12}+\frac{1}{2}\varsigma'(-2,t)-\frac{1}{2}\varsigma'(-2)\right]dt$$

Since $\int\limits_0^x\varsigma'(-2,u)\,du = \frac{B_4 - B_4(x)}{36}+\frac{1}{3}\varsigma'(-3,x)-\frac{1}{3}\varsigma'(-3)$ we obtain $\int\limits_0^1\varsigma'(-2,u)\,du = 0$

and hence we get

$$\int\limits_0^1 t\varsigma'(-1,t)dt = \frac{1}{2}\varsigma'(-2)$$

We already know from (F.8b) in Volume VI that

$$\varsigma'(-2) = -\frac{\varsigma(3)}{4\pi^2}$$

and we therefore get the same result as above.

We now refer back to (4.4.229hi)

$$\int\limits_0^1\varsigma'(-1,t)\sin 2n\pi t\,dt = \frac{1}{8\pi n^2}$$

and make the summation

(4.4.229z)   $$\sum_{n=1}^\infty\frac{1}{n}\int\limits_0^1\varsigma'(-1,t)\sin 2n\pi t\,dt = \frac{\varsigma(3)}{8\pi}$$

We now assume that we may write

$$\sum_{n=1}^\infty\frac{1}{n}\int\limits_0^1\varsigma'(-1,t)\sin 2n\pi t\,dt = \int\limits_0^1\varsigma'(-1,t)\sum_{n=1}^\infty\frac{\sin 2n\pi t}{n}\,dt$$

and that we may legitimately employ (7.5) from Volume V for $0 < t < 1$

$$\sum_{n=1}^{\infty} \frac{\sin 2n\pi t}{n} = \frac{\pi}{2}(1-2t)$$

notwithstanding that the Fourier series is not valid at the integration end points. Making these assumptions gives us

$$\sum_{n=1}^{\infty} \frac{1}{n} \int_0^1 \varsigma'(-1,t) \sin 2n\pi t \, dt = \frac{\pi}{2} \int_0^1 \varsigma'(-1,t)(1-2t) \, dt$$

$$= \frac{\pi}{2} \int_0^1 \varsigma'(-1,t) \, dt - \pi \int_0^1 t \varsigma'(-1,t) \, dt$$

$$= -\pi \int_0^1 t \varsigma'(-1,t) \, dt$$

where we have used (4.3.145). Then using the result (4.4.229y) just obtained above we obtain

$$= -\frac{\pi}{2} \varsigma'(-2)$$

and this brings us back to (4.4.229z). In the same Eulerian fashion we could perhaps multiply (4.4.229hii) by $1/n$ and then perform the same summation technique to obtain the integral $\int_0^1 \varsigma'(-1,t) \log \sin \pi t \, dt$.

And now for something completely different!

## MORE IDENTITIES INVOLVING HARMONIC NUMBERS AND POLYLOGARITHMS

**Theorem 4.11:**

$$(4.4.230) \qquad \sum_{n=1}^{\infty} \frac{H_n^{(r)}}{n^q} = \varsigma(r)\varsigma(q) - \frac{(-1)^{r-1}}{(r-1)!} \int_0^1 \frac{\log^{r-1} x \, Li_q(x)}{1-x} \, dx \qquad \text{, for } q, r \geq 2$$

**Proof:**

This is a well-known relation and can be found for example in [126, p.157]. The following proof is due to Freitas [69a].

In (3.86) we gave an elementary derivation of

$$\int_0^1 x^{k-1} \log^{r-1} x \, dx = (-1)^{r-1} \frac{(r-1)!}{k^r} \qquad (k \geq 1, r \geq 1)$$

and therefore completing the summation we get

$$H_n^{(r)} = \sum_{k=1}^n \frac{1}{k^r} = \frac{(-1)^{r-1}}{(r-1)!} \sum_{k=1}^n \int_0^1 x^{k-1} \log^{r-1} x \, dx$$

(4.4.230a) $$\qquad H_n^{(r)} = \frac{(-1)^{r-1}}{(r-1)!} \int_0^1 \log^{r-1}(x) \frac{1-x^n}{1-x} \, dx$$

Hence upon making a further summation we have

$$\sum_{n=1}^\infty \frac{H_n^{(r)}}{n^q} = \frac{(-1)^{r-1}}{(r-1)!} \int_0^1 \frac{\log^{r-1} x}{1-x} \sum_{n=1}^\infty \frac{1-x^n}{n^q} \, dx$$

$$= \frac{(-1)^{r-1}}{(r-1)!} \int_0^1 \frac{\log^{r-1} x}{1-x} \left[ \varsigma(q) - Li_q(x) \right] dx \qquad , q \geq 2$$

$$= \varsigma(r)\varsigma(q) - \frac{(-1)^{r-1}}{(r-1)!} \int_0^1 \frac{\log^{r-1} x}{1-x} Li_q(x) \, dx$$

where the last step follows from (4.4.51b)

(4.4.231) $$\qquad \int_0^1 \frac{\log^{r-1} x}{1-x} \, dx = (-1)^{r-1} \varsigma(r)\Gamma(r) \qquad , r \geq 2$$

(alternatively, the above equation can be easily derived using (4.4.230a) and taking the limit as $n \rightarrow \infty$).

For example, with $r = q = 2$ we obtain

(4.4.232) $$\qquad \int_0^1 \frac{\log^2 x}{1-x} Li_2(x) \, dx = \sum_{n=1}^\infty \frac{H_n^{(2)}}{n^2} - \varsigma^2(2)$$

By definition we have for $\text{Re}(q) > 1$

$$\sum_{n=1}^{\infty} \frac{H_n^{(p)}}{n^q} = \sum_{n=1}^{\infty} \frac{1}{n^q} \sum_{k=1}^{n} \frac{1}{k^p}$$

and using (3.23), provided $\mathrm{Re}\,(p) > 1$, this becomes

$$= \sum_{n=1}^{\infty} \frac{1}{n^p} \sum_{k=n}^{\infty} \frac{1}{k^q}$$

$$= \sum_{n=1}^{\infty} \frac{1}{n^p} \left[ \sum_{k=1}^{\infty} \frac{1}{k^q} - \sum_{k=1}^{n-1} \frac{1}{k^q} \right]$$

$$= \varsigma(p)\varsigma(q) - \sum_{n=1}^{\infty} \frac{1}{n^p} \sum_{k=1}^{n-1} \frac{1}{k^q}$$

$$= \varsigma(p)\varsigma(q) - \sum_{n=1}^{\infty} \frac{1}{n^p} \left[ \sum_{k=1}^{n} \frac{1}{k^q} - \frac{1}{k^q} \right]$$

$$= \varsigma(p)\varsigma(q) - \sum_{n=1}^{\infty} \frac{H_n^{(q)}}{n^p} + \sum_{n=1}^{\infty} \frac{1}{n^{p+q}}$$

We therefore have the symmetrical result (this is generalised in (4.4.247c))

(4.4.232a) $$\sum_{n=1}^{\infty} \frac{H_n^{(p)}}{n^q} + \sum_{n=1}^{\infty} \frac{H_n^{(q)}}{n^p} = \varsigma(p)\varsigma(q) + \varsigma(p+q)$$

and in particular we get

(4.4.232b) $$\sum_{n=1}^{\infty} \frac{H_n^{(p)}}{n^p} = \frac{1}{2} \left[ \varsigma^2(p) + \varsigma(2p) \right]$$

From (4.4.232) we obtain

(4.4.232b) $$\int_0^1 \frac{\log^2 x}{1-x} Li_2(x)\, dx = \frac{1}{2} \left[ \varsigma(4) - \varsigma^2(2) \right]$$

As another example, with $r = q = 3$ we obtain

(4.4.232c) $$\int_0^1 \frac{\log^2 x}{1-x} Li_3(x)\, dx = \varsigma^2(3) - \varsigma(6)$$

because

(4.4.233)
$$\sum_{n=1}^{\infty} \frac{H_n^{(3)}}{n^3} = \frac{1}{2}\varsigma^2(3) + \frac{1}{2}\varsigma(6)$$

Using (4.4.230) we have

$$\sum_{n=1}^{\infty} \frac{H_n^{(r)}}{n^q} = \varsigma(r)\varsigma(q) - \frac{(-1)^{r-1}}{(r-1)!} \int_0^1 \frac{\log^{r-1} x \, Li_q(x)}{1-x} dx$$

$$\sum_{n=1}^{\infty} \frac{H_n^{(p)}}{n^q} = \varsigma(p)\varsigma(q) - \frac{(-1)^{p-1}}{(r-1)!} \int_0^1 \frac{\log^{p-1} x \, Li_q(x)}{1-x} dx$$

From (4.4.232a) we get

$$\sum_{n=1}^{\infty} \frac{H_n^{(p)}}{n^q} + \sum_{n=1}^{\infty} \frac{H_n^{(q)}}{n^p} = 2\varsigma(p)\varsigma(q) - \frac{(-1)^{p-1}}{(p-1)!} \int_0^1 \frac{\log^{p-1} x \, Li_q(x)}{1-x} dx - \frac{(-1)^{q-1}}{(q-1)!} \int_0^1 \frac{\log^{q-1} x \, Li_p(x)}{1-x} dx$$

$$= \varsigma(p)\varsigma(q) + \varsigma(p+q)$$

and

(4.4.233a)
$$\frac{1}{2}\Big[\varsigma^2(p) - \varsigma(2p)\Big] = \frac{(-1)^{p-1}}{(p-1)!} \int_0^1 \frac{\log^{p-1} x \, Li_p(x)}{1-x} dx$$

The above analysis may be generalised as follows. As before we have

$$\sum_{k=1}^{n} \frac{t^k}{k^r} = \frac{(-1)^{r-1}}{(r-1)!} \sum_{k=1}^{n} \int_0^t x^{k-1} \log^{r-1} x \, dx$$

$$= \frac{(-1)^{r-1}}{(r-1)!} \int_0^t \log^{r-1}(x) \frac{1-x^n}{1-x} dx$$

Hence upon making a further summation we have

$$\sum_{n=1}^{\infty} \frac{1}{n^q} \sum_{k=1}^{n} \frac{t^k}{k^r} = \frac{(-1)^{r-1}}{(r-1)!} \int_0^t \frac{\log^{r-1} x}{1-x} \sum_{n=1}^{\infty} \frac{1-x^n}{n^q} dx$$

(4.4.233b)
$$= \varsigma(q) \frac{(-1)^{r-1}}{(r-1)!} \int_0^t \frac{\log^{r-1} x}{1-x} dx - \frac{(-1)^{r-1}}{(r-1)!} \int_0^t \frac{\log^{r-1} x \, Li_q(x)}{1-x} dx$$

We have

(4.4.233bi) $\int_0^t \dfrac{\log x}{1-x}\,dx = Li_2(1-t) - \varsigma(2) = -\log t \log(1-t) - Li_2(t)$

With integration by parts we get

$$\int_0^t \frac{\log^2 x}{1-x}\,dx = -\log^2 t \log(1-t) + 2\int_0^t \frac{\log x \log(1-x)}{x}\,dx$$

$$\int_0^t \frac{\log x \log(1-x)}{x}\,dx = -\log t\, Li_2(t) + \int_0^t \frac{Li_2(x)}{x}\,dx = -\log t\, Li_2(t) + Li_3(t)$$

and therefore we have

(4.4.233bii) $\int_0^t \dfrac{\log^2 x}{1-x}\,dx = -\log^2 t \log(1-t) - 2\log t\, Li_2(t) + 2Li_3(t)$

With a little more work we find that (as per (3.223))

(4.4.233biii) $\int_0^t \dfrac{\log^3 x}{1-x}\,dx = -\log^3 t \log(1-t) - 3Li_2(t)\log^2 t + 6Li_3(t)\log t - 6Li_4(t)$

From this limited data we may conjecture that

$$\int_0^t \frac{\log^r x}{1-x}\,dx = -\log^r t \log(1-t) + \sum_{k=1}^r c_k Li_{k+1}(t)\log^{r-k} t$$

and the coefficients may be easily determined by differentiation. The final result is

(4.4.233c) $\int_0^t \dfrac{\log^r x}{1-x}\,dx = -\log^r t \log(1-t) + \sum_{k=1}^r (-1)^k r(r-1)...(r-k+1) Li_{k+1}(t)\log^{r-k} t$

From (4.4.231) we have

$$\int_0^t \frac{\log^{r-1} x}{1-x}\,dx = (-1)^{r-1}\varsigma(r)\Gamma(r) \quad \text{for } r \geq 2$$

in agreement with (4.4.233c).

Now referring back to (4.4.233b) we obtain

(4.4.233d) $\displaystyle\sum_{n=1}^{\infty}\frac{1}{n^q}\sum_{k=1}^{n}\frac{t^k}{k^r}=$

$$\varsigma(q)\frac{(-1)^r}{(r-1)!}\log^{r-1}t\log(1-t)+\varsigma(q)\frac{(-1)^{r-1}}{(r-1)!}\sum_{k=1}^{r-1}(-1)^k(r-1)(r-2)...(r-k)Li_{k+1}(t)\log^{r-k-1}t$$

$$-\frac{(-1)^{r-1}}{(r-1)!}\int_0^t\frac{\log^{r-1}x\,Li_q(x)}{1-x}dx$$

When $t=1$ it is clear that (4.4.233d) reduces to (4.4.229).

For example, with $r=q=2$ we get from (4.4.233b)

$$\sum_{n=1}^{\infty}\frac{1}{n^2}\sum_{k=1}^{n}\frac{t^k}{k^2}=-\varsigma(2)\int_0^t\frac{\log x}{1-x}dx+\int_0^t\frac{\log x\,Li_2(x)}{1-x}dx$$

We therefore obtain

(4.4.233e) $\displaystyle\int_0^t\frac{\log x\,Li_2(x)}{1-x}dx=\sum_{n=1}^{\infty}\frac{1}{n^2}\sum_{k=1}^{n}\frac{t^k}{k^2}-\varsigma^2(2)+\varsigma(2)Li_2(1-t)$

The above formula is used in (3.217a).

We may extend (4.4.233b) as follows

$$\sum_{n=1}^{\infty}\frac{u^n}{n^q}\sum_{k=1}^{n}\frac{t^k}{k^r}=\frac{(-1)^{r-1}}{(r-1)!}\int_0^t\frac{\log^{r-1}x}{1-x}\sum_{n=1}^{\infty}\frac{u^n-(ux)^n}{n^q}dx$$

and therefore we have

(4.4.233ei)

$$\sum_{n=1}^{\infty}\frac{u^n}{n^q}\sum_{k=1}^{n}\frac{t^k}{k^r}=\varsigma(q)\frac{(-1)^{r-1}Li_q(u)}{(r-1)!}\int_0^t\frac{\log^{r-1}x}{1-x}dx-\frac{(-1)^{r-1}}{(r-1)!}\int_0^t\frac{\log^{r-1}x\,Li_q(ux)}{1-x}dx$$

Coffey indicated the following analysis in [45d].

From (4.1.6) we have

$$H_n^{(1)} = \int_0^1 \frac{1-(1-x)^n}{x} \, dx = \int_0^1 \frac{1-x^n}{1-x} \, dx$$

and using (4.4.28) we have $\quad \dfrac{1}{n^p} = \dfrac{1}{\Gamma(p)} \displaystyle\int_0^\infty u^{p-1} e^{-nu} \, du$

Therefore we get for $p > 1$

$$\sum_{n=1}^\infty \frac{H_n^{(1)}}{n^p} = \sum_{n=1}^\infty \frac{1}{n^p} \int_0^1 \frac{1-x^n}{1-x} \, dx$$

$$= \sum_{n=1}^\infty \frac{1}{\Gamma(p)} \int_0^\infty u^{p-1} e^{-nu} \, du \int_0^1 \frac{1-x^n}{1-x} \, dx$$

$$= \frac{1}{\Gamma(p)} \int_0^\infty u^{p-1} \, du \int_0^1 \sum_{n=1}^\infty e^{-nu} \frac{1-x^n}{1-x} \, dx$$

As a geometric series we have

$$\sum_{n=1}^\infty e^{-nu} \left(1-x^n\right) = \sum_{n=1}^\infty e^{-nu} - \sum_{n=1}^\infty \left(xe^{-u}\right)^n$$

$$= \frac{1}{e^u-1} - \frac{x}{e^u-x} = \frac{(1-x)e^u}{(e^u-1)(e^u-x)}$$

and hence we have the double integral

(4.4.233f) $\qquad \displaystyle\sum_{n=1}^\infty \frac{H_n^{(1)}}{n^p} = \frac{1}{\Gamma(p)} \int_0^\infty \frac{u^{p-1}e^u}{e^u-1} \, du \int_0^1 \frac{dx}{e^u-x}$

More generally, by the same method, Coffey [45d] showed that

(4.4.233g) $\qquad \displaystyle\sum_{n=1}^\infty \frac{1}{n^p} \big[\gamma + \psi(kn+1)\big] = \frac{1}{\Gamma(p)} \int_0^\infty \frac{u^{p-1}e^u}{e^u-1} \, du \int_0^1 \left(\frac{1-x^k}{1-x}\right) \frac{dx}{e^u-x^k}$

(4.4.233h) $\qquad \displaystyle\sum_{n=1}^\infty \frac{(-1)^n}{n^p} \big[\gamma + \psi(kn+1)\big] = -\frac{1}{\Gamma(p)} \int_0^\infty \frac{u^{p-1}e^u}{e^u+1} \, du \int_0^1 \left(\frac{1-x^k}{1-x}\right) \frac{dx}{e^u+x^k}$

We have the elementary integral

$$\int_0^1 \frac{dx}{e^u - x} = -\log(e^u - x)\Big|_0^1 = \log\frac{e^u}{e^u - 1} = \log\frac{1}{1 - e^{-u}}$$

and therefore we get

$$\sum_{n=1}^{\infty} \frac{H_n^{(1)}}{n^p} = \frac{1}{\Gamma(p)} \int_0^{\infty} \frac{u^p e^u}{e^u - 1}\, du - \frac{1}{\Gamma(p)} \int_0^{\infty} \frac{u^{p-1} e^u \log(e^u - 1)}{e^u - 1}\, du$$

or alternatively

(4.4.233i) $$\sum_{n=1}^{\infty} \frac{H_n^{(1)}}{n^p} = \frac{1}{\Gamma(p)} \int_0^{\infty} \frac{u^{p-1}}{1 - e^{-u}} \log\frac{1}{1 - e^{-u}}\, du$$

A change of variables $v = e^{-u}$ gives us

(4.4.233j) $$\sum_{n=1}^{\infty} \frac{H_n^{(1)}}{n^p} = \frac{(-1)^p}{\Gamma(p)} \int_0^1 \frac{\log^{p-1} v}{v(1 - v)} \log(1 - v)\, dv$$

$$= \frac{(-1)^p}{\Gamma(p)} \int_0^1 \frac{\log^{p-1} v}{v} \log(1 - v)\, dv + \frac{(-1)^p}{\Gamma(p)} \int_0^1 \frac{\log^{p-1} v}{1 - v} \log(1 - v)\, dv$$

We have using integration by parts

$$\int_0^1 \frac{\log^{p-1} v}{v} \log(1 - v)\, dv = \frac{1}{p} \log^p v \log(1 - v)\Big|_0^1 + \frac{1}{p} \int_0^1 \frac{\log^p v}{1 - v}\, dv$$

$$= \frac{1}{p} \int_0^1 \frac{\log^p v}{1 - v}\, dv$$

and using (4.4.231) this becomes for $p \geq 1$

$$= \frac{1}{p}(-1)^p \varsigma(p+1)\Gamma(p+1)$$

Therefore we obtain

(4.4.233k) $$\sum_{n=1}^{\infty} \frac{H_n^{(1)}}{n^p} = (-1)^p \varsigma(p+1)\Gamma(p) + \frac{(-1)^p}{\Gamma(p)} \int_0^1 \frac{\log^{p-1} v}{1 - v} \log(1 - v)\, dv$$

and we could employ integration by parts to evaluate the above integral (albeit the required integration becomes increasingly complex).

Borwein et al. [30b] have recently shown that (see also (E.56))

$$(4.4.233l) \qquad -\int_0^1 \frac{1-t^\alpha}{1-t} \log t \, dt = \frac{1}{2}\Big[\big(\psi(1+\alpha)-\psi(1)\big)^2 + \psi'(1) - \psi'(1+\alpha)\Big]$$

and, therefore with $\alpha = n$, making use of (8.57a) in Volume V

$$\psi(1+n) - \psi(1) = H_n^{(1)} \qquad\qquad \psi'(1+n) - \psi'(1) = -H_n^{(2)}$$

we obtain

$$(4.4.233m) \qquad -\int_0^1 \frac{1-t^n}{1-t} \log t \, dt = -\int_0^1 \frac{1-(1-x)^n}{x}\log(1-x)\, dx = \frac{1}{2}\Big[H_n^{(2)} + \big[H_n^{(1)}\big]^2\Big]$$

Employing the same method as above we see that

$$\frac{1}{2}\sum_{n=1}^\infty \frac{\Big[H_n^{(2)} + \big[H_n^{(1)}\big]^2\Big]}{n^p} = -\sum_{n=1}^\infty \frac{1}{n^p}\int_0^1 \frac{1-x^n}{1-x}\log x \, dx$$

$$= -\sum_{n=1}^\infty \frac{1}{\Gamma(p)}\int_0^\infty u^{p-1}e^{-nu}\,du \int_0^1 \frac{1-x^n}{1-x}\log x \, dx$$

$$= -\frac{1}{\Gamma(p)}\int_0^\infty \frac{u^{p-1}e^u}{e^u-1}\,du \int_0^1 \frac{\log x\,dx}{e^u-x}$$

We have

$$\int \frac{\log x}{A-x}\,dx = -\log x \log\left(1-\frac{x}{A}\right) - Li_2\left(\frac{x}{A}\right)$$

and therefore

$$\int_0^1 \frac{\log x}{A-x}\,dx = -Li_2\left(\frac{1}{A}\right)$$

Hence we obtain

(4.4.233m) $\qquad \dfrac{1}{2}\displaystyle\sum_{n=1}^{\infty}\dfrac{\left[H_n^{(2)}+\left[H_n^{(1)}\right]^2\right]}{n^p}=\dfrac{1}{\Gamma(p)}\displaystyle\int_0^{\infty}\dfrac{u^{p-1}e^u}{e^u-1}Li_2\left(e^{-u}\right)du$

A change of variables $v=e^{-u}$ gives us

(4.4.233o) $\qquad \displaystyle\sum_{n=1}^{\infty}\dfrac{H_n^{(2)}}{n^p}+\displaystyle\sum_{n=1}^{\infty}\dfrac{\left[H_n^{(1)}\right]^2}{n^p}=2\dfrac{(-1)^p}{\Gamma(p)}\displaystyle\int_0^1\dfrac{\log^{p-1}v}{v(1-v)}Li_2(v)\,dv$

Reference to (4.4.229) gives us

(4.4.233p) $\qquad \displaystyle\sum_{n=1}^{\infty}\dfrac{H_n^{(p)}}{n^2}=\varsigma(p)\varsigma(2)-\dfrac{(-1)^{p-1}}{(p-1)!}\displaystyle\int_0^1\dfrac{\log^{p-1}x\,Li_2(x)}{1-x}\,dx$

and adding (4.4.233o) and (4.4.233p) together we obtain

$$\sum_{n=1}^{\infty}\dfrac{H_n^{(2)}}{n^p}+\sum_{n=1}^{\infty}\dfrac{H_n^{(p)}}{n^2}+\sum_{n=1}^{\infty}\dfrac{\left[H_n^{(1)}\right]^2}{n^p}=\varsigma(p)\varsigma(2)+\dfrac{(-1)^p}{\Gamma(p)}\int_0^1\dfrac{\log^{p-1}v}{v(1-v)}Li_2(v)\,dv$$

Using (4.4.232a) this becomes

(4.4.233q) $\qquad \varsigma(p+2)+\displaystyle\sum_{n=1}^{\infty}\dfrac{\left[H_n^{(1)}\right]^2}{n^p}=\dfrac{(-1)^p}{\Gamma(p)}\displaystyle\int_0^1\dfrac{\log^{p-1}v}{v(1-v)}Li_2(v)\,dv$

Differentiating (4.4.233l) with respect to $\alpha$ we get

(4.4.233r) $\qquad \displaystyle\int_0^1\dfrac{t^{\alpha}}{1-t}\log^2 t\,dt=\dfrac{1}{2}\Big[2\big(\psi(1+\alpha)-\psi(1)\big)\psi'(1+\alpha)-\psi''(1+\alpha)\Big]$

and, therefore with $\alpha=n$, we obtain

(4.4.233s) $\qquad \displaystyle\int_0^1\dfrac{t^n}{1-t}\log^2 t\,dt=\dfrac{1}{2}\Big[2\big(\psi(1+n)-\psi(1)\big)\psi'(1+n)-\psi''(1+n)\Big]$

We have from (8.57a) in Volume V

$$\psi'(1+n)=\varsigma(2)-H_n^{(2)}\qquad\qquad \psi''(1+n)=2H_n^{(3)}-2\varsigma(3)$$

and hence we have

$$(4.4.233t) \qquad \int_0^1 \frac{t^n}{1-t} \log^2 t \, dt = H_n^{(1)} \left[ \varsigma(2) - H_n^{(2)} \right] + \varsigma(3) - H_n^{(3)}$$

The Wolfram Integrator also kindly furnishes the integral (see also (4.4.233c) above)

$$\int_0^x \frac{\log^2 t}{1-t} \, dt = 2 Li_3(x) - \log^2 x \log(1-x) - 2 \log x \, Li_2(x)$$

and hence

$$\int_0^1 \frac{\log^2 t}{1-t} \, dt = 2 \varsigma(3)$$

Therefore we obtain

$$\int_0^1 \frac{1-t^n}{1-t} \log^2 t \, dt = \varsigma(3) + H_n^{(3)} - H_n^{(1)} \left[ \varsigma(2) - H_n^{(2)} \right]$$

Coffey's summation method then gives us

$$\sum_{n=1}^\infty \frac{\varsigma(3) + H_n^{(3)} - H_n^{(1)} \left[ \varsigma(2) - H_n^{(2)} \right]}{n^p} = \frac{1}{\Gamma(p)} \int_0^\infty \frac{u^{p-1} e^u}{e^u - 1} \, du \int_0^1 \frac{\log^2 x \, dx}{e^u - x}$$

We have

$$\int \frac{\log^2 x}{A - x} \, dx = -\log^2 x \log\left(1 - \frac{x}{A}\right) - 2 Li_2\left(\frac{x}{A}\right) \log x + 2 Li_3\left(\frac{x}{A}\right)$$

and therefore

$$\int_0^1 \frac{\log^2 x}{A - x} \, dx = 2 Li_3\left(\frac{1}{A}\right)$$

We then have

$$(4.4.233u) \quad \varsigma(3)\varsigma(p) + \sum_{n=1}^\infty \frac{H_n^{(3)}}{n^p} - \varsigma(2) \sum_{n=1}^\infty \frac{H_n^{(1)}}{n^p} + \sum_{n=1}^\infty \frac{H_n^{(1)} H_n^{(2)}}{n^p} = \frac{2}{\Gamma(p)} \int_0^\infty \frac{u^{p-1} e^u}{e^u - 1} Li_3\left(e^{-u}\right) du$$

A change of variables $v = e^{-u}$ gives us

$$(4.4.233v) \quad \varsigma(3)\varsigma(p) + \sum_{n=1}^{\infty}\frac{H_n^{(3)}}{n^p} - \varsigma(2)\sum_{n=1}^{\infty}\frac{H_n^{(1)}}{n^p} + \sum_{n=1}^{\infty}\frac{H_n^{(1)}H_n^{(2)}}{n^p} = \frac{2(-1)^p}{\Gamma(p)}\int_0^{\infty}\frac{\log^{p-1}v\,Li_3(v)}{v(1-v)}dv$$

It is clear that the range of identities may be extended by successively differentiating (4.4.233l).

With reference to (4.4.233t) we have

$$\int_0^1 \frac{t^n}{1-t}\log^2 t\,dt = \frac{t^{n+1}}{n+1}\frac{\log^2 t}{1-t}\bigg|_0^1 - \frac{1}{n+1}\int_0^1 \frac{t^{n+1}\left(2(1-t)\log t/t + \log^2 t\right)}{(1-t)^2}dt$$

$$= -\frac{1}{n+1}\int_0^1 \frac{t^{n+1}\left(2(1-t)\log t/t + \log^2 t\right)}{(1-t)^2}dt$$

Using L'Hôpital's rule we can see that the integrand is bounded at $t = 0,1$ and hence throughout the interval $[0,1]$; the integral is therefore finite and we see that

$$\lim_{n\to\infty}\int_0^1 \frac{t^n}{1-t}\log^2 t\,dt = 0$$

The above limit is also apparent from (4.4.234a). This therefore implies from (4.4.233t) that

$$(4.4.233w) \qquad \lim_{n\to\infty} H_n^{(1)}\left[\varsigma(2) - H_n^{(2)}\right] = 0$$

Since $\lim_{n\to\infty} H_n^{(1)}$ does not exist, the above limit is non-trivial. We have by reference to the definition of Euler's constant $\gamma$

$$\lim_{n\to\infty}\left[\log n - H_n^{(1)}\right]\left[\varsigma(2) - H_n^{(2)}\right] = 0$$

and therefore we deduce that

$$(4.4.233x) \qquad \lim_{n\to\infty}\log n\left[\varsigma(2) - H_n^{(2)}\right] = 0$$

This limit is employed in (E.33i) in Volume VI.

We have from (4.4.238aa)

$$\int_0^x \frac{t^n}{1-t} \log t \, dt = \sum_{k=1}^n \frac{x^k}{k^2} - \log x \sum_{k=1}^n \frac{x^k}{k} + Li_2(1-x) - \varsigma(2)$$

and the following integral appears in G&R [74, p.535]

(4.4.233y) $\quad \int_0^1 \frac{t^{\mu-1}}{1-t} \log t \, dt = -\psi'(\mu) = -\varsigma(2,\mu)$

where $\varsigma(s,a)$ is the Hurwitz zeta function defined by $\varsigma(s,a) = \sum_{k=0}^\infty \frac{1}{(k+a)^s}$

We have

$$\varsigma(s,n) = \sum_{k=0}^\infty \frac{1}{(k+n)^s} = \varsigma(s) - H_{n-1}^{(s)} = \varsigma(s) - H_n^{(s)} + \frac{1}{n^s}$$

and therefore we obtain the same result as before.

$$\int_0^1 \frac{t^{n-1}}{1-t} \log t \, dt = H_n^{(2)} - \frac{1}{n^2} - \varsigma(2)$$

Differentiating (4.4.233y) under the integral we have

$$\int_0^1 \frac{t^{\mu-1}}{1-t} \log^2 t \, dt = -\psi''(\mu) = -\varsigma'(2,\mu)$$

Following the Cauchy series product method employed in (3.32) and using Adamchik's formula (3.18) we obtain

$$\frac{1}{1-x} \sum_{n=1}^\infty \left[ \frac{H_n^{(2)}}{n} + \frac{H_n^{(1)}}{n^2} \right] x^n = \sum_{n=0}^\infty x^n \sum_{k=1}^n \left[ \frac{H_k^{(2)}}{k} + \frac{H_k^{(1)}}{k^2} \right]$$

$$= \sum_{n=1}^\infty \left[ H_n^{(3)} + H_n^{(1)} H_n^{(2)} \right] x^n$$

Multiplying across by $\log x / x$ and integrating we obtain

$$\sum_{n=1}^\infty \left[ \frac{H_n^{(2)}}{n} + \frac{H_n^{(1)}}{n^2} \right] \int_0^t \frac{x^{n-1}}{1-x} \log x \, dx = \sum_{n=1}^\infty \left[ H_n^{(3)} + H_n^{(1)} H_n^{(2)} \right] \int_0^t x^{n-1} \log x \, dx$$

and this becomes

$$\sum_{n=1}^{\infty}\left[\frac{H_n^{(2)}}{n}+\frac{H_n^{(1)}}{n^2}\right]\left[\sum_{k=1}^{n-1}\frac{t^k}{k^2}-\log t\sum_{k=1}^{n-1}\frac{t^k}{k}+Li_2(1-t)\right]=\sum_{n=1}^{\infty}\left[H_n^{(3)}+H_n^{(1)}H_n^{(2)}\right]\left[-\frac{1}{n^2}+\frac{\log t}{n}\right]t^n$$

We have with $t=1$

$$\sum_{n=1}^{\infty}\left[\frac{H_n^{(2)}}{n}+\frac{H_n^{(1)}}{n^2}\right]\left[H_n^{(2)}-\frac{1}{n^2}-\varsigma(2)\right]=-\sum_{n=1}^{\infty}\frac{\left[H_n^{(3)}+H_n^{(1)}H_n^{(2)}\right]}{n^2}$$

This may be rearranged as follows

$$\sum_{n=1}^{\infty}\frac{H_n^{(2)}}{n}\left[H_n^{(2)}-\frac{1}{n^2}-\varsigma(2)\right]+\sum_{n=1}^{\infty}\frac{H_n^{(1)}H_n^{(2)}}{n^2}-\sum_{n=1}^{\infty}\frac{H_n^{(1)}}{n^4}-\varsigma(2)\sum_{n=1}^{\infty}\frac{H_n^{(1)}}{n^2}=-\sum_{n=1}^{\infty}\frac{\left[H_n^{(3)}+H_n^{(1)}H_n^{(2)}\right]}{n^2}$$

and using $\sum_{n=1}^{\infty}\frac{H_n^{(1)}}{n^4}=-\varsigma(2)\varsigma(3)+3\varsigma(5)$ and $\sum_{n=1}^{\infty}\frac{H_n^{(1)}}{n^2}=2\varsigma(3)$ this is simplified to

$$\sum_{n=1}^{\infty}\frac{H_n^{(2)}}{n}\left[H_n^{(2)}-\varsigma(2)\right]=\sum_{n=1}^{\infty}\frac{H_n^{(2)}}{n^3}-\sum_{n=1}^{\infty}\frac{H_n^{(3)}}{n^2}-2\sum_{n=1}^{\infty}\frac{H_n^{(1)}H_n^{(2)}}{n^2}+\varsigma(2)\varsigma(3)+3\varsigma(5)$$

We already know $\sum_{n=1}^{\infty}\frac{H_n^{(2)}}{n^3}$, $\sum_{n=1}^{\infty}\frac{H_n^{(3)}}{n^2}$ and $\sum_{n=1}^{\infty}\frac{H_n^{(1)}H_n^{(2)}}{n^2}$ and hence we may determine $\sum_{n=1}^{\infty}\frac{H_n^{(2)}}{n}\left[H_n^{(2)}-\varsigma(2)\right]$.

A similar analysis may also be performed using Adamchik's formula (3.19)

$$\sum_{k=1}^{n}\frac{\left(H_k^{(1)}\right)^2}{k}+\sum_{k=1}^{n}\frac{H_k^{(2)}}{k}=\sum_{k=1}^{n}\frac{1}{k}\left[\left(H_k^{(1)}\right)^2+H_k^{(2)}\right]=\frac{1}{3}\left(H_n^{(1)}\right)^3+H_n^{(1)}H_n^{(2)}+\frac{2}{3}H_n^{(3)}$$

**Theorem 4.12:**

(4.4.234)     $$\sum_{n=1}^{\infty}\frac{1}{2^n}\sum_{k=1}^{n}\binom{n}{k}\frac{H_k^{(3)}}{k^3}=\varsigma^2(3)+\varsigma(6)$$

**Proof:**

From (3.55) we know that

$$2Li_3(x) = \sum_{n=1}^{\infty} \frac{1}{2^n} \sum_{k=1}^{n} \binom{n}{k} \frac{x^k}{k^3}$$

and therefore we have

$$\int_0^1 \frac{\log^2 x \, Li_3(x)}{1-x} \, dx = \frac{1}{2} \sum_{n=1}^{\infty} \frac{1}{2^n} \sum_{k=1}^{n} \binom{n}{k} \frac{1}{k^3} \int_0^1 \frac{x^k \log^2 x}{1-x} \, dx$$

The following indefinite integral may be verified by differentiation (the general form was suggested to me by inspection of the result obtained by using the Wolfram Integrator for $k = 4$)

(4.4.234a) $\qquad \displaystyle\int \frac{x^k \log^2 x}{1-x} \, dx =$

$$-2\sum_{j=1}^{k} \frac{x^j}{j^3} + 2\log x \sum_{j=1}^{k} \frac{x^j}{j^2} - \log^2 x \sum_{j=1}^{k} \frac{x^j}{j} - \log(1-x)\log^2 x - 2\log x \, Li_2(x) + 2Li_3(x)$$

and therefore we get

(4.4.235) $\qquad \displaystyle\int_0^1 \frac{x^k \log^2 x}{1-x} \, dx = -2\sum_{j=1}^{k} \frac{1}{j^3} + 2Li_3(1) = -2H_k^{(3)} + 2\varsigma(3)$

(this integral is contained in G&R [74, p.538]). Accordingly we obtain

$$\int_0^1 \frac{\log^2 x \, Li_3(x)}{1-x} \, dx = \sum_{n=1}^{\infty} \frac{1}{2^n} \sum_{k=1}^{n} \binom{n}{k} \frac{1}{k^3} \left[ \varsigma(3) - H_k^{(3)} \right]$$

(4.4.236) $\qquad\qquad = 2\varsigma^2(3) - \sum_{n=1}^{\infty} \frac{1}{2^n} \sum_{k=1}^{n} \binom{n}{k} \frac{H_k^{(3)}}{k^3}$

Equating (4.4.232c) and (4.4.236) we have

$$\sum_{n=1}^{\infty} \frac{1}{2^n} \sum_{k=1}^{n} \binom{n}{k} \frac{H_k^{(3)}}{k^3} = \varsigma^2(3) + \varsigma(6)$$

We now consider the following identity in the limit as $k \rightarrow \infty$

$$\int_0^t \frac{x^k \log^2 x}{1-x} \, dx = -2\sum_{j=1}^{k} \frac{t^j}{j^3} + 2\log t \sum_{j=1}^{k} \frac{t^j}{j^2} - \log^2 t \sum_{j=1}^{k} \frac{t^j}{j} - \log(1-t)\log^2 t - 2\log t \, Li_2(t) + 2Li_3(t)$$

It is clear that the integral on the left-hand side will vanish as $k \to \infty$ provided $0 \le t < 1$ and the right-hand side will also approach zero as a result of the series definitions of the various polylogarithms.

Now let us multiply the above equation by $A^k$ and complete the summation: we obtain for the left-hand side

$$\sum_{k=1}^{\infty} A^k \int_0^t \frac{x^k \log^2 x}{1-x} \, dx = \int_0^t \sum_{k=1}^{\infty} \frac{A^k x^k \log^2 x}{1-x} \, dx = \int_0^t \frac{Ax \log^2 x}{(1-x)(1-Ax)} \, dx$$

and for the right-hand side we have

$$= -2 \sum_{k=1}^{\infty} A^k \sum_{j=1}^{k} \frac{t^j}{j^3} + 2 \log t \sum_{k=1}^{\infty} A^k \sum_{j=1}^{k} \frac{t^j}{j^2} - \log^2 t \sum_{k=1}^{\infty} A^k \sum_{j=1}^{k} \frac{t^j}{j}$$

$$+ \frac{A}{1-A} \Big[ -\log(1-t) \log^2 t - 2 \log t \, Li_2(t) + 2 Li_3(t) \Big]$$

We now use (3.23) to rearrange this as

$$-2 \sum_{k=1}^{\infty} \frac{t^k}{k^3} \sum_{j=k}^{\infty} A^j + 2 \log t \sum_{k=1}^{\infty} \frac{t^k}{k^2} \sum_{j=k}^{\infty} A^j - \log^2 t \sum_{k=1}^{\infty} \frac{t^k}{k} \sum_{j=k}^{\infty} A^j$$

$$+ \frac{A}{1-A} \Big[ -\log(1-t) \log^2 t - 2 \log t \, Li_2(t) + 2 Li_3(t) \Big]$$

Since $\sum_{j=k}^{\infty} A^j = \frac{A^k}{1-A}$ where $|A| < 1$, we have

$$= \frac{1}{1-A} \left[ -2 \sum_{k=1}^{\infty} \frac{t^k A^k}{k^3} + 2 \log t \sum_{k=1}^{\infty} \frac{t^k A^k}{k^2} - \log^2 t \sum_{k=1}^{\infty} \frac{t^k A^k}{k} \right]$$

$$+ \frac{A}{1-A} \Big[ -\log(1-t) \log^2 t - 2 \log t \, Li_2(t) + 2 Li_3(t) \Big]$$

Hence we have

$$\int_0^t \frac{Ax \log^2 x}{(1-x)(1-Ax)} \, dx =$$

$$\frac{1}{1-A}\Big[\log^2 t\left\{-A\log(1-t)+\log(1-At)\right\}-2\log t\left\{ALi_2(t)-Li_2(At)\right\}+2\left\{ALi_3(t)-Li_3(At)\right\}\Big]$$

It should be noted that this could have been derived more directly by using partial fractions because

$$\int_0^t \frac{Ax\log^2 x}{(1-x)(1-Ax)}\,dx = \frac{A}{1-A}\left[\int_0^t \frac{\log^2 x}{(1-x)}\,dx - \int_0^t \frac{\log^2 x}{(1-Ax)}\,dx\right]$$

From (4.4.235) we have

$$\int_0^1 \frac{x^k \log^2 x}{1-x}\,dx = -2H_k^{(3)}+2\varsigma(3)$$

Completing the summation we get

$$\sum_{k=1}^{\infty}\frac{1}{k^s}\int_0^1 \frac{x^k \log^2 x}{1-x}\,dx = -2\sum_{k=1}^{\infty}\frac{H_k^{(3)}}{k^s}+2\varsigma(3)\varsigma(s)$$

or equivalently

(4.4.237) $$\int_0^1 \frac{Li_s(x)\log^2 x}{1-x}\,dx = 2\varsigma(3)\varsigma(s)-2\sum_{k=1}^{\infty}\frac{H_k^{(3)}}{k^s}$$

This corresponds with (4.4.229) with $r = 3$.

From [101aa] we have

$$H_k^{(r)} = \frac{1}{k^r}+\frac{(-1)^{r+1}\psi^{(r)}(k)}{\Gamma(r)}+\varsigma(r)$$

where $\psi^{(n)}(x)=\dfrac{d^{n+1}}{dx^{n+1}}\log\Gamma(x)=\dfrac{d^n}{dx^n}\psi(x)$ and $k,r$ are positive integers $\neq 1$.

Substituting this in (4.4.234), for example, and using (3.55) we may obtain a series involving $\sum_{n=1}^{\infty}\dfrac{1}{2^n}\sum_{k=1}^{n}\binom{n}{k}\dfrac{\psi^{(3)}(k)}{k^3}$.

Using (4.4.234a) and completing the summation we obtain

(4.4.237a) $$\int_0^x \frac{Li_p(t)\log^2 t}{1-t}\,dt = -2\sum_{k=1}^{\infty}\frac{1}{k^p}\sum_{j=1}^{k}\frac{x^j}{j^3}+2\log x\sum_{k=1}^{\infty}\frac{1}{k^p}\sum_{j=1}^{k}\frac{x^j}{j^2}-\log^2 x\sum_{k=1}^{\infty}\frac{1}{k^p}\sum_{j=1}^{k}\frac{x^j}{j}$$

$$-\zeta(p)\log(1-x)\log^2 x - 2\zeta(p)\log x \, Li_2(x) + 2\zeta(p)Li_3(x)$$

With $x = 1$ we get

(4.4.237b) $$\int\limits_0^1 \frac{Li_p(t)\log^2 t}{1-t}\,dt = 2\zeta(3)\zeta(p) - 2\sum_{k=1}^{\infty}\frac{1}{k^p}\sum_{j=1}^{k}\frac{1}{j^3} = 2\zeta(3)\zeta(p) - 2\sum_{k=1}^{\infty}\frac{H_k^{(3)}}{k^p}$$

**Theorem 4.13:**

(4.4.238) $$\sum_{n=1}^{\infty}\frac{1}{2^{n+1}}\sum_{k=1}^{n}\binom{n}{k}\frac{H_k^{(2)}}{k^4} = \zeta(2)\zeta(4) + \zeta^2(3) - \frac{25}{12}\zeta(6)$$

**Proof:**

Employing the same methodology as before we have

$$\int\limits_0^1 \frac{\log x \, Li_4(x)}{1-x}\,dx = \frac{1}{2}\sum_{n=1}^{\infty}\frac{1}{2^n}\sum_{k=1}^{n}\binom{n}{k}\frac{1}{k^4}\int\limits_0^1 \frac{x^k \log x}{1-x}\,dx$$

As in the previous theorem, the following indefinite integral may be verified by differentiation (the general form was again suggested to me by using the Wolfram Integrator for $k = 4$ )

(4.4.238aa) $$\int \frac{x^k \log x}{1-x}\,dx = \sum_{j=1}^{k}\frac{x^j}{j^2} - \log x\sum_{j=1}^{k}\frac{x^j}{j} + Li_2(1-x)$$

Therefore we have

(4.4.238a) $$\int\limits_0^1 \frac{x^k \log x}{1-x}\,dx = H_k^{(2)} - \zeta(2)$$

(this integral is contained in G&R [74, p.535]). Accordingly we obtain

$$\int\limits_0^1 \frac{\log x \, Li_4(x)}{1-x}\,dx = \frac{1}{2}\sum_{n=1}^{\infty}\frac{1}{2^n}\sum_{k=1}^{n}\binom{n}{k}\frac{1}{k^4}\left[H_k^{(2)} - \zeta(2)\right]$$

$$= \frac{1}{2}\sum_{n=1}^{\infty}\frac{1}{2^n}\sum_{k=1}^{n}\binom{n}{k}\frac{H_k^{(2)}}{k^4} - \zeta(2)\zeta(4)$$

From Freitas [69a] we have

$$(4.4.238b) \qquad \int_0^1 \frac{\log x \, Li_4(x)}{1-x} dx = \varsigma^2(3) - \frac{25}{12}\varsigma(6)$$

and hence we obtain

$$(4.4.239) \qquad \sum_{n=1}^{\infty} \frac{1}{2^{n+1}} \sum_{k=1}^{n} \binom{n}{k} \frac{H_k^{(2)}}{k^4} = \varsigma(2)\varsigma(4) + \varsigma^2(3) - \frac{25}{12}\varsigma(6)$$

In a similar manner we have

$$\int_0^1 \frac{\log x \, Li_2(x)}{1-x} dx = \frac{1}{2} \sum_{n=1}^{\infty} \frac{1}{2^n} \sum_{k=1}^{n} \binom{n}{k} \frac{1}{k^2} \left[ H_k^{(2)} - \varsigma(2) \right]$$

$$= \frac{1}{2} \sum_{n=1}^{\infty} \frac{1}{2^n} \sum_{k=1}^{n} \binom{n}{k} \frac{H_k^{(2)}}{k^2} - \varsigma^2(2)$$

From Freitas [69a] we have (see also (3.211e) of Volume I)

$$(4.4.239a) \qquad \int_0^1 \frac{\log x \, Li_2(x)}{1-x} dx = -\frac{3}{4}\varsigma(4)$$

and hence we obtain

$$(4.4.240) \qquad \sum_{n=1}^{\infty} \frac{1}{2^{n+1}} \sum_{k=1}^{n} \binom{n}{k} \frac{H_k^{(2)}}{k^2} = \varsigma^2(2) - \frac{3}{4}\varsigma(4)$$

Using the same techniques as before we have

$$\int_0^1 \frac{\log x \, Li_3(x)}{1-x} dx = \frac{1}{2} \sum_{n=1}^{\infty} \frac{1}{2^n} \sum_{k=1}^{n} \binom{n}{k} \frac{1}{k^3} \left[ H_k^{(2)} - \varsigma(2) \right]$$

$$= \frac{1}{2} \sum_{n=1}^{\infty} \frac{1}{2^n} \sum_{k=1}^{n} \binom{n}{k} \frac{H_k^{(2)}}{k^3} - \varsigma(2)\varsigma(3)$$

Using Freitas [69a] again we have

$$\int_0^1 \frac{\log x \, Li_3(x)}{1-x} dx = 2\varsigma(2)\varsigma(3) - \frac{9}{2}\varsigma(5)$$

and hence we obtain

$$(4.4.241) \qquad \sum_{n=1}^{\infty} \frac{1}{2^{n+1}} \sum_{k=1}^{n} \binom{n}{k} \frac{H_k^{(2)}}{k^3} = 2\varsigma(2)\varsigma(3) - \frac{9}{2}\varsigma(5)$$

These formulae are collected below for ease of reference

$$\sum_{n=1}^{\infty} \frac{1}{2^{n+1}} \sum_{k=1}^{n} \binom{n}{k} \frac{H_k^{(2)}}{k^2} = \varsigma^2(2) - \frac{3}{4}\varsigma(4)$$

$$\sum_{n=1}^{\infty} \frac{1}{2^{n+1}} \sum_{k=1}^{n} \binom{n}{k} \frac{H_k^{(2)}}{k^3} = 2\varsigma(2)\varsigma(3) - \frac{9}{2}\varsigma(5)$$

$$\sum_{n=1}^{\infty} \frac{1}{2^{n+1}} \sum_{k=1}^{n} \binom{n}{k} \frac{H_k^{(2)}}{k^4} = \varsigma(2)\varsigma(4) + \varsigma^2(3) - \frac{25}{12}\varsigma(6)$$

$$\sum_{n=1}^{\infty} \frac{1}{2^{n+1}} \sum_{k=1}^{n} \binom{n}{k} \frac{H_k^{(3)}}{k^3} = \frac{1}{2}\varsigma^2(3) + \frac{1}{2}\varsigma(6)$$

It is clear that the list of formulae may be easily extended. Structurally similar identities may also be obtained in the same manner using the definition of the polylogarithm in equation (4.4.45).

$$Li_{s+1}(x) = \sum_{n=1}^{\infty} \frac{1}{n2^n} \sum_{k=1}^{n} \binom{n}{k} \frac{x^k}{k^s}$$

Let us now consider the identity in (4.4.85)

$$\frac{x}{s-1} \sum_{n=0}^{\infty} \frac{1}{n+1} \sum_{k=0}^{n} \binom{n}{k} \frac{(-1)^k x^k}{(k+1)^{s-1}} = \frac{\log x}{s-1} Li_{s-1}(x) + Li_s(x)$$

Multiplying both sides by $\log x/(1-x)$ in the case where $s = 4$ and integrating, we obtain

$$\frac{1}{3} \sum_{n=0}^{\infty} \frac{1}{n+1} \sum_{k=0}^{n} \binom{n}{k} \frac{(-1)^k}{(k+1)^3} \int_0^1 \frac{x^{k+1} \log x}{1-x} dx = \frac{1}{3} \int_0^1 \frac{\log^2 x \, Li_3(x)}{1-x} dx + \int_0^1 \frac{\log x \, Li_4(x)}{1-x} dx$$

and this becomes using (4.4.238a), (4.4.236a) and (4.4.238b)

$$\frac{1}{3} \sum_{n=0}^{\infty} \frac{1}{n+1} \sum_{k=0}^{n} \binom{n}{k} \frac{(-1)^k}{(k+1)^3} \left[ H_{k+1}^{(2)} - \varsigma(2) \right] = \frac{1}{3}\left(\varsigma^2(3) - \varsigma(6)\right) + \left(\varsigma^2(3) - \frac{25}{12}\varsigma(6)\right)$$

Using (4.4.85a) we obtain

$$(4.4.242) \quad \frac{1}{3}\sum_{n=0}^{\infty}\frac{1}{n+1}\sum_{k=0}^{n}\binom{n}{k}\frac{(-1)^k H_{k+1}^{(2)}}{(k+1)^3} = \frac{4}{3}\varsigma^2(3) - \frac{29}{12}\varsigma(6) + \varsigma(2)\varsigma(4)$$

In (4.2.28) we showed that

$$\sum_{k=0}^{n}\binom{n}{k}\frac{(-1)^k}{(1+k)^3} = \frac{1}{n+1}\left\{\frac{1}{2}\left(H_{n+1}^{(1)}\right)^2 + \frac{1}{2}H_{n+1}^{(2)}\right\}$$

Therefore, using the binomial inversion formula [75, p.192]

$$a_n = \sum_{k=0}^{n}(-1)^k\binom{n}{k}b_k \quad \Leftrightarrow \quad b_n = \sum_{k=0}^{n}(-1)^k\binom{n}{k}a_k$$

we obtain

$$\sum_{k=0}^{n}(-1)^k\binom{n}{k}\frac{1}{k+1}\left\{\frac{1}{2}\left(H_{k+1}^{(1)}\right)^2 + \frac{1}{2}H_{k+1}^{(2)}\right\} = \frac{1}{(1+n)^3}$$

and hence

$$(4.4.243) \quad \sum_{n=0}^{\infty}\sum_{k=0}^{n}(-1)^k\binom{n}{k}\frac{1}{k+1}\left\{\frac{1}{2}\left(H_{k+1}^{(1)}\right)^2 + \frac{1}{2}H_{k+1}^{(2)}\right\} = \sum_{n=0}^{\infty}\frac{1}{(1+n)^3} = \varsigma(3)$$

Jumping ahead to the next theorem, we have

$$\int_0^1 x^k \log^2(1-x)\,dx = \frac{2}{k+1}\left[\frac{1}{2}\left(H_{k+1}^{(1)}\right)^2 + \frac{1}{2}H_{k+1}^{(2)}\right]$$

and substituting this in the left-hand side of (4.4.243) we obtain

$$\frac{1}{2}\sum_{n=0}^{\infty}\sum_{k=0}^{n}(-1)^k\binom{n}{k}\int_0^1 x^k \log^2(1-x)\,dx = \frac{1}{2}\sum_{n=0}^{\infty}\int_0^1 (1-x)^n \log^2(1-x)\,dx$$

$$= \frac{1}{2}\int_0^1 \frac{\log^2(1-x)}{x}\,dx$$

We therefore have the well-known result

$$\varsigma(3) = \frac{1}{2}\int_0^1 \frac{\log^2(1-x)}{x}\,dx$$

The above result can of course be obtained much more directly using integration by parts

$$\int \frac{\log^2(1-x)}{x}\,dx = \log x \log^2(1-x) + 2\int \frac{\log x \log(1-x)}{1-x}\,dx$$

and from equations (3.39) and (3.40) we showed previously that

$$\frac{d}{dx}Li_2(1-x) = \frac{\log x}{1-x} \qquad \text{and} \qquad \frac{d}{dx}Li_3(1-x) = -\frac{Li_2(1-x)}{1-x}$$

$$\int \frac{\log(1-x)\log x\,dx}{1-x} = \log(1-x)Li_2(1-x) + \int \frac{Li_2(1-x)dx}{1-x}$$

$$= \log(1-x)Li_2(1-x) - Li_3(1-x)$$

Therefore we have

$$\int \frac{\log^2(1-x)}{x}\,dx = \log x \log^2(1-x) + 2\log(1-x)Li_2(1-x) - 2Li_3(1-x)$$

We may generalise (4.4.243) as follows

$$\sum_{n=0}^{\infty} t^n \sum_{k=0}^{n}(-1)^k \binom{n}{k}\frac{1}{k+1}\left\{\frac{1}{2}(H_{k+1}^{(1)})^2 + \frac{1}{2}H_{k+1}^{(2)}\right\} = \sum_{n=0}^{\infty}\frac{t^n}{(1+n)^3}$$

and therefore we have

$$\frac{1}{2}\sum_{n=0}^{\infty}\int_0^1 t^n(1-x)^n \log^2(1-x)\,dx = \frac{1}{2}\int_0^1 \frac{\log^2(1-x)}{1-t(1-x)}\,dx$$

Hence we obtain

(4.4.244) $\qquad \dfrac{1}{2}\displaystyle\int_0^1 \frac{\log^2(1-x)}{1-t(1-x)}\,dx = \sum_{n=0}^{\infty}\frac{t^n}{(1+n)^3} = \frac{1}{t}Li_3(t)$

and we have seen a similar result previously in (4.4.155k). In particular for $t = 1/2$, we have

(4.4.244a) $$\int_0^1 \frac{\log^2(1-x)}{1+x}\,dx = 2Li_3(1/2)$$

Using the Wolfram Integrator we get

$$\int_0^t \frac{\log^2(1-x)}{1+x}\,dx = \log\left(\frac{1+t}{2}\right)\log^2(1-t) + 2Li_2\left(\frac{1-t}{2}\right)\log(1-t) - 2Li_3\left(\frac{1-t}{2}\right) + 2Li_3(1/2)$$

and we note that (4.4.244a) results when $t = 1$.

Differentiating (4.4.244) we obtain

(4.4.244b) $$\int_0^1 \frac{(1-x)\log^2(1-x)}{\left[1-t(1-x)\right]^2}\,dx = \frac{2}{t^2}\left[Li_2(t) - Li_3(t)\right]$$

Integrating (4.4.244) we obtain

$$\frac{1}{2}\int_0^u dt \int_0^1 \frac{\log^2(1-x)}{1-t(1-x)}\,dx = Li_4(u)$$

and changing the order of integration we have

(4.4.245) $$-\frac{1}{2}\int_0^1 \frac{\log^2(1-x)\log\left[1-u(1-x)\right]}{1-x}\,dx = Li_4(u)$$

Putting $u = 1$ we obtain

(4.4.245a) $$\int_0^1 \frac{\log x \log^2(1-x)}{1-x}\,dx = -2\varsigma(4)$$

which is similar to the integral evaluated by Rutledge and Douglass [116aa] in (4.4.167qa)

(4.4.245b) $$\int_0^1 \frac{\log x \log^2(1-x)}{x}\,dx = -\frac{1}{2}\varsigma(4)$$

With $u = 1/2$ we obtain

(4.4.245c)    $-\dfrac{1}{2}\displaystyle\int_0^1 \dfrac{\log^2(1-x)\log\left[(1+x)/2\right]}{1-x}\,dx = Li_4(1/2)$

Putting $u=1$ in (4.4.245) we get

(4.4.245d)    $-\dfrac{1}{2}\displaystyle\int_0^1 \dfrac{\log^2(1-x)\log x}{1-x}\,dx = Li_4(1) = \varsigma(4)$

Differentiating (4.4.245) we obtain

(4.4.245e)    $\displaystyle\int_0^1 \dfrac{\log^2(1-x)}{1-u(1-x)}\,dx = \dfrac{2}{u}Li_3(u)$

and therefore we have

(4.4.245f)    $\displaystyle\int_0^1 \dfrac{\log^2(1-x)}{x}\,dx = 2Li_3(1) = 2\varsigma(3)$

(4.4.245g)    $\displaystyle\int_0^1 \dfrac{\log^2(1-x)}{1+x}\,dx = 2Li_3(1/2)$

Putting $y=1-x$ in (4.4.245e) we get

(4.4.245h)    $\displaystyle\int_0^1 \dfrac{\log^2 y}{1-uy}\,dy = \dfrac{2}{u}Li_3(u)$

and making the substitution $t=-\log y$ this becomes

(4.4.245i)    $\displaystyle\int_0^\infty \dfrac{t^2}{e^t-u}\,dt = \dfrac{2}{u}Li_3(u)$

and we recognise this integral as being equivalent to (4.4.25).

Making a change of variables $t=1-x$ in (4.4.245) we obtain

(4.4.245j)    $\displaystyle\int_0^1 \dfrac{\log^2 t\,\log\left[1-ut\right]}{t}\,dt = -2Li_4(u)$

and with $u=1/2$ we have

(4.4.245k) $$\int_0^1 \frac{\log^2 t \log[1-t/2]}{t} dt = -2 Li_4(1/2)$$

The integral (4.4.245k) was regarded with some importance by Kölbig in his 1986 paper [91a] (and appears to have been unknown at that time): matters have really progressed since then because the Wolfram Integrator can now evaluate the more general integral (4.4.245j) with the minimum of effort.

Integration by parts gives us

$$\int_0^1 \frac{\log^2 t \log[1-t/2]}{t} dt = \frac{1}{3} \log^3 t \log[1-t/2]\Big|_0^1 + \frac{1}{6} \int_0^1 \frac{\log^3 t}{[1-t/2]} dt$$

$$= \frac{1}{6} \int_0^1 \frac{\log^3 t}{[1-t/2]} dt$$

A further integration by parts gives us

$$\int \frac{\log^3 t}{[1-t/2]} dt = -12 \left[ \frac{1}{6} \log^3 t \log[1-t/2] + \frac{1}{2} \log^2 t \, Li_3(t/2) - \log t \, Li_3(t/2) + Li_4(t/2) \right]$$

and hence

$$\int_0^1 \frac{\log^3 t}{[1-t/2]} dt = -12 Li_4(1/2)$$

and (4.4.245k) immediately follows. With the substitution $t = 1-u$ we get

$$\int_0^1 \frac{\log^3 (1-u)t}{1+u} dt = -6 Li_4(1/2)$$

Dividing (4.4.245j) by $u$ and integrating we get

(4.4.245l) $$\int_0^1 \frac{du}{u} \int_0^1 \frac{\log^2 t \log[1-ut]}{t} dx = -2 Li_5(1)$$

Changing the order of integration results in

(4.4.245m) $$\int_0^1 \frac{\log^2 t}{t} dt \int_0^1 \log[1-ut] \frac{du}{u} = -2 Li_5(1)$$

and we have

$$\int_0^1 \log[1-ut]\frac{du}{u} = -\sum_{n=1}^{\infty}\int_0^1 \frac{t^n u^{n-1}}{n}\,du = -\sum_{n=1}^{\infty}\frac{t^n}{n^2} = -Li_2(t)$$

We therefore obtain

(4.4.245n)     $$\int_0^1 \frac{\log^2 t\, Li_2(t)}{t}\,dt = 2\varsigma(5)$$

We used the following integral earlier in this theorem

$$\int \frac{x^k \log x}{1-x}\,dx = \sum_{j=1}^{k}\frac{x^j}{j^2} - \log x\sum_{j=1}^{k}\frac{x^j}{j} + Li_2(1-x)$$

and hence we have

$$\int_0^x \frac{x^k \log x}{1-x}\,dx = \sum_{j=1}^{k}\frac{x^j}{j^2} - \log x\sum_{j=1}^{k}\frac{x^j}{j} + Li_2(1-x) - Li_2(1)$$

As $k \to \infty$ the integral approaches zero for $|x| < 1$ and, from the right-hand side, we deduce Euler's identity for the dilogarithm

$$\sum_{j=1}^{\infty}\frac{x^j}{j^2} - \log x\sum_{j=1}^{\infty}\frac{x^j}{j} + Li_2(1-x) - Li_2(1) = 0$$

$$Li_2(x) + \log x\log(1-x) + Li_2(1-x) - \varsigma(2) = 0$$

Using (4.4.238aa) and completing the summation we obtain

(4.4.245o)     $$\int_0^x \frac{Li_p(x)\log x}{1-x}\,dx = \sum_{k=1}^{\infty}\frac{1}{k^p}\sum_{j=1}^{k}\frac{x^j}{j^2} - \log x\sum_{k=1}^{\infty}\frac{1}{k^p}\sum_{j=1}^{k}\frac{x^j}{j} + \varsigma(p)\,Li_2(1-x) - \varsigma(p)\varsigma(2)$$

**Theorem 4.14:**

(4.4.246)     $$\int_0^1 x^n \log^2(1-x)\,dx = \frac{1}{n+1}\left[\left(H_{n+1}^{(1)}\right)^2 + H_{n+1}^{(2)}\right] = \frac{2}{n+1}\left[H_{n+1}^{(2)} + \sum_{k=1}^{n}\frac{H_k^{(1)}}{k+1}\right]$$

**Proof:**

Using integration by parts we have

$$\int x^{k-1} \log(1-x)\,dx = -\frac{x^k}{k} \log(1-x) + \frac{1}{k} \int \frac{x^k}{1-x}\,dx$$

We provided an elementary derivation of the following integral in (4.4.148)

$$\int \frac{x^k}{1-x}\,dx = -\log(1-x) - \sum_{p=1}^{k} \frac{x^p}{p}$$

and we therefore have

$$k \int x^{k-1} \log(1-x)\,dx = -(1-x^k) \log(1-x) - \sum_{p=1}^{k} \frac{x^p}{p}$$

Hence we obtain

$$-\int_0^1 x^{k-1} \log(1-x)\,dx = \frac{H_k}{k}$$

Therefore we have

$$-\sum_{k=1}^{n+1} \int_0^1 x^{k-1} \log(1-x)\,dx = \sum_{k=1}^{n+1} \frac{H_k}{k}$$

We now recall Adamchik's identity (3.17)

$$\sum_{k=1}^{n+1} \frac{H_k^{(1)}}{k} = \frac{1}{2} \left( H_{n+1}^{(1)} \right)^2 + \frac{1}{2} H_{n+1}^{(2)}$$

and obtain

$$-\int_0^1 \frac{(1-x^{n+1})}{1-x} \log(1-x)\,dx = \frac{1}{2} \left( H_{n+1}^{(1)} \right)^2 + \frac{1}{2} H_{n+1}^{(2)}$$

With integration by parts we have

$$-\int_0^1 (1-x^{n+1}) \frac{\log(1-x)}{1-x}\,dx = \frac{(1-x^{n+1})}{2} \log^2(1-x)\Big|_0^1 + \frac{(n+1)}{2} \int_0^1 x^n \log^2(1-x)\,dx$$

$$= \frac{(n+1)}{2} \int_0^1 x^n \log^2(1-x)\,dx$$

Therefore we deduce (see a different proof in (4.4.155h))

$$\int_0^1 x^n \log^2(1-x)\,dx = \frac{2}{n+1}\left[\frac{1}{2}\left(H_{n+1}^{(1)}\right)^2 + \frac{1}{2}H_{n+1}^{(2)}\right]$$

This identity may also be obtained using the second derivative of the beta function.

Now, with simple algebra, we have

$$\sum_{k=1}^n \frac{H_k^{(1)}}{k+1} = \sum_{k=1}^n \frac{H_{k+1}^{(1)}}{k+1} - \sum_{k=1}^n \frac{1}{(k+1)^2}$$

$$= \sum_{k=1}^{n+1} \frac{H_k^{(1)}}{k} - H_{n+1}^{(2)}$$

$$\sum_{k=1}^n \frac{H_k^{(1)}}{k+1} + H_{n+1}^{(2)} = \sum_{k=1}^{n+1} \frac{H_k^{(1)}}{k} = \frac{1}{2}\left(H_{n+1}^{(1)}\right)^2 + \frac{1}{2}H_{n+1}^{(2)}$$

Hence we have

$$\int_0^1 x^n \log^2(1-x)\,dx = \frac{2}{n+1}\left[H_{n+1}^{(2)} + \sum_{k=1}^n \frac{H_k^{(1)}}{k+1}\right]$$

This result was obtained by De Doelder [55] in 1991 (and this corrects the misprint contained in the recent paper by Freitas [69a, p.8]).

From (4.4.246) we have

$$I_{n-1} = \int_0^1 x^{n-1} \log^2(1-x)\,dx = \frac{1}{n}\left[\left(H_n^{(1)}\right)^2 + H_n^{(2)}\right] = \frac{2}{n}\left[H_n^{(2)} + \sum_{k=1}^{n-1} \frac{H_k^{(1)}}{k+1}\right]$$

Let us consider the series $\sum_{n=1}^\infty \frac{I_{n-1}}{n}$. Using the familiar representation $\frac{1}{n} = \int_0^1 t^{n-1}dt$ we have

$$\sum_{n=1}^{\infty} \frac{I_{n-1}}{n} = \sum_{n=1}^{\infty} \int_0^1 t^{n-1} \, dt \int_0^1 x^{n-1} \log^2(1-x) \, dx$$

$$= \int_0^1 \sum_{n=1}^{\infty} t^{n-1} x^{n-1} \int_0^1 \log^2(1-x) \, dx \, dt$$

$$= \int_0^1 \int_0^1 \frac{\log^2(1-x)}{1-xt} \, dx \, dt$$

$$= \int_0^1 \frac{\log^3(1-x)}{x} \, dx$$

For the right-hand side we get

$$\sum_{n=1}^{\infty} \frac{1}{n^2} \left[ \left( H_n^{(1)} \right)^2 + H_n^{(2)} \right] = \sum_{n=1}^{\infty} \frac{\left( H_n^{(1)} \right)^2}{n^2} + \sum_{n=1}^{\infty} \frac{H_n^{(2)}}{n^2}$$

and these Euler sums have already been evaluated in (4.4.168) and (4.4.43) as

$$\sum_{n=1}^{\infty} \frac{\left( H_n^{(1)} \right)^2}{n^2} = \frac{17}{4} \varsigma(4) \qquad \sum_{n=1}^{\infty} \frac{H_n^{(2)}}{n^2} = \frac{7}{4} \varsigma(4)$$

Therefore we get the well-known result

$$\int_0^1 \frac{\log^3(1-x)}{x} \, dx = 6\varsigma(4)$$

which is a particular case of (3.108a).

From (3.82) we have

$$\frac{1}{2} \log^2(1-x) + Li_2(x) = \sum_{k=1}^{\infty} \frac{H_k}{k} x^k$$

and hence upon integration we get

$$\frac{1}{2} \int_0^1 x^{n-1} \log^2(1-x) \, dx + \int_0^1 x^{n-1} Li_2(x) \, dx = \int_0^1 \sum_{k=1}^{\infty} \frac{H_k}{k} x^{k+n-1} \, dx$$

This results in

(4.4.246a)        $\dfrac{1}{2n}\left[\left(H_n^{(1)}\right)^2 + H_n^{(2)}\right] + \displaystyle\sum_{k=1}^{\infty}\dfrac{1}{k^2(k+n)} = \sum_{k=1}^{\infty}\dfrac{H_k}{k(k+n)}$

Employing partial fractions we have

$$\sum_{k=1}^{\infty}\frac{1}{k^2(k+n)} = \frac{1}{n}\sum_{k=1}^{\infty}\frac{1}{k^2} - \frac{1}{n}\sum_{k=1}^{\infty}\frac{1}{k(k+n)}$$

$$= \frac{1}{n}\sum_{k=1}^{\infty}\frac{1}{k^2} - \frac{1}{n^2}\sum_{k=1}^{\infty}\left(\frac{1}{k} - \frac{1}{(k+n)}\right)$$

$$= \frac{1}{n}\sum_{k=1}^{\infty}\frac{1}{k^2} - \frac{H_n}{n^2}$$

$$= \frac{1}{n}\varsigma(2) - \frac{H_n}{n^2}$$

Therefore we have

$$\sum_{k=1}^{\infty}\frac{H_k}{k(k+n)} = \frac{1}{2n}\left[\left(H_n^{(1)}\right)^2 + H_n^{(2)}\right] + \frac{1}{n}\varsigma(2) - \frac{H_n}{n^2}$$

In their paper, "When is 0.999...equal to 1?", Pemantle and Schneider [105aa] show by means of the summation package Sigma in the computer algebra system Mathematica that

$$\sum_{k=1}^{N}\frac{H_k}{k(k+n)} = \frac{1}{2}\frac{\left(H_n^{(1)}\right)^2}{n} - \frac{H_n^{(1)}}{n^2} + \frac{1}{2}\frac{H_n^{(2)}}{n} + \frac{H_N^{(2)}}{n} - \frac{(nH_N^{(1)}-1)}{n^2}\sum_{j=1}^{n}\frac{1}{N+j} - \frac{1}{n}\sum_{i=1}^{n}\frac{1}{i}\sum_{j=1}^{i}\frac{1}{N+j}$$

and hence as $N \to \infty$ they obtained the result above. In turn this was employed to prove that

$$\sum_{j,k=1}^{\infty}\frac{H_j(H_{k+1}-1)}{jk(k+1)(j+k)} = -4\varsigma(2) - 2\varsigma(3) + 4\varsigma(2)\varsigma(3) + 2\varsigma(5)$$

$$= 0.999222...$$

(and the result explains the title of their paper).

Shortly thereafter in 2005, this result was also proved by Panholzer and Prodinger [105a] without requiring computer algebra.

We also have by completing the summation of (4.4.246)

$$\sum_{n=1}^{\infty} \frac{1}{n^s} \int_0^1 x^n \frac{\log^2(1-x)}{x}\, dx = \sum_{n=1}^{\infty} \frac{1}{n^{s+1}} \left[ \left( H_n^{(1)} \right)^2 + H_n^{(2)} \right]$$

and therefore

(4.4.246b)
$$\int_0^1 Li_s(x) \frac{\log^2(1-x)}{x}\, dx = \sum_{n=1}^{\infty} \frac{1}{n^{s+1}} \left[ \left( H_n^{(1)} \right)^2 + H_n^{(2)} \right]$$

We also see that

$$\sum_{n=1}^{\infty} (-1)^{n+1} \int_0^1 x^n \frac{\log^2(1-x)}{x}\, dx = \sum_{n=1}^{\infty} \frac{(-1)^{n+1}}{n} \left[ \left( H_n^{(1)} \right)^2 + H_n^{(2)} \right]$$

and hence

$$\int_0^1 \frac{\log^2(1-x)}{1+x}\, dx = \sum_{n=1}^{\infty} \frac{(-1)^{n+1}}{n} \left[ \left( H_n^{(1)} \right)^2 + H_n^{(2)} \right]$$

The Wolfram Integrator gives us

(4.4.246c)

$$\int \frac{\log^2(1-x)}{1+x}\, dx = \log\left[ (1+x)/2 \right] \log^2(1-x) + 2Li_2\left[ (1-x)/2 \right] \log(1-x) - 2Li_3\left[ (1-x)/2 \right]$$

and thus we have

(4.4.246d)
$$\int_0^1 \frac{\log^2(1-x)}{1+x}\, dx = 2Li_3(1/2) = \sum_{n=1}^{\infty} \frac{(-1)^n}{n} \left[ \left( H_n^{(1)} \right)^2 + H_n^{(2)} \right]$$

We note that $\sum_{n=1}^{\infty} (-1)^n \frac{\left( H_n^{(1)} \right)^2}{n}$ may be obtained from (3.110a).

We saw in (4.4.232a) that

$$\sum_{n=1}^{\infty} \frac{H_n^{(p)}}{n^q} + \sum_{n=1}^{\infty} \frac{H_n^{(q)}}{n^p} = \varsigma(p)\varsigma(q) + \varsigma(p+q)$$

Not long ago (April 2006) I realised that this identity could be easily generalised as follows. We have

$$\sum_{n=1}^{\infty} \frac{H_n^{(p)}}{n^q} x^n = \sum_{n=1}^{\infty} \frac{x^n}{n^q} \sum_{k=1}^{n} \frac{1}{k^p} = \sum_{n=1}^{\infty} \frac{1}{n^p} \sum_{k=n}^{\infty} \frac{x^k}{k^q}$$

and for $\mathrm{Re}(p) > 1$ this is equal to

$$= \sum_{n=1}^{\infty} \frac{1}{n^p} \left[ Li_q(x) - \sum_{k=1}^{n-1} \frac{x^k}{k^q} \right] = \varsigma(p) Li_q(x) - \sum_{n=1}^{\infty} \frac{1}{n^p} \left[ \sum_{k=1}^{n} \frac{x^k}{k^q} - \frac{x^n}{n^q} \right]$$

and we end up with

(4.4.247a)     $$\sum_{n=1}^{\infty} \frac{H_n^{(p)}}{n^q} x^n - \varsigma(p) Li_q(x) + \sum_{n=1}^{\infty} \frac{1}{n^p} \sum_{k=1}^{n} \frac{x^k}{k^q} = Li_{p+q}(x)$$

For example we have for $-1 \le x < 1$

$$Li_3(x) = \sum_{n=1}^{\infty} \frac{H_n^{(2)}}{n} x^n - \varsigma(2) Li_1(x) + \sum_{n=1}^{\infty} \frac{1}{n^2} \sum_{k=1}^{n} \frac{x^k}{k}$$

(4.4.247b)     $$= \sum_{n=1}^{\infty} \frac{H_n^{(2)}}{n} x^n + \log(1-x)\varsigma(2) + \sum_{n=1}^{\infty} \frac{1}{n^2} \sum_{k=1}^{n} \frac{x^k}{k}$$

For $(p,q) = (2,2)$ and $(3,1)$ respectively we get

(4.4.247ci)     $$Li_4(x) = \sum_{n=1}^{\infty} \frac{H_n^{(2)}}{n^2} x^n - \varsigma(2) Li_2(x) + \sum_{n=1}^{\infty} \frac{1}{n^2} \sum_{k=1}^{n} \frac{x^k}{k^2}$$

and (also refer to (4.4.168i))

(4.4.247cii)     $$Li_4(x) = \sum_{n=1}^{\infty} \frac{H_n^{(3)}}{n} x^n + \varsigma(3) \log(1-x) + \sum_{n=1}^{\infty} \frac{1}{n^3} \sum_{k=1}^{n} \frac{x^k}{k}$$

Accordingly, by a rather elementary procedure, we have demonstrated that

(4.4.247ciii)

$$\sum_{n=1}^{\infty} \frac{H_n^{(2)}}{n^2} x^n - \varsigma(2) Li_2(x) + \sum_{n=1}^{\infty} \frac{1}{n^2} \sum_{k=1}^{n} \frac{x^k}{k^2} = \sum_{n=1}^{\infty} \frac{H_n^{(3)}}{n} x^n + \varsigma(3) \log(1-x) + \sum_{n=1}^{\infty} \frac{1}{n^3} \sum_{k=1}^{n} \frac{x^k}{k}$$

Similarly, with $(p,q) = (2,3), \ (3,2)$ and $(4,1)$ we have

(4.4.247di) $\qquad Li_5(x) = \sum_{n=1}^{\infty} \dfrac{H_n^{(2)}}{n^3} x^n - \varsigma(2)Li_3(x) + \sum_{n=1}^{\infty} \dfrac{1}{n^2} \sum_{k=1}^{n} \dfrac{x^k}{k^3}$

(4.4.247dii) $\qquad Li_5(x) = \sum_{n=1}^{\infty} \dfrac{H_n^{(3)}}{n^2} x^n - \varsigma(3)Li_2(x) + \sum_{n=1}^{\infty} \dfrac{1}{n^3} \sum_{k=1}^{n} \dfrac{x^k}{k^2}$

(4.4.247diii) $\qquad Li_5(x) = \sum_{n=1}^{\infty} \dfrac{H_n^{(4)}}{n} x^n + \varsigma(4)\log(1-x) + \sum_{n=1}^{\infty} \dfrac{1}{n^4} \sum_{k=1}^{n} \dfrac{x^k}{k}$

We may note that, for example, dividing (4.4.247ci) by $x$ and integrating results in (4.4.247di).

In the same way we see that

$$\sum_{n=1}^{\infty} \frac{H_n^{(p)}}{n^q} (-1)^n x^n = \sum_{n=1}^{\infty} (-1)^n \frac{x^n}{n^q} \sum_{k=1}^{n} \frac{1}{k^p} = \sum_{n=1}^{\infty} \frac{1}{n^p} \sum_{k=n}^{\infty} (-1)^k \frac{x^k}{k^q}$$

and for $\operatorname{Re}(p) > 1$ this is equal to

$$= \sum_{n=1}^{\infty} \frac{1}{n^p} \left[ Li_q(-x) - \sum_{k=1}^{n-1} (-1)^k \frac{x^k}{k^q} \right]$$

$$= \varsigma(p)Li_q(-x) - \sum_{n=1}^{\infty} \frac{1}{n^p} \left[ \sum_{k=1}^{n} (-1)^k \frac{x^k}{k^q} - (-1)^n \frac{x^n}{n^q} \right]$$

and we end up with

(4.4.247e) $\qquad \sum_{n=1}^{\infty} (-1)^n \dfrac{H_n^{(p)}}{n^q} x^n - \varsigma(p)Li_q(-x) + \sum_{n=1}^{\infty} \dfrac{1}{n^p} \sum_{k=1}^{n} (-1)^k \dfrac{x^k}{k^q} = Li_{p+q}(-x)$

(letting $x \to -x$ in (4.4.247c) would of course also suffice).

We also have

$$\sum_{n=1}^{\infty} \frac{x^n}{n^q} \sum_{k=1}^{n} \frac{y^k}{k^p} = \sum_{n=1}^{\infty} \frac{y^n}{n^p} \sum_{k=n}^{\infty} \frac{x^k}{k^q}$$

and for $\operatorname{Re}(p) > 1$ this is equal to

$$= \sum_{n=1}^{\infty} \frac{y^n}{n^p} \left[ Li_q(x) - \sum_{k=1}^{n-1} \frac{x^k}{k^q} \right] = Li_p(y) Li_q(x) - \sum_{n=1}^{\infty} \frac{y^n}{n^p} \left[ \sum_{k=1}^{n} \frac{x^k}{k^q} - \frac{x^n}{n^q} \right]$$

and we end up with

(4.4.247f) $\qquad \sum_{n=1}^{\infty} \frac{x^n}{n^q} \sum_{k=1}^{n} \frac{y^k}{k^p} - Li_p(y) Li_q(x) + \sum_{n=1}^{\infty} \frac{y^n}{n^p} \sum_{k=1}^{n} \frac{x^k}{k^q} = Li_{p+q}(xy)$

With $y = -1$ we obtain for $p > 0$

(4.4.247fi) $\qquad \sum_{n=1}^{\infty} \frac{x^n}{n^q} \sum_{k=1}^{n} \frac{(-1)^k}{k^p} + \varsigma_a(p) Li_q(x) + \sum_{n=1}^{\infty} \frac{(-1)^n}{n^p} \sum_{k=1}^{n} \frac{x^k}{k^q} = Li_{p+q}(-x)$

A specific example is shown below (again in full detail)

$$\sum_{n=1}^{\infty} \frac{H_n^{(3)}}{n 2^n} = \sum_{n=1}^{\infty} \frac{1}{n 2^n} \sum_{k=1}^{n} \frac{1}{k^3} = \sum_{n=1}^{\infty} \frac{1}{n^3} \sum_{k=n}^{\infty} \frac{1}{k 2^k}$$

$$= \sum_{n=1}^{\infty} \frac{1}{n^3} \left[ \log 2 - \sum_{k=1}^{n-1} \frac{1}{k 2^k} \right]$$

$$= \varsigma(3) \log 2 - \sum_{n=1}^{\infty} \frac{1}{n^3} \left[ \sum_{k=1}^{n} \frac{1}{k 2^k} - \frac{1}{n 2^n} \right]$$

$$= \varsigma(3) \log 2 - \sum_{n=1}^{\infty} \frac{1}{n^3} \sum_{k=1}^{n} \frac{1}{k 2^k} + \sum_{n=1}^{\infty} \frac{1}{2^n n^4}$$

$$= \varsigma(3) \log 2 - \sum_{n=1}^{\infty} \frac{1}{n^3} \sum_{k=1}^{n} \frac{1}{k 2^k} + Li_4(1/2)$$

In the reverse direction we see that

$$\sum_{n=1}^{\infty} \frac{1}{n^3} \sum_{k=1}^{n} \frac{1}{k 2^k} = \sum_{n=1}^{\infty} \frac{1}{n 2^n} \sum_{k=n}^{\infty} \frac{1}{k^3} = \sum_{n=1}^{\infty} \frac{1}{n 2^n} \left[ \varsigma(3) - H_{n-1}^{(3)} \right]$$

$$= \sum_{n=1}^{\infty} \frac{1}{n 2^n} \left[ \varsigma(3) - H_{n-1}^{(3)} \right]$$

**Remark:**

Consider the series

$$(4.4.247g) \quad S_2(-x,s) = \sum_{n=1}^{\infty} \frac{1}{n^2 2^n} \sum_{k=1}^{n} \binom{n}{k} \frac{(-1)^k x^k}{k^s}$$

**Comment:**

We have

$$\sum_{n=1}^{\infty} \frac{S_n^a}{n^2 2^n} = \sum_{n=1}^{\infty} \frac{1}{n^2 2^n} \frac{(-1)^s nx}{s\Gamma(s)} \int_0^1 \log^s t (1-xt)^{n-1} dt$$

$$= \frac{(-1)^s x}{2\Gamma(s+1)} \int_0^1 du \int_0^1 \log^s t \sum_{n=1}^{\infty} u^{n-1} \frac{(1-xt)^{n-1}}{2^{n-1}} dt$$

$$= \frac{(-1)^s x}{2\Gamma(s+1)} \int_0^1 du \int_0^1 \frac{\log^s t}{1 - \frac{u}{2}(1-xt)} dt$$

$$= \frac{(-1)^s x}{2\Gamma(s+1)} \int_0^1 \log^s t \, dt \int_0^1 \frac{1}{1 - \frac{u}{2}(1-xt)} du$$

$$= \frac{(-1)^s x}{\Gamma(s+1)} \int_0^1 \log^s t \, dt \left\{ -\frac{1}{1-xt} \log\left[ 1 - \frac{u}{2}(1-xt) \right] \Big|_0^1 \right\}$$

$$= \frac{(-1)^s x}{\Gamma(s+1)} \int_0^1 \frac{\log^s t \log\left[ (1+xt)/2 \right] xt}{1-xt} dt$$

$$= \frac{(-1)^s x}{\Gamma(s+1)} \int_0^1 \frac{\log^s t \log(1+xt)}{1-xt} dt - \frac{(-1)^s x \log 2}{\Gamma(s+1)} \int_0^1 \frac{\log^s t}{1-xt} dt$$

$$= \frac{(-1)^s x}{\Gamma(s+1)} \int_0^1 \frac{\log^s t \log(1+xt)}{1-xt} dt - \log 2 \, Li_{s+1}(x)$$

According to Freitas [6a], integrals of the above type (with $x=1$) were considered in 1981 by Gastmans and Troost [69b] but I have not yet had access to their paper. Reference may also be made profitably to Rutledge and Douglass [116aa] and [116ab] and Kölbig's paper on Nielsen's generalised polylogarithms [91a].

In [118a] Sasvári gave an elementary proof of Binet's formula for the gamma function

(4.4.248) $$\Gamma(x+1) = \left(\frac{x}{e}\right)^x \sqrt{2\pi x} \exp \vartheta(x)$$

where $$\vartheta(x) = \int_0^\infty \left(\frac{1}{e^t - 1} - \frac{1}{t} + \frac{1}{2}\right) \frac{e^{-xt}}{t} dt$$

Hence we have upon taking logarithms

(4.4.248a) $$\log \Gamma(x+1) = \left(x + \frac{1}{2}\right) \log x - x + \log \sqrt{2\pi} + \int_0^\infty \left(\frac{1}{e^t - 1} - \frac{1}{t} + \frac{1}{2}\right) \frac{e^{-xt}}{t} dt$$

We have the generating function for the Bernoulli numbers in (A.2) in Volume V

(4.4.248b) $$\frac{t}{e^t - 1} = \sum_{n=0}^\infty B_n \frac{t^n}{n!} = 1 - \frac{1}{2} t + \sum_{n=1}^\infty B_{2n} \frac{t^{2n}}{(2n)!} \quad , (|t| < 2\pi)$$

and, in an Eulerian fashion, we therefore obtain (ignoring subtleties such as convergence!)

$$\vartheta(x) \approx \int_0^\infty \sum_{n=1}^\infty B_{2n} \frac{t^{2n-1}}{(2n)!} \frac{e^{-xt}}{t} dt$$

$$= \sum_{n=1}^\infty \frac{B_{2n}}{(2n)!} \int_0^\infty t^{2n-2} e^{-xt} dt$$

From (4.4.57c) we have

$$\int_0^\infty t^{2n-2} e^{-xt} dt = \frac{(2n-2)!}{x^{2n-1}}$$

and hence we obtain

$$\vartheta(x) \approx \sum_{n=1}^\infty \frac{B_{2n}}{2n(2n-1)x^{2n-1}}$$

$$= \frac{1}{12x} - \frac{1}{360x^3} + \frac{1}{1260x^5} - \dots$$

We therefore have an asymptotic formula for $\Gamma(x+1)$

$$(4.4.249) \qquad \Gamma(x+1) \cong \left(\frac{x}{e}\right)^x \sqrt{2\pi x} \cdot \exp\left[\sum_{n=1}^{\infty} \frac{B_{2n}}{2n(2n-1)x^{2n-1}}\right]$$

$$\cong \left(\frac{x}{e}\right)^x \sqrt{2\pi x} \cdot \exp\left[\frac{1}{12x} - \frac{1}{360x^3} + \frac{1}{1260x^5} - \ldots\right]$$

The above formula is given by Bressoud in [34, p.302] and the result is quite extraordinary given that the interval of convergence of the generating function for the Bernoulli numbers is $|t| < 2\pi$ and yet we have integrated the resulting expression over the range $[0, \infty)$. Sasvári helps to explain the unusual outcome by referring to an inequality given by Póyla and Szegö [108a], namely

$$\sum_{n=1}^{2N} B_{2n} \frac{t^{2n-1}}{(2n)!} < \frac{1}{e^t - 1} - \frac{1}{t} + \frac{1}{2} < \sum_{n=1}^{2N+1} B_{2n} \frac{t^{2n-1}}{(2n)!}$$

and therefore

$$\int_0^{\infty} \sum_{n=1}^{2N} B_{2n} \frac{t^{2n-1}}{(2n)!} \frac{e^{-xt}}{t} dt \;<\; \int_0^{\infty} \left(\frac{1}{e^t - 1} - \frac{1}{t} + \frac{1}{2}\right) \frac{e^{-xt}}{t} dt \;<\; \int_0^{\infty} \sum_{n=1}^{2N+1} B_{2n} \frac{t^{2n-1}}{(2n)!} \frac{e^{-xt}}{t} dt$$

Hence we have

$$\sum_{n=1}^{2N} \frac{B_{2n}}{2n(2n-1)x^{2n-1}} \;<\; \int_0^{\infty} \left(\frac{1}{e^t - 1} - \frac{1}{t} + \frac{1}{2}\right) \frac{e^{-xt}}{t} dt \;<\; \sum_{n=1}^{2N+1} \frac{B_{2n}}{2n(2n-1)x^{2n-1}}$$

Integrating (4.4.248a) results in

$$(4.4.248c) \qquad \int_0^u \log\Gamma(x+1)dx = \frac{1}{2}u(u+1)\log u - \frac{3}{4}u^2 - \frac{1}{2}u + \frac{u}{2}\log 2\pi$$

$$+ \int_0^{\infty} \left(\frac{1}{e^t - 1} - \frac{1}{t} + \frac{1}{2}\right) \frac{1 - e^{-ut}}{t^2} dt$$

A simple derivation of the following formula is contained in [43a] and [25, p.217]

$$(4.4.250) \qquad \log x - \psi(x+1) = \int_0^{\infty} \left(\frac{1}{e^t - 1} - \frac{1}{t}\right) e^{-xt} dt$$

and this may also be obtained by differentiating (4.4.248a). Differentiating (4.4.250) we obtain

$$\frac{1}{x} - \psi'(x+1) = \int\limits_0^\infty \left( \frac{t}{e^t - 1} - 1 \right) e^{-xt} dt$$

With $x = 1$ we have

$$1 - \psi'(2) = \int\limits_0^\infty \left( \frac{t}{e^t - 1} - 1 \right) e^{-t} dt$$

We have the following equation in [126, p.15]

$$\log x - \psi(x) = \int\limits_0^x \left( \frac{1}{1 - e^{-t}} - \frac{1}{t} \right) e^{-xt} dt$$

$$= \int\limits_0^x \left( \frac{e^t}{e^t - 1} - \frac{1}{t} \right) e^{-xt} dt$$

$$= \int\limits_0^x \left( 1 + \frac{1}{e^t - 1} - \frac{1}{t} \right) e^{-xt} dt$$

(4.4.250a)
$$= \int\limits_0^x \left[ \frac{1}{2} + \sum_{n=1}^\infty B_{2n} \frac{t^{2n-1}}{(2n)!} \right] e^{-xt} dt$$

Using (A.16) this becomes

$$= \int\limits_0^x \left( \frac{1}{2} + \frac{1}{2} \coth\left[ \frac{t}{2} \right] - \frac{1}{t} \right) e^{-xt} dt$$

and therefore we get

$$\log x - \psi(x) = \frac{1}{2x} \left[ 1 - e^{-x^2} \right] + \int\limits_0^x \left( \frac{1}{2} \coth\left[ \frac{t}{2} \right] - \frac{1}{t} \right) e^{-xt} dt$$

Using (4.4.57b) we have

$$\int\limits_0^y u^m e^{-\mu u} du = \frac{m!}{\mu^{m+1}} - e^{-\mu y} \sum_{k=0}^m \frac{m! \, y^k}{k! \, \mu^{m-k+1}}$$

and accordingly we have

$$\int_0^x t^{2n-1} e^{-xt} dt = \frac{(2n-1)!}{x^{2n}} - e^{-x^2} \sum_{k=0}^{2n-1} \frac{(2n-1)!}{k! \, x^{2n-2k}}$$

Therefore we formally obtain from (4.4.250a)

$$\log x - \psi(x) = \frac{1}{2} \int_0^x e^{-xt} dt + \sum_{n=1}^{\infty} \frac{B_{2n}}{2nx^{2n}} \left[ 1 - e^{-x^2} \sum_{k=0}^{2n-1} \frac{x^{2k}}{k!} \right]$$

(4.4.251)
$$= \frac{1}{2x} - \frac{e^{-x^2}}{2x} + \sum_{n=1}^{\infty} \frac{B_{2n}}{2nx^{2n}} \left[ 1 - e^{-x^2} \sum_{k=0}^{2n-1} \frac{x^{2k}}{k!} \right]$$

Since $\displaystyle \lim_{n \to \infty} \left[ 1 - e^{-x^2} \sum_{k=0}^{2n-1} \frac{x^{2k}}{k!} \right] = 0$, it is clear that this greatly assists the convergence of the series. On the other hand, as mentioned in Appendix A of Volume VI, $|B_{2n}|$ increases very rapidly (for example, $B_{32}$ is greater than the number of people living on the planet Earth). The net effect of these opposing forces is that (4.4.251) is convergent for a certain range of $x$. Letting $x \to \sqrt{x}$ we have

$$\log \sqrt{x} - \psi(\sqrt{x}) = \frac{1}{2\sqrt{x}} - \frac{e^{-x}}{2\sqrt{x}} + e^{-x} \sum_{n=1}^{\infty} \frac{B_{2n}}{2nx^n} \left[ e^x - \sum_{k=0}^{2n-1} \frac{x^k}{k!} \right]$$

Using Taylor's theorem we have

$$e^x - \sum_{k=0}^{2n-1} \frac{x^k}{k!} = \frac{x^{2n}}{(2n)!} e^{\xi_n} \quad \text{where } 0 < \xi_n < x$$

and this gives us a factor of $(2n)!$ in the denominator to assist the convergence. Applying Euler's formula (1.7) for $B_{2n}$ we deduce

$$\log \sqrt{x} - \psi(\sqrt{x}) = \frac{1}{2\sqrt{x}} - \frac{e^{-x}}{2\sqrt{x}} + e^{-x} \sum_{n=1}^{\infty} \frac{(-1)^{n+1} (2n)! \, \zeta(2n)}{n(2\pi)^{2n} x^n} \cdot \frac{x^{2n}}{(2n)!} e^{\xi_n}$$

$$= \frac{1}{2\sqrt{x}} - \frac{e^{-x}}{2\sqrt{x}} + e^{-x} \sum_{n=1}^{\infty} \frac{(-1)^{n+1} \zeta(2n) e^{\xi_n} x^n}{n(2\pi)^{2n}}$$

If $(x_n)$ is a decreasing sequence of strictly positive numbers, and $\lim_{n \to \infty} x_n = 0$, then

$\sum (-1)^n x_n$ is convergent by the alternating series test [17a, p.308]. From (A.10) we know

that $\lim_{n \to \infty} \varsigma(n) = 1$ and therefore the above series is convergent provided $\sqrt{x} \le 2\pi$.

However, we need to consider the nature of the $e^{\xi_n}$ factor and I therefore have strong reservations as to the validity of this analysis.

With $x = 1$ we have

(4.4.252) $\qquad -\psi(1) = \dfrac{1}{2} - \dfrac{e^{-1}}{2} + \sum_{n=1}^{\infty} \dfrac{B_{2n}}{2n} \left[ 1 - e^{-1} \sum_{k=0}^{2n-1} \dfrac{1}{k!} \right]$

and, since $\psi(1) = -\gamma$, we have an expression for $\gamma$ involving $e$ and $\pi$

(4.4.252a) $\qquad \gamma = \dfrac{1}{2}\left[ 1 - \dfrac{1}{e} \right] + \sum_{n=1}^{\infty} \dfrac{(-1)^{n+1}(2n-1)!\,\varsigma(2n)}{(2\pi)^{2n}} \left[ 1 - \dfrac{1}{e} \sum_{k=0}^{2n-1} \dfrac{1}{k!} \right]$

Equation (4.4.252a) could be used to determine a numerical approximation for $\gamma$ and, due to the fact that we have an alternating series, it may be used to quantify the error in the approximation. It should be noted that I have not checked the validity of interchanging the order of the integration and summation in the above analysis.

The following formula was derived by Knopp [90, p.527] (see also [126, p.6])

(4.4.252b) $\qquad \gamma = \dfrac{1}{2} + \sum_{n=1}^{N} \dfrac{B_{2n}}{2n} - (2N+1)! \int_{0}^{\infty} \dfrac{Q_{2N+1}(x)}{x^{2N+2}}\, dx$

where the functions $Q_N(x)$ are defined by

$\qquad Q_N(x) = x - \dfrac{1}{2} \qquad\qquad (n=1; 0 < x < 1)$

$\qquad Q_N(x) = \dfrac{B_n(x)}{n!} \qquad\qquad (n \in N_0; 0 \le x \le 1)$

In [43a] Chao-Ping Chen and Feng Qi employ the following asymptotic expansion

$\qquad \log x - \psi(x) = \dfrac{1}{2x} + \dfrac{1}{12x^2} + O(x^{-4})$

whereas the Wolfram website records the following formula

$$\log x - \psi(x) \cong \frac{1}{2x} + \sum_{n=1}^{\infty} \frac{B_{2n}}{2nx^{2n}}$$

and both are in accordance with (4.4.251). In [126. p.22], Srivastava and Choi record the following asymptotic formula

$$\log x - \psi(x) \cong \frac{1}{2x} + \sum_{n=1}^{\infty} \frac{(-1)^{n-1} B_n}{2nx^{2n}}$$

and I believe that the difference arises because the authors have used the different definition of the Bernoulli polynomials contained in Whittaker and Watson [135, p.125].

$$f(x) = \log \sqrt{x} - \psi(\sqrt{x}) - \frac{1}{2\sqrt{x}} + \frac{e^{-x}}{2\sqrt{x}} = e^{-x} \sum_{n=1}^{\infty} \frac{B_{2n}}{2nx^n} \left[ e^x - \sum_{k=0}^{2n-1} \frac{x^k}{k!} \right]$$

We have the integral form of Taylor's theorem [66, p.698]

$$F(x) = F(0) + \frac{x}{1!} F(0) + \dots + \frac{x^{n-1}}{(n-1)!} F^{(n-1)}(0) + R_n$$

where $R_n = \dfrac{1}{(n-1)!} \displaystyle\int_0^x (x-t)^{n-1} F^{(n)}(t)\, dt$. Therefore we have with $F(x) = e^x$

$$e^x - \sum_{k=0}^{2n-1} \frac{x^k}{k!} = \frac{1}{(2n-1)!} \int_0^x (x-t)^{2n-1} e^t\, dt = \frac{1}{(2n-1)!} \int_0^x \sum_{k=0}^{2n-1} \binom{2n-1}{k} (-1)^k x^{2n-1-k} t^k e^t\, dt$$

Therefore we may write

$$f(x) = \log \sqrt{x} - \psi(\sqrt{x}) - \frac{1}{2\sqrt{x}} + \frac{e^{-x}}{2\sqrt{x}} = e^{-x} \sum_{n=1}^{\infty} \frac{B_{2n}}{(2n)! x^n} \int_0^x (x-t)^{2n-1} e^t\, dt$$

The incomplete gamma function is defined as

$$\Gamma(a, x) = \int_x^{\infty} t^{a-1} e^{-t}\, dt$$

and making the substitution $t = u^2$ we have

$$\Gamma(a,x) = \int\limits_{x^2}^{\infty} u^{2a-1} e^{-u^2} du$$

We also have

$$\Gamma(n,x) = (n-1)! \, e^{-x} \sum_{k=0}^{n-1} \frac{x^k}{k!}$$

and therefore

$$\Gamma(2n,x^2) = (2n-1)! \, e^{-x^2} \sum_{k=0}^{2n-1} \frac{x^{2k}}{k!}$$

$$\log x - \psi(x) = \frac{1}{2x} - \frac{e^{-x^2}}{2x} + \sum_{n=1}^{\infty} \frac{B_{2n}}{x^{2n}} \left[ \frac{1}{2n} - \frac{\Gamma(2n,x^2)}{(2n)!} \right]$$

Math-ph/0408036 [abs, ps, pdf, other]

1997.
On the Rapid Computation of Various Polylogarithmic Constants

$$\sum_{k=1}^{\infty}\frac{1}{k^2}=\frac{\pi^2}{6} \ \text{ and } \ \sin x = x\prod_{k=1}^{\infty}\left(1-\frac{x^2}{k^2\pi^2}\right).$$ Math. Mag.69, 122-125, 1996.

*Corrected Proof, Available online 12 February 2007.*

Donal F. Connon
Elmhurst
Dundle Road
Matfield
Kent TN12 7HD
dconnon@btopenworld.com